\documentclass[10pt]{article}

\usepackage{tikz-cd}

\usepackage{tikz}
\usetikzlibrary{matrix,arrows,decorations.pathmorphing,shapes.geometric}
\usepackage[percent]{overpic} 
\usepackage{etex}
\usepackage{subcaption}
\usepackage[utf8]{inputenc}
\usepackage{dsfont}
\usepackage{a4wide}
\usepackage{graphicx}
\usepackage{wrapfig}
\usepackage[hmarginratio={1:1},     % equal left and right margins
vmarginratio={1:1},     % equal top and bottom margins
textwidth=17cm,        % new text width
textheight=23cm,
heightrounded,]{geometry}
\usepackage{bbm}

\usepackage{enumitem} 
\renewcommand{\theenumii}{\arabic{enumii}}

\renewcommand{\labelenumi}{(\roman{enumi})}

\usepackage{amsmath,accents}
\usepackage{amsthm}
\usepackage{amssymb}

\usepackage{mathrsfs, mathtools}

\usepackage{tikz}
\usepackage{textcomp}
\usepackage{pdfpages}
\usepackage{tikz-cd}
\usepackage{tikz-3dplot}
\usepackage{pgfplots}
\usetikzlibrary{calc}
\pgfplotsset{width=7cm,compat=1.8}
\usetikzlibrary{decorations.pathreplacing}
\usetikzlibrary{decorations.markings}
\usetikzlibrary{patterns}
\usepgflibrary{shapes.geometric}

\tikzset{
	string/.style={draw=#1, postaction={decorate}, decoration={markings,mark=at position .51 with {\arrow[color=#1]{>}}}},
	costring/.style={draw=#1, postaction={decorate}, decoration={markings,mark=at position .51 with {\arrow[draw=#1]{<}}}},
	ostring/.style={draw=#1, postaction={decorate}, decoration={markings,mark=at position .47 with {\arrow[draw=#1]{>}}}},
	ustring/.style={draw=#1, postaction={decorate}, decoration={markings,mark=at position .56 with {\arrow[draw=#1]{>}}}},
	oostring/.style={draw=#1, postaction={decorate}, decoration={markings,mark=at position .43 with {\arrow[draw=#1]{>}}}},
	uustring/.style={draw=#1, postaction={decorate}, decoration={markings,mark=at position .59 with {\arrow[draw=#1]{>}}}},
	directed/.style={string=blue!50!black}, 
	odirected/.style={ostring=blue!50!black}, 
	udirected/.style={ustring=blue!50!black}, 
	oodirected/.style={oostring=blue!50!black}, 
	uudirected/.style={uustring=blue!50!black},     
	redirected/.style={costring= blue!50!black},
	redirectedgreen/.style={costring= green!50!black},
	directedgreen/.style={string= green!50!black},
	redirectedlightgreen/.style={costring= green!65!black},
	directedlightgreen/.style={string= green!65!black},
	redirectedred/.style={costring= red!50!black},
	directedred/.style={string= red!50!black}%
}

\tikzset{-dot-/.style={decoration={
			markings,
			mark=at position 0.5 with {\fill circle (1.875pt);}},postaction={decorate}}}

\tikzset{
	Fdot/.style={circle, draw, fill, inner sep=0pt}, 
	Odot/.style={circle, draw, inner sep=0.1pt, minimum size=0.1cm}
}

\usepackage{import}
\usepackage{calc}
\newcommand{\pdfpic}[2][1.0]{
	\begin{tikzpicture}[baseline={([yshift=-.5ex]current bounding box.center)}]
	\node at (0,0) {\def\svgscale{#1} \import{}{#2} };
	\end{tikzpicture}
}

\def\nicedashedcolourscheme{\shadedraw[top color=blue!22, bottom color=blue!22, draw=gray, dashed]}
\def\nicedashedpalecolourscheme{\shadedraw[top color=blue!12, bottom color=blue!12, draw=gray, dashed]}
\def\nicehalfpalecolourscheme{\shadedraw[top color=blue!22, bottom color=blue!22, draw=white]}
\def\nicenotpalecolourscheme{\shadedraw[top color=blue!32, bottom color=blue!32, draw=white]}
\def\nicecolourscheme{\shadedraw[top color=blue!22, bottom color=blue!22, draw=blue!22]}
\def\nicepalecolourscheme{\shadedraw[top color=blue!12, bottom color=blue!12, draw=white]}
\def\nicenocolourscheme{\shadedraw[top color=gray!2, bottom color=gray!25, draw=white]}
\def\nicereallynocolourscheme{\shadedraw[top color=white!2, bottom color=white!25, draw=white]}
\def\boringcolourscheme{\draw[fill=blue!20, dashed]}

\newcommand\tikzzbox[1]
{#1}

\usepackage[all,cmtip]{xy}

\let\emph\undefined
\newcommand{\emph}[1]{\textsl{#1}}
\let\itshape\undefined
\let\itshape\slshape

\definecolor{Myblue}{rgb}{0,0,0.6}  
\usepackage[colorlinks,citecolor=Myblue,linkcolor=Myblue,urlcolor=Myblue,pdfpagemode=None]{hyperref}

\numberwithin{equation}{section}
\usepackage{mathtools}
\numberwithin{equation}{section}
\numberwithin{figure}{section}

\newtheoremstyle{style1}% name of the style to be used
{13pt}% measure of space to leave above the theorem. E.g.: 3pt
{13pt}% measure of space to leave below the theorem. E.g.: 3pt
{}% name of font to use in the body of the theorem
{}% measure of space to indent
{\normalfont\bfseries}% name of head font
{.}% punctuation between head and body
{.5em}% space after theorem head; " " = normal interword space
{}

\theoremstyle{style1}

\newtheorem{definition}{Definition}[section]
\newtheorem{example}[definition]{Example}
\newtheorem{remark}[definition]{Remark}

\newcommand{\catf}[1]{{\mathrm{#1}}}
\newcommand{\red}[1]{\color{red} #1 \color{black}}
\usepackage{tocloft}
\let\P\undefined
\newcommand{\P}{\mathcal{P}}
\newtheoremstyle{style2}% name of the style to be used
{13pt}% measure of space to leave above the theorem. E.g.: 3pt
{13pt}% measure of space to leave below the theorem. E.g.: 3pt
{\slshape}% name of font to use in the body of the theorem
{}% measure of space to indent
{\normalfont\bfseries}% name of head font
{.}% punctuation between head and body
{.5em}% space after theorem head; " " = normal interword space
{}

\newtheorem{lemma}[definition]{Lemma}
\newtheorem{theorem}[definition]{Theorem}
\newtheorem{proposition}[definition]{Proposition}

\def\hocolim{\mathrm{hocolim}}

\newcommand\arxiv[2]      {\href{https://arXiv.org/abs/#1}{#2}}
\newcommand\doi[2]        {\href{https://dx.doi.org/#1}{#2}}

\usepackage{tikz}
\usetikzlibrary{matrix,arrows,decorations.pathmorphing,shapes.geometric}
\usepackage{tikz-cd}

\newcommand*{\tarrow}[2][]{\arrow[Rrightarrow, #1]{#2}\arrow[dash, shorten >= 0.5pt, #1]{#2}}
\newcommand*{\qarrow}[2][]{\arrow[RRightarrow, #1]{#2}\arrow[equal, double distance = 0.25pt, shorten >= 1.28pt, #1]{#2}}

\newcommand{\Woike}[1]{\color{blue} LW: #1\color{black}}

\usepackage{enumitem}

\newcommand{\Grpd}{\catf{Grpd}}

\usepackage{needspace}
\newcommand{\spaceplease}{\needspace{5\baselineskip}}

\newcommand{\SO}{\operatorname{SO}}
\newcommand{\Or}{\operatorname{O}}

\newcommand{\C}{\mathds{C}}
\newcommand{\N}{\mathds{N}}
\newcommand{\R}{\mathds{R}} 
\newcommand{\Z}{\mathds{Z}}
\newcommand{\K}{\Bbbk}

\newcommand{\Ca}{\mathcal{C}}
\newcommand{\Red}{\catf{Red}}
\newcommand{\cen}{\catf{C}}
\newcommand{\tft}{\catf{TFT}}
\newcommand{\gtft}{G\text{-}\catf{TFT}}
\newcommand{\PBun}{\catf{PBun}}

\newcommand{\cat}[1]{\mathcal{#1}}

\newcommand{\colim}{\operatorname{colim}}
\newcommand{\Aut}{\operatorname{Aut}}
\newcommand{\End}{\operatorname{End}}
\newcommand{\Hom}{\operatorname{Hom}}
\newcommand{\im}{\operatorname{im}}
\newcommand{\id}{\operatorname{id}}
\newcommand{\Sym}{{\mathrm{Sym}}}
\newcommand{\ev}{\operatorname{ev}}
\newcommand{\tr}{\operatorname{tr}}
\newcommand{\EW}{\operatorname{EW}}
\newcommand{\btimes}{\mathbin{\square}} % operator for tensor product in a bicategory
\renewcommand{\Box}{\btimes}
\newcommand{\relpro}{\operatornamewithlimits{\btimes}}

\newcommand{\Alg}{\catf{Alg}}
\newcommand{\sAlg}{\catf{Alg}^{\operatorname{s}}}
\newcommand{\fvs}{\catf{FinVect}}
\newcommand{\vs}{\catf{Vect}}
\newcommand{\Tvs}{{2\catf{Vect}}}
\newcommand{\TTvs}{{3\catf{Vect}}}
\newcommand{\U}{\text{\normalfont U}}

\newcommand{\Cob}{{\catf{Cob}}}
\newcommand{\DBord}{{\catf{Bord}_{3,2}^{\operatorname{def}}}}
\newcommand{\Op}{{\catf{Op}_\mathfrak{C}}}
\newcommand{\SymS}{{\catf{Sym}_\mathfrak{C}}}
\newcommand{\Cat}{\catf{Cat}}
\newcommand{\ssFrob}{\catf{ssFrob}}
\let\to\undefined
\newcommand{\to}{\longrightarrow}
\let\mapsto\undefined
\newcommand{\mapsto}{\longmapsto}
\newcommand{\TVectBun}{2\catf{VecBunGrpd}}
\let\FinVect\undefined
\newcommand{\FinVect}{\catf{Vect}}
\newcommand{\Mod}{\catf{Mod}}
\newcommand{\Par}{\catf{Par}}
\newcommand{\opp}{\text{opp}}
\newcommand{\D}{\catf{D}}
\let\Bar\undefined
\newcommand{\Bar}{\textrm{B}\hspace{0.07em}}
\newcommand{\Orb}{\operatorname{Orb}}
\newcommand{\CYCat}{\catf{CY}\catf{cat}^{\operatorname{s}}}
\newcommand{\CYTCat}{\catf{CY}\catf{2cat}}
\newcommand{\FAlg}{\catf{FAlg}}
\newcommand{\mFus}{\catf{mFus}}
\newcommand{\msFus}{\catf{msFus}}
\newcommand{\tric}{\mathcal{T}}%symbol for generic Gray category with duals; alternative: \catf{C}
\newcommand{\ealg}{\mathcal{E}}%construction of 3-category of E_1-algebras 
\newcommand{\EC}{\ealg(\tric)}
\newcommand{\zz}{\mathcal{Z}}%construction of 3-category of E_1-algebras 
\newcommand\orb[1]{{#1}_{\textrm{orb}}}
\newcommand{\one}{\mathbbm{1}}
\newcommand{\be}{\begin{equation}}
\newcommand{\ee}{\end{equation}}
\newcommand{\defi}{\stackrel{\textrm{def}}{=}}
\newcommand{\tev}{\widetilde{\operatorname{ev}}}
\newcommand{\coev}{\operatorname{coev}}
\newcommand{\tcoev}{\widetilde{\operatorname{coev}}}
\newcommand{\dX}{{}^\vee\hspace{-1.8pt}X}
\newcommand{\Xd}{X^\vee}
\def\lra{\longrightarrow}
\def\lmt{\longmapsto}
\newcommand*{\longhookleftarrow}{\ensuremath{\leftarrow\joinrel\relbar\joinrel\rhook}}
\newcommand*{\longhookrightarrow}{\ensuremath{\lhook\joinrel\relbar\joinrel\rightarrow}}
\newcommand*{\twoheadlongrightarrow}{\ensuremath{\relbar\joinrel\twoheadrightarrow}}

\DeclareMathSymbol{\Phiit}{\mathalpha}{letters}{"08}\newcommand{\bflux}{\mathit{\normalfont\Upphi}} \newcommand{\vekphi}{{\vek{\Phiit}}}% kursive girghiche Buchstaben
\DeclareMathSymbol{\Psiit}{\mathalpha}{letters}{"09}\newcommand{\vekpsi}{{\vek{\Psiit}}}
\DeclareMathSymbol{\Sigmait}{\mathalpha}{letters}{"06}
\DeclareMathSymbol{\Xiit}{\mathalpha}{letters}{"04}
\DeclareMathSymbol{\Piit}{\mathalpha}{letters}{"05}\let\Pi\undefined\newcommand{\Pi}{\Piit}
\DeclareMathSymbol{\Gammait}{\mathalpha}{letters}{"00}
\DeclareMathSymbol{\Omegait}{\mathalpha}{letters}{"0A}\let\Omega\undefined\newcommand{\Omega}{\Omegait}
\DeclareMathSymbol{\Upsilonit}{\mathalpha}{letters}{"07}
\DeclareMathSymbol{\Thetait}{\mathalpha}{letters}{"02}
\DeclareMathSymbol{\Lambdait}{\mathalpha}{letters}{"03}\let\Lambda\undefined\newcommand{\Lambda}{\Lambdait}

\let\Phi\undefined\newcommand{\Phi}{\Phiit}
\let\Sigma\undefined\newcommand{\Sigma}{\Sigmait}
\let\Psi\undefined\newcommand{\Psi}{\Psiit}
\let\Gamma\undefined\newcommand{\Gamma}{\Gammait}

\definecolor{Blue}  {rgb} {0.282352,0.239215,0.803921}
\definecolor{Green} {rgb} {0.133333,0.545098,0.133333}
\definecolor{Red}   {rgb} {0.803921,0.000000,0.000000}
\definecolor{Violet}{rgb} {0.580392,0.000000,0.827450}
\definecolor{C1}{RGB} {115,2,127}
\definecolor{C2}{RGB} {230,0,255}
\definecolor{C3}{RGB} {255,0,23}
\definecolor{C4}{RGB} {3,0,112}
\definecolor{C5}{RGB} {170,46,2}

\newcommand{\cbbb}[1]{{\color{Blue}{#1}}}
\newcommand{\cggg}[1]{{\color{Green}{#1}}}
\newcommand{\crrr}[1]{{\color{Red}{#1}}}
\newcounter{jfc}

\newcommand{\LM}[1] {\addtocounter{jfc}{1}\marginpar{\small \cggg{cmt~\thejfc~by~L:}
		{}}{ ~\\ \phantom{~~~} ( {\it \cggg{#1}} ) \\ }}
\newcommand{\NC}[1] {\addtocounter{jfc}{1}\marginpar{\small \cbbb{cmt~\thejfc~by~N:}
		{}}{ ~\\ \phantom{~~~} ( {\it \cbbb{#1}} ) \\ }}
\makeatletter
\g@addto@macro\bfseries{\boldmath}
\makeatother
\begin{document}

	\vspace{5mm}
	
	\begin{center}
		\textbf{\LARGE{Orbifold completion of 3-categories}}\\
		\vspace{1cm}
		{\large Nils Carqueville}$^*$ \ \ and \ \ {\large Lukas Müller}$^\vee$
		\\[0.7cm]
		\normalsize{\texttt{\href{mailto:nils.carqueville@univie.ac.at}{nils.carqueville@univie.ac.at}}} \quad   %
		\normalsize{\texttt{\href{mailto:lmueller@perimeterinstitute.ca}{lmueller@perimeterinstitute.ca}}}
		\\[0.3cm]  %
		\hspace{-1.2cm} {\normalsize\slshape $^\#$Universit\"at Wien, Fakult\"at f\"ur Physik, Boltzmanngasse 5, 1090 Wien, \"{O}sterreich}
		\\[0.1cm]
		\hspace{-1.2cm} {\normalsize\slshape $^\vee$Perimeter Institute, 31 Caroline Street North,
			N2L 2Y5 Waterloo, Canada}
	\end{center}

\vfill 

\begin{abstract}
	We develop a general theory of 3-dimensional ``orbifold completion'', to describe (generalised) orbifolds of topological quantum field theories as well as all their defects. 
	Given a semistrict 3-category~$\tric$ with adjoints for all 1- and 2-morphisms (more precisely, a Gray category with duals), we construct the 3-category $\orb{\tric}$ as a Morita category of certain $E_1$-algebras in~$\tric$ which encode triangulation invariance. 
	We prove that in $\orb{\tric}$ again all 1- and 2-morphisms have adjoints, that it contains~$\tric$ as a full subcategory, and we 
	%arXiv_v2: 
		%argue
		 argue, but do not prove, 
	that it satisfies a universal property which implies $\orb{(\orb{\tric})} \cong \orb{\tric}$. 
	This is a categorification of the work in \cite{cr1210.6363}. 
	
	Orbifold completion by design allows us to lift the orbifold construction from closed TQFT to the much richer world of defect TQFTs. 
	We illustrate this by constructing a universal 3-dimensional state sum model with all defects from first principles, and we explain how recent work on defects between Witt equivalent Reshetikhin--Turaev theories naturally appears as a special case of orbifold completion. 
\end{abstract}

\vfill 

\tableofcontents

\newpage

\section{Introduction and summary}
\label{sec:introduction}

Topological quantum field theories are symmetric monoidal functors on bordism categories. 
The latter may or may not come with labelled stratifications, and the former may or may not take values in higher categories. 
Among the ordinary (1-)categories of labelled stratified bordisms, the oriented case has so far received the most attention, and leads to $n$-dimensional \textsl{defect TQFTs} $\zz\colon \catf{Bord}_{n,n-1}^{\operatorname{def}}(\mathds{D}) \lra \mathcal C$, cf.\ \cite{dkr1107.0495, CRS1}. 
Here~$\mathcal C$ is some prescribed target category, e.\,g.\ of (super) vector spaces, and~$\mathds{D}$ consists of prescribed label sets~$D_j$ for $j$-dimensional strata as well as adjacency rules on how they can meet. 
A sketch of a defect morphism for $n=2$ is 
\be 
\label{eq:DefectBordIntro}
\vcenter{\hbox{
\begin{overpic}[scale=0.5]{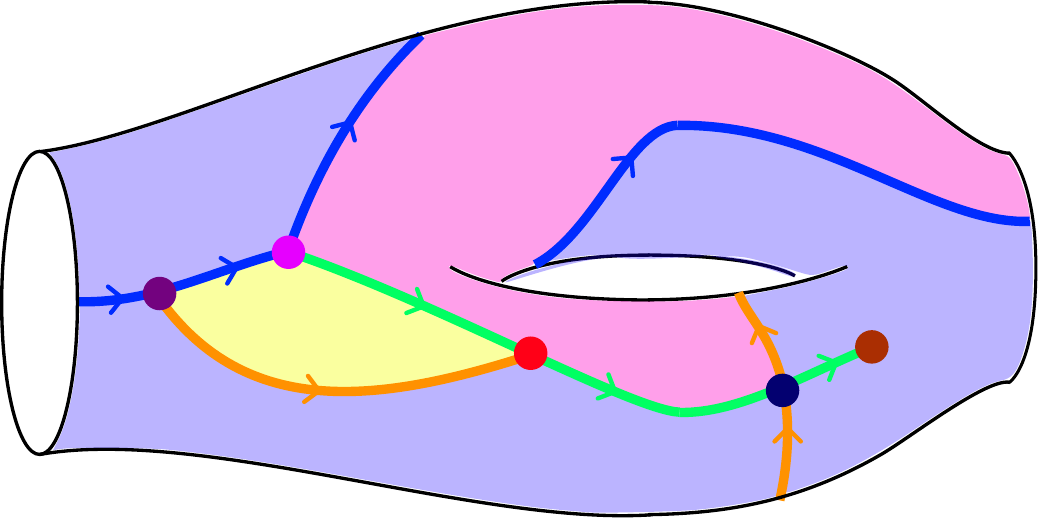}
\put(50,5){$u_1$}
\put(92,22){$u_1$}
\put(15,34){$u_1$} 
\put(32,15){$u_2$} 
\put(45,36){$u_3$}
\put(45,36){$u_3$}
\put(8,22.5){$X_1$}
\put(30,8){$X_2$} 
\put(20,26){$X_3$}
\put(31.5,30){$X_{10}$}
\put(38,23){$X_4$} 
\put(63,12){$X_5$}
\put(76,16.5){$X_6$}
\put(75,38){$X_7$}
\put(77,5){$X_8$}
\put(66,17){$X_9$}
\put(12,18){$\varphi_1$}
\put(27,21.5){$\varphi_2$}
\put(51.5,18){$\varphi_3$}
\put(70.5,8.5){$\varphi_4$}
\put(84,19){$\varphi_5$} 
\put(13.5,23.5){\small $\color{C1}+$}
\put(29.5,25.5){\small$\color{C2}-$}
\put(49.9,11.5){\small$\color{C3}+$} 
\put(77,10){\small$\color{C4}+$}
\put(86,14.5){\small$\color{C5}-$}
\end{overpic}}}
\ee 
where the labels $u_i\in D_2$ correspond to ``bulk theories'' or closed TQFTs, $X_j\in D_1$ describe ``line defects'', and $\varphi_k\in D_0$ are ``point defects''. 
By restricting~$\zz$ to trivially stratified bordisms with only a single label, one recovers the \textsl{closed} TQFTs of \cite{AtiyahTQFT}. 
We refer to \cite{2dDefectTQFTLectureNotes} for a review of 2-dimensional defect theories, and note that recently topological defects have been used to discuss ``non-invertible symmetries'' of full QFTs (see~\cite{snow} for a recent review). 

This notion of generalised symmetries is intimately related to the theory of \textsl{(generalised) orbifolds} of \cite{ffrs0909.5013, cr1210.6363, CRS1}, where the qualifier ``generalised'' covers all non-invertible defects; if instead of orientations one considers framings, one arrives at the ``condensation monads'' of \cite{GaiottoJohnsonFreyd}. 
In arbitrary dimension~$n$, the \textsl{orbifold construction} $\zz \lmt \zz_{\mathcal A}$ produces (new) closed TQFTs~$\zz_{\mathcal A}$ from (old) defect TQFTs~$\zz$, with state sum models and gaugings of (finite) symmetry group actions as particular examples.\footnote{The name ``(generalised) orbifold'' draws from the case of sigma models whose target manifolds~$X$ come with a $G$-action; then the orbifold of the sigma model~$X$ should be a sigma model whose target is the orbifold $X/G$.} 
Besides~$\zz$, the construction needs a collection of special defect labels~$\mathcal A$ as input. 
Then in a nutshell, $\zz_{\mathcal A}$ is evaluated on any bordism by choosing a triangulation, labelling the Poincar\'{e} dual stratification with~$\mathcal A$ (the possibility of which imposes a first set of conditions on these labels), evaluating with~$\zz$, and taking a colimit. 
Independence of the choice of triangulation imposes further defining conditions on the \textsl{orbifold datum}~$\mathcal A$, see \cite[Sect.\,3]{CRS1} for a detailed discussion. 

It is generally expected that an orbifold datum~$\mathcal A$ is an $E_1$-algebra with extra structure in the $n$-category associated to the defect TQFT~$\zz$. 
For $n\in\{2,3\}$, this has been made precise in \cite{cr1210.6363, CMS, CRS1}. 
For $n=2$, $\mathcal A$ amounts to a $\Delta$-separable symmetric Frobenius algebra with defining relations 
\be 
\label{eq:DssFrob}
\tikzzbox{%
	%%%%%%%%%%%%%%%%%%%%%%%%%%%%
	\begin{tikzpicture}[very thick,scale=0.53,color=green!50!black, baseline=0.59cm]
	\draw[-dot-] (3,0) .. controls +(0,1) and +(0,1) .. (2,0);
	\draw[-dot-] (2.5,0.75) .. controls +(0,1) and +(0,1) .. (3.5,0.75);
	\draw (3.5,0.75) -- (3.5,0); 
	\draw (3,1.5) -- (3,2.25); 
	\end{tikzpicture} 
	%%%%%%%%%%%%%%%%%%%%%%%%%%%% 
}%
=
\tikzzbox{%
	%%%%%%%%%%%%%%%%%%%%%%%%%%%%
	\begin{tikzpicture}[very thick,scale=0.53,color=green!50!black, baseline=0.59cm]
	\draw[-dot-] (3,0) .. controls +(0,1) and +(0,1) .. (2,0);
	\draw[-dot-] (2.5,0.75) .. controls +(0,1) and +(0,1) .. (1.5,0.75);
	\draw (1.5,0.75) -- (1.5,0); 
	\draw (2,1.5) -- (2,2.25); 
	\end{tikzpicture} 
	%%%%%%%%%%%%%%%%%%%%%%%%%%%% 
}%
\, , \quad
\tikzzbox{%
	%%%%%%%%%%%%%%%%%%%%%%
	\begin{tikzpicture}[very thick,scale=0.33,color=green!50!black, baseline]
	\draw (-0.5,-0.5) node[Odot] (unit) {}; 
	\fill (0,0.6) circle (5.0pt) node (meet) {};
	\draw (unit) .. controls +(0,0.5) and +(-0.5,-0.5) .. (0,0.6);
	\draw (0,-1.5) -- (0,1.5); 
	\end{tikzpicture} 
	%%%%%%%%%%%%%%%%%%%%%% 
}%
=
\tikzzbox{%
	%%%%%%%%%%%%%%%%%%%%%%%%%%%%
	\begin{tikzpicture}[very thick,scale=0.33,color=green!50!black, baseline]
	\draw (0,-1.5) -- (0,1.5); 
	\end{tikzpicture} 
	%%%%%%%%%%%%%%%%%%%%%%%%%%%% 
}%
=
\tikzzbox{%
	%%%%%%%%%%%%%%%%%%%%%%
	\begin{tikzpicture}[very thick,scale=0.33,color=green!50!black, baseline]
	\draw (0.5,-0.5) node[Odot] (unit) {}; 
	\fill (0,0.6) circle (5.0pt) node (meet) {};
	\draw (unit) .. controls +(0,0.5) and +(0.5,-0.5) .. (0,0.6);
	\draw (0,-1.5) -- (0,1.5); 
	\end{tikzpicture} 
	%%%%%%%%%%%%%%%%%%%%% 
}%
\, , \quad
\tikzzbox{%
	%%%%%%%%%%%%%%%%%%%
	\begin{tikzpicture}[very thick,scale=0.53,color=green!50!black, baseline=-0.59cm, rotate=180]
	\draw[-dot-] (3,0) .. controls +(0,1) and +(0,1) .. (2,0);
	\draw[-dot-] (2.5,0.75) .. controls +(0,1) and +(0,1) .. (1.5,0.75);
	\draw (1.5,0.75) -- (1.5,0); 
	\draw (2,1.5) -- (2,2.25); 
	\end{tikzpicture} 
	%%%%%%%%%%%%%%%%%%% 
}%
=
\tikzzbox{%
	%%%%%%%%%%%%%%%%%%%
	\begin{tikzpicture}[very thick,scale=0.53,color=green!50!black, baseline=-0.59cm, rotate=180]
	\draw[-dot-] (3,0) .. controls +(0,1) and +(0,1) .. (2,0);
	\draw[-dot-] (2.5,0.75) .. controls +(0,1) and +(0,1) .. (3.5,0.75);
	\draw (3.5,0.75) -- (3.5,0); 
	\draw (3,1.5) -- (3,2.25); 
	\end{tikzpicture} 
	%%%%%%%%%%%%%%%%%%% 
}%
\, , \quad
\tikzzbox{%
	%%%%%%%%%%%%%%%%%%%
	\begin{tikzpicture}[very thick,scale=0.33,color=green!50!black, baseline=0, rotate=180]
	\draw (0.5,-0.5) node[Odot] (unit) {}; 
	\fill (0,0.6) circle (5.0pt) node (meet) {};
	\draw (unit) .. controls +(0,0.5) and +(0.5,-0.5) .. (0,0.6);
	\draw (0,-1.5) -- (0,1.5); 
	\end{tikzpicture} 
	%%%%%%%%%%%%%%%%%%% 
}%
=
\tikzzbox{%
	%%%%%%%%%%%%%%%%%%%
	\begin{tikzpicture}[very thick,scale=0.33,color=green!50!black, baseline=0, rotate=180]
	\draw (0,-1.5) -- (0,1.5); 
	\end{tikzpicture} 
	%%%%%%%%%%%%%%%%%%% 
}%
=
\tikzzbox{%
	%%%%%%%%%%%%%%%%%%%
	\begin{tikzpicture}[very thick,scale=0.33,color=green!50!black, baseline=0cm, rotate=180]
	\draw (-0.5,-0.5) node[Odot] (unit) {}; 
	\fill (0,0.6) circle (5.0pt) node (meet) {};
	\draw (unit) .. controls +(0,0.5) and +(-0.5,-0.5) .. (0,0.6);
	\draw (0,-1.5) -- (0,1.5); 
	\end{tikzpicture} 
	%%%%%%%%%%%%%%%%%%% 
}%
\, , \quad 
\tikzzbox{%
	%%%%%%%%%%%%%%%%%%%%%%
	\begin{tikzpicture}[very thick,scale=0.33,color=green!50!black, baseline=0cm]
	\draw[-dot-] (0,0) .. controls +(0,-1) and +(0,-1) .. (-1,0);
	\draw[-dot-] (1,0) .. controls +(0,1) and +(0,1) .. (0,0);
	\draw (-1,0) -- (-1,1.5); 
	\draw (1,0) -- (1,-1.5); 
	\draw (0.5,0.8) -- (0.5,1.5); 
	\draw (-0.5,-0.8) -- (-0.5,-1.5); 
	\end{tikzpicture}
	%%%%%%%%%%%%%%%%%%%%%%
}%
=
\tikzzbox{%
	%%%%%%%%%%%%%%%%%%%%%%
	\begin{tikzpicture}[very thick,scale=0.33,color=green!50!black, baseline=0cm]
	\draw[-dot-] (0,1.5) .. controls +(0,-1) and +(0,-1) .. (1,1.5);
	\draw[-dot-] (0,-1.5) .. controls +(0,1) and +(0,1) .. (1,-1.5);
	\draw (0.5,-0.8) -- (0.5,0.8); 
	\end{tikzpicture}
	%%%%%%%%%%%%%%%%%%%%%%
}%
=
\tikzzbox{%
	%%%%%%%%%%%%%%%%%%%%%%
	\begin{tikzpicture}[very thick,scale=0.33,color=green!50!black, baseline=0cm]
	\draw[-dot-] (0,0) .. controls +(0,1) and +(0,1) .. (-1,0);
	\draw[-dot-] (1,0) .. controls +(0,-1) and +(0,-1) .. (0,0);
	\draw (-1,0) -- (-1,-1.5); 
	\draw (1,0) -- (1,1.5); 
	\draw (0.5,-0.8) -- (0.5,-1.5); 
	\draw (-0.5,0.8) -- (-0.5,1.5); 
	\end{tikzpicture}
	%%%%%%%%%%%%%%%%%%%%%%
}%
\, , \quad 
%%%%%%%%%%%%%%%%%%%%%%
\begin{tikzpicture}[very thick,scale=0.33,color=green!50!black, baseline=-0.1cm]
\draw[-dot-] (0,0) .. controls +(0,-1) and +(0,-1) .. (1,0);
\draw[-dot-] (0,0) .. controls +(0,1) and +(0,1) .. (1,0);
\draw (0.5,-0.8) -- (0.5,-1.5); 
\draw (0.5,0.8) -- (0.5,1.5); 
\end{tikzpicture}
%%%%%%%%%%%%%%%%%%%%%%
\, = \, 
%%%%%%%%%%%%%%%%%%%%%%
\begin{tikzpicture}[very thick,scale=0.33,color=green!50!black, baseline=-0.1cm]
\draw (0.5,-1.5) -- (0.5,1.5); 
\end{tikzpicture}
%%%%%%%%%%%%%%%%%%%%%%
\, , \quad 
%%%%%%%%%%%%%%%%%%%%%%
\begin{tikzpicture}[very thick,scale=0.33,color=green!50!black, baseline=0cm]
\draw[-dot-] (0,0) .. controls +(0,1) and +(0,1) .. (-1,0);
\draw[directedgreen, color=green!50!black] (1,0) .. controls +(0,-1) and +(0,-1) .. (0,0);
\draw (-1,0) -- (-1,-1.5); 
\draw (1,0) -- (1,1.5); 
\draw (-0.5,1.2) node[Odot] (end) {}; 
\draw (-0.5,0.8) -- (end); 
\end{tikzpicture}
%%%%%%%%%%%%%%%%%%%%%%
= 
%%%%%%%%%%%%%%%%%%%%%%
\begin{tikzpicture}[very thick,scale=0.33,color=green!50!black, baseline=0cm]
\draw[redirectedgreen, color=green!50!black] (0,0) .. controls +(0,-1) and +(0,-1) .. (-1,0);
\draw[-dot-] (1,0) .. controls +(0,1) and +(0,1) .. (0,0);
\draw (-1,0) -- (-1,1.5); 
\draw (1,0) -- (1,-1.5); 
\draw (0.5,1.2) node[Odot] (end) {}; 
\draw (0.5,0.8) -- (end); 
\end{tikzpicture}
%%%%%%%%%%%%%%%%%%%%%%
\ee 
in the pivotal 2-category~$\mathcal{B}_\zz$ associated to~$\zz$ (cf.\ \cite{dkr1107.0495}). 
For $n=3$, an orbifold datum~$\mathcal A$ is a generalisation of spherical fusion categories internal to the 3-category with coherent adjoints~$\tric_{\zz}$ (more precisely, a ``Gray category with duals''\footnote{The name is slightly unfortunate, because it does not only involve the existence of duals (or adjoints) of morphisms, but also coherent identifications between left and right adjoints.} as defined in \cite{BMS}, cf.\ Section~\ref{subsec:3categories}) constructed in \cite{CMS}; we give the precise definition of 3-dimensional orbifold data in Definition~\ref{def:OrbifoldDatum}. 

It is also expected that by considering~$\zz_{\mathcal A}$ for all orbifold data~$\mathcal A$ at once, the orbifold construction should lift to a map $\zz \lmt \zz_{\textrm{orb}}$ from defect TQFTs to \textsl{defect} TQFTs 
%arXiv_v2: 
	%(as opposed to closed TQFTs~$\zz_{\mathcal A}$ as above). 
%arXiv_v3: 
	 (as opposed to merely closed TQFTs~$\zz_{\mathcal A}$ as above whose domain are plain bordisms, while the domain of $\zz_{\textrm{orb}}$ are stratified and decorated bordisms). 
This amounts to identifying all defects between the closed TQFTs~$\zz_{\mathcal A}$, all defects between them, etc. 
Algebraically, this leads to a higher Morita category of the algebras~$\mathcal A$ with compatibility with adjoints (coming from orientation-reversal of strata) built in. 
For $n=2$, this was made precise in \cite{cr1210.6363} by associating to any pivotal 2-category~$\mathcal B$ its \textsl{orbifold completion}, i.\,e.\ the 2-category $\mathcal B_{\textrm{orb}}$ of $\Delta$-separable symmetric Frobenius algebras~$\mathcal A$ in~$\mathcal B$, with 1- and 2-morphisms given by bimodules and bimodule maps. 
Then one shows that $\mathcal B_{\textrm{orb}}$ has a natural pivotal structure, and that it is complete in the sense that $(\mathcal B_{\textrm{orb}})_{\textrm{orb}} \cong \orb{\mathcal B}$. 
This can be thought of as an ``oriented'' variant of 
%arXiv_v2: 
	%2-idempotent completion 
	 2-condensation completion
where the symmetry condition on $\mathcal A \in \mathcal B_{\textrm{orb}}$ is the property that distinguishes it from the 
%arXiv_v2: 
	%2-idempotent completion 
	 2-condensation completion
construction of \cite{DouglasReutter2018, GaiottoJohnsonFreyd} in the context of \textsl{framed} TQFTs; the symmetry condition can also be thought of as a potential obstruction for ``gauging'' the generalised symmetry~$\mathcal A$. 
Via application to TQFTs, the results of \cite{cr1210.6363} lead to the discovery of new relations between closed TQFTs and singularity theory 
%arXiv_v3: 
	%\cite{CRCR, OEReck}
	 \cite{CRCR, NRC1, OEReck}
%arXiv_v2: 
	(which are not known to have analogues in the ``framed'' setting of condensation monads), 
and a general description of ``non-abelian quantum symmetry'' \cite{BCP2}. 
%arXiv_v2: 
	Quite generally in any dimension $n\geqslant 2$, ``oriented'' orbifold completion is expected to be richer than ``framed'' condensation completion, as illustrated e.\,g.\ in the context of fully extended TQFTs by the fact that $\SO(n)$-homotopy fixed point structures are considerably more intricate than full dualisability by itself. 

\medskip 

In the present paper, we develop a 3-categorical version of orbifold completion, and we study some of its properties and applications. 
More precisely, given a Gray category with duals~$\tric$, we construct a 3-category $\orb{\tric}$ whose objects are orbifold data~$\mathcal A$ in~$\tric$, i.\,e.\ $E_1$-algebras with additional properties listed in Figure~\ref{fig:OrbifoldDatumAxioms} (using the graphical calculus reviewed in Section~\ref{sec:Background}). 
These properties categorify the relations in~\eqref{eq:DssFrob} and were shown in \cite{CRS1} to encode invariance under 3-dimensional oriented Pachner moves. 
In particular, $\mathcal A$ involves a 1-endomorphism $A\colon a\lra a$, a multiplication 2-morphism~$\mu$, and an associator~$\alpha$: 
\be 
%%%%%%%%%%%%%%%%%%%%%% 
\tikzzbox{
}%%popende
%%%%%%%%%%%%%%%%%%%%%% 
\, .
\ee 
Among the defining relations in Figure~\ref{fig:OrbifoldDatumAxioms} is the condition that quantum dimensions of~$\mu$ are identities, as well as pentagon and invertibility axioms like  
\be 
%%%%%%%%%%%%%%%%%%%%%% 
\tikzzbox{%
}%%popende
%%%%%%%%%%%%%%%%%%%%%% 
\, . 
\ee 

1-morphisms in $\orb{\tric}$ are bimodules~$\mathcal M$ over the underlying $E_1$-algebras (cf.\ Definition~\ref{def:1MorphismsEAlg}) that satisfy ``coloured'' versions of the conditions for orbifold data~$\mathcal A$ such as 
\be 
%%%%%%%%%%%%%%%%%%%%%% 
\tikzzbox{%
}%%popende
%%%%%%%%%%%%%%%%%%%%%% 
\, . \ee 
We give all details in Definition~\ref{def:Corb}, which also explicitly presents 2- and 3-morphisms of $\orb{\tric}$ as those in the Morita 3-category of $E_1$-algebras in~$\tric$ (reviewed in Section~\ref{sub:EofC}) which are compatible with all the orbifold structure. 
These conditions had previously appeared in \cite{MuleRunk, CMRSS1}, here we explain how they naturally arise from the perspective of higher representation theory. 

\medskip 

Our first main result is that $\orb{\tric}$ has the desired duality properties: 

\begin{theorem}
	\label{thm:MainTheorem}
	Let~$\tric$ be a Gray category with duals (satisfying a mild condition on the existence of certain (co)limits, cf.\ Section~\ref{subsec:3catOrbData}). 
	Then $\orb{\tric}$ is a 3-category that admits adjoints for all 1- and 2-morphisms. 
\end{theorem}

The main ingredients of the proof are a description of the composition of 1-morphisms in terms of the splitting of 
%arXiv_v2:
	%2-idempotents
	 a condensation 2-monad 
(Definition~\ref{def:2IdempotentSplit} and Proposition~\ref{prop:RelProd}) as well as explicit constructions of all adjoints (Lemma~\ref{lem:Adj2Mor} and Proposition~\ref{prop:Adj1Mor}). 
The verification of coherence axioms is considerably facilitated by the graphical calculus of \cite{BMS}. 

To fully justify the term orbifold \textsl{completion} also in three dimensions, one expects an equivalence of 3-categories $\orb{(\orb{\tric})} \cong \orb{\tric}$. 
We leave this for future work, but we do propose (in Section~\ref{subsubsec:UniversalProperty}) a universal property for orbifold completion in arbitrary dimension, and we prove that the 2-dimensional orbifold completion $\mathcal B \lmt \orb{\mathcal B}$ of \cite{cr1210.6363} does satisfy this property (cf.\ Proposition~\ref{prop:BorbSatisfiesUniversalProperty}). 
The completion property $\orb{(\orb{\mathcal B})} \cong \orb{\mathcal B}$ is a direct corollary of this more conceptual description, and we expect similar relations in higher dimensions. 

\medskip 

The applications of Theorem~\ref{thm:MainTheorem} that we cover in the present paper (in Sections~\ref{sec:ExamppleCompletion} and~\ref{sec:TQFTs}) concern 3-dimensional state sum models with defects, and more generally Reshetikhin--Turaev theory. 
To summarise these results, we shall first briefly discuss that our original motivation for 3-dimensional orbifold completion was to construct 4-dimensional state sum models from first principles, and re-discover the 4-manifold invariants of Douglas and Reutter \cite{DouglasReutter2018} as part of a closed orbifold TQFT. 
This application will appear in joint work with Vincentas Mulevi\v{c}ius in \cite{LukasNilsVincentas}. 

A key idea, in part inspired by \cite{cr1210.6363, CRS3, GaiottoJohnsonFreyd}, is that the $n$-dimensional state sum model with defects~$\zz^{\textrm{ss}}_n$ for some rigid symmetric monoidal target category~$\mathcal C$ can be constructed by iterative orbifold completions. 
Indeed, $\zz^{\textrm{ss}}_1$ by definition is just evaluation of string diagrams in $\mathcal C_1 := \mathcal C$. 
For the inductive step, assume that we have constructed a symmetric monoidal $(n-1)$-category~$\mathcal C_{n-1}$. 
Hence we can consider the delooping $n$-category $\Bar\mathcal C_{n-1}$ and define the sets $D^{\textrm{triv}_n}_j$ of ``trivial'' $j$-dimensional defect labels as the sets of  $(n-j)$-morphisms in $\Bar\mathcal C_{n-1}$. 
Using an iterative construction, one obtains the ``trivial'' defect TQFT  $\zz^{\textrm{triv}}_n\colon \catf{Bord}_n^{\operatorname{def}}(\mathds{D}^{\textrm{triv}_n}) \lra \mathcal C$. 
The crucial step is then to define $\mathcal C_n := E(\Bar\mathcal C_{n-1})_{\textrm{orb}}$, where $E(-)$ is the ``Euler completion'' of \cite{CRS1}.\footnote{For $n=2$, the difference between $(-)_{\textrm{orb}}$ and $E(-)_{\textrm{orb}}$ is that the latter has \textsl{all} {separable} symmetric Frobenius algebras in~$\mathcal C$ as objects, not only {$\Delta$-separable} ones. We review the 3-dimensional Euler completion of \cite[Sect.\,4.2]{CRS1} in Section~\ref{subsubsec:3dStateSumModels}.} 

Taking the sets $D^{\textrm{ss}_n}_j$ of $j$-dimensional defect labels to consist of the $(n-j)$-morphisms of~$\mathcal C_n$, we can define the \textsl{$n$-dimensional defect state sum model} 
\be 
\zz^{\textrm{ss}}_n\colon \catf{Bord}_n^{\operatorname{def}}(\mathds{D}^{\textrm{ss}_n}) \lra \mathcal C
\ee 
by interpreting $\mathcal C_n$-labelled stratified bordisms as mere string diagrams in~$\mathcal C$, since after all~$\mathcal C_n$ is a Morita $n$-category of the delooping for a Morita $(n-1)$-category of a delooping of \dots\ $\mathcal C$.\footnote{There are several variations of this construction. For example, instead of starting with a symmetric monoidal 1-category~$\mathcal C$, we can start with a symmetric monoidal $k$-category for some $k<n$, and then apply $E(\Bar(-)_{\textrm{orb}})^{n-k}$.}

To construct 4-dimensional state sum models along the lines just sketched, the last and main steps in the iteration are to construct the orbifold completion of the specific Morita 3-category~$\tric'$ of orbifold data in the delooping of the 2-category of 2-dimensional state sum models, and to identify orbifold data in the delooping of $E(\tric'_{\textrm{orb}})$. 
In Section~\ref{subsec:StateSumModels} we study the 3-category~$\tric'$ in detail, and in \cite{LukasNilsVincentas} we will show how the spherical fusion 2-categories~$\mathcal S$ of \cite{DouglasReutter2018} naturally give rise to 4-dimensional orbifold data in $\Bar E(\tric'_{\textrm{orb}})$ -- in complete analogy to how spherical fusion 1-categories give rise to the 3-dimensional orbifold data underlying Turaev--Viro--Barrett--Westbury theory, cf.\ \cite[Sect.\,4]{CRS3}. 
Building on the work of Douglas and Reutter it is then straightforward to check that restricting our defect state sum model~$\zz^{\textrm{ss}}_4$ to trivially stratified and~$\mathcal S$-labelled bordisms, one obtains a closed TQFT whose 4-manifold invariants are those of \cite{DouglasReutter2018}. 

\medskip 

We conclude this section with a brief summary of the applications treated in this paper. 
Using the above notation, Section~\ref{subsec:StateSumModels} provides a detailed study of the 3-category $ E(\Bar\ssFrob)_{\textrm{orb}}$, i.\,e.\ the (Euler completion of the) orbifold completion of the delooping of the 2-category $\ssFrob = E(\Bar\vs)_{\textrm{orb}}$ of separable symmetric Frobenius algebras in $\vs$. 
%arXiv_v2: 
	(Here and below $\vs$ denotes the category of vector spaces over some fixed algebraically closed field of characteristic~0). 
In particular, we construct (cf.\ Proposition~\ref{prop:EWor}) a pivotal equivalence between $\ssFrob$ and a 2-category $\CYCat$ of semisimple Calabi--Yau categories, and we prove: 

\smallskip 
\noindent
\textbf{Theorem \ref{thm:OrbOfTriv}. }
The 3-category $\catf{sFus}$ of spherical fusion categories, bimodule categories with trace, bimodule functors, 
and bimodule natural transformations defined in~\cite{Bimodtrace} is a subcategory of $\orb{E(\Bar\CYCat)}$, or equivalently of $\orb{E(\Bar\ssFrob)}$. 
\smallskip 

A direct consequence of this is a new proof of the fact, first proved in \cite{DuaTen}, that 
%arXiv_v3: 
	%spherical fusion categories 
	 fusion categories admitting a spherical structure 
are 3-dualisable in the 3-category $\Alg_1(2\vs)$. 
%arXiv_v3: 
	Note that it is widely conjectured that every fusion category admits a spherical structure. 

In Section~\ref{sec:TQFTs} we briefly review 3-dimensional defect TQFTs and then apply our algebraic theory of orbifold completion to them to show: 

\smallskip 
\noindent
\textbf{Theorem \ref{thm:OrbifoldDefectTQFT}. }
Let $\cat{Z}\colon \DBord (\mathds{D})\to \vs$ be a defect TQFT. 
The orbifold construction gives a defect TQFT $\cat{Z}_{\operatorname{orb}}\colon \DBord (\mathds{D}_{\operatorname{orb}})\to \vs$. 
\smallskip 

Applying this to the special case of the trivial defect TQFT, we obtain the 3-dimensional defect state sum model~$\zz^{\textrm{ss}}_3$ along the lines outlined for arbitrary dimension above. 
It restricts to Turaev--Viro--Barrett--Westbury theory on trivially stratified bordisms, and we expect~$\zz^{\textrm{ss}}_3$ to also generalise the work of \cite{Meusburger3dDefectStateSumModels} to arbitrary defect bordisms. 

To make the connection to general TQFTs of Reshetikhin--Turaev type we fix any modular fusion category~$\mathcal M$ and consider the 2-category $\ssFrob(\mathcal M)$. 
In \cite{KMRS} a connection was made between the general analysis \cite{fsv1203.4568} of semisimple defect theories in three dimensions and the general orbifold theory of \cite{CRS1, CRS3}. 
A key step in \cite{KMRS} was the clever ad hoc introduction of the notion of \textsl{Frobenius algebras~$F$ over a pair $(A,B)$ of commutative $\Delta$-separable Frobenius algebras} to describe surface defects. 
In Section~\ref{subsec:DomainWallsForRT} we show that this notion as well as the related higher maps naturally occur as morphisms in the orbifold completion of $\Bar\ssFrob(\mathcal M)$: 

\smallskip 
\noindent
\textbf{Theorem \ref{thm:OrbOfssFrobM}. }
Let~$\mathcal M$ be a modular fusion category. 
The 3-category in which
\begin{itemize}
	\item 
	objects are commutative $\Delta$-separable Frobenius algebras in $\cat{M}$, 
	\item 
	1-morphisms from $A$ to $B$ are $\Delta$-separable symmetric Frobenius algebras $F$ over $(A,B)$, 
	\item 
	2-morphisms from $F$ to $G$ are $G$-$F$-bimodules $M$ over $(A,B)$, and
	\item 
	3-morphisms are bimodule maps 
\end{itemize} 
is a subcategory of $\orb{(\Bar\Delta\ssFrob(\mathcal{M}))}$. 
\smallskip 

Combining this with Theorem~\ref{thm:OrbifoldDefectTQFT} above we recover the defect TQFT of \cite{KMRS} from the universal orbifold construction. 

\medskip 

\noindent
\textbf{Acknowledgements. }
We thank 
	Theo Johnson-Freyd, 
	Christopher Lieberum, 
	Vincentas Mulevi\v{c}ius, 
	David Reutter, 
	Claudia Scheimbauer,
	and Lukas Woike 
for helpful discussions and comments. 
N.\,C.\ is supported by the DFG Heisenberg Programme. 
L.\,M.\ gratefully acknowledges support by the Max Planck Institute for Mathematics in Bonn,
where part of this work was carried out,
and the Simons Collaboration on Global Categorical Symmetries. Research at Perimeter Institute is supported in part by the Government of Canada through the Department of Innovation, Science and Economic Development and by the Province of Ontario through the Ministry of Colleges and Universities. The Perimeter Institute is in the Haldimand Tract, land promised to the Six Nations.

\section{Background}
\label{sec:Background}

In this section we collect some background on 2- and 3-categories with adjoints and their graphical calculus. 
In Section~\ref{subsec:2categories} we review pivotal 2-categories~$\mathcal B$ as well as their orbifold completions, i.\,e.\ the representation theory of $\Delta$-separable symmetric Frobenius algebras internal to~$\mathcal B$. 
In Section~\ref{subsec:3categories} we discuss a semistrict version of 3-categories with compatible adjoints, namely the ``Gray categories with duals'' introduced and studied in \cite{BMS}, and we recall some of the coherence results on these structures.

\subsection{Pivotal 2-categories and 2-dimensional orbifold completion}
\label{subsec:2categories}

For the basic definitions and theory of 2-categories (which we do not assume to be strict) we refer to \cite{JohnsonYauBook}. 
We use ``$\cdot$'' or mere concatenation to denote vertical composition, and ``$\circ$'' for horizontal composition. 
Our conventions for the graphical calculus are to read diagrams from bottom to top and from right to left. 
Hence for objects $u,v,w$ as well as suitably composable 1-morphisms $X,Y,Z,P,Q$ and 2-morphisms $\varphi,\psi,\zeta$ in some 2-category, we have 
\be 
\tikzzbox{%
	%%%%%%%%%%%%%%%%%%%%%% 
	\begin{tikzpicture}[very thick,scale=1.0,color=blue!50!black, baseline=1.2cm]
	\coordinate (d1) at (-3,0);
	\coordinate (d2) at (+3,0);
	\coordinate (u1) at (-3,3);
	\coordinate (u2) at (+3,3);
	\coordinate (s) at ($(-0.5,-0.2)$);
	\coordinate (b1) at ($(d1)+(s)$);
	\coordinate (b2) at ($(d2)+(s)$);
	\coordinate (t1) at ($(u1)+(s)$);
	\coordinate (t2) at ($(u2)+(s)$);
	%
	% colouring front: 
	\fill [orange!40!white, opacity=0.7] (b1) -- (b2) -- (t2) -- (t1); 
	%
	% strings and labels front: 
	\draw[very thick] ($(d1)+(2.25,0)+(s)$) -- ($(d1)+(2.25,3)+(s)$);
	\fill ($(d1)+(2.25,1.5)+(s)$) circle (2.5pt) node[left] {{\small $\zeta$}};
	\fill ($(d1)+(2.3,2.5)+(s)$) circle (0pt) node[left] {{\small $Q$}};
	\fill ($(d1)+(2.3,0.5)+(s)$) circle (0pt) node[left] {{\small $P$}};
	\draw[very thick] ($(d1)+(3.75,0)+(s)$) -- ($(d1)+(3.75,3)+(s)$);
	\fill ($(d1)+(3.7,1.5)+(s)$) circle (0pt) node[right] {{\small $Y$}};
	\fill ($(d1)+(3.75,1)+(s)$) circle (2.5pt) node[right] {{\small $\varphi$}};
	\fill ($(d1)+(3.75,2)+(s)$) circle (2.5pt) node[right] {{\small $\psi$}};
	\fill ($(d1)+(3.7,2.5)+(s)$) circle (0pt) node[right] {{\small $Z$}};
	\fill ($(d1)+(3.7,0.5)+(s)$) circle (0pt) node[right] {{\small $X$}};
	\fill[red!80!black] ($(d1)+(5.75,0.5)+(s)$) circle (0pt) node {{\scriptsize $u$}};
	\fill[red!80!black] ($(d1)+(3,0.5)+(s)$) circle (0pt) node {{\scriptsize $v$}};
	\fill[red!80!black] ($(d1)+(0.25,0.5)+(s)$) circle (0pt) node {{\scriptsize $w$}};
	%
	% red and black lines front: 
	\draw[thin, black] (t1) -- (t2); 
	\draw[thin, black] (b1) -- (b2); 
	\draw[thin, black] (b1) -- (t1); 
	\draw[thin, black] (b2) -- (t2);
	\end{tikzpicture}
	%%%%%%%%%%%%%%%%%%%%%% 
} 
= \zeta \circ (\psi \cdot \varphi) \colon P\circ X \lra Q \circ Z \, 
\ee 

A 1-morphism $X\colon u\lra v$ has \textsl{left} and \textsl{right adjoints} if there are ${}^\vee\!X \colon v\lra u$ and $X^\vee \colon v\lra u$ together with 2-morphisms 
\be 
\label{eq:AdjMaps}
\ev_X  = 
\tikzzbox{%
	%%%%%%%%%%%%%%%%%%%%%% 
	\begin{tikzpicture}[very thick,scale=1.0,color=blue!50!black, baseline=0.5cm]
	\coordinate (X) at (0.5,0);
	\coordinate (Xd) at (-0.5,0);
	\coordinate (d1) at (-1,0);
	\coordinate (d2) at (+1,0);
	\coordinate (u1) at (-1,1.25);
	\coordinate (u2) at (+1,1.25);
	%
	% colouring: 
	\fill [orange!40!white, opacity=0.7] (d1) -- (d2) -- (u2) -- (u1); 
	\draw[thin] (d1) -- (d2) -- (u2) -- (u1) -- (d1); 
	%
	% strings: 
	\draw[directed] (X) .. controls +(0,1) and +(0,1) .. (Xd);
	%
	% labels: 
	\fill (X) circle (0pt) node[below] {{\small $X\vphantom{X^\vee}$}};
	\fill (Xd) circle (0pt) node[below] {{\small ${}^\vee\!X$}};
	\fill[red!80!black] (0,0.3) circle (0pt) node {{\scriptsize $v$}};
	\fill[red!80!black] (0.8,1) circle (0pt) node {{\scriptsize $u$}};
	\end{tikzpicture}
	%%%%%%%%%%%%%%%%%%%%%% 
}
\, , \quad 
\coev_X = 
\tikzzbox{%
	%%%%%%%%%%%%%%%%%%%%%% 
	\begin{tikzpicture}[very thick,scale=1.0,color=blue!50!black, baseline=0.5cm]
	\coordinate (X) at (-0.5,1.25);
	\coordinate (Xd) at (0.5,1.25);
	\coordinate (d1) at (-1,0);
	\coordinate (d2) at (+1,0);
	\coordinate (u1) at (-1,1.25);
	\coordinate (u2) at (+1,1.25);
	%
	% colouring: 
	\fill [orange!40!white, opacity=0.7] (d1) -- (d2) -- (u2) -- (u1); 
	\draw[thin] (d1) -- (d2) -- (u2) -- (u1) -- (d1);
	%
	% strings: 
	\draw[redirected] (X) .. controls +(0,-1) and +(0,-1) .. (Xd);
	%
	% labels: 
	\fill (X) circle (0pt) node[above] {{\small $X\vphantom{X^\vee}$}};
	\fill (Xd) circle (0pt) node[above] {{\small ${}^\vee\!X$}};
	\fill[red!80!black] (0,0.9) circle (0pt) node {{\scriptsize $u$}};
	\fill[red!80!black] (0.8,0.2) circle (0pt) node {{\scriptsize $v$}};
	\end{tikzpicture}
	%%%%%%%%%%%%%%%%%%%%%% 
}
\quad\textrm{and}\quad 
\tev_X  = 
\tikzzbox{%
	%%%%%%%%%%%%%%%%%%%%%% 
	\begin{tikzpicture}[very thick,scale=1.0,color=blue!50!black, baseline=0.5cm]
	\coordinate (X) at (0.5,0);
	\coordinate (Xd) at (-0.5,0);
	\coordinate (d1) at (-1,0);
	\coordinate (d2) at (+1,0);
	\coordinate (u1) at (-1,1.25);
	\coordinate (u2) at (+1,1.25);
	%
	% colouring: 
	\fill [orange!40!white, opacity=0.7] (d1) -- (d2) -- (u2) -- (u1); 
	\draw[thin] (d1) -- (d2) -- (u2) -- (u1) -- (d1);
	%
	% strings: 
	\draw[redirected] (X) .. controls +(0,1) and +(0,1) .. (Xd);
	%
	% labels: 
	\fill (Xd) circle (0pt) node[below] {{\small $X\vphantom{X^\vee}$}};
	\fill (X) circle (0pt) node[below] {{\small $X^\vee$}};
	\fill[red!80!black] (0,0.3) circle (0pt) node {{\scriptsize $u$}};
	\fill[red!80!black] (0.8,1) circle (0pt) node {{\scriptsize $v$}};
	\end{tikzpicture}
	%%%%%%%%%%%%%%%%%%%%%% 
}
\, , \quad 
\tcoev_X = 
\tikzzbox{%
	%%%%%%%%%%%%%%%%%%%%%% 
	\begin{tikzpicture}[very thick,scale=1.0,color=blue!50!black, baseline=0.5cm]
	\coordinate (X) at (-0.5,1.25);
	\coordinate (Xd) at (0.5,1.25);
	\coordinate (d1) at (-1,0);
	\coordinate (d2) at (+1,0);
	\coordinate (u1) at (-1,1.25);
	\coordinate (u2) at (+1,1.25);
	%
	% colouring: 
	\fill [orange!40!white, opacity=0.7] (d1) -- (d2) -- (u2) -- (u1); 
	\draw[thin] (d1) -- (d2) -- (u2) -- (u1) -- (d1);
	%
	% strings: 
	\draw[directed] (X) .. controls +(0,-1) and +(0,-1) .. (Xd);
	%
	% labels: 
	\fill (Xd) circle (0pt) node[above] {{\small $X\vphantom{X^\vee}$}};
	\fill (X) circle (0pt) node[above] {{\small $X^\vee$}};
	\fill[red!80!black] (0,0.9) circle (0pt) node {{\scriptsize $v$}};
	\fill[red!80!black] (0.8,0.2) circle (0pt) node {{\scriptsize $u$}};
	\end{tikzpicture}
	%%%%%%%%%%%%%%%%%%%%%% 
} 
\, , 
\ee 
respectively, such that the \textsl{Zorro moves} 
%arXiv_v3: 
	(also known as \textsl{snake} or \textsl{zig-zag identities}) 
\be 
\label{eq:ZorroMoves}
\tikzzbox{%
	%%%%%%%%%%%%%%%%%%%%%% 
	\begin{tikzpicture}[very thick,scale=1.0,color=blue!50!black, baseline=0cm]
	\coordinate (A) at (-1,1.25);
	\coordinate (A2) at (1,-1.25);
	\coordinate (d1) at (-1.5,-1.25);
	\coordinate (d2) at (+1.5,-1.25);
	\coordinate (u1) at (-1.5,1.25);
	\coordinate (u2) at (+1.5,1.25);
	%
	% colouring: 
	\fill [orange!40!white, opacity=0.7] (d1) -- (d2) -- (u2) -- (u1); 
	\draw[thin] (d1) -- (d2) -- (u2) -- (u1) -- (d1);
	%
	% strings: 
	\draw[directed] (0,0) .. controls +(0,-1) and +(0,-1) .. (-1,0);
	\draw[directed] (1,0) .. controls +(0,1) and +(0,1) .. (0,0);
	\draw (-1,0) -- (A); 
	\draw (1,0) -- (A2); 
	%
	% labels: 
	\fill ($(A)+(0.1,0)$) circle (0pt) node[below left] {{\small $X$}};
	\fill ($(A2)+(-0.1,0)$) circle (0pt) node[above right] {{\small $X$}};
	\fill[red!80!black] (-0.5,0) circle (0pt) node {{\scriptsize $u$}};
	\fill[red!80!black] (+0.5,0) circle (0pt) node {{\scriptsize $v$}};
	\end{tikzpicture}
	%%%%%%%%%%%%%%%%%%%%%% 
}
= 
\tikzzbox{%
	%%%%%%%%%%%%%%%%%%%%%% 
	\begin{tikzpicture}[very thick,scale=1.0,color=blue!50!black, baseline=0cm]
	\coordinate (A) at (0,1.25);
	\coordinate (A2) at (0,-1.25);
	\coordinate (d1) at (-1,-1.25);
	\coordinate (d2) at (+1,-1.25);
	\coordinate (u1) at (-1,1.25);
	\coordinate (u2) at (+1,1.25);
	%
	% colouring: 
	\fill [orange!40!white, opacity=0.7] (d1) -- (d2) -- (u2) -- (u1); 
	\draw[thin] (d1) -- (d2) -- (u2) -- (u1) -- (d1);
	%
	% strings: 
	\draw (A) -- (A2); 
	%
	% labels: 
	\fill ($(A2)+(-0.1,0)$) circle (0pt) node[above right] {{\small $X$}};
	\fill[red!80!black] (-0.5,0) circle (0pt) node {{\scriptsize $v$}};
	\fill[red!80!black] (+0.5,0) circle (0pt) node {{\scriptsize $u$}};
	\end{tikzpicture}
	%%%%%%%%%%%%%%%%%%%%%% 
}
\, , \quad 
\tikzzbox{%
	%%%%%%%%%%%%%%%%%%%%%% 
	\begin{tikzpicture}[very thick,scale=1.0,color=blue!50!black, baseline=0cm]
	\coordinate (A) at (1,1.25);
	\coordinate (A2) at (-1,-1.25);
	\coordinate (d1) at (-1.5,-1.25);
	\coordinate (d2) at (+1.5,-1.25);
	\coordinate (u1) at (-1.5,1.25);
	\coordinate (u2) at (+1.5,1.25);
	%
	% colouring: 
	\fill [orange!40!white, opacity=0.7] (d1) -- (d2) -- (u2) -- (u1); 
	\draw[thin] (d1) -- (d2) -- (u2) -- (u1) -- (d1);
	%
	% strings: 
	\draw[directed] (0,0) .. controls +(0,1) and +(0,1) .. (-1,0);
	\draw[directed] (1,0) .. controls +(0,-1) and +(0,-1) .. (0,0);
	\draw (-1,0) -- (A2); 
	\draw (1,0) -- (A); 
	%
	% labels: 
	\fill ($(A)+(-0.1,0)$) circle (0pt) node[below right] {{\small $\dX$}};
	\fill ($(A2)+(0.1,0)$) circle (0pt) node[above left] {{\small $\dX$}};
	\fill[red!80!black] (-0.5,0) circle (0pt) node {{\scriptsize $v$}};
	\fill[red!80!black] (+0.5,0) circle (0pt) node {{\scriptsize $u$}};
	\end{tikzpicture}
	%%%%%%%%%%%%%%%%%%%%%% 
}
= 
\tikzzbox{%
	%%%%%%%%%%%%%%%%%%%%%% 
	\begin{tikzpicture}[very thick,scale=1.0,color=blue!50!black, baseline=0cm]
	\coordinate (A) at (0,1.25);
	\coordinate (A2) at (0,-1.25);
	\coordinate (d1) at (-1,-1.25);
	\coordinate (d2) at (+1,-1.25);
	\coordinate (u1) at (-1,1.25);
	\coordinate (u2) at (+1,1.25);
	%
	% colouring: 
	\fill [orange!40!white, opacity=0.7] (d1) -- (d2) -- (u2) -- (u1); 
	\draw[thin] (d1) -- (d2) -- (u2) -- (u1) -- (d1);
	%
	% strings: 
	\draw (A) -- (A2); 
	%
	% labels: 
	\fill ($(A2)+(-0.1,0)$) circle (0pt) node[above right] {{\small $\dX$}};
	\fill[red!80!black] (-0.5,0) circle (0pt) node {{\scriptsize $u$}};
	\fill[red!80!black] (+0.5,0) circle (0pt) node {{\scriptsize $v$}};
	\end{tikzpicture}
	%%%%%%%%%%%%%%%%%%%%%% 
}
\ee 
are satisfied for the left adjoint, and analogously for the right adjoint. 
If they exist, left and right adjoints are unique up to 
%arXiv_v3: 
	a 
unique 2-isomorphism compatible with the adjunction maps~\eqref{eq:AdjMaps}. 

\begin{definition}%[Section 2.1 \cite{cr1210.6363}]
	\begin{enumerate}
		\item 
		A \emph{pivotal 2-category} is a 2-category~$\mathcal B$ together with the choice of right adjunction data $(\Xd, \tev_X, \tcoev_X)$ as well as left adjunction data $(\Xd, \ev_X, \coev_X)$ (with the same underlying 1-morphisms $\dX=\Xd$) for all 1-morphisms~$X$ in~$\mathcal B$ such that
		\be
		\label{eq:pivotality}
		\tikzzbox{
			%%%%%%%%%%%%%%%%%%%%%%
			\begin{tikzpicture}[very thick,scale=0.85,color=blue!50!black, baseline=0cm]
			\coordinate (A) at (-1,1.25);
			\coordinate (A2) at (1,-1.45);
			\coordinate (d1) at (-1.5,-1.45);
			\coordinate (d2) at (+1.5,-1.45);
			\coordinate (u1) at (-1.5,1.25);
			\coordinate (u2) at (+1.5,1.25);
			%
			% colouring: 
			\fill [orange!40!white, opacity=0.7] (d1) -- (d2) -- (u2) -- (u1); 
			\draw[thin] (d1) -- (d2) -- (u2) -- (u1) -- (d1);
			\fill (A) circle (0pt) node[below right] {{\small $Z^\vee$}};
			\fill (A2) circle (0pt) node[above left] {{\small $X^\vee$}};
			\fill (0,0) circle (3.0pt) node[left] {{\small $\xi$}};
			\draw[redirected] (0,0) .. controls +(0,-1) and +(0,-1) .. (-1,0);
			\draw[redirected] (1,0) .. controls +(0,1) and +(0,1) .. (0,0);
			\draw (-1,0) -- (A); 
			\draw (1,0) -- (A2); 
			\end{tikzpicture}
			%%%%%%%%%%%%%%%%%%%%%%
		}
		\! = \!
		\tikzzbox{
			%%%%%%%%%%%%%%%%%%%%%%
			\begin{tikzpicture}[very thick,scale=0.85,color=blue!50!black, baseline=0cm]
			\coordinate (A) at (1,1.25);
			\coordinate (A2) at (-1,-1.45);
			\coordinate (d1) at (-1.5,-1.45);
			\coordinate (d2) at (+1.5,-1.45);
			\coordinate (u1) at (-1.5,1.25);
			\coordinate (u2) at (+1.5,1.25);
			%
			% colouring: 
			\fill [orange!40!white, opacity=0.7] (d1) -- (d2) -- (u2) -- (u1); 
			\draw[thin] (d1) -- (d2) -- (u2) -- (u1) -- (d1);
			\fill (A) circle (0pt) node[below left] {{\small ${}^{\vee}\!Z$}};
			\fill (A2) circle (0pt) node[above right] {{\small $\dX$}};
			\draw[directed] (0,0) .. controls +(0,1) and +(0,1) .. (-1,0);
			\draw[directed] (1,0) .. controls +(0,-1) and +(0,-1) .. (0,0);
			\fill (0,0) circle (3.0pt) node[right] {{\small $\xi$}};
			\draw (-1,0) -- (A2); 
			\draw (1,0) -- (A); 
			\end{tikzpicture}
			%%%%%%%%%%%%%%%%%%%%%%
		}
		\, , \quad
		\tikzzbox{
			%%%%%%%%%%%%%%%%%%%%%%
			\begin{tikzpicture}[very thick,scale=0.65,color=blue!50!black, baseline=-0.2cm, rotate=180]
			\coordinate (X) at (3,-1.25);
			\coordinate (Y) at (2,-1.25);
			\coordinate (XY) at (-1,2.25);
			\coordinate (d1) at (-1.5,-1.25);
			\coordinate (d2) at (+3.5,-1.25);
			\coordinate (u1) at (-1.5,2.25);
			\coordinate (u2) at (+3.5,2.25);
			%
			% colouring: 
			\fill [orange!40!white, opacity=0.7] (d1) -- (d2) -- (u2) -- (u1); 
			\draw[thin] (d1) -- (d2) -- (u2) -- (u1) -- (d1);
			\fill (X) circle (0pt) node[below right] {{\small $\!\! X^\vee$}};
			\fill (Y) circle (0pt) node[below right] {{\small $\!Y^\vee$}};
			\fill ($(XY)+(0,0.2)$) circle (0pt) node[above left] {{\small $(Y\otimes X)^\vee\!{}$}};
			\draw[redirected] (1,0) .. controls +(0,1) and +(0,1) .. (2,0);
			\draw[redirected] (0,0) .. controls +(0,2) and +(0,2) .. (3,0);
			\draw[redirected, ultra thick] (-1,0) .. controls +(0,-1) and +(0,-1) .. (0.5,0);
			\draw (2,0) -- (Y);
			\draw (3,0) -- (X);
			\draw[dotted] (0,0) -- (1,0);
			\draw[ultra thick] (-1,0) -- (XY);
			\end{tikzpicture}
			%%%%%%%%%%%%%%%%%%%%%%
		}
		\! =  \!
		\tikzzbox{
			%%%%%%%%%%%%%%%%%%%%%%
			\begin{tikzpicture}[very thick,scale=0.65,color=blue!50!black, baseline=-0.2cm, rotate=180]
			\coordinate (X) at (-1,-1.25);
			\coordinate (Y) at (-2,-1.25);
			\coordinate (XY) at (2,2.25);
			\coordinate (d1) at (-2.5,-1.25);
			\coordinate (d2) at (+2.5,-1.25);
			\coordinate (u1) at (-2.5,2.25);
			\coordinate (u2) at (+2.5,2.25);
			%
			% colouring: 
			\fill [orange!40!white, opacity=0.7] (d1) -- (d2) -- (u2) -- (u1); 
			\draw[thin] (d1) -- (d2) -- (u2) -- (u1) -- (d1);
			\fill (X) circle (0pt) node[below left] {{\small $\; \dX\!\!$}};
			\fill (Y) circle (0pt) node[below left] {{\small $\; {}^\vee Y\!\!$}};
			\fill ($(XY)+(0,0.2)$) circle (0pt) node[above right] {{\small $\!{}^\vee\hspace{-1.4pt}(Y\otimes X)$}};
			\draw[redirected] (0,0) .. controls +(0,1) and +(0,1) .. (-1,0);
			\draw[redirected] (1,0) .. controls +(0,2) and +(0,2) .. (-2,0);
			\draw[redirected, ultra thick] (2,0) .. controls +(0,-1) and +(0,-1) .. (0.5,0);
			\draw (-1,0) -- (X);
			\draw (-2,0) -- (Y);
			\draw[dotted] (0,0) -- (1,0);
			\draw[ultra thick] (2,0) -- (XY);
			\end{tikzpicture}
			%%%%%%%%%%%%%%%%%%%%%% 
		}
		\ee 
		for all 2-morphisms $\xi\colon Z\lra X$ and all composable 1-morphisms $X,Y$ in~$\mathcal B$. 
		\item 
		The \textsl{left} and \textsl{right traces} of a 2-endomorphism~$\chi$ of $X\colon u \lra v$ in a pivotal 2-category are 
		\be 
		\tr_{\textrm{l}}(\chi) := 
		\tikzzbox
		{
			%%%%%%%%%%%%%%%%%%%%%%
			\begin{tikzpicture}[very thick,scale=0.7,color=blue!50!black, baseline=-0.1cm]
			\fill [orange!40!white, opacity=0.7] (-1.7,-1.7) -- (-1.7,1.7) -- (1.7,1.7) -- (1.7,-1.7); 
			\draw[thin] (-1.7,-1.7) -- (-1.7,1.7) -- (1.7,1.7) -- (1.7,-1.7) -- (-1.7,-1.7);
%			\fill[orange!40!white, opacity=0.7] (0,0) circle (1.7);
%			\draw[thin] (0,0) circle (1.7);
			\fill (-135:1.32) circle (0pt) node[red!80!black] {{\scriptsize$u$}};
			\fill (-135:0.55) circle (0pt) node[red!80!black] {{\scriptsize$v$}};
			\draw (0,0) circle (0.95);
			\fill (45:1.32) circle (0pt) node {{\small$X$}};
			\draw[<-, very thick] (0.100,-0.95) -- (-0.101,-0.95) node[above] {}; 
			\draw[<-, very thick] (-0.100,0.95) -- (0.101,0.95) node[below] {}; 
			\fill (0:0.95) circle (3.5pt) node[left] {{\small$\chi$}};
			\end{tikzpicture} 
			%%%%%%%%%%%%%%%%%%%%%% 
		} 
		\in \End(1_u)
		\, , \quad 
		\tr_{\textrm{r}}(\chi) := 
		\tikzzbox
		{
			%%%%%%%%%%%%%%%%%%%%%%
			\begin{tikzpicture}[very thick,scale=0.7,color=blue!50!black, baseline=-0.1cm]
			\fill [orange!40!white, opacity=0.7] (-1.7,-1.7) -- (-1.7,1.7) -- (1.7,1.7) -- (1.7,-1.7); 
			\draw[thin] (-1.7,-1.7) -- (-1.7,1.7) -- (1.7,1.7) -- (1.7,-1.7) -- (-1.7,-1.7);
%			\fill[orange!40!white, opacity=0.7] (0,0) circle (1.7);
%			\draw[thin] (0,0) circle (1.7);
			\fill (-45:1.32) circle (0pt) node[red!80!black] {{\scriptsize$v$}};
			\fill (-45:0.55) circle (0pt) node[red!80!black] {{\scriptsize$u$}};
			\draw (0,0) circle (0.95);
			\fill (135:1.32) circle (0pt) node {{\small$X$}};
			\draw[->, very thick] (0.100,-0.95) -- (-0.101,-0.95) node[above] {}; 
			\draw[->, very thick] (-0.100,0.95) -- (0.101,0.95) node[below] {}; 
			\fill (180:0.95) circle (3.5pt) node[right] {{\small$\chi$}};
			\end{tikzpicture} 
			%%%%%%%%%%%%%%%%%%%%%% 
		}
		\in \End(1_v)
		\, , 
		\ee
		and the \textsl{left} and \textsl{right quantum dimensions of~$X$} are $\dim_{\textrm{l}}(X) := \tr_{\textrm{l}}(1_X)$ 
		and 
		$\dim_{\textrm{r}}(X) := \tr_{\textrm{r}}(1_X)$. 
		\item 
		A 2-functor $F\colon \mathcal B \lra \mathcal B'$ between pivotal 2-categories is \textsl{pivotal} if for all 1-morphisms~$X$ in~$\mathcal B$ the diagram
		\be
		\label{eq:ConditionPivotal}
		\begin{tikzcd}[column sep=3em, row sep=3em]
		{F}(\Xd) \ar[r,"="] \ar[d, "\cong",swap] &  {F}({}^\vee\!X) \ar[d, "\cong"] 
		\\ 
		{F}(X)^\vee \ar[r,"=",swap] &  {}^\vee\!{F}(X)
		\end{tikzcd}
		\ee 
		commutes, where the vertical arrows are induced by the fact that 2-functors preserve adjoints, and the uniqueness property of adjoints. 
	\end{enumerate}
\end{definition} 

Examples of pivotal 2-categories are naturally obtained from 2-dimensional defect TQFTs, see \cite{dkr1107.0495} or the review \cite{2dDefectTQFTLectureNotes}: objects are labels for 2-strata (or ``bulk theories''), 1-morphisms are lists of composable labels for 1-strata (or 
%arXiv_v3: 
	%(fusion products of)
``line defects''), while 2-morphisms and their compositions are extracted from the TQFT itself. 
Adjunctions and pivotality arise from orientation reversal, which is strictly involutive. 
In particular, 2-dimensional defect state sum models give rise to a pivotal structure on the 2-category of separable symmetric Frobenius algebras, or equivalently of semisimple Calabi--Yau categories, as we review in Section~\ref{subsec:StateSumModels} below. 
    
It is not true that every pivotal 2-functor that is an equivalence also admits a pivotal inverse.
To show this we give a general Whitehead-style theorem identifying conditions for a pivotal inverse to exist, for which we first  introduce the following:  

\begin{definition}
	\label{def:PivotalEquivalence}
	Let $\mathcal{B}$ be a pivotal 2-category with idempotent complete morphism categories. 
	A 1-morphism $X\colon u\to v$ is a \emph{pivotal equivalence} if $X^\vee\colon v\to u$ is a weak inverse to~$X$ whose witnessing 2-isomorphisms to and from the identity are given in terms of the left and right adjunction data for $X^\vee$. 
\end{definition}

Note that a 1-morphism is a pivotal equivalence if and only if its quantum dimensions are trivial. 

\begin{proposition}[Whitehead theorem for pivotal 2-categories]\label{Prop: Whitehead}
	Let $\mathcal{F}\colon \cat{B}\to \cat{B}'$ be a pivotal 2-functor. 
	There exists a pivotal functor $\mathcal{F}^{-1}\colon \cat{B}' \to \cat{B}$ such that $\mathcal{F}\circ \mathcal{F}^{-1}\cong \id_{\cat{B}'} $ and $\mathcal{F}^{-1}\circ \mathcal{F}\cong \id_{\cat{B}} $ 
	%arXiv_v2: 
		%if 
		 iff
	every object $b'\in \mathcal{B}'$ is pivotally equivalent to $\mathcal{F}(b)$ for some $b\in \mathcal{B}$, and $\mathcal{F}$ is fully faithful on morphism 
	%arXiv_v3:  
		%category.
	categories.  
\end{proposition} 
\begin{proof}
	It is standard that we can construct an inverse functor to $\mathcal{F}$ as follows: For every object $b'\in \mathcal{B}$ we pick an object $\mathcal{F}^{-1}(b')\in \mathcal{B}$ together with an equivalence $\gamma_{b'}\colon \mathcal{F}(\mathcal{F}^{-1}(b')) \to b'$. 
	The value on a 1-morphism $f\colon a' \to b'$ can be constructed by choosing a morphism $\mathcal{F}^{-1}(f)\colon \mathcal{F}^{-1}(a')\to \mathcal{F}^{-1}(b')$ such that $\mathcal{F}(\mathcal{F}^{-1}(f))$ is 2-isomorphic to $\gamma_{b'}^{-1} \circ f \circ \gamma_{a'}$. 
	The value on 2-morphisms can then be constructed by using that $\mathcal{F}$ is fully faithful on Hom categories; this is however not important for the proof.
	To conclude that $\mathcal{F}^{-1}$ is pivotal we can use that $\mathcal{F}$ is pivotal if we know that the adjunction data on $\mathcal{F}\mathcal{F}^{-1}(f^\vee)\cong \gamma_{a'}^{-1}f^\vee \gamma_{b'}$ with respect to $\mathcal{F}\mathcal{F}^{-1}(f)\cong \gamma_{b'}^{-1} \circ f \circ \gamma_{a'}$ are those coming from $\mathcal{B}'$. 
	The adjoint coming from $\mathcal{B}$ is $\gamma_{a'}^{\vee} \circ f^\vee \circ \gamma_{b'}^{-1\vee}$. 
	We see that the relation we want to check holds exactly when $\gamma_{-}$ is a pivotal equivalence. 
	But by assumption we can choose all the $\gamma_{-}$ to be pivotal. 
	This finishes the proof. 
\end{proof}

\medskip 

A \textsl{Frobenius algebra (on an object~$a$ of a 2-category~$\mathcal B$)} is a 1-morphism $A\in\mathcal B(a,a)$ together with (co)multiplication and (co)unit 2-morphisms that satisfy 
\be
\label{eq:FrobeniusAlgebra}
\tikzzbox{%
	%%%%%%%%%%%%%%%%%%%%%%%%%%%%
	\begin{tikzpicture}[very thick,scale=0.53,color=green!50!black, baseline=0.59cm]
	\draw[-dot-] (3,0) .. controls +(0,1) and +(0,1) .. (2,0);
	\draw[-dot-] (2.5,0.75) .. controls +(0,1) and +(0,1) .. (3.5,0.75);
	\draw (3.5,0.75) -- (3.5,0); 
	\draw (3,1.5) -- (3,2.25); 
	\end{tikzpicture} 
	%%%%%%%%%%%%%%%%%%%%%%%%%%%% 
}%
=
\tikzzbox{%
	%%%%%%%%%%%%%%%%%%%%%%%%%%%%
	\begin{tikzpicture}[very thick,scale=0.53,color=green!50!black, baseline=0.59cm]
	\draw[-dot-] (3,0) .. controls +(0,1) and +(0,1) .. (2,0);
	\draw[-dot-] (2.5,0.75) .. controls +(0,1) and +(0,1) .. (1.5,0.75);
	\draw (1.5,0.75) -- (1.5,0); 
	\draw (2,1.5) -- (2,2.25); 
	\end{tikzpicture} 
	%%%%%%%%%%%%%%%%%%%%%%%%%%%% 
}%
\, , \quad
\tikzzbox{%
	%%%%%%%%%%%%%%%%%%%%%%
	\begin{tikzpicture}[very thick,scale=0.4,color=green!50!black, baseline]
	\draw (-0.5,-0.5) node[Odot] (unit) {}; 
	\fill (0,0.6) circle (5.0pt) node (meet) {};
	\draw (unit) .. controls +(0,0.5) and +(-0.5,-0.5) .. (0,0.6);
	\draw (0,-1.5) -- (0,1.5); 
	\end{tikzpicture} 
	%%%%%%%%%%%%%%%%%%%%%% 
}%
=
\tikzzbox{%
	%%%%%%%%%%%%%%%%%%%%%%%%%%%%
	\begin{tikzpicture}[very thick,scale=0.4,color=green!50!black, baseline]
	\draw (0,-1.5) -- (0,1.5); 
	\end{tikzpicture} 
	%%%%%%%%%%%%%%%%%%%%%%%%%%%% 
}%
=
\tikzzbox{%
	%%%%%%%%%%%%%%%%%%%%%%
	\begin{tikzpicture}[very thick,scale=0.4,color=green!50!black, baseline]
	\draw (0.5,-0.5) node[Odot] (unit) {}; 
	\fill (0,0.6) circle (5.0pt) node (meet) {};
	\draw (unit) .. controls +(0,0.5) and +(0.5,-0.5) .. (0,0.6);
	\draw (0,-1.5) -- (0,1.5); 
	\end{tikzpicture} 
	%%%%%%%%%%%%%%%%%%%%% 
}%
\, , \quad
\tikzzbox{%
	%%%%%%%%%%%%%%%%%%%
	\begin{tikzpicture}[very thick,scale=0.53,color=green!50!black, baseline=-0.59cm, rotate=180]
	\draw[-dot-] (3,0) .. controls +(0,1) and +(0,1) .. (2,0);
	\draw[-dot-] (2.5,0.75) .. controls +(0,1) and +(0,1) .. (1.5,0.75);
	\draw (1.5,0.75) -- (1.5,0); 
	\draw (2,1.5) -- (2,2.25); 
	\end{tikzpicture} 
	%%%%%%%%%%%%%%%%%%% 
}%
=
\tikzzbox{%
	%%%%%%%%%%%%%%%%%%%
	\begin{tikzpicture}[very thick,scale=0.53,color=green!50!black, baseline=-0.59cm, rotate=180]
	\draw[-dot-] (3,0) .. controls +(0,1) and +(0,1) .. (2,0);
	\draw[-dot-] (2.5,0.75) .. controls +(0,1) and +(0,1) .. (3.5,0.75);
	\draw (3.5,0.75) -- (3.5,0); 
	\draw (3,1.5) -- (3,2.25); 
	\end{tikzpicture} 
	%%%%%%%%%%%%%%%%%%% 
}%
\, , 
\quad
\tikzzbox{%
	%%%%%%%%%%%%%%%%%%%
	\begin{tikzpicture}[very thick,scale=0.4,color=green!50!black, baseline=0, rotate=180]
	\draw (0.5,-0.5) node[Odot] (unit) {}; 
	\fill (0,0.6) circle (5.0pt) node (meet) {};
	\draw (unit) .. controls +(0,0.5) and +(0.5,-0.5) .. (0,0.6);
	\draw (0,-1.5) -- (0,1.5); 
	\end{tikzpicture} 
	%%%%%%%%%%%%%%%%%%% 
}%
=
\tikzzbox{%
	%%%%%%%%%%%%%%%%%%%
	\begin{tikzpicture}[very thick,scale=0.4,color=green!50!black, baseline=0, rotate=180]
	\draw (0,-1.5) -- (0,1.5); 
	\end{tikzpicture} 
	%%%%%%%%%%%%%%%%%%% 
}%
=
\tikzzbox{%
	%%%%%%%%%%%%%%%%%%%
	\begin{tikzpicture}[very thick,scale=0.4,color=green!50!black, baseline=0cm, rotate=180]
	\draw (-0.5,-0.5) node[Odot] (unit) {}; 
	\fill (0,0.6) circle (5.0pt) node (meet) {};
	\draw (unit) .. controls +(0,0.5) and +(-0.5,-0.5) .. (0,0.6);
	\draw (0,-1.5) -- (0,1.5); 
	\end{tikzpicture} 
	%%%%%%%%%%%%%%%%%%% 
}%
\, , \quad 
\tikzzbox{%
	%%%%%%%%%%%%%%%%%%%%%%
	\begin{tikzpicture}[very thick,scale=0.4,color=green!50!black, baseline=0cm]
	\draw[-dot-] (0,0) .. controls +(0,-1) and +(0,-1) .. (-1,0);
	\draw[-dot-] (1,0) .. controls +(0,1) and +(0,1) .. (0,0);
	\draw (-1,0) -- (-1,1.5); 
	\draw (1,0) -- (1,-1.5); 
	\draw (0.5,0.8) -- (0.5,1.5); 
	\draw (-0.5,-0.8) -- (-0.5,-1.5); 
	\end{tikzpicture}
	%%%%%%%%%%%%%%%%%%%%%%
}%
=
\tikzzbox{%
	%%%%%%%%%%%%%%%%%%%%%%
	\begin{tikzpicture}[very thick,scale=0.4,color=green!50!black, baseline=0cm]
	\draw[-dot-] (0,1.5) .. controls +(0,-1) and +(0,-1) .. (1,1.5);
	\draw[-dot-] (0,-1.5) .. controls +(0,1) and +(0,1) .. (1,-1.5);
	\draw (0.5,-0.8) -- (0.5,0.8); 
	\end{tikzpicture}
	%%%%%%%%%%%%%%%%%%%%%%
}%
=
\tikzzbox{%
	%%%%%%%%%%%%%%%%%%%%%%
	\begin{tikzpicture}[very thick,scale=0.4,color=green!50!black, baseline=0cm]
	\draw[-dot-] (0,0) .. controls +(0,1) and +(0,1) .. (-1,0);
	\draw[-dot-] (1,0) .. controls +(0,-1) and +(0,-1) .. (0,0);
	\draw (-1,0) -- (-1,-1.5); 
	\draw (1,0) -- (1,1.5); 
	\draw (0.5,-0.8) -- (0.5,-1.5); 
	\draw (-0.5,0.8) -- (-0.5,1.5); 
	\end{tikzpicture}
	%%%%%%%%%%%%%%%%%%%%%%
}%
\, , 
\ee
where here and below we show no labels and colourings if they can be suppressed at no relevant cost. 
The algebra is \textsl{separable} if there exists a section for the multiplication as a bimodule map, and \textsl{$\Delta$-separable} if this section is the comultiplication, i.\,e.\ if 
\be 
%%%%%%%%%%%%%%%%%%%%%%
\begin{tikzpicture}[very thick,scale=0.4,color=green!50!black, baseline=-0.1cm]
\draw[-dot-] (0,0) .. controls +(0,-1) and +(0,-1) .. (1,0);
\draw[-dot-] (0,0) .. controls +(0,1) and +(0,1) .. (1,0);
\draw (0.5,-0.8) -- (0.5,-1.5); 
\draw (0.5,0.8) -- (0.5,1.5); 
\end{tikzpicture}
%%%%%%%%%%%%%%%%%%%%%%
\, = \, 
%%%%%%%%%%%%%%%%%%%%%%
\begin{tikzpicture}[very thick,scale=0.4,color=green!50!black, baseline=-0.1cm]
\draw (0.5,-1.5) -- (0.5,1.5); 
\end{tikzpicture}
%%%%%%%%%%%%%%%%%%%%%%
\, . 
\ee 
If~$\mathcal B$ is pivotal, then a Frobenius algebra is \textsl{symmetric} if 
\be 
\label{eq:Symmetry}
%%%%%%%%%%%%%%%%%%%%%%
\begin{tikzpicture}[very thick,scale=0.4,color=green!50!black, baseline=0cm]
\draw[-dot-] (0,0) .. controls +(0,1) and +(0,1) .. (-1,0);
\draw[directedgreen, color=green!50!black] (1,0) .. controls +(0,-1) and +(0,-1) .. (0,0);
\draw (-1,0) -- (-1,-1.5); 
\draw (1,0) -- (1,1.5); 
\draw (-0.5,1.2) node[Odot] (end) {}; 
\draw (-0.5,0.8) -- (end); 
\end{tikzpicture}
%%%%%%%%%%%%%%%%%%%%%%
= 
%%%%%%%%%%%%%%%%%%%%%%
\begin{tikzpicture}[very thick,scale=0.4,color=green!50!black, baseline=0cm]
\draw[redirectedgreen, color=green!50!black] (0,0) .. controls +(0,-1) and +(0,-1) .. (-1,0);
\draw[-dot-] (1,0) .. controls +(0,1) and +(0,1) .. (0,0);
\draw (-1,0) -- (-1,1.5); 
\draw (1,0) -- (1,-1.5); 
\draw (0.5,1.2) node[Odot] (end) {}; 
\draw (0.5,0.8) -- (end); 
\end{tikzpicture}
%%%%%%%%%%%%%%%%%%%%%%
\, . 
\ee 
Such algebras are the objects of the 2-dimensional version (see \cite[Def.\,4.1\,\&\,5.1]{cr1210.6363}) of the 3-dimensional construction which is the main topic of the present paper: 

\begin{definition}
	\label{def:2dOrbifoldCompletion}
	Let~$\mathcal B$ be a pivotal 2-category. 
	Its \textsl{orbifold completion} is the 2-category~$\orb{\mathcal B}$ 
	whose objects are $\Delta$-separable symmetric Frobenius algebras, 
	whose 1-morphisms are bimodules over the underlying algebras, 
	whose 2-morphisms are bimodule maps, 
	whose horizontal composition is the relative tensor product over the intermediate algebra, 
	and whose unit 1-morphisms are algebras viewed as bimodules over themselves. 
\end{definition}

The relative tensor product $X\otimes_A Y$ of composable 1-morphisms in $\orb{\mathcal B}$ can be computed by splitting an idempotent, see e.\,g.\ \cite[Lem.\,2.3]{cr1210.6363}: 
\be 
\label{eq:XAY}
X\otimes_A Y = \operatorname{Im}(p_{X,A,Y})  \, , 
\quad 
p_{X,A,Y} := 
%%%%%%%%%%%%%%%%%%%%%%
\begin{tikzpicture}[very thick,scale=0.75,color=green!50!black, baseline]
\draw[color=blue!50!black] (-1,-1) -- (-1,1); 
\draw[color=blue!50!black] (1,-1) -- (1,1); 
\fill[green!50!black] (-0.1,0.1) circle (0pt) node {{\small $A$}};
\fill[blue!50!black] (-1.3,-0.8) circle (0pt) node {{\small $X$}};
\fill[blue!50!black] (1.3,-0.8) circle (0pt) node {{\small $Y$}};
\fill[color=green!50!black] (-1,0.6) circle (2.5pt) node (meet) {};
\fill[color=green!50!black] (1,0.6) circle (2.5pt) node (meet) {};
\draw[-dot-, color=green!50!black] (0.35,-0.0) .. controls +(0,-0.5) and +(0,-0.5) .. (-0.35,-0.0);
\draw[color=green!50!black] (0.35,-0.0) .. controls +(0,0.25) and +(-0.25,-0.25) .. (1,0.6);
\draw[color=green!50!black] (-0.35,-0.0) .. controls +(0,0.25) and +(0.25,-0.25) .. (-1,0.6);
\draw[color=green!50!black] (0,-0.75) node[Odot] (down) {}; 
\draw[color=green!50!black] (down) -- (0,-0.35); 
\end{tikzpicture} 
%%%%%%%%%%%%%%%%%%%%%%  
\, . 
\ee 
Orbifold completion itself is also an idempotent operation, see \cite[Prop.\,4.2\,\&\,4.7]{cr1210.6363}: 

\begin{theorem}
	Let~$\mathcal B$ be a pivotal 2-category such that the idempotents $p_{X,A,Y}$ in~\eqref{eq:XAY} split. 
	\begin{enumerate}
		\item 
		$\orb{\mathcal B}$ has a pivotal structure with
		\be 
		\ev_X = \!
		%%%%%%%%%%%%%%%%%%%%%%
		\begin{tikzpicture}[very thick,scale=1.0,color=blue!50!black, baseline=.8cm]
		\draw[line width=0pt] 
		(1.75,1.75) node[line width=0pt, color=green!50!black] (A) {{\small$A\vphantom{\Xd }$}}
		(1,0) node[line width=0pt] (D) {{\small$X\vphantom{\Xd }$}}
		(0,0) node[line width=0pt] (s) {{\small$\Xd $}}; 
		\draw[directed] (D) .. controls +(0,1.5) and +(0,1.5) .. (s);
		
		\draw[color=green!50!black] (1.25,0.55) .. controls +(0.0,0.25) and +(0.25,-0.15) .. (0.86,0.95);
		\draw[-dot-, color=green!50!black] (1.25,0.55) .. controls +(0,-0.5) and +(0,-0.5) .. (1.75,0.55);
		
		\draw[color=green!50!black] (1.5,-0.1) node[Odot] (unit) {}; 
		\draw[color=green!50!black] (1.5,0.15) -- (unit);
		
		\fill[color=green!50!black] (0.86,0.95) circle (2pt) node (meet) {};
		
		\draw[color=green!50!black] (1.75,0.55) -- (A);
		\end{tikzpicture}
		%%%%%%%%%%%%%%%%%%%%%%
		\circ \iota_{X^\vee,B,X}
		\, , \quad
		\coev_X =  \pi_{X,A,X^\vee} \circ 
		%%%%%%%%%%%%%%%%%%%%%%
		\begin{tikzpicture}[very thick,scale=1.0,color=blue!50!black, baseline=-.8cm,rotate=180]
		\draw[line width=0pt] 
		(3.21,1.85) node[line width=0pt, color=green!50!black] (B) {{\small$B\vphantom{\dX }$}}
		(3,0) node[line width=0pt] (D) {{\small$X\vphantom{\Xd }$}}
		(2,0) node[line width=0pt] (s) {{\small$\Xd $}}; 
		\draw[redirected] (D) .. controls +(0,1.5) and +(0,1.5) .. (s);
		
		\fill[color=green!50!black] (2.91,0.85) circle (2pt) node (meet) {};
		
		\draw[color=green!50!black] (2.91,0.85) .. controls +(0.2,0.25) and +(0,-0.75) .. (B);
		
		\end{tikzpicture}
		%%%%%%%%%%%%%%%%%%%%%%
		\ee 
		for all $X\in\orb{\mathcal B}(A,B)$, where $\iota_{X^\vee,B,X}\colon X^\vee \otimes_B X \lra X^\vee \otimes X$ and $\pi_{X,A,X^\vee}\colon X\otimes X^\vee \lra X\otimes_A X^\vee$ are part of the splitting maps for $p_{X^\vee,B,X}$ and $p_{X,A,X^\vee}$, respectively, and analogously for the right adjoint. 
		\item 
		The embedding 2-functor $\iota_{\mathcal B}\colon \mathcal B \lra \orb{\mathcal B}$, $u\lmt 1_u$, is pivotal. 
		\item 
		$\iota_{\orb{\mathcal B}} \colon \orb{\mathcal B} \lra \orb{(\orb{\mathcal B})}$ is an equivalence. 
	\end{enumerate}
\end{theorem}

\begin{remark}
	Orbifold completion is the purely algebraic part of the generalised orbifold construction for \textsl{oriented} 2-dimensional TQFTs of \cite{cr1210.6363} (which was worked out for arbitrary dimension in \cite{CRS1}). 
	In this context the defining relations \eqref{eq:FrobeniusAlgebra}--\eqref{eq:Symmetry} of objects in $\orb{\mathcal B}$ precisely encode invariance under oriented Pachner moves between triangulations of 2-dimensional bordisms. 
	By dropping the symmetry condition~\eqref{eq:Symmetry} on objects, one obtains another idempotent operation introduced as ``equivariant completion'' $\mathcal B_{\textrm{eq}}$ in \cite{cr1210.6363}, which does not need~$\mathcal B$ to be pivotal. 
	As worked out in detail in \cite{FrauenbergerMasterThesis}, $\mathcal B_{\textrm{eq}}$ is equivalent to the condensation completion of \cite{GaiottoJohnsonFreyd} restricted to 2-categories. 
	The latter theory is inspired by the cobordism hypothesis for fully extended \textsl{framed} TQFTs. 
\end{remark}

\subsection{3-categories with adjoints} 
\label{subsec:3categories}

For the general definition and theory of (not necessarily strict) 3-categories we refer to \cite{GPS, Gurskibook}. 
Every 3-category is equivalent to a \textsl{Gray category}, i.\,e.\ to a category enriched over the symmetric monoidal category of strict 2-categories and strict 2-functors with the Gray tensor product. 
As reviewed extensively in \cite[Sect.\,3]{BMS} and more concisely in \cite[Sect.\,3.1.2]{CMS}, a Gray category~$\tric$ is a 3-categorical structure whose Hom 2-categories (with horizontal and vertical compositions ``$\circ$'' and~``$\cdot$'', respectively) are strict, as are the composition 2-functors
\be 
\label{eq:BoxComp}
u\btimes (-) \colon \tric(c,a) \lra \tric(c,b) \, , 
	\quad 
	(-)\btimes u \colon \tric(b,c) \lra \tric(a,c) 
\ee 
for all 1-morphisms $u\in\tric(a,b)$. 
The only not necessarily strict aspect of~$\tric$ is the \textsl{tensorator}, i.\,e.\ a natural family of 3-isomorphisms
\be 
\sigma_{X,Y} \colon \big( X \btimes 1_{u'} \big) \circ \big( 1_v \btimes Y \big) 
	\stackrel{\cong}{\lra}
	\big( 1_{v'} \btimes Y \big) \circ \big( X \btimes 1_u \big) 
\ee 
for all $\btimes$-composable 2-morphisms $Y\colon u\lra u'$ and $X\colon v\lra v'$. 

The tensorator is subject to axioms 
%arXiv_v3: 
	%which 
	 that
are manifest in the graphical calculus developed in \cite{TrimbleSurfaceDiagrams, BMS} which we discuss next. 
Our conventions deviate from those in \cite{BMS}, but they are in line with those in \cite[Sect.\,3]{CMS}, to which we refer for more details and further illustrations. 
Our graphical conventions for the compositions ``$\circ$'' and~``$\cdot$''are as in Section~\ref{subsec:2categories}. 
The additional composition ``$\btimes$'' is read ``from front to back''. 
Hence for example 
\be 
\label{eq:3dExample}
\tikzzbox{\begin{tikzpicture}[thick,scale=2.321,color=blue!50!black, baseline=0.3cm, >=stealth, 
	style={x={(-0.6cm,-0.4cm)},y={(1cm,-0.2cm)},z={(0cm,0.9cm)}}]
	%: where to put leftmost T-line: 
	\pgfmathsetmacro{\yy}{0.2}
	\coordinate (P) at (0.5, \yy, 0);
	\coordinate (R) at (0.625, 0.5 + \yy/2, 0);
	\coordinate (L) at (0.5, 0, 0);
	\coordinate (R1) at (0.25, 1, 0);
	\coordinate (R2) at (0.5, 1, 0);
	\coordinate (R3) at (0.75, 1, 0);
	% top vertices: 
	\coordinate (Pt) at (0.5, \yy, 1);
	\coordinate (Rt) at (0.375, 0.5 + \yy/2, 1);
	\coordinate (Lt) at (0.5, 0, 1);
	\coordinate (R1t) at (0.25, 1, 1);
	\coordinate (R2t) at (0.5, 1, 1);
	\coordinate (R3t) at (0.75, 1, 1);
	%alpha: 
	\coordinate (alpha) at (0.5, 0.5, 0.5);
	%
	%T-line: 
	%
	% A-planes: 
	\fill [orange!80,opacity=0.545] (L) -- (P) -- (alpha) -- (Pt) -- (Lt);
	\fill [orange!80,opacity=0.545] (Pt) -- (Rt) -- (alpha);
	\fill [orange!80,opacity=0.545] (Rt) -- (R1t) -- (R1) -- (P) -- (alpha);
	\fill [orange!80,opacity=0.545] (Rt) -- (R2t) -- (R2) -- (R) -- (alpha);
	%T-line: 
	\draw[string=green!30!black, ultra thick] (alpha) -- (Rt);
	% A-label (back): 
	% T-label (back):
	\fill[color=blue!50!black] (0.5,0.77,0.77) circle (0pt) node[left] (0up) { {\scriptsize$B$} };
	\fill [orange!80,opacity=0.545] (Pt) -- (R3t) -- (R3) -- (R) -- (alpha);
	% A-plane (front)
	\fill [orange!80,opacity=0.545] (P) -- (R) -- (alpha);
	%
	% T-lines: 
	\draw[string=blue!50!black, ultra thick] (P) -- (alpha);
	\draw[string=blue!50!black, ultra thick] (R) -- (alpha);
	\draw[string=blue!50!black, ultra thick] (alpha) -- (Pt);
	%
	% labels: 
	\fill[color=blue!50!black] (alpha) circle (1.2pt) node[left] (0up) { {\scriptsize$\varphi$} };
	%
	% labels: 
	\fill[red!80!black] (0.5,0.35,-0.17) circle (0pt) node {{\scriptsize $\hspace{7em}u$}};
	\fill[red!80!black] (0.5,0.35,-0.02) circle (0pt) node {{\scriptsize $\hspace{1em}v$}};
	\fill[red!80!black] (0.5,0.35,0.4) circle (0pt) node {{\scriptsize $\hspace{-4em}w$}};
	\fill[red!80!black] (0.5,0.35,-0.07) circle (0pt) node {{\scriptsize $\hspace{9.5em}s$}};
	\fill[red!80!black] (0.5,0.35,0.05) circle (0pt) node {{\scriptsize $\hspace{12em}r$}};
	\fill[red!80!black] (0.5,0.35,0.95) circle (0pt) node {{\scriptsize $\hspace{3em}t$}};
	\fill[black] (0.5,0.35,-0.15) circle (0pt) node {{\scriptsize $\hspace{0em}a$}};
	\fill[black] (0.5,0.35,-0.19) circle (0pt) node {{\scriptsize $\hspace{9.3em}b$}};
	\fill[black] (0.5,0.35,-0.08) circle (0pt) node {{\scriptsize $\hspace{11.8em}c$}};
	\fill[black] (0.5,0.35,1.05) circle (0pt) node {{\scriptsize $\hspace{10em}d$}};
	%
	% T labels: 
	\fill[color=blue!50!black] (0.5,0.35,0.25) circle (0pt) node[left] (0up) { {\scriptsize$X\,$} };
	\fill[color=blue!50!black] (0.5,0.8,0.21) circle (0pt) node[left] (0up) { {\scriptsize$\hspace{-3em}Y$} };
	\fill[color=blue!50!black] (0.5,0.4,0.71) circle (0pt) node[left] (0up) { {\scriptsize$A\,$} };
	%
	% black boundaries: 
	\draw [black,opacity=1, very thin] (Pt) -- (Lt) -- (L) -- (P);
	\draw [black,opacity=1, very thin] (Pt) -- (Rt);
	\draw [black,opacity=1, very thin] (Rt) -- (R1t) -- (R1) -- (P);
	\draw [black,opacity=1, very thin] (Rt) -- (R2t) -- (R2) -- (R);
	\draw [black,opacity=1, very thin] (Pt) -- (R3t) -- (R3) -- (R);
	\draw [black,opacity=1, very thin] (P) -- (R);
	\end{tikzpicture}} 
\quad 
\textrm{represents a 3-morphism } \varphi\colon X \circ \big(1_r \btimes Y \big) \lra A \circ \big( B \btimes 1_u \big) 
\ee
where $u\in \tric(a,b)$, $v\in \tric(a,c)$, $X\colon r \btimes v\lra w$, etc. 
Alternatively, we may only display the source and target of a 3-morphism, in which case~$\varphi$ in~\eqref{eq:3dExample} becomes
\be 
\begin{tikzcd}[column sep=4em, row sep=-1em]
%%%%%%%%%%%%%%%%%%%%%%
\tikzzbox{\begin{tikzpicture}[ultra thick,scale=0.7, baseline=-0.08cm]
	\coordinate (s0) at (0,-1);
	\coordinate (s1) at (0,1);
	\coordinate (v) at (-2,0);
	\coordinate (t) at (-3,0);
	\coordinate (d) at (-1,-0.5);
	\coordinate (u) at (-1,0.5);
	\coordinate (q) at (0,0);
	%
	% lines: 
	\draw[color=orange!80,opacity=0.545] (s1) -- (v) -- (t);
	\draw[color=orange!80,opacity=0.545] (s0) -- (v);
	\draw[color=orange!80,opacity=0.545] (q) -- (d);
	%
	% vertices: 
	\fill[color=blue!50!black] (v) circle (3pt) node[color=blue!50!black, below, font=\footnotesize] {$X$};
	\fill[color=blue!50!black] (d) circle (3pt) node[color=blue!50!black, below, font=\footnotesize] {$Y$};
	\fill[color=red!80!black] (q) circle (0pt) node[color=red!80!black, right, font=\footnotesize] {$s$};
	\fill[color=red!80!black] (s0) circle (0pt) node[color=red!80!black, right, font=\footnotesize] {$u$};
	\fill[color=red!80!black] (s1) circle (0pt) node[color=red!80!black, right, font=\footnotesize] {$r$};
	\fill[color=red!80!black] (t) circle (0pt) node[color=red!80!black, left, font=\footnotesize] {$w$};
	\fill[color=red!80!black] (-1.3,-0.2) circle (0pt) node[color=red!80!black, font=\footnotesize] {$v$};
	\end{tikzpicture}}%%popende 
%%%%%%%%%%%%%%%%%%%%%% 
\ar[r, "\varphi"] 
&
%%%%%%%%%%%%%%%%%%%%%%
\tikzzbox{\begin{tikzpicture}[ultra thick,scale=0.7, baseline=-0.08cm]
	\coordinate (s0) at (0,-1);
	\coordinate (s1) at (0,1);
	\coordinate (v) at (-2,0);
	\coordinate (t) at (-3,0);
	\coordinate (d) at (-1,0.5);
	\coordinate (u) at (-1,0.5);
	\coordinate (q) at (0,0);
	%
	% lines: 
	\draw[color=orange!80,opacity=0.545] (s1) -- (v) -- (t);
	\draw[color=orange!80,opacity=0.545] (s0) -- (v);
	\draw[color=orange!80,opacity=0.545] (q) -- (d);
	%
	% vertices: 
	\fill[color=blue!50!black] (v) circle (3pt) node[color=blue!50!black, above, font=\footnotesize] {$A$};
	\fill[color=blue!50!black] (d) circle (3pt) node[color=blue!50!black, above, font=\footnotesize] {$B$};
	\fill[color=red!80!black] (q) circle (0pt) node[color=red!80!black, right, font=\footnotesize] {$s$};
	\fill[color=red!80!black] (s0) circle (0pt) node[color=red!80!black, right, font=\footnotesize] {$u$};
	\fill[color=red!80!black] (s1) circle (0pt) node[color=red!80!black, right, font=\footnotesize] {$r$};
	\fill[color=red!80!black] (t) circle (0pt) node[color=red!80!black, left, font=\footnotesize] {$w$};
	\fill[color=red!80!black] (-1.3,0) circle (0pt) node[color=red!80!black, font=\footnotesize] {$t$};
	\end{tikzpicture}}%%popende 
%%%%%%%%%%%%%%%%%%%%%% 
\end{tikzcd}
\, .
\ee 

\medskip 

The notion of \textsl{Gray category with duals} was introduced and 
%arXiv_v3: 
	%extensively 
studied in \cite{BMS}. 
It involves a Gray category~$\tric$ together with strict pivotal structures for its Hom-2-categories and composition 2-functors~\eqref{eq:BoxComp}, as well as adjoints for 1-morphisms as follows: 
for every $u\in\tric(a,b)$ there is a choice of $u^\# \in \tric(b,a)$ with $u^{\#\#} = u$ and $\coev_u \colon 1_b \lra u\btimes u^\#$ as well as an invertible \textsl{cusp isomorphism} (or \textsl{triangulator}, or \textsl{Zorro movie})
\be 
c_u = 
\tikzzbox{%
	%%%%%%%%%%%%%%%%%%%%%% 
	\begin{tikzpicture}[thick,scale=0.5,color=black, baseline=1.0cm]
	\coordinate (p1) at (0,0);
	\coordinate (p2) at (1.5,-0.5);
	\coordinate (p2s) at (4,-0.5);
	\coordinate (p3) at (1.5,0.5);
	\coordinate (p4) at (3,1);
	\coordinate (p5) at (1.5,1.5);
	\coordinate (p6) at (-1,1.5);
	\coordinate (u1) at (0,2.5);
	\coordinate (u2) at (1.5,2);
	\coordinate (u2s) at (4,3.5);
	\coordinate (u3) at (1.5,3);
	\coordinate (u4) at (3,3.5);
	\coordinate (u5) at (1.5,4);
	\coordinate (u6) at (-1,3.5);
	%
	% triple fill: 
	\fill [orange!90!white, opacity=1] 
	(p4) .. controls +(0,0.25) and +(1,0) ..  (p5) 
	-- (p5) -- (0.79,1.5) -- (1.5,2.5) 
	-- (1.5,2.5) .. controls +(0,0) and +(0,0.5) ..  (p4);
	%
	% double fill: 
	\fill [orange!30!white, opacity=0.8] 
	(1.5,2.5) .. controls +(0,0) and +(0,0.5) ..  (p4)
	-- 
	(p4) .. controls +(0,-0.25) and +(1,0) ..  (p3)
	-- (p3) .. controls +(-1,0) and +(0,0.25) ..  (p1)
	-- 
	(p1) .. controls +(0,0.5) and +(0,0) ..  (1.5,2.5);
	%
	% red line middle:  
	\draw[thin] (p4) .. controls +(0,0.25) and +(1,0) ..  (p5) -- (p6); 
	\draw[thin] (p1) .. controls +(0,0.25) and +(-1,0) ..  (p3)
	-- (p3) .. controls +(1,0) and +(0,-0.25) ..  (p4); 
	%
	% black auxiliary lines: 
	\draw[thin] (p2s) --  (u2s); 
	\draw[thin] (p6) --  (u6); 
	\fill [orange!20!white, opacity=0.8] 
	(p2s) -- (p2)
	-- (p2) .. controls +(-1,0) and +(0,-0.25) ..  (p1) 
	-- (p1) .. controls +(0,0.5) and +(0,0) ..  (1.5,2.5)
	-- (1.5,2.5) .. controls +(0,0) and +(0,0.5) .. (p4)
	-- (p4) .. controls +(0,0.25) and +(-1,0) ..  (p5) 
	-- (p5) -- (p6) -- (u6) -- (u2s)
	;
	%
	%%%%%%%%%%%%%%%%%%%%%%%%%%%5
	% near cusp: 
	\draw[thin, dotted] (p4) .. controls +(0,0.5) and +(0,0) ..  (1.5,2.5); 
	\draw[thin] (p1) .. controls +(0,0.5) and +(0,0) ..  (1.5,2.5); 
	%%%%%%%%%%%%%%%%%%%%%%%%%%%
	%
	%
	% red boundaries: 
	\draw[thin] (u6) -- (u2s); 
	\draw[thin] (p1) .. controls +(0,-0.25) and +(-1,0) ..  (p2) -- (p2s); 
	\draw[thin] (0.79,1.5) -- (p6); 
	%
	%%%%%%%%%%%%%%%%%%%%%%%%%%%
	% labels: 
	\fill[red!80!black] (3.5,0) circle (0pt) node {{\scriptsize $u$}};
	\fill (1.5,2.5) circle (2.5pt) node[above] {{\small $c_u$}};
	\fill (3.25,1.8) circle (0pt) node {{\scriptsize $\coev_u$}};
	\fill (-0.1,1) circle (0pt) node {{\scriptsize $\ev_u$}};
	%%%%%%%%%%%%%%%%%%%%%%%%%%%
	%
	\end{tikzpicture}
	%%%%%%%%%%%%%%%%%%%%%% 
}   
\quad 
\text{where } \ev_u := \big( \coev_{u^\#} \big)^{\!\vee} \, , 
\ee 
and $(-)^\vee$ is part of the pivotal structure on $\tric(a,a)$. 
For a detailed discussion of the axioms we refer to \cite[Sect.\,4]{BMS}, and to \cite[Sect.\,3.2.2]{CMS} for a concise summary using our graphical conventions. 

The adjunction 2-morphisms $\ev_u, \coev_u$ together with~$c_u$ witness~$u^\#$ as the left adjoint of~$u$. 
But~$u^\#$ is also a right adjoint via $\ev_{u^\#}, \coev_{u^\#}$. 
Using this with the tensorator and cusp isomorphisms, one obtains canonical 3-isomorphisms
\be 
\begin{tikzcd}[column sep=4em, row sep=-1em]
\Phi_{Y,X} \colon Y^\# \circ X^\# = 
%%%%%%%%%%%%%%%%%%%%%%
\tikzzbox{\begin{tikzpicture}[ultra thick,scale=0.45,baseline=-0.1cm]
	\coordinate (ml) at (-1.5,2);
	\coordinate (coev) at (0.5,2);
	\coordinate (ev) at (-0.5,0);
	\coordinate (coev2) at (2,0);
	\coordinate (ev2) at (1,-2);
	\coordinate (mr) at (2.5,-2);
	\coordinate (ac2) at ($(coev)-0.5*(coev)+0.5*(ev)$);
	\coordinate (ac) at ($(coev2)-0.5*(coev2)+0.5*(ev2)$);
	\coordinate (B) at (2.5,-1);
	\coordinate (B2) at (2.5,1);
	\coordinate (B3) at (2.5,-2.3);
	\coordinate (Bcoev) at (0.4,1.1);
	\coordinate (Bev) at (-1.1,-2.3);
	%
	% lines: 
	\draw[color=orange!80,opacity=0.545] (ml) -- (coev) -- (ev) -- (coev2) -- (ev2)-- (mr);
	%
	% vertices: 
	\fill[color=blue!50!black] (ac) circle (4pt) node[color=blue!50!black, left, font=\footnotesize] {$X$};
	\fill[color=blue!50!black] (ac2) circle (4pt) node[color=blue!50!black, left, font=\footnotesize] {$Y$};
	%
	% labels: 
	\fill[color=red!80!black] (mr) circle (0pt) node[color=red!80!black, font=\footnotesize, right] {$u^\#$};
	\fill[color=red!80!black] (ml) circle (0pt) node[color=red!80!black, font=\footnotesize, left] {$w^\#$};
	\fill[color=red!80!black] (1,0.2) circle (0pt) node[color=red!80!black, font=\footnotesize] {$v^\#$};
	\end{tikzpicture}}%%popende 
%%%%%%%%%%%%%%%%%%%%%%  
\ar[r, "\cong"] 
&
\tikzzbox{\begin{tikzpicture}[ultra thick,scale=0.45,baseline=-0.1cm]
	\coordinate (ml) at (-1.5,2);
	\coordinate (coev) at (1,2);
	\coordinate (ev) at (-1,-2);
	\coordinate (mr) at (2.5,-2);
	\coordinate (ac2) at ($(coev)-0.3*(coev)+0.3*(ev)$);
	\coordinate (ac) at ($(coev)-0.7*(coev)+0.7*(ev)$);
	\coordinate (B) at (2.5,-0.8);
	\coordinate (B2) at (2.5,0.8);
	\coordinate (B3) at (2.5,-2.3);
	\coordinate (Bcoev) at (0.5,1.35);
	\coordinate (Bev) at (-1.5,-2.3);
	%
	% lines: 
	\draw[color=orange!80,opacity=0.545] (ml) -- (coev) -- (ev) -- (mr);
	%
	%	\draw[color=green!70!black,opacity=0.545] (ac2) -- (Bcoev) .. controls +(-1,-0.2) and +(0,1) .. (Bev) -- (B3);
	%	%
	% vertices: 
	\fill[color=blue!50!black] (ac) circle (4pt) node[color=blue!50!black, left, font=\footnotesize] {$X$};
	\fill[color=blue!50!black] (ac2) circle (4pt) node[color=blue!50!black, left, font=\footnotesize] {$Y$};
	%
	% labels: 
	\fill[color=red!80!black] (mr) circle (0pt) node[color=red!80!black, font=\footnotesize, right] {$u^\#$};
	\fill[color=red!80!black] (ml) circle (0pt) node[color=red!80!black, font=\footnotesize, left] {$w^\#$};
	\fill[color=red!80!black] (0.4,-0.3) circle (0pt) node[color=red!80!black, font=\footnotesize] {$v$};
	\end{tikzpicture}}%%popende 
%%%%%%%%%%%%%%%%%%%%%%  
= (X\circ Y)^{\#}
\end{tikzcd}
\ee 
for composable 2-morphisms $Y\colon v\lra w$ and $X\colon u\lra v$ in a Gray category with duals, cf.\ \cite[Fig.\,23\,c)]{BMS}. 
Moreover, there is also a canonical 3-isomorphism 
\be 
\begin{tikzcd}[column sep=4em, row sep=-1em]
\Theta_X \colon X^{\#\#} = 
%%%%%%%%%%%%%%%%%%%%%%  
\tikzzbox{\begin{tikzpicture}[ultra thick,scale=0.45,baseline=-0.1cm]
	\coordinate (ml) at (-3,2);
	\coordinate (coev1) at (2.5,0);
	\coordinate (ev1) at (-1.5,0);
	\coordinate (coev2) at (1.5,0);
	\coordinate (ev2) at (-2.5,0);
	\coordinate (mr) at (3,-2);
	\coordinate (X) at (0,0);
	%
	% lines: 
	\draw[color=orange!80,opacity=0.545] (ml) .. controls +(5,0) and +(-0.5,1) .. (coev1)
	-- (coev1) .. controls +(-1,-1.5) and +(1,-1.5) .. (ev1) -- (coev2)
	-- (coev2) .. controls +(-1,1.5) and +(1,1.5) .. (ev2)
	-- (ev2) .. controls +(0.5,-1) and +(-5,0) .. (mr);
	%
	%	\draw[color=green!70!black,opacity=0.545] (ac2) -- (Bcoev) .. controls +(-1,-0.2) and +(0,1) .. (Bev) -- (B3);
	%	%
	% vertices: 
	\fill[color=blue!50!black] (X) circle (4pt) node[color=blue!50!black, above, font=\footnotesize] {$X$};
	%
	% labels: 
	\fill[color=red!80!black] (mr) circle (0pt) node[color=red!80!black, font=\footnotesize, right] {$u$};
	\fill[color=red!80!black] (ml) circle (0pt) node[color=red!80!black, font=\footnotesize, left] {$v$};
	\end{tikzpicture}}%%popende 
%%%%%%%%%%%%%%%%%%%%%%  
\ar[r, "\cong"] 
&
%%%%%%%%%%%%%%%%%%%%%%  
\tikzzbox{\begin{tikzpicture}[ultra thick,scale=0.45,baseline=-0.1cm]
	\coordinate (ml) at (-1.5,0);
	\coordinate (mr) at (1.5,0);
	\coordinate (X) at (0,0);
	%
	% lines: 
	\draw[color=orange!80,opacity=0.545] (ml) -- (mr);
	%
	% vertices: 
	\fill[color=blue!50!black] (X) circle (4pt) node[color=blue!50!black, above, font=\footnotesize] {$X$};
	%
	% labels: 
	\fill[color=red!80!black] (mr) circle (0pt) node[color=red!80!black, font=\footnotesize, right] {$u$};
	\fill[color=red!80!black] (ml) circle (0pt) node[color=red!80!black, font=\footnotesize, left] {$v$};
	\end{tikzpicture}}%%popende 
%%%%%%%%%%%%%%%%%%%%%%  
= X
\end{tikzcd}
\ee 
which is constructed from cusp, tensorator and adjunction 3-morphisms, cf.\ \cite[Eq.\,(24)\,\&\,Fig.\,32]{BMS}. 
These maps are compatible in the sense that 
\be 
\label{eq:ThetaThetaTheta}
\begin{tikzcd}[column sep=5em, row sep=3em]
Y^{\#\#} \circ X^{\#\#} 
\ar[r, "\Phi_{Y^\#,X^\#}"] 
\ar[d, "\Theta_Y \circ \Theta_X", swap] 
& 
\big(X^\# \circ Y^\# \big)^\# 
\ar[d, "\Phi_{X,Y}^\#"]
\\ 
Y\circ X 
\ar[r, "\Theta_{Y\circ X}^{-1}", swap] 
& 
\big( Y\circ X \big)^{\#\#}
\end{tikzcd} 
\ee 
commutes, which means that $\Theta_Y\circ\Theta_X$ and~$\Theta_{Y\circ X}$ are equal up to tensorator and cusp isomorphisms. 

\medskip 

Given a monoidal 2-category~$\mathcal C$, we denote its \textsl{delooping} 3-category by $\operatorname{B}\!\mathcal C$, which by definition has only a single object whose End 2-category is~$\mathcal C$. 
Clearly $\operatorname{B}\!\mathcal C$ is a Gray category with duals if~$\mathcal C$ is a pivotal monoidal 2-category with compatible duals for objects. 

Further examples of Gray categories with duals are naturally obtained from 3-dimensional defect TQFTs, see \cite{CMS}: 
objects are labels for 3-strata, 1- and 2-morphisms are lists of labels for 2- and 1-strata, respectively, and only 3-morphisms are constructed from the TQFT itself. 
In particular, 3-dimensional defect state sum models, i.\,e.\ Turaev--Viro--Barrett--Westbury models with arbitrary defects, should provide 3-categories that are equivalent to Gray categories with duals. 
Work in this direction was carried out in \cite{CRS1, CRS2, CRS3, Meusburger3dDefectStateSumModels}, and can be concluded with the main results of the present paper, see Section~\ref{sec:TQFTs} below.

\section{Morita categories of ${E_1}$-algebras}
\label{sub:EofC}

Let $\tric$ be a 3-category such that all Hom 2-categories admit finite sifted 2-colimits that commute with composition. 
(We refer to Appendix~\ref{app:2colimits} for some background on 2-colimits; a diagram 2-category~$\mathcal D$ is \textsl{sifted} if 2-colimits of shape~$\mathcal D$ commute with finite products in Cat.)
In this section we define a new 3-category $\ealg(\tric)$ which locally is a version of the 3-category $\Alg_1(\mathcal{B})$ of $E_1$-algebras in a monoidal 2-category~$\mathcal B$, see e.\,g.\ \cite{TheoClaudia}, in the sense that if~$\tric$ has a single object~$*$, then $\ealg(\tric) = \Alg_1(\tric(*,*))$.  

\begin{definition}[Objects]
	\label{def:ObjEAlg}
	An \textsl{object~$\mathcal A$ of $\ealg(\tric)$} consists of an object $a\in \tric$ and a unital $E_1$-algebra
	\be 
	%%%%%%%%%%%%%%%%%%%%%% 
	\tikzzbox{
}%%popende
	%%%%%%%%%%%%%%%%%%%%%% 
	\, .  
	\ee 
\end{definition}   

%arXiv_v3: 
	Throughout we use green to colour surfaces corresponding to 1-morphisms that are part of the structure of $E_1$-algebras (and hence also the orbifold data introduced in Section~\ref{sec:OrbComp3cat} below). 
	Similarly, we shall use colours with orange or red in the graphical calculus for surfaces related to (bi-)modules over such algebras: 

\begin{definition}[1-morphisms]
	\label{def:1MorphismsEAlg}
	Let $\mathcal A = (a,A,\mu_A, \alpha_A, u_A, u^{\textrm{l}}_A, u^{\textrm{r}}_A)$ and $\mathcal B = (b,B,\mu_B, \alpha_B, u_B, u^{\textrm{l}}_B, u^{\textrm{r}}_B)$ be objects of $\ealg(\tric)$. 
	A \textsl{1-morphism~$\mathcal M \colon \mathcal A \lra \mathcal B$} consists of a 1-morphism $M\colon a\to b$
	together with 2-morphisms 
	\be 
	%%%%%%%%%%%%%%%%%%%%%% 
	\tikzzbox{
}%%popende
	%%%%%%%%%%%%%%%%%%%%%% 
	\ee 
	in~$\tric$ such that 
	\be 
	\label{eq:23MovesFor1Morphisms1}
	\begin{tikzcd}[column sep=5em, row sep=-2em]
	%%%%%%%%%%%%%%%%%%%%%%
	\tikzzbox{%
}%%popende 
	%%%%%%%%%%%%%%%%%%%%%%   
	\arrow{ur}[swap]{\alpha_B}
	&
	\end{tikzcd}
	\ee 
	
	\be 
	\label{eq:23MovesFor1Morphisms2}
	\begin{tikzcd}[column sep=5em, row sep=-2em]
	%%%%%%%%%%%%%%%%%%%%%%
	\tikzzbox{%
}%%popende 
	%%%%%%%%%%%%%%%%%%%%%%   
	\arrow{ur}[swap]{\alpha^{\textrm{l}}_M}
	&
	\end{tikzcd}
	\ee 
	
	\be 
	\label{eq:23MovesFor1Morphisms3}
	\begin{tikzcd}[column sep=5em, row sep=-2em]
	%%%%%%%%%%%%%%%%%%%%%%
	\tikzzbox{%
}%%popende 
	%%%%%%%%%%%%%%%%%%%%%%   
	\arrow{ur}[swap]{\alpha^{\textrm{m}}_M}
	&
	\end{tikzcd}
	\ee 
	
	\be 
	\label{eq:23MovesFor1Morphisms4}
	\begin{tikzcd}[column sep=5em, row sep=-2em]
	%%%%%%%%%%%%%%%%%%%%%%
	\tikzzbox{%
}%%popende 
	%%%%%%%%%%%%%%%%%%%%%%   
	\arrow{ur}[swap]{\alpha^{\textrm{r}}_M}
	&
	\end{tikzcd}
	\ee 
	and 
	\be 
	\begin{tikzcd}[column sep=5em, row sep=-1em]
	%%%%%%%%%%%%%%%%%%%%%%
	\tikzzbox{%
}%%popende 
	%%%%%%%%%%%%%%%%%%%%%% 
	\arrow{ru}[swap]{u^{\textrm{l}}_A}
	\arrow{uu}{\alpha^{\textrm{r}}_M}
	& 
	\end{tikzcd}
	\qquad \quad 
	\begin{tikzcd}[column sep=5em, row sep=-1em]
	%%%%%%%%%%%%%%%%%%%%%%
	\tikzzbox{%
}%%popende 
	%%%%%%%%%%%%%%%%%%%%%% 
	\arrow{ru}[swap]{u^{\textrm{l}}_B}
	\arrow{uu}{\alpha^{\textrm{l}}_M}
	& 
	\end{tikzcd}
	\ee 
	
	\be 
	\begin{tikzcd}[column sep=5em, row sep=-1em]
	%%%%%%%%%%%%%%%%%%%%%%
	\tikzzbox{%
}%%popende 
	%%%%%%%%%%%%%%%%%%%%%% 
	\arrow{ru}[swap]{u^{\textrm{r}}_M}
	& 
	\end{tikzcd}
	\qquad \quad 
	\begin{tikzcd}[column sep=5em, row sep=-1em]
	%%%%%%%%%%%%%%%%%%%%%%
	\tikzzbox{%
}%%popende 
	%%%%%%%%%%%%%%%%%%%%%% 
	\arrow{ru}[swap]{u^{\textrm{l}}_M}
	\arrow{uu}[swap]{\alpha^{\textrm{m}}_M}
	& 
	\end{tikzcd}
	\ee  
	commute. 
\end{definition}

\begin{definition}[2-morphisms]
	\label{def:2MorphismsEAlg}
	A \emph{2-morphism} 
	\be 
	\mathcal F \colon 
	\mathcal M 
	= 
	\big( M, \triangleright_M, \triangleleft_M, u_M^{\textrm{l}}, u_M^{\textrm{r}}, \alpha_M^{\textrm{l}},\alpha_M^{\textrm{m}}, \alpha_M^{\textrm{r}} \big) 
	\lra 
	\mathcal M'
	= 
	\big( M', \triangleright_{M'}, \triangleleft_{M'}, u_{M'}^{\textrm{l}}, u_{M'}^{\textrm{r}}, \alpha_{M'}^{\textrm{l}},\alpha_{M'}^{\textrm{m}}, \alpha_{M'}^{\textrm{r}} \big) 
	\ee 
	between 1-morphisms $\mathcal M, \mathcal M' \colon \mathcal A \lra \mathcal B$ in $\ealg(\tric)$ 
	consists of 2- and 3-morphisms 
	\be 
	\label{eq:2MorphismStructureMaps}
	%%%%%%%%%%%%%%%%%%%%%% 
	\tikzzbox{
}%%popende
	%%%%%%%%%%%%%%%%%%%%%% 
	\, . 
	\end{align}
\end{definition}

For later use we note that the inverses of~\eqref{eq:2MorphismStructureMaps} are presented as follows: 
\begin{align}
%\triangleright_F^{-1} 
%= 
%%%%%%%%%%%%%%%%%%%%%% 
\tikzzbox{%
}%%popende
	%%%%%%%%%%%%%%%%%%%%%% 
	\, . 
	\ee 
\end{definition} 

\medskip 

The 2-category structure on $\ealg(\tric)(\mathcal A, \mathcal B) $ is induced from
the composition in~$\tric$ in a straightforward way. 
The composition of 1-morphisms is given in terms of relative box products. 
Here is where our assumptions on the existence of certain 2-colimits enters, as we discuss next. 
Our guide for the description below is the exposition of tricategories in \cite[Sect.\,A.4]{GregorDiss}. 

Let $\mathcal A, \mathcal B, \mathcal C \in \ealg(\tric)$. 
The definition of the composition $\Box_{\mathcal B} \colon \ealg(\tric)(\mathcal B ,\mathcal C) \times \ealg(\tric)(\mathcal A,\mathcal B) \to \ealg(\tric)(\mathcal A,\mathcal C) $ can be split into two steps, where we refer to Appendix~\ref{app:2colimits} for details on 2-colimits:
\begin{itemize}
	\item 
	First we define a 2-functor $ \Box^\Delta_{\mathcal B} \colon \ealg(\tric)(\mathcal B,\mathcal C) \times \mathcal{E}(\cat{C})(\mathcal A,\mathcal B) \to \ealg(\tric)(\mathcal A,\mathcal C)^{\tau_{2}\Delta} $ sending a pair of 
	objects $(\mathcal M, \mathcal N)$ to the truncated simplicial object
	\begin{equation}
	\begin{tikzcd}
	M\Box B \Box B \Box N 
	\ar[r] \ar[r, shift left=2.5] \ar[r, shift right=2.5] 
	& 
	M \Box B \Box N 
	\ar[r, shift left]\ar[r, shift right] 
	& 
	M \Box N \, .
	\end{tikzcd}
	\end{equation} 
	The coherences involved in the diagram are part of the bimodule structure of~$\mathcal M$ and~$\mathcal N$. 
	Note that we do not require to include brackets since we assume that $\tric$ is a Gray category and also that this becomes in a canonical way a diagram in $\ealg(\tric)(\mathcal A, \mathcal C)$. 
	\item 
	We define $\Box_{\mathcal B} $ as the composition of $\Box^{\Delta}_{\mathcal B}$ with the 2-colimit 2-functor 
	$\operatorname{colim}\colon \ealg(\tric)(\mathcal A,\mathcal C)^{\tau_{2}\Delta} 
	\to \ealg(\tric)(\mathcal A,\mathcal C)$. 
	Since $\tau_{2}\Delta$ is sifted this colimit exists and can be computed in the Hom 2-categories of $\tric$.    
\end{itemize} 
For 1-morphism $\mathcal M\colon \mathcal B \lra \mathcal C$ and $\mathcal N\colon \mathcal A \lra \mathcal B$, we denote their composition $\mathcal M \btimes_{\mathcal B} \mathcal N$ in $\ealg(\tric)$ graphically as 
\be 
\label{eq:RelativeGraphicalNotation}
%%%%%%%%%%%%%%%%%%%%%%
\tikzzbox{\begin{tikzpicture}[ultra thick,scale=0.75, baseline=-0.1cm]
	\coordinate (r0) at (0,-0.75);
	\coordinate (r1) at (0,0.75);
	\coordinate (l0) at (-3,-0.75);
	\coordinate (l1) at (-3,0.75);
	\coordinate (acm) at (-1.5,-0.75);
	\coordinate (acn) at (-1.5,0.75);
	%
	% area: 
	\fill[color=orange!50,opacity=0.35] ($(l0)+(0,-0.8pt)$) -- ($(r0)+(0,-0.8pt)$) -- ($(r1)+(0,0.8pt)$) -- ($(l1)+(0,0.8pt)$);
	%
	% lines: 
	\draw[color=orange!50,opacity=0.545] (l0) -- (r0);
	\draw[color=orange!50,opacity=0.545] (l1) -- (r1);
	%
	%	\draw[blue!50!green] ($(acm)+(0,-0.8pt)$) -- ($(acn)+(0,0.8pt)$);
	%
	% vertices: 
	\fill[color=red!30!black] (r0) circle (0pt) node[color=red!80!black, right, font=\footnotesize] {$N$};
	\fill[color=red!30!black] (r1) circle (0pt) node[color=red!80!black, right, font=\footnotesize] {$M$};
	\fill[color=red!30!black] (l0) circle (0pt) node[color=red!80!black, left, font=\footnotesize] {$N$};
	\fill[color=red!30!black] (l1) circle (0pt) node[color=red!80!black, left, font=\footnotesize] {$M$};
	\fill[color=red!30!black] ($(l1)-0.5*(l1)+0.5*(l0)+(-0.05,-0.05)$) circle (0pt) node[color=red!80!black, left, font=\footnotesize] {$\displaystyle{\relpro_{\mathcal B}}$};
	\fill[color=red!30!black] ($(l1)-0.5*(l1)+0.5*(l0)+(3.75,-0.05)$) circle (0pt) node[color=red!80!black, left, font=\footnotesize] {$\displaystyle{\relpro_{\mathcal B}}$};
	\end{tikzpicture}}%%popende 
%%%%%%%%%%%%%%%%%%%%%%   
\;\; = \;\;
%%%%%%%%%%%%%%%%%%%%%%
\tikzzbox{\begin{tikzpicture}[ultra thick,scale=0.75, baseline=-0.65cm]
	\coordinate (r0) at (0,-0.75);
	\coordinate (r1) at (0,0.75);
	\coordinate (l0) at (-3,-0.75);
	\coordinate (l1) at (-3,0.75);
	%
	% lines: 
	\draw[color=orange!50,opacity=0.545] (l0) -- (r0);
	% vertices: 
	\fill[color=red!30!black] (r0) circle (0pt) node[color=red!80!black, right, font=\footnotesize] {$M \relpro_{\mathcal B} N$};
	\fill[color=red!30!black] (l0) circle (0pt) node[color=red!80!black, left, font=\footnotesize] {$M \relpro_{\mathcal B} N$};
	\end{tikzpicture}}%%popende 
%%%%%%%%%%%%%%%%%%%%%%   
\, . 
\ee 

It is useful to reformulate the definition of $M \relpro_{\mathcal B} N$ in terms of a universal property involving balancings, cf.\ \cite{Cooke2019} for the $\Bbbk$-linear case. 
Recall that a 2-morphism $F\in\Hom_{\tric}(M\btimes N, T)$ is \textsl{balanced} if it comes together with a 3-morphism 
\be 
\label{eq:BalancingMap}
f = 
%%%%%%%%%%%%%%%%%%%%%% 
\tikzzbox{\begin{tikzpicture}[thick,scale=2.321,color=blue!50!black, baseline=1.4cm, >=stealth, 
	style={x={(-0.6cm,-0.4cm)},y={(1cm,-0.2cm)},z={(0cm,0.9cm)}}]
	%: where to put leftmost T-line: 
	\pgfmathsetmacro{\yy}{0.2}
	\coordinate (P) at (0.5, \yy, 0);
	\coordinate (L) at (0.5, 0, 0);
	\coordinate (R1) at (0.25, 1, 0);
	\coordinate (R2) at (0.5, 1, 0);
	\coordinate (R3) at (0.75, 1, 0);
	\coordinate (R) at ($(R3)+0.5*(P)-0.5*(R3)$);
	% top vertices: 
	\coordinate (Pt) at (0.5, \yy, 1);
	\coordinate (Rt) at (0.625, 0.5 + \yy/2, 1);
	\coordinate (Lt) at (0.5, 0, 1);
	\coordinate (R1t) at (0.25, 1, 1);
	\coordinate (R2t) at (0.5, 1, 1);
	\coordinate (R3t) at (0.75, 1, 1);
	\coordinate (aP) at (0.5, \yy, 1);
	\coordinate (aR) at (0.625, 0.5 + \yy/2, 1);
	\coordinate (aL) at (0.5, 0, 1);
	\coordinate (aR1) at (0.25, 1, 1);
	\coordinate (aR2) at (0.5, 1, 1);
	\coordinate (aR3) at (0.75, 1, 1);
	% top vertices: 
	\coordinate (aPt) at (0.5, \yy, 2);
	\coordinate (aRt) at (0.375, 0.5 + \yy/2, 2);
	\coordinate (aLt) at (0.5, 0, 2);
	\coordinate (aR1t) at (0.25, 1, 2);
	\coordinate (aR2t) at (0.5, 1, 2);
	\coordinate (aR3t) at (0.75, 1, 2);
	\fill [purple!80,opacity=0.2] (P) -- (aPt) -- (aLt) -- (L);
	\fill [orange!80,opacity=0.545] (P) -- (aPt) -- (aR1t) -- (R1);
	\fill [green!50,opacity=0.545] (R) -- (Pt) -- (aRt) -- (aR2t) -- (R2);
	%
	% triangle back
	\fill[red!80!black] (0.5, 0.38, 1.4) circle (0pt) node[right] (0up) {{\scriptsize$\triangleright_{M}$}};
	%
	% action line above: 
	\draw[color=red!80!black, ultra thick, rounded corners=0.5mm] (aRt) -- (Pt);
	\fill [orange!80,opacity=0.545] (P) -- (aPt) -- (aR3t) -- (R3);
	%
	% action line below: 
	\draw[blue!80!black, ultra thick] (P) -- (aPt);
	%
	% F-line: 
	\draw[color=red!80!black, ultra thick, rounded corners=0.5mm] (R) -- (Pt);
	%
	% T-label: 
	\fill[color=purple!50!black] (0.55, \yy, 0.15) circle (0pt) node[left] (0up) { {\scriptsize$T$} };
	\fill[blue!80!black] (0.45,0.2,1.5) circle (0pt) node[left] (0up) { {\scriptsize$F$} };
	\fill[blue!80!black] (0.45,0.2,0.5) circle (0pt) node[left] (0up) { {\scriptsize$F$} };
	\fill[red!80!black] (0.5,0.89,0.0) circle (0pt) node[left] (0up) { {\scriptsize$N$} };
	\fill[color=blue!60!black] (0.5,1.04,0.13) circle (0pt) node[left] (0up) { {\scriptsize$B$} };
	\fill[color=red!80!black] (0.5,1.19,0.28) circle (0pt) node[left] (0up) { {\scriptsize$M$} };
	%
	% vertices: 
	\fill[color=black] (Pt) circle (1.2pt) node[left] (0up) { {\scriptsize$f$} };
	%
	% triangle front
	\fill[red!80!black] (0.5, 0.38, 0.3) circle (0pt) node[right] (0up) {{\scriptsize$\triangleleft_{N}$}};
	%
	% black boundaries: 
	\draw [black,opacity=1, very thin] (L) -- (aLt) -- (aPt) -- (aR3t) -- (R3) -- (P) -- cycle;
	\draw [black,opacity=1, very thin] (aPt) -- (aR1t) -- (R1) -- (P);
	\draw [black,opacity=1, very thin] (R) -- (R2) -- (aR2t) -- (aRt);
	\end{tikzpicture}}%%popende
%%%%%%%%%%%%%%%%%%%%%% 
\ee 
that satisfies a pentagon compatibility condition with the associators of $\mathcal M,\mathcal B,\mathcal N$. 
Given another such balanced 2-morphism $(G,g)$, a \textsl{balanced} 3-morphism $(F,f) \lra (G,g)$ is a 3-morphism $\zeta\colon F\lra G$ such that 
\be 
%%%%%%%%%%%%%%%%%%%%%% 
\tikzzbox{\begin{tikzpicture}[thick,scale=2.321,color=blue!50!black, baseline=1.4cm, >=stealth, 
	style={x={(-0.6cm,-0.4cm)},y={(1cm,-0.2cm)},z={(0cm,0.9cm)}}]
	%: where to put leftmost T-line: 
	\pgfmathsetmacro{\yy}{0.2}
	\coordinate (P) at (0.5, \yy, 0);
	\coordinate (L) at (0.5, 0, 0);
	\coordinate (R1) at (0.25, 1, 0);
	\coordinate (R2) at (0.5, 1, 0);
	\coordinate (R3) at (0.75, 1, 0);
	\coordinate (R) at ($(R3)+0.5*(P)-0.5*(R3)$);
	\coordinate (v) at ($(P)+0.6*(Pt)-0.6*(P)$);
	%
	% top vertices: 
	\coordinate (Pt) at (0.5, \yy, 1);
	\coordinate (Rt) at (0.625, 0.5 + \yy/2, 1);
	\coordinate (Lt) at (0.5, 0, 1);
	\coordinate (R1t) at (0.25, 1, 1);
	\coordinate (R2t) at (0.5, 1, 1);
	\coordinate (R3t) at (0.75, 1, 1);
	\coordinate (aP) at (0.5, \yy, 1);
	\coordinate (aR) at (0.625, 0.5 + \yy/2, 1);
	\coordinate (aL) at (0.5, 0, 1);
	\coordinate (aR1) at (0.25, 1, 1);
	\coordinate (aR2) at (0.5, 1, 1);
	\coordinate (aR3) at (0.75, 1, 1);
	% top vertices: 
	\coordinate (aPt) at (0.5, \yy, 2);
	\coordinate (aRt) at (0.375, 0.5 + \yy/2, 2);
	\coordinate (aLt) at (0.5, 0, 2);
	\coordinate (aR1t) at (0.25, 1, 2);
	\coordinate (aR2t) at (0.5, 1, 2);
	\coordinate (aR3t) at (0.75, 1, 2);
	\fill [purple!80,opacity=0.2] (P) -- (aPt) -- (aLt) -- (L);
	\fill [orange!80,opacity=0.545] (P) -- (aPt) -- (aR1t) -- (R1);
	\fill [green!50,opacity=0.545] (R) -- (Pt) -- (aRt) -- (aR2t) -- (R2);
	%
	% triangle back
	\fill[red!80!black] (0.5, 0.38, 1.4) circle (0pt) node[right] (0up) {{\scriptsize$\triangleright_{M}$}};
	%
	% action line above: 
	\draw[color=red!80!black, ultra thick, rounded corners=0.5mm] (aRt) -- (Pt);
	\fill [orange!80,opacity=0.545] (P) -- (aPt) -- (aR3t) -- (R3);
	%
	% action line below: 
	\draw[blue!80!black, ultra thick] (P) -- (aPt);
	%
	% F-line: 
	\draw[color=red!80!black, ultra thick, rounded corners=0.5mm] (R) -- (Pt);
	%
	% T-label: 
	\fill[color=purple!50!black] (0.55, \yy, 0.15) circle (0pt) node[left] (0up) { {\scriptsize$T$} };
	\fill[blue!80!black] (0.45,0.2,0.3) circle (0pt) node[left] (0up) { {\scriptsize$F$} };
	\fill[blue!80!black] (0.45,0.2,1.7) circle (0pt) node[left] (0up) { {\scriptsize$G$} };
	\fill[red!80!black] (0.5,0.89,0.0) circle (0pt) node[left] (0up) { {\scriptsize$N$} };
	\fill[color=blue!60!black] (0.5,1.04,0.13) circle (0pt) node[left] (0up) { {\scriptsize$B$} };
	\fill[color=red!80!black] (0.5,1.19,0.28) circle (0pt) node[left] (0up) { {\scriptsize$M$} };
	%
	% vertices: 
	\fill[color=black] (Pt) circle (1.2pt) node[left] (0up) { {\scriptsize$g$} };
	\fill[color=black] (v) circle (1.2pt) node[left] (0up) { {\scriptsize$\zeta$} };
	%
	% triangle front
	\fill[red!80!black] (0.5, 0.38, 0.3) circle (0pt) node[right] (0up) {{\scriptsize$\triangleleft_{N}$}};
	%
	% black boundaries: 
	\draw [black,opacity=1, very thin] (L) -- (aLt) -- (aPt) -- (aR3t) -- (R3) -- (P) -- cycle;
	\draw [black,opacity=1, very thin] (aPt) -- (aR1t) -- (R1) -- (P);
	\draw [black,opacity=1, very thin] (R) -- (R2) -- (aR2t) -- (aRt);
	\end{tikzpicture}}%%popende
%%%%%%%%%%%%%%%%%%%%%% 
=
%%%%%%%%%%%%%%%%%%%%%% 
\tikzzbox{\begin{tikzpicture}[thick,scale=2.321,color=blue!50!black, baseline=1.4cm, >=stealth, 
	style={x={(-0.6cm,-0.4cm)},y={(1cm,-0.2cm)},z={(0cm,0.9cm)}}]
	%: where to put leftmost T-line: 
	\pgfmathsetmacro{\yy}{0.2}
	\coordinate (P) at (0.5, \yy, 0);
	\coordinate (L) at (0.5, 0, 0);
	\coordinate (R1) at (0.25, 1, 0);
	\coordinate (R2) at (0.5, 1, 0);
	\coordinate (R3) at (0.75, 1, 0);
	\coordinate (R) at ($(R3)+0.5*(P)-0.5*(R3)$);
	\coordinate (v) at ($(aPt)+0.6*(Pt)-0.6*(aPt)$);
	%
	% top vertices: 
	\coordinate (Pt) at (0.5, \yy, 1);
	\coordinate (Rt) at (0.625, 0.5 + \yy/2, 1);
	\coordinate (Lt) at (0.5, 0, 1);
	\coordinate (R1t) at (0.25, 1, 1);
	\coordinate (R2t) at (0.5, 1, 1);
	\coordinate (R3t) at (0.75, 1, 1);
	\coordinate (aP) at (0.5, \yy, 1);
	\coordinate (aR) at (0.625, 0.5 + \yy/2, 1);
	\coordinate (aL) at (0.5, 0, 1);
	\coordinate (aR1) at (0.25, 1, 1);
	\coordinate (aR2) at (0.5, 1, 1);
	\coordinate (aR3) at (0.75, 1, 1);
	% top vertices: 
	\coordinate (aPt) at (0.5, \yy, 2);
	\coordinate (aRt) at (0.375, 0.5 + \yy/2, 2);
	\coordinate (aLt) at (0.5, 0, 2);
	\coordinate (aR1t) at (0.25, 1, 2);
	\coordinate (aR2t) at (0.5, 1, 2);
	\coordinate (aR3t) at (0.75, 1, 2);
	\fill [purple!80,opacity=0.2] (P) -- (aPt) -- (aLt) -- (L);
	\fill [orange!80,opacity=0.545] (P) -- (aPt) -- (aR1t) -- (R1);
	\fill [green!50,opacity=0.545] (R) -- (Pt) -- (aRt) -- (aR2t) -- (R2);
	%
	% triangle back
	\fill[red!80!black] (0.5, 0.38, 1.4) circle (0pt) node[right] (0up) {{\scriptsize$\triangleright_{M}$}};
	%
	% action line above: 
	\draw[color=red!80!black, ultra thick, rounded corners=0.5mm] (aRt) -- (Pt);
	\fill [orange!80,opacity=0.545] (P) -- (aPt) -- (aR3t) -- (R3);
	%
	% action line below: 
	\draw[blue!80!black, ultra thick] (P) -- (aPt);
	%
	% F-line: 
	\draw[color=red!80!black, ultra thick, rounded corners=0.5mm] (R) -- (Pt);
	%
	% T-label: 
	\fill[color=purple!50!black] (0.55, \yy, 0.15) circle (0pt) node[left] (0up) { {\scriptsize$T$} };
	\fill[blue!80!black] (0.45,0.2,0.3) circle (0pt) node[left] (0up) { {\scriptsize$F$} };
	\fill[blue!80!black] (0.45,0.2,1.7) circle (0pt) node[left] (0up) { {\scriptsize$G$} };
	\fill[red!80!black] (0.5,0.89,0.0) circle (0pt) node[left] (0up) { {\scriptsize$N$} };
	\fill[color=blue!60!black] (0.5,1.04,0.13) circle (0pt) node[left] (0up) { {\scriptsize$B$} };
	\fill[color=red!80!black] (0.5,1.19,0.28) circle (0pt) node[left] (0up) { {\scriptsize$M$} };
	%
	% vertices: 
	\fill[color=black] (Pt) circle (1.2pt) node[left] (0up) { {\scriptsize$f$} };
	\fill[color=black] (v) circle (1.2pt) node[left] (0up) { {\scriptsize$\zeta$} };
	%
	% triangle front
	\fill[red!80!black] (0.5, 0.38, 0.3) circle (0pt) node[right] (0up) {{\scriptsize$\triangleleft_{N}$}};
	%
	% black boundaries: 
	\draw [black,opacity=1, very thin] (L) -- (aLt) -- (aPt) -- (aR3t) -- (R3) -- (P) -- cycle;
	\draw [black,opacity=1, very thin] (aPt) -- (aR1t) -- (R1) -- (P);
	\draw [black,opacity=1, very thin] (R) -- (R2) -- (aR2t) -- (aRt);
	\end{tikzpicture}}%%popende
%%%%%%%%%%%%%%%%%%%%%% 
\, . 
\ee 
Balanced 2-morphisms $M\btimes N\lra T$ and their 3-morphisms form a category $\operatorname{Bal}(M\btimes N,T)$. 
The relative product $M \relpro_{\mathcal B} N$ by definition comes with a 2-morphism $\Pi\in \operatorname{Bal}(M\btimes N, M \relpro_{\mathcal B} N)$, and one finds that the universal property says that 
\begin{align}
\Pi^* \colon \Hom_\tric \big( M\relpro_{\mathcal A} N, T \big) 
& 
\lra \textrm{Bal}\big( M\btimes N, T \big) 
\nonumber
\\ 
\varphi 
& 
\lmt \varphi \circ \Pi 
\end{align}
is an equivalence. 

\medskip 

The definition of the composition in terms of colimits allows us to construct
coherence isomorphisms by using the universal property. 
This will automatically ensure that the necessary coherence conditions are satisfied. 
We begin with the associator, which is given by 
\begin{equation}\label{eq:Defass}
\begin{tikzcd}[column sep=1em, row sep=3em]
& \EC(\mathcal C,\mathcal D)\times \EC(\mathcal B,\mathcal C) \times \EC(\mathcal A,\mathcal B) \ar[ld] \ar[rd]  & \\ 
\EC(\mathcal C,\mathcal D)^{\tau_2\Delta}\times \EC(\mathcal A,\mathcal C) \ar[rd] \ar[d, "\colim_{\tau_2 \Delta}",swap] & & \EC(\mathcal B,\mathcal D)^{\tau_2\Delta}\times \EC(\mathcal A,\mathcal B) \ar[ld] \ar[d, "\colim_{\tau_2 \Delta}"] \\ 
\EC(\mathcal C,\mathcal D)\times \EC(\mathcal A,\mathcal C) \ar[d] & \EC(\mathcal A,\mathcal D)^{\tau_2\Delta\times \tau_2\Delta} \ar[d, "\operatorname{colim}_{\tau_2\Delta\times \tau_2\Delta}"] \ar[rd, Rightarrow] & \EC(\mathcal B,\mathcal D)\times \EC(\mathcal A,\mathcal B) \ar[d]  \\ 
\EC(\mathcal A,\mathcal D)^{\tau_2\Delta} \ar[r, "\colim_{\tau_2 \Delta}",swap] \ar[ru, Rightarrow] & \EC(\mathcal A,\mathcal D)  & \EC(\mathcal A,\mathcal D)^{\tau_2\Delta} \ar[l, "\colim_{\tau_2 \Delta}"]
\end{tikzcd}
\end{equation}   
where the upper rhombus commutes strictly,\footnote{This is the case because the composition in $\tric$ is strict. Otherwise its associator would enter here.} and the lower subdiagrams commute because of the compatibility of the colimit with composition. 
The \textsl{unit 1-morphisms} are given by the 2-functor 
\begin{align}
1_{\mathcal A} \colon * &  \to \EC(\mathcal A,\mathcal A) \nonumber 
\\
* &  
\longmapsto \big(A,\mu_A, \mu_A, u^{\textrm{l}}_A, u^{\textrm{r}}_A, \alpha_A,\alpha_A,\alpha_A\big) \, . 
\end{align}
Next we construct the \textsl{unitor 2-morphisms} $l_{\mathcal A,\mathcal B}$ and $r_{\mathcal A,\mathcal B}$. For this consider
the natural transformation 
\begin{equation}
\begin{tikzcd}[column sep=1em, row sep=3em]
& \EC(\mathcal{B},\mathcal{B})\times \EC(\mathcal{A},\mathcal{B}) \ar[rd] \ar[d, Rightarrow, "{l}"] & \\ 
\EC(\mathcal{A},\mathcal{B}) \ar[ru, "1_\mathcal{B}\times \id_{\EC(\mathcal{A},\mathcal{B})}"] \ar[rr, "c_{-}", swap] & \  & \EC(\mathcal{A},\mathcal{B})^{\tau_2\Delta} 
\end{tikzcd}
\end{equation} 
where $c_{-}$ denotes the functor sending an object to the constant diagram, and $\tilde{l}$ is 
the natural transformation with component at $\mathcal{N}\in \EC(\mathcal{A},\mathcal{B})$ given by
\begin{equation}
\begin{tikzcd}
B\Box B \Box B \Box N \ar[dd, "\triangleright",swap] \ar[r]\ar[r, shift left=2.5] \ar[r, shift right=2.5] & B \Box B \Box N  \ar[dd, "\triangleright"] \ar[r, shift left]\ar[r, shift right] &  \ar[dd, "\triangleright"] B \Box N  \\ 
\\
N \ar[r]\ar[r, shift left=2.5] \ar[r, shift right=2.5] &  N \ar[r, shift left]\ar[r, shift right] &  N
\end{tikzcd}
\end{equation} 
where the vertical arrows denote the action of $B$ on $N$ and the coherence data for the action can be used to fill the squares. 
It is straightforward to see that the map $\triangleright \colon B\Box N \to N $ defines a universal cocone, which implies that the 2-colimit of the upper diagram is $N$ and hence that the map induced by $\tilde{l}$ on the level of 2-colimits is an equivalence. 
This allows us to define $l$ by
\begin{equation}
\begin{tikzcd}[column sep=1em, row sep=3em]
& \EC(\mathcal{B},\mathcal{B})\times \EC(\mathcal{A},\mathcal{B}) \ar[rd] \ar[d, Rightarrow, "{l}"] & 
\\ 
\EC(\mathcal{A},\mathcal{B}) \ar[rrd, "\id", bend right=15,swap] \ar[ru, "1_\mathcal{B}\times \id_{\EC(\mathcal{A},\mathcal{B})}"] \ar[rr, "c_{-}" near start, swap]  & \ 
\ar[d, Rightarrow, shorten >=2ex]  & \EC(\mathcal{A},\mathcal{B})^{\tau_2\Delta}  \ar[d, "\colim"] 
\\ 
& \ & \EC(\mathcal{A},\mathcal{B})
\end{tikzcd}
\end{equation} 
where the lower triangle is filled by a natural isomorphism which exists because the 2-colimit of
the constant diagram is the identity (this is not true for 2-colimits over constant diagrams in general). 
The unitor~$r$ is defined in a completely analogous way.  

To give the component of the \textsl{pentagonator} at  $(\mathcal{M},\mathcal{N},\mathcal{O},\mathcal{P})\in \EC(\mathcal{D},\mathcal{E})\times \EC(\mathcal{C},\mathcal{D})\times \EC(\mathcal{B},\mathcal{C})\times \EC(\mathcal{A},\mathcal{B})$, we use the following notational convention: 
if we write a composite like $ \mathcal{M} \relpro_{\mathcal D} \mathcal{N} \relpro_{\mathcal{C}} \mathcal{O}$ without any brackets we mean the colimit over the truncated bisimplicial object featuring also in~\eqref{eq:Defass}. 
We adopt analogous conventions for situations involving more bimodules. 
In this notation the associator is then the map
\begin{align}
\alpha_{\mathcal{M},\mathcal{N},\mathcal{O}}\colon (\mathcal{M} \relpro_{\mathcal D} \mathcal{N}) \relpro_{\mathcal{C}} \mathcal{O} \longrightarrow 
\mathcal{M} \relpro_{\mathcal D} \mathcal{N} \relpro_{\mathcal{C}} \mathcal{O} \longrightarrow 
\mathcal{M} \relpro_{\mathcal D} (\mathcal{N} \relpro_{\mathcal{C}} \mathcal{O})  
\end{align}   
where the morphisms are induced from the compatibility of composition with colimits and the Fubini theorem for colimits, cf.\ Appendix~\ref{app:2colimits}. 
The pentagonator is now given by (suppressing the symbols relative box products) 
\begin{equation}
\begin{tikzcd}[row sep=3em]
& & (\mathcal{M}  \mathcal{N})  (\mathcal{O}   \mathcal{P}) \ar[rrd, bend left=10] \ar[rd] \ar[dd,Rightarrow] & & \\ 
((\mathcal{M}  \mathcal{N})  \mathcal{O})  \mathcal{P}\ar[rd] \ar[rddd] \ar[rru, bend left=10]\ar[r] & (\mathcal{M}  \mathcal{N})  \mathcal{O}   \mathcal{P} \ar[ru] \ar[rd] \ar[d,Rightarrow] & & \mathcal{M}  \mathcal{N}  (\mathcal{O}   \mathcal{P}) \ar[r]  \ar[d,Rightarrow]& \mathcal{M}  (\mathcal{N}  (\mathcal{O}   \mathcal{P})) \\ 
& (\mathcal{M}  \mathcal{N}  \mathcal{O} )  \mathcal{P} \ar[dd] \ar[r] \ar[dr,Rightarrow] & \mathcal{M}  \mathcal{N}  \mathcal{O}   \mathcal{P} \ar[ru] \ar[r] \ar[d] & \mathcal{M}  (\mathcal{N}  \mathcal{O}   \mathcal{P})\ar[ru] \ar[dl,Rightarrow] & \\ 
& & \mathcal{M}  (\mathcal{N}  \mathcal{O})   \mathcal{P} \ar[rd]  & & \\ 
& (\mathcal{M}  (\mathcal{N}  \mathcal{O}))   \mathcal{P} \ar[rr]\ar[ru] & & \mathcal{M}  (\mathcal{N}  \mathcal{O} )  \mathcal{P}) \ar[ruuu] \ar[uu]
\end{tikzcd}
\end{equation}
where the triangles commute by definition and the squares commute because of the coherence of the Fubini theorem. 

Finally, there are three modifications $\lambda, \mu, \rho$ encoding the compatibility between the unit constraints $l,r$ and the associator~$\alpha$. 
The one most complicated to construct is $\mu$, the other two follow directly from the universal property of 2-colimits. 
The component of~$\mu$ at a pair of 1-morphisms
$(\mathcal{M},\mathcal{N})\in \EC(\mathcal{B},\mathcal{C}) \times \EC(\mathcal{A},\mathcal{B})$ is a 2-morphism 
\begin{equation}
\begin{tikzcd}[column sep=1em, row sep=3em]
(\mathcal{M} \relpro_\mathcal{B} \mathcal{B})  \relpro_\mathcal{B} \mathcal{N} \ar[rr, "\alpha"] \ar[rd, "r", swap] & \ar[d, Rightarrow, "\mu_{\mathcal{M},\mathcal{N}}"] & \mathcal{M} \relpro_\mathcal{B} (\mathcal{B}  \relpro_\mathcal{B} \mathcal{N}) \ar[ld, "l"] \\ 
& \mathcal{M}   \relpro_\mathcal{B} \mathcal{N}
\end{tikzcd}
\end{equation} 
which we define to be    
\begin{equation}
\label{eq:MuModification}
\begin{tikzcd}[column sep=1em, row sep=3em]
(\mathcal{M} \relpro_\mathcal{B} \mathcal{B})  \relpro_\mathcal{B} \mathcal{N} \ar[rr, "\alpha"] \ar[rd] \ar[rddd, "r", swap, bend right] &  & \mathcal{M} \relpro_\mathcal{B} (\mathcal{B}  \relpro_\mathcal{B} \mathcal{N}) \ar[lddd, "l", bend left] \ar[ld] \\ 
& \mathcal{M} \relpro_\mathcal{B} \mathcal{B}  \relpro_\mathcal{B} \mathcal{N} \ar[ld, Rightarrow, shorten >=30, shorten <=0] \ar[rd, Leftarrow, shorten >=30, shorten <=0] \ar[dd, bend right, "r'", swap] \ar[dd, bend left, "l'"] & \\ 
\ & & \ar[ll, "{\mu'}", Rightarrow, shorten >=68, shorten <=68] \ \\
& \mathcal{M}   \relpro_\mathcal{B} \mathcal{N}
\end{tikzcd}
\end{equation}where $r',l'$ are the maps induced by acting with $B$ to the left and right in the truncated bisimplical diagram composed with the projection to the relative tensor product. 
Since the projection is balanced we get an induced natural transformation between the two cocones which induces $\mu'$. 
The upper triangle in~\eqref{eq:MuModification} commutes strictly, the outer ones by the universal property of 2-colimits. 

This concludes the construction of the 3-category $\ealg(\tric)$.

\section{Orbifold completion of 3-categories}
\label{sec:OrbComp3cat}

In this section we introduce a categorification of the 2-dimensional orbifold completion of \cite{cr1210.6363}, as briefly reviewed in Section~\ref{subsec:2categories}. 
In Section~\ref{subsec:3catOrbData} we start out with a Gray category with duals~$\tric$ satisfying a mild finiteness condition, and construct a 3-category $\orb{\tric}$ as a subcategory of the 3-category $\ealg(\tric)$ of Section~\ref{sub:EofC}.
The defining conditions on the data of $\orb{\tric}$ are such that they encode invariance under oriented Pachner moves when used as labels for stratified 3-bordisms. 
Composition of 1-morphisms in $\orb{\tric}$ is given by the relative product, and we show how to compute it by splitting 
%arXiv_v2: 
	%higher idempotents
	 higher idempotents, called condensation monads in~\cite{GaiottoJohnsonFreyd}. 
Then in Section~\ref{subsubsec:Adjoints} we prove that 1- and 2-morphisms in $\orb{\tric}$ have compatible adjoints, while in Section~\ref{subsubsec:UniversalProperty} we identify a universal property for 2-dimensional orbifold completion, and outline its generalisation to arbitrary dimension.

\subsection{3-category of orbifold data}
\label{subsec:3catOrbData}

Let~$\tric$ be a Gray category with duals such that all Hom 2-categories admit finite sifted 2-limits and 2-colimits that commute with composition. 
The subcategory $\orb{\tric}$ of $\ealg(\tric)$ is defined by imposing further conditions on its objects and morphisms which as we will show below ensure the existence of adjoints and allows us to compute the relative box product by splitting a ``higher idempotent".

\medskip 

The conditions we want to impose on $\mathcal A \in \ealg(\tric)$ are those which make it a ``special orbifold datum'' in the sense of \cite{CRS1}, where the notion was introduced and studied as an algebraic structure that encodes invariance under Pachner moves in arbitrary dimension~$n$. Here we are interested in $n=3$, see \cite[Def.\,4.2]{CRS1}. 

\begin{definition}
	\label{def:OrbifoldDatum}
	An object $\mathcal A = (a,A,\mu,\alpha,u,u^{\textrm{l}},u^{\textrm{r}}) \in \ealg(\tric)$ is a \textsl{(3-dimensional, special) orbifold datum} in~$\tric$ if the conditions \eqref{eq:O1}--\eqref{eq:O8} in Figure~\ref{fig:OrbifoldDatumAxioms} are satisfied, where 
	\begin{align}\label{eq:alpha}
	\alpha' 
	:= 
	%%%%%%%%%%%%%%%%%%%%%%
	\tikzzbox{
}%%popende 
	%%%%%%%%%%%%%%%%%%%%%% 
	\, 
	\end{align}
are constructed from the associator $\alpha$ and the adjunction data for $\mu$. 
\end{definition}

\begin{figure}[!ht]
	\captionsetup[subfigure]{labelformat=empty}
	\centering
	\vspace{-60pt}
	\begin{subfigure}[b]{0.5\textwidth}
		\centering
		\be 
		\tag{O1}
		\label{eq:O1}
		%%%%%%%%%%%%%%%%%%%%%% 
		\tikzzbox{%
}%%popende
		%%%%%%%%%%%%%%%%%%%%%% 
		\ee
		%		\label{eq:O8}
	\end{subfigure}%\hfill\raisebox{0em}{(O8)}
	\vspace{0pt}
	\caption{Defining conditions on orbifold data $(a,A,\mu,\alpha,u,u^{\textrm{l}},u^{\textrm{r}})$, with $\overline\alpha := \alpha^{-1}$. 
		%arXiv_v3:
			Stripes indicate $\#$-adjoints.}
	\label{fig:OrbifoldDatumAxioms}
\end{figure}

We note that every $\mathcal A \in \ealg(\tric)$ satisfies the conditions \eqref{eq:O1}--\eqref{eq:O3} which simply express the pentagon axiom and the invertibility of~$\alpha$. 
The other conditions \eqref{eq:O4}--\eqref{eq:O8} are proper constraints on~$\mathcal A$ (but none of them involve the unit~$u$ or the unitors~$u^{\textrm{l}}, u^{\textrm{r}}$). 

%arXiv_v2: 
	\begin{remark}\label{Rem: Connection to seperable algebras}
		A \emph{rigid algebra} in a monoidal bicategory $\mathcal{B}$ is an $E_1$-algebra $(A,\mu,\alpha,u,u^{\textrm{l}},u^{\textrm{r}})$ in $\mathcal{B}$ where the multiplication $\mu$ admits a right adjoint $\mu^R\colon A  \to A \Box A $ which is equipped with a compatible structure of an $A$-$A$-bimodule morphism, see for example~\cite[Def.\,2.1.1]{Decoppet_Rig} for the precise definition. 
		A \emph{separable algebra} is a rigid algebra equipped with a compatible section $s^\mu \colon \id_M \to \mu \circ \mu^R $ of the counit $\epsilon^R$~\cite[Def.\,2.1.7]{Decoppet_Rig}. 
		%arXiv_v3: 
			(Such a separable algebra categorifies the notion of ``separable algebra'' in a monoidal 1-category, to which it reduces in the homotopy category of~$\mathcal B$.)
		  
		Every orbifold datum $(a,A,\mu,\alpha,u,u^{\textrm{l}},u^{\textrm{r}})$ is, in particular, a separable algebra in $\End_\tric(a)$, where the bimodule structure for $\mu^{\#}$ consists of the isomorphisms $\alpha'$ and $\alpha''$ and the section is $\coev_{\mu^{\#}}$. That this is a section follows from Equation~\eqref{eq:O8}. 
		We want to highlight that orbifold data are different from separable algebras. 
		The former strongly depend on the identification of left and right duals in a Gray category with adjoints, i.\,e.\ they are sensitive to a higher pivotal structure, whereas separable algebras can be defined without such an identification.             
	\end{remark}

\begin{remark}
	\begin{enumerate}[label={(\roman*)}]
		\item 
		\label{item:OrbifoldTQFTs}
		The original motivation to introduce $n$-dimensional orbifold data for arbitrary~$n$ (generalising the case $n=2$ of \cite{cr1210.6363}) was as follows (see \cite{CRS1} for details). 
		Given an $n$-dimensional defect TQFT~$\zz$, i.\,e.\ a symmetric monoidal functor on a category of stratified bordisms whose $j$-strata are labelled with elements in prescribed sets~$D_j$, one can construct a closed TQFT~$\zz_{\mathcal A}$ from defect labels $\mathcal A_j \in D_j$ if they satisfy certain conditions. 
		To wit, $\zz_{\mathcal A}$ is evaluated on any given (non-stratified) bordism by first choosing a ``good'' stratification, decorating its $j$-strata with~$\mathcal A_j$, evaluating with~$\zz$ and taking a colimit over all good stratifications. 
		The latter can be taken to be Poincar\'{e} duals of oriented triangulations (or the more convenient ``admissible skeleta'' of \cite{CMRSS1} for $n=3$). 
		For this procedure to be well-defined for $n=3$, we view the labels~$\mathcal A_j$ as $(3-j)$-morphisms in the 3-category~$\tric_{\zz}$ associated to~$\zz$ in \cite{CMS} and impose the conditions of Definition~\ref{def:OrbifoldDatum} with the identification $(a,A,\mu,\alpha) = (\mathcal A_3, \mathcal A_2, \mathcal A_1, \mathcal A_0)$. 
		\item 
		The definition of orbifold datum can be broadened to include additional data $\phi \in \textrm{Aut}(1_{1_a})$ and $\psi \in \textrm{Aut}(1_A)$. 
		This is often useful for applications, e.\,g.\ to Reshetikhin--Turaev theory in \cite{CRS3, CMRSS1, CMRSS2,MuleRunk, MuleRunk2, Mule1}. 
		As explained in \cite[Sect.\,4.2]{CRS3}, the technical complications involving $\phi,\psi$ are not relevant for the development of the general theory, hence for now we disregard them. 
		(In particular, the additional data $\phi,\psi$ are precisely included by passing to the ``Euler completion'' reviewed in Section~\ref{sec:ExamppleCompletion}.)
	\end{enumerate}
	\label{rem:OrbifoldTQFTs}
\end{remark} 

The orbifold procedure summarised in Remark~\ref{rem:OrbifoldTQFTs}\ref{item:OrbifoldTQFTs} produces ordinary closed TQFTs~$\zz_{\mathcal A}$ for any given defect TQFT~$\zz$. 
It is natural to allow for all defects between them and combine all~$\zz_{\mathcal A}$ into a single defect TQFT~$\orb{\zz}$. 
This was done in the 2-dimensional case in \cite{cr1210.6363}. 
Using the analysis of decomposition invariance of \cite{CRS1, CMRSS1} we are let to the following enhancement of the 3-category~$\ealg(\tric)$ of Section~\ref{sub:EofC}: 

\begin{definition}
	\label{def:Corb}
	The \textsl{orbifold completion} $\orb{\tric}$ of~$\tric$ is the sub-3-category of~$\ealg(\tric)$ 
	\begin{itemize}
		\item 
		whose objects are orbifold data, i.\,e.\ satisfy the conditions \eqref{eq:O1}--\eqref{eq:O8} in Figure~\ref{fig:OrbifoldDatumAxioms}, 
		\item 
		whose 1-morphisms 
		$
		\mathcal M 
		= 
		( M, \triangleright_M, \triangleleft_M, u_M^{\textrm{l}}, u_M^{\textrm{r}}, \alpha_M^{\textrm{l}},\alpha_M^{\textrm{m}}, \alpha_M^{\textrm{r}}) 
		$ 
		satisfy constraints analogous to those in~\eqref{eq:O1}--\eqref{eq:O7}, where some~$\alpha$ are replaced by $\alpha_M^{\textrm{l}},\alpha_M^{\textrm{m}}$ or~$\alpha_M^{\textrm{r}}$ and some~$\mu$ are replaced by~$\triangleright_M$ or~$\triangleleft_M$, as well as 
		\be 
		%%%%%%%%%%%%%%%%%%%%%% 
		\tikzzbox{\begin{tikzpicture}[thick,scale=1.961,color=blue!50!black, baseline=0.8cm, >=stealth, 
			style={x={(-0.6cm,-0.4cm)},y={(1cm,-0.2cm)},z={(0cm,0.9cm)}}]
			\coordinate (Lb) at (0, 0, 0);
			\coordinate (Rb) at (0, 1.5, 0);
			\coordinate (Rt) at (0, 1.5, 1);
			\coordinate (Lt) at (0, 0, 1);
			\coordinate (d) at (0, 0.45, 0.5);
			\coordinate (b) at (0, 1.05, 0.5);
			\fill [orange!80,opacity=0.545] (Lb) -- (Rb) -- (Rt) -- (Lt);
			\fill[color=red!80!black] (d) circle (0pt) node[left] (0up) { {\scriptsize$\triangleright_M$} };
			\fill[inner color=green!30!white,outer color=green!55!white, very thick, rounded corners=0.5mm] (d) .. controls +(0,0,0.4) and +(0,0,0.4) .. (b) -- (b) .. controls +(0,0,-0.4) and +(0,0,-0.4) .. (d);
			\draw[color=red!80!black, ultra thick, rounded corners=0.5mm, postaction={decorate}, decoration={markings,mark=at position .51 with {\arrow[draw=red!80!black]{>}}}] (d) .. controls +(0,0,0.4) and +(0,0,0.4) .. (b);
			\draw[color=red!80!black, ultra thick, rounded corners=0.5mm, postaction={decorate}, decoration={markings,mark=at position .53 with {\arrow[draw=red!80!black]{<}}}] (d) .. controls +(0,0,-0.4) and +(0,0,-0.4) .. (b);
			\fill[color=red!80!black] (0, 1.35, 0.15) circle (0pt) node (0up) { {\scriptsize$M$} };
			\fill[color=green!50!black] (0, 0.85, 0.4) circle (0pt) node (0up) { {\scriptsize$A$} };
			\fill[color=red!80!black, opacity=0.2] (0, 0.65, 0.6) circle (0pt) node (0up) { {\scriptsize$M$} };
			\end{tikzpicture}}%%popende
		%%%%%%%%%%%%%%%%%%%%%% 
		\, 
		=
		\,
		%%%%%%%%%%%%%%%%%%%%%% 
		\tikzzbox{\begin{tikzpicture}[thick,scale=1.961,color=blue!50!black, baseline=0.8cm, >=stealth, 
			style={x={(-0.6cm,-0.4cm)},y={(1cm,-0.2cm)},z={(0cm,0.9cm)}}]
			\coordinate (Lb) at (0, 0, 0);
			\coordinate (Rb) at (0, 1.5, 0);
			\coordinate (Rt) at (0, 1.5, 1);
			\coordinate (Lt) at (0, 0, 1);
			\coordinate (d) at (0, 0.45, 0.5);
			\coordinate (b) at (0, 1.05, 0.5);
			\fill [orange!80,opacity=0.545] (Lb) -- (Rb) -- (Rt) -- (Lt);
			\fill[color=red!80!black] (b) circle (0pt) node[right] (0up) { {\scriptsize$\triangleright_M$} };
			\fill[inner color=green!30!white,outer color=green!55!white, very thick, rounded corners=0.5mm] (d) .. controls +(0,0,0.4) and +(0,0,0.4) .. (b) -- (b) .. controls +(0,0,-0.4) and +(0,0,-0.4) .. (d);
			\draw[color=red!80!black, ultra thick, rounded corners=0.5mm, postaction={decorate}, decoration={markings,mark=at position .51 with {\arrow[draw=red!80!black]{<}}}] (d) .. controls +(0,0,0.4) and +(0,0,0.4) .. (b);
			\draw[color=red!80!black, ultra thick, rounded corners=0.5mm, postaction={decorate}, decoration={markings,mark=at position .53 with {\arrow[draw=red!80!black]{>}}}] (d) .. controls +(0,0,-0.4) and +(0,0,-0.4) .. (b);
			\fill[color=red!80!black] (0, 1.35, 0.15) circle (0pt) node (0up) { {\scriptsize$M$} };
			\fill[color=green!50!black] (0, 0.85, 0.4) circle (0pt) node (0up) { {\scriptsize$A$} };
			\fill[color=red!80!black, opacity=0.2] (0, 0.65, 0.6) circle (0pt) node (0up) { {\scriptsize$M$} };
			\end{tikzpicture}}%%popende
		%%%%%%%%%%%%%%%%%%%%%% 
		\, 
		=
		\, 
		%%%%%%%%%%%%%%%%%%%%%% 
		\tikzzbox{\begin{tikzpicture}[thick,scale=1.961,color=blue!50!black, baseline=0.8cm, >=stealth, 
			style={x={(-0.6cm,-0.4cm)},y={(0.5cm,-0.2cm)},z={(0cm,0.9cm)}}]
			\coordinate (Lb) at (0, 0, 0);
			\coordinate (Rb) at (0, 1.5, 0);
			\coordinate (Rt) at (0, 1.5, 1);
			\coordinate (Lt) at (0, 0, 1);
			\coordinate (d) at (0, 0.45, 0.5);
			\coordinate (b) at (0, 1.05, 0.5);
			\fill [orange!80,opacity=0.545] (Lb) -- (Rb) -- (Rt) -- (Lt);
			\fill[color=red!80!black] (0, 1.25, 0.15) circle (0pt) node (0up) { {\scriptsize$M$} };
			\end{tikzpicture}}%%popende
		%%%%%%%%%%%%%%%%%%%%%% 
		\,
		=
		\, 
		%%%%%%%%%%%%%%%%%%%%%% 
		\tikzzbox{\begin{tikzpicture}[thick,scale=1.961,color=blue!50!black, baseline=0.8cm, >=stealth, 
			style={x={(-0.6cm,-0.4cm)},y={(1cm,-0.2cm)},z={(0cm,0.9cm)}}]
			\coordinate (Lb) at (0, 0, 0);
			\coordinate (Rb) at (0, 1.5, 0);
			\coordinate (Rt) at (0, 1.5, 1);
			\coordinate (Lt) at (0, 0, 1);
			\coordinate (d) at (0, 0.45, 0.5);
			\coordinate (b) at (0, 1.05, 0.5);
			\fill[inner color=green!30!white,outer color=green!55!white, very thick, rounded corners=0.5mm] (d) .. controls +(0,0,0.4) and +(0,0,0.4) .. (b) -- (b) .. controls +(0,0,-0.4) and +(0,0,-0.4) .. (d);
			\fill [orange!80,opacity=0.545] (Lb) -- (Rb) -- (Rt) -- (Lt);
			\fill[color=red!80!black] (d) circle (0pt) node[left] (0up) { {\scriptsize$\triangleleft_M$} };
			\draw[color=red!80!black, ultra thick, rounded corners=0.5mm, postaction={decorate}, decoration={markings,mark=at position .51 with {\arrow[draw=red!80!black]{>}}}] (d) .. controls +(0,0,0.4) and +(0,0,0.4) .. (b);
			\draw[color=red!80!black, ultra thick, rounded corners=0.5mm, postaction={decorate}, decoration={markings,mark=at position .53 with {\arrow[draw=red!80!black]{<}}}] (d) .. controls +(0,0,-0.4) and +(0,0,-0.4) .. (b);
			\fill[color=red!80!black] (0, 1.35, 0.15) circle (0pt) node (0up) { {\scriptsize$M$} };
			\fill[color=green!50!black, opacity=0.4] (0, 0.85, 0.4) circle (0pt) node (0up) { {\scriptsize$B$} };
			\fill[color=red!80!black] (0, 0.65, 0.6) circle (0pt) node (0up) { {\scriptsize$M$} };
			\end{tikzpicture}}%%popende
		%%%%%%%%%%%%%%%%%%%%%% 
		\, 
		=
		\,
		%%%%%%%%%%%%%%%%%%%%%% 
		\tikzzbox{\begin{tikzpicture}[thick,scale=1.961,color=blue!50!black, baseline=0.8cm, >=stealth, 
			style={x={(-0.6cm,-0.4cm)},y={(1cm,-0.2cm)},z={(0cm,0.9cm)}}]
			\coordinate (Lb) at (0, 0, 0);
			\coordinate (Rb) at (0, 1.5, 0);
			\coordinate (Rt) at (0, 1.5, 1);
			\coordinate (Lt) at (0, 0, 1);
			\coordinate (d) at (0, 0.45, 0.5);
			\coordinate (b) at (0, 1.05, 0.5);
			\fill[inner color=green!30!white,outer color=green!55!white, very thick, rounded corners=0.5mm] (d) .. controls +(0,0,0.4) and +(0,0,0.4) .. (b) -- (b) .. controls +(0,0,-0.4) and +(0,0,-0.4) .. (d);
			\fill [orange!80,opacity=0.545] (Lb) -- (Rb) -- (Rt) -- (Lt);
			\fill[color=red!80!black] (b) circle (0pt) node[right] (0up) { {\scriptsize$\triangleleft_M$} };
			\draw[color=red!80!black, ultra thick, rounded corners=0.5mm, postaction={decorate}, decoration={markings,mark=at position .51 with {\arrow[draw=red!80!black]{<}}}] (d) .. controls +(0,0,0.4) and +(0,0,0.4) .. (b);
			\draw[color=red!80!black, ultra thick, rounded corners=0.5mm, postaction={decorate}, decoration={markings,mark=at position .53 with {\arrow[draw=red!80!black]{>}}}] (d) .. controls +(0,0,-0.4) and +(0,0,-0.4) .. (b);
			\fill[color=red!80!black] (0, 1.35, 0.15) circle (0pt) node (0up) { {\scriptsize$M$} };
			\fill[color=green!50!black, opacity=0.4] (0, 0.85, 0.4) circle (0pt) node (0up) { {\scriptsize$B$} };
			\fill[color=red!80!black] (0, 0.65, 0.6) circle (0pt) node (0up) { {\scriptsize$M$} };
			\end{tikzpicture}}%%popende
		%%%%%%%%%%%%%%%%%%%%%% 
		\, ,
		\ee 
		\item 
		whose 2-morphisms satisfy the conditions \eqref{eq:T1}--\eqref{eq:T7} in Figures~\ref{fig:2MorphismsInCorb} and~\ref{fig:2MorphismsInCorbPart2}, 
		\item 
		and whose 3-morphisms $\xi\colon (F,\triangleleft_F, \triangleright_F) \lra (G,\triangleleft_{G}, \triangleright_{G})$ need not satisfy any other conditions but those in~\eqref{eq:3MorphismsInEAlg}. 
	\end{itemize}
\end{definition}

We note that the conditions on 1-morphisms~$\mathcal M$ in~$\orb{\tric}$ analogous to~\eqref{eq:O1}--\eqref{eq:O3} are precisely those in~\eqref{eq:23MovesFor1Morphisms1}--\eqref{eq:23MovesFor1Morphisms4} together with the invertibility of $\alpha_M^{\textrm{l}},\alpha_M^{\textrm{m}}, \alpha_M^{\textrm{r}}$. 
The conditions analogous to~\eqref{eq:O4}--\eqref{eq:O7} are indeed obtained in complete analogy; for example, one of the eight relations corresponding to~\eqref{eq:O6} is 
\be 
%%%%%%%%%%%%%%%%%%%%%% 
\tikzzbox{\begin{tikzpicture}[thick,scale=2.321,color=blue!50!black, baseline=1.2cm, >=stealth, 
	style={x={(-0.6cm,-0.4cm)},y={(1cm,-0.2cm)},z={(0cm,0.9cm)}}]
	%: where to put leftmost T-line: 
	\pgfmathsetmacro{\yy}{0.2}
	\coordinate (P) at (0.5, \yy, 0);
	\coordinate (R) at (0.375, 0.5 + \yy/2, 0);
	\coordinate (L) at (0.5, 0, 0);
	\coordinate (R1) at (0.25, 1, 0);
	\coordinate (R2) at (0.5, 1, 0);
	\coordinate (R3) at (0.5 + 2/12, \yy + -2*0.45/3, 0);
	% top vertices: 
	\coordinate (Pt) at (0.5, \yy, 1);
	\coordinate (Rt) at (0.625, 0.5 + \yy/2, 1);
	\coordinate (Lt) at (0.5, 0, 1);
	\coordinate (R1t) at (0.25, 1, 1);
	\coordinate (R2t) at (0.5, 1, 1);
	\coordinate (R3t) at (0.5 + 2/12, \yy + -2*0.45/3, 1);
	%alpha: 
	\coordinate (alpha) at (0.5, 0.5, 0.5);
	\coordinate (aP) at (0.5, \yy, 1);
	\coordinate (aR) at (0.625, 0.5 + \yy/2, 1);
	\coordinate (aL) at (0.5, 0, 1);
	\coordinate (aR1) at (0.25, 1, 1);
	\coordinate (aR2) at (0.5, 1, 1);
	\coordinate (aR3) at (0.5 + 2/12, \yy + -2*0.45/3, 1);
	% top vertices: 
	\coordinate (aPt) at (0.5, \yy, 2);
	\coordinate (aRt) at (0.375, 0.5 + \yy/2, 2);
	\coordinate (aLt) at (0.5, 0, 2);
	\coordinate (aR1t) at (0.25, 1, 2);
	\coordinate (aR2t) at (0.5, 1, 2);
	\coordinate (aR3t) at (0.5 + 2/12, \yy + -2*0.45/3, 2);
	%alpha: 
	\coordinate (aalpha) at (0.5, 0.5, 1.5);
	%
	% A-planes: 
	\fill [orange!80,opacity=0.545] (L) -- (P) -- (alpha) -- (Pt) -- (aalpha) -- (aPt) -- (aLt);
	\fill [orange!80,opacity=0.545] (R1) -- (R) -- (alpha) -- (Pt) -- (aalpha) -- (aRt) -- (aR1t);
	%%%%%%%%%%%%%%
	\fill [green!50,opacity=0.545] (aP) -- (aalpha) -- (aR) -- (alpha);
	\fill[color=blue!60!black] (0.5,0.55,0.84) circle (0pt) node[left] (0up) { {\scriptsize$A$} };
	%%%%%%%%%%%%%%%
	\draw[color=red!80!black, ultra thick, rounded corners=0.5mm, postaction={decorate}, decoration={markings,mark=at position .37 with {\arrow[draw=red!80!black]{>}}}] (alpha) -- (aP) -- (aalpha); 
	\fill [green!50,opacity=0.545] (R3) -- (P) -- (alpha) -- (Rt) -- (aalpha) -- (aPt) -- (aR3t);
	\fill [orange!80,opacity=0.545] (P) -- (alpha) -- (R);
	\fill [orange!80,opacity=0.545] (aPt) -- (aalpha) -- (aRt);
	\fill [green!50,opacity=0.545] (R2) -- (R) -- (alpha) -- (Rt) -- (aalpha) -- (aRt) -- (aR2t);
	\fill[color=red!80!black] (0.5,0.9,0.3) circle (0pt) node[left] (0up) { {\scriptsize$\triangleright_{M}$} };
	\fill[color=red!80!black] (0.5,0.55,0.19) circle (0pt) node[left] (0up) { {\scriptsize$M$} };
	\draw[string=red!80!black, ultra thick] (R) -- (alpha);
	%
	% T-lines: 
	\draw[color=red!80!black, ultra thick, rounded corners=0.5mm, postaction={decorate}, decoration={markings,mark=at position .51 with {\arrow[draw=red!80!black]{<}}}] (P) -- (alpha);
	%
	% labels: 
	\fill[color=green!30!black] (alpha) circle (1.2pt) node[left] (0up) { {\scriptsize$(\alpha^{\textrm{r}}_M)''$} };
	% A labels: 
	\fill[color=red!80!black] (0.5,0.2,0.15) circle (0pt) node[left] (0up) { {\scriptsize$M$} };
	\fill[color=blue!60!black] (0.6, -0.06, 0.08) circle (0pt) node (0up) { {\scriptsize$A$} };
	\fill[color=blue!60!black] (0.5,1.04,0.13) circle (0pt) node[left] (0up) { {\scriptsize$A$} };
	\fill[color=red!80!black] (0.5,1.19,0.28) circle (0pt) node[left] (0up) { {\scriptsize$M$} };
	\fill[color=red!80!black] (0.5,0.58,1.95) circle (0pt) node[left] (0up) { {\scriptsize$M$} };
	% T labels: 
	\fill[color=red!80!black] (0.5,0.37,0.25) circle (0pt) node[left] (0up) { {\scriptsize$\triangleright_{M}$} };
	\fill[color=red!80!black] (0.5,0.4,0.71) circle (0pt) node[left] (0up) { {\scriptsize$\triangleright_{M}$} };
	%
	%
	%T-line: 
	\draw[color=red!80!black, ultra thick, rounded corners=0.5mm, postaction={decorate}, decoration={markings,mark=at position .51 with {\arrow[draw=red!80!black]{>}}}] (aalpha) -- (aRt);
	%
	% T-lines: 
	\draw[color=green!30!black, ultra thick, rounded corners=0.5mm, postaction={decorate}, decoration={markings,mark=at position .57 with {\arrow[draw=green!30!black]{<}}}] (alpha) -- (aR) -- (aalpha); 
	\draw[color=red!80!black, ultra thick, rounded corners=0.5mm, postaction={decorate}, decoration={markings,mark=at position .51 with {\arrow[draw=red!80!black]{<}}}] (aalpha) -- (aPt);
	%
	% labels: 
	\fill[color=green!30!black] (aalpha) circle (1.2pt) node[left] (0up) { {\scriptsize$(\bar\alpha^{\textrm{r}}_M)''$} }; 
	\fill[color=green!30!black] (0.5,0.73,1.07) circle (0pt) node[left] (0up) { {\scriptsize$\mu$} };
	%
	% black boundaries: 
	\draw [black,opacity=1, very thin] (aPt) -- (aLt) -- (L) -- (P);
	\draw [black,opacity=1, very thin] (aPt) -- (aR1t) -- (R1) -- (R);
	\draw [black,opacity=1, very thin] (aRt) -- (aR2t) -- (R2) -- (R);
	\draw [black,opacity=1, very thin] (P) -- (R3) -- (aR3t) -- (aPt);
	\draw [black,opacity=1, very thin] (P) -- (R);
	%
	% black auxiliary boundaries: 
	\draw [black,opacity=0.4141, densely dotted, semithick] (Lt) -- (Pt);
	\draw [black,opacity=0.4141, densely dotted, semithick] (Pt) -- (0.375, 0.5 + \yy/2, 1) -- (R1t);
	\draw [black,opacity=0.4141, densely dotted, semithick] (Rt) -- (R2t);
	\draw [black,opacity=0.4141, densely dotted, semithick] (Rt) -- (R3t);
	\draw [black,opacity=0.4141, densely dotted, semithick] (Rt) -- (Pt);
	\end{tikzpicture}}%%popende
%%%%%%%%%%%%%%%%%%%%%% 
= 
%%%%%%%%%%%%%%%%%%%%%% 
\tikzzbox{\begin{tikzpicture}[thick,scale=2.321,color=blue!50!black, baseline=1.2cm, >=stealth, 
	style={x={(-0.6cm,-0.4cm)},y={(1cm,-0.2cm)},z={(0cm,0.9cm)}}]
	%: where to put leftmost T-line: 
	\pgfmathsetmacro{\yy}{0.2}
	\coordinate (P) at (0.5, \yy, 0);
	\coordinate (R) at (0.375, 0.5 + \yy/2, 0);
	\coordinate (L) at (0.5, 0, 0);
	\coordinate (R1) at (0.25, 1, 0);
	\coordinate (R2) at (0.5, 1, 0);
	\coordinate (R3) at (0.5 + 2/12, \yy + -2*0.45/3, 0); %(0.75, -0.25, 0);
	% top vertices: 
	\coordinate (Pt) at (0.5, \yy, 1);
	\coordinate (Rt) at (0.625, 0.5 + \yy/2, 1);
	\coordinate (Lt) at (0.5, 0, 1);
	\coordinate (R1t) at (0.25, 1, 1);
	\coordinate (R2t) at (0.5, 1, 1);
	\coordinate (R3t) at (0.5 + 2/12, \yy + -2*0.45/3, 1);
	\coordinate (aP) at (0.5, \yy, 1);
	\coordinate (aR) at (0.625, 0.5 + \yy/2, 1);
	\coordinate (aL) at (0.5, 0, 1);
	\coordinate (aR1) at (0.25, 1, 1);
	\coordinate (aR2) at (0.5, 1, 1);
	\coordinate (aR3) at (0.5 + 2/12, \yy + -2*0.45/3, 1);
	% top vertices: 
	\coordinate (aPt) at (0.5, \yy, 2);
	\coordinate (aRt) at (0.375, 0.5 + \yy/2, 2);
	\coordinate (aLt) at (0.5, 0, 2);
	\coordinate (aR1t) at (0.25, 1, 2);
	\coordinate (aR2t) at (0.5, 1, 2);
	\coordinate (aR3t) at (0.5 + 2/12, \yy + -2*0.45/3, 2);
	\fill [orange!80,opacity=0.545] (P) -- (aPt) -- (aLt) -- (L);
	\fill [orange!80,opacity=0.545] (P) -- (aPt) -- (aR1t) -- (R1);
	\fill [green!50,opacity=0.545] (R) -- (aRt) -- (aR2t) -- (R2);
	\fill[red!80!black] (0.33, 0.5 + \yy/2, 0.94) circle (0pt) node[left] (0up) { {\scriptsize$\triangleright_{M}$} };
	\draw[color=red!80!black, ultra thick, rounded corners=0.5mm, postaction={decorate}, decoration={markings,mark=at position .51 with {\arrow[draw=red!80!black]{>}}}] (R) -- (aRt);
	\fill[color=red!80!black] (0.5,0.56,0.2) circle (0pt) node[left] (0up) { {\scriptsize$M$} };
	\fill [green!50,opacity=0.545] (P) -- (aPt) -- (aR3t) -- (R3);
	\draw[color=red!80!black, ultra thick, rounded corners=0.5mm, postaction={decorate}, decoration={markings,mark=at position .51 with {\arrow[draw=red!80!black]{<}}}] (P) -- (aPt);
	\fill[color=red!80!black] (0.5, \yy, 0.94) circle (0pt) node[right] (0up) { {\scriptsize$\triangleright_{M}$} };
	\fill[color=red!80!black] (0.5,0.2,0.15) circle (0pt) node[left] (0up) { {\scriptsize$M$} };
	\fill[color=blue!60!black] (0.6, -0.06, 0.08) circle (0pt) node (0up) { {\scriptsize$A$} };
	\fill[color=blue!60!black] (0.5,1.04,0.13) circle (0pt) node[left] (0up) { {\scriptsize$A$} };
	\fill[color=red!80!black] (0.5,1.19,0.28) circle (0pt) node[left] (0up) { {\scriptsize$M$} };
	%
	% black boundaries: 
	\draw [black,opacity=1, very thin] (L) -- (aLt) -- (aPt) -- (aR3t) -- (R3) -- (P) -- cycle;
	\draw [black,opacity=1, very thin] (aPt) -- (aR1t) -- (R1) -- (P);
	\draw [black,opacity=1, very thin] (R) -- (R2) -- (aR2t) -- (aRt);
	%
	% black auxiliary boundaries: 
	\draw [black,opacity=0.4141, densely dotted, semithick] (Lt) -- (Pt) -- (R1t);
	\draw [black,opacity=0.4141, densely dotted, semithick] (Pt) -- (R3t);
	\draw [black,opacity=0.4141, densely dotted, semithick] (R2t) -- (0.375, 0.5 + \yy/2, 1);
	\end{tikzpicture}}%%popende
%%%%%%%%%%%%%%%%%%%%%%
\, . 
\ee 
In Proposition~\ref{prop:1MorphismsOrb} we give a simpler condition which implies that the additional conditions are satisfied.  

We also observe that the conditions \eqref{eq:T1}--\eqref{eq:T5} on 2-morphisms~$\mathcal F$ already hold in~$\ealg(\tric)$, see Definition~\ref{def:2MorphismsEAlg}, while the conditions~\eqref{eq:T6} and~\eqref{eq:T7} are proper constraints on~$\mathcal F$ to be in~$\orb{\tric}$. 
1-morphisms in $\orb{\tric}$ are closed under composition, given by relative products $\mathcal M \relpro_{\mathcal B} \mathcal N$ as discussed around~\eqref{eq:RelativeGraphicalNotation}, thanks to the definition in terms of a universal property: 
the structure maps of $\mathcal M \relpro_{\mathcal B} \mathcal N$ with respect to the actions of~$\mathcal A$ and~$\mathcal C$ are induced from those of~$\mathcal N$ and~$\mathcal M$, respectively, and are only spectators in the construction of $\mathcal M \relpro_{\mathcal B} \mathcal N$. 

\begin{figure}[!ht]
	\captionsetup[subfigure]{labelformat=empty}
	\centering
%	\vspace{-60pt}
	\hspace{7pt}
	\begin{subfigure}[b]{0.5\textwidth}
		\centering
		\be 
		\tag{T1}
		\label{eq:T1}
		%%%%%%%%%%%%%%%%%%%%%% 
		\tikzzbox{
}%%popende
		%%%%%%%%%%%%%%%%%%%%%%  
		\ee
	\end{subfigure}%\hfill\raisebox{0em}{(O8)}
	\vspace{0pt}
	\caption{Defining conditions on 2-morphisms $\mathcal F\colon \mathcal M \lra \mathcal M'$ in $\tric_{\textrm{orb}}$ (where $\mathcal M, \mathcal M'\colon \mathcal A \lra \mathcal B$), part 1/2.}
	\label{fig:2MorphismsInCorb}
\end{figure}

\begin{figure}[!ht]
	\captionsetup[subfigure]{labelformat=empty}
	\centering
	\vspace{-5pt}
	%	\hspace{-60pt}
	\begin{subfigure}[b]{1.0\textwidth}
		\centering
		\be
		\tag{T6}
		\label{eq:T6}
		%%%%%%%%%%%%%%%%%%%%%% 
		\tikzzbox{
}%%popende
		%%%%%%%%%%%%%%%%%%%%%%
		\ee
	\end{subfigure}%\hfill\raisebox{0em}{(O8)}
%	\\
%	%\vspace{-5pt}
%	%	\hspace{-60pt}
%	\begin{subfigure}[b]{1.0\textwidth}
%		\centering
%		%		\label{eq:O8}
%	\end{subfigure}
	\vspace{0pt}
	\caption{Defining conditions on 2-morphisms $\mathcal F\colon \mathcal M \lra \mathcal M'$ in $\tric_{\textrm{orb}}$ (where $\mathcal M, \mathcal M'\colon \mathcal A \lra \mathcal B$), part 2/2.}
	\label{fig:2MorphismsInCorbPart2}
\end{figure}

\begin{remark}
	The definition of~$\orb{\tric}$ is such that if $\tric = \tric_{\zz}$ is the 3-category associated to a defect TQFT~$\zz$ as in Remark~\ref{rem:OrbifoldTQFTs}\ref{item:OrbifoldTQFTs}, then the $k$-morphisms in $\orb{(\tric_{\zz})}$ can be used as the $(3-k)$-dimensional defects of the orbifold defect TQFT $\orb{\zz}$. 
	Moreover, one then finds 
	$\orb{(\tric_{\zz})} \cong \tric_{\orb{\zz}}$. 
	This will be discussed in more detail in Section~\ref{sec:TQFTs}. 
\end{remark}

The conditions \eqref{eq:T1}--\eqref{eq:T7} on 1-morphisms are needed to ensure that~$\orb{\tric}$ has adjoints, as we will see in Section~\ref{subsec:properties}. 
This gives a second, purely algebraic motivation for Definition~\ref{def:Corb}. 

\medskip

It is straightforward to check that an orbifold datum~$\mathcal A$ in~$\tric$ induces a Frobenius algebra in the homotopy 2-category of~$\tric$. 
Indeed, there are ``Frobeniusator'' 3-isomorphisms 
\be 
\begin{tikzcd}[column sep=5em, row sep=2em]
%%%%%%%%%%%%%%%%%%%%%%
\tikzzbox{\begin{tikzpicture}[ultra thick,scale=0.75]
	\coordinate (r0) at (0,-0.75);
	\coordinate (r1) at (0,0.75);
	\coordinate (l0) at (-3,-0.75);
	\coordinate (l1) at (-3,0.75);
	\coordinate (m) at (-2,-0.75);
	\coordinate (md) at (-1,0.75);
	%
	% lines: 
	\draw[color=green!70!black,opacity=0.545] (l0) -- (r0);
	\draw[color=green!70!black,opacity=0.545] (l1) -- (r1);
	\draw[color=green!70!black,opacity=0.545] (m) -- (md);
	%
	% vertices: 
	\fill[color=green!30!black] (r0) circle (0pt) node[color=green!30!black, right, font=\footnotesize] {$A$};
	\fill[color=green!30!black] (r1) circle (0pt) node[color=green!30!black, right, font=\footnotesize] {$A$};
	\fill[color=green!30!black] (l0) circle (0pt) node[color=green!30!black, left, font=\footnotesize] {$A$};
	\fill[color=green!30!black] (l1) circle (0pt) node[color=green!30!black, left, font=\footnotesize] {$A$};
	\fill[color=green!30!black] (m) circle (3pt) node[color=green!30!black, below, font=\footnotesize] {$\mu_A$};
	\fill[color=green!30!black] (md) circle (3pt) node[color=green!30!black, above, font=\footnotesize] {$\mu^\vee_A$};
	\end{tikzpicture}}%%popende 
%%%%%%%%%%%%%%%%%%%%%%   
\arrow{r}{\overline{\alpha}''}
\arrow{r}[swap]{\cong}
& 
%%%%%%%%%%%%%%%%%%%%%%
\tikzzbox{\begin{tikzpicture}[ultra thick,scale=0.75]
	\coordinate (r0) at (0,-0.75);
	\coordinate (r1) at (0,0.75);
	\coordinate (l0) at (-3,-0.75);
	\coordinate (l1) at (-3,0.75);
	\coordinate (m) at (-1,0);
	\coordinate (md) at (-2,0);
	%
	% lines: 
	\draw[color=green!70!black,opacity=0.545] (r0) -- (m);
	\draw[color=green!70!black,opacity=0.545] (r1) -- (m);
	\draw[color=green!70!black,opacity=0.545] (l0) -- (md);
	\draw[color=green!70!black,opacity=0.545] (l1) -- (md);
	\draw[color=green!70!black,opacity=0.545] (m) -- (md);
	%
	% vertices: 
	\fill[color=green!30!black] (r0) circle (0pt) node[color=green!30!black, right, font=\footnotesize] {$A$};
	\fill[color=green!30!black] (r1) circle (0pt) node[color=green!30!black, right, font=\footnotesize] {$A$};
	\fill[color=green!30!black] (l0) circle (0pt) node[color=green!30!black, left, font=\footnotesize] {$A$};
	\fill[color=green!30!black] (l1) circle (0pt) node[color=green!30!black, left, font=\footnotesize] {$A$};
	\fill[color=green!30!black] (m) circle (3pt) node[color=green!30!black, above, font=\footnotesize] {$\mu_A$};
	\fill[color=green!30!black] (md) circle (3pt) node[color=green!30!black, below, font=\footnotesize] {$\mu^\vee_A$};
	\end{tikzpicture}}%%popende 
%%%%%%%%%%%%%%%%%%%%%%
\arrow{r}{\alpha''}
\arrow{r}[swap]{\cong}
& 
%%%%%%%%%%%%%%%%%%%%%%
\tikzzbox{\begin{tikzpicture}[ultra thick,scale=0.75]
	\coordinate (r0) at (0,-0.75);
	\coordinate (r1) at (0,0.75);
	\coordinate (l0) at (-3,-0.75);
	\coordinate (l1) at (-3,0.75);
	\coordinate (m) at (-2,0.75);
	\coordinate (md) at (-1,-0.75);
	%
	% lines: 
	\draw[color=green!70!black,opacity=0.545] (l0) -- (r0);
	\draw[color=green!70!black,opacity=0.545] (l1) -- (r1);
	\draw[color=green!70!black,opacity=0.545] (m) -- (md);
	%
	% vertices: 
	\fill[color=green!30!black] (r0) circle (0pt) node[color=green!30!black, right, font=\footnotesize] {$A$};
	\fill[color=green!30!black] (r1) circle (0pt) node[color=green!30!black, right, font=\footnotesize] {$A$};
	\fill[color=green!30!black] (l0) circle (0pt) node[color=green!30!black, left, font=\footnotesize] {$A$};
	\fill[color=green!30!black] (l1) circle (0pt) node[color=green!30!black, left, font=\footnotesize] {$A$};
	\fill[color=green!30!black] (m) circle (3pt) node[color=green!30!black, above, font=\footnotesize] {$\mu_A$};
	\fill[color=green!30!black] (md) circle (3pt) node[color=green!30!black, below, font=\footnotesize] {$\mu^\vee_A$};
	\end{tikzpicture}}%%popende 
%%%%%%%%%%%%%%%%%%%%%%    
\, . 
\end{tikzcd}
\ee 
With this one finds that every $\mathcal B$-$\mathcal A$-bimodule~$\mathcal M$, i.\,e.\ 1-morphism $\mathcal M \colon \mathcal A \lra \mathcal B$ in $\ealg(\tric)$, also has the structure of a $\mathcal B$-$\mathcal A$-bicomodule with coactions 
\be 
%%%%%%%%%%%%%%%%%%%%%%
\tikzzbox{\begin{tikzpicture}[ultra thick,scale=0.75,baseline=-0.1cm]
	\coordinate (ml) at (-3,0);
	\coordinate (mr) at (0,0);
	\coordinate (ac) at (-1.5,0);
	\coordinate (A) at (-3,-1);
	%
	% lines: 
	\draw[color=orange!80,opacity=0.545] (ml) -- (mr);
	\draw[color=green!70!black,opacity=0.545] (ac) -- (A);
	%
	% vertices: 
	\fill[color=green!30!black] (mr) circle (0pt) node[color=red!80!black, right, font=\footnotesize] {$M$};
	\fill[color=green!30!black] (ml) circle (0pt) node[color=red!80!black, left, font=\footnotesize] {$M$};
	\fill[color=green!30!black] (A) circle (0pt) node[color=green!30!black, left, font=\footnotesize] {$A$};
	green!30!black
	\fill[color=red!80!black] (ac) circle (3pt) node[color=red!80!black, above, font=\footnotesize] {${\triangleright_M}^\dagger$};
	\end{tikzpicture}}%%popende 
%%%%%%%%%%%%%%%%%%%%%%  
:= 
%%%%%%%%%%%%%%%%%%%%%%
\tikzzbox{\begin{tikzpicture}[ultra thick,scale=0.75,baseline=-0.1cm]
	\coordinate (ml) at (-3,0);
	\coordinate (mr) at (0,0);
	\coordinate (ac) at (-1.5,0);
	\coordinate (A) at (-3,-1);
	\coordinate (md) at (-1,-0.5);
	\coordinate (u) at (-0.3,-0.5);
	%
	% lines: 
	\draw[color=orange!80,opacity=0.545] (ml) -- (mr);
	\draw[color=green!70!black,opacity=0.545] (ac) -- (md);
	\draw[color=green!70!black,opacity=0.545] (md) -- (u);
	\draw[color=green!70!black,opacity=0.545] (md) -- (A);
	%
	% vertices: 
	\fill[color=green!30!black] (mr) circle (0pt) node[color=red!80!black, right, font=\footnotesize] {$M$};
	\fill[color=green!30!black] (ml) circle (0pt) node[color=red!80!black, left, font=\footnotesize] {$M$};
	\fill[color=green!30!black] (A) circle (0pt) node[color=green!30!black, left, font=\footnotesize] {$A$};
	green!30!black
	\fill[color=red!80!black] (ac) circle (3pt) node[color=red!80!black, above, font=\footnotesize] {$\triangleright_M$};
	\fill[color=green!30!black] (md) circle (3pt) node[color=green!30!black, below, font=\footnotesize] {$\mu_A^\vee$};
	\fill[color=green!80!black] (u) circle (3pt) node[color=green!80!black, below, font=\footnotesize] {$u_A$};
	\end{tikzpicture}}%%popende 
%%%%%%%%%%%%%%%%%%%%%%  
\, , \quad 
%%%%%%%%%%%%%%%%%%%%%%
\tikzzbox{\begin{tikzpicture}[ultra thick,scale=0.75,baseline=-0.1cm]
	\coordinate (ml) at (-3,0);
	\coordinate (mr) at (0,0);
	\coordinate (ac) at (-1.5,0);
	\coordinate (A) at (-3,1);
	%
	% lines: 
	\draw[color=orange!80,opacity=0.545] (ml) -- (mr);
	\draw[color=green!70!black,opacity=0.545] (ac) -- (A);
	%
	% vertices: 
	\fill[color=green!30!black] (mr) circle (0pt) node[color=red!80!black, right, font=\footnotesize] {$M$};
	\fill[color=green!30!black] (ml) circle (0pt) node[color=red!80!black, left, font=\footnotesize] {$M$};
	\fill[color=green!30!black] (A) circle (0pt) node[color=green!30!black, left, font=\footnotesize] {$B$};
	green!30!black
	\fill[color=red!80!black] (ac) circle (3pt) node[color=red!80!black, below, font=\footnotesize] {${\triangleleft_M}^\dagger$};
	\end{tikzpicture}}%%popende 
%%%%%%%%%%%%%%%%%%%%%%  
:= 
%%%%%%%%%%%%%%%%%%%%%%
\tikzzbox{\begin{tikzpicture}[ultra thick,scale=0.75,baseline=-0.1cm]
	\coordinate (ml) at (-3,0);
	\coordinate (mr) at (0,0);
	\coordinate (ac) at (-1.5,0);
	\coordinate (A) at (-3,1);
	\coordinate (md) at (-1,0.5);
	\coordinate (u) at (-0.3,0.5);
	%
	% lines: 
	\draw[color=orange!80,opacity=0.545] (ml) -- (mr);
	\draw[color=green!70!black,opacity=0.545] (ac) -- (md);
	\draw[color=green!70!black,opacity=0.545] (md) -- (u);
	\draw[color=green!70!black,opacity=0.545] (md) -- (A);
	%
	% vertices: 
	\fill[color=green!30!black] (mr) circle (0pt) node[color=red!80!black, right, font=\footnotesize] {$M$};
	\fill[color=green!30!black] (ml) circle (0pt) node[color=red!80!black, left, font=\footnotesize] {$M$};
	\fill[color=green!30!black] (A) circle (0pt) node[color=green!30!black, left, font=\footnotesize] {$B$};
	green!30!black
	\fill[color=red!80!black] (ac) circle (3pt) node[color=red!80!black, below, font=\footnotesize] {$\triangleleft_M$};
	\fill[color=green!30!black] (md) circle (3pt) node[color=green!30!black, above, font=\footnotesize] {$\mu_B^\vee$};
	\fill[color=green!80!black] (u) circle (3pt) node[color=green!80!black, above, font=\footnotesize] {$u_B$};
	\end{tikzpicture}}%%popende 
%%%%%%%%%%%%%%%%%%%%%%  
\, . 
\ee 
Using that $\mu_A^\vee$, $\mu_B^\vee$ are adjoint to $\mu_A$ and $\mu_B$, respectively, we can equip the 2-morphisms ${\triangleright_M}^\dagger$ and ${\triangleleft_M}^\dagger$ with the structure of both a left and right adjoint to $\triangleright_M$ and $\triangleleft_M$, respectively. The adjunction data for $\triangleright_M$ is given by
\begin{equation}
\begin{tikzcd}[column sep=4em, row sep=2.5em]
%%%%%%%%%%%%%%%%%%%%%%
\tikzzbox{\begin{tikzpicture}[ultra thick,scale=0.6,baseline=-0.1cm]
	\coordinate (ml) at (-3,0);
	\coordinate (mr) at (0,0);
	\coordinate (ac) at (-1.5,0);
	\coordinate (ac2) at (-2.5,0);
	\coordinate (md) at (-1,-0.5);
	\coordinate (u) at (-0.3,-0.5);
	%
	% lines: 
	\draw[color=orange!80,opacity=0.545] (ml) -- (mr);
	\draw[color=green!70!black,opacity=0.545] (ac) -- (md);
	\draw[color=green!70!black,opacity=0.545] (md) -- (u);
	\draw[color=green!70!black,opacity=0.545] (md) .. controls +(-1,0) and +(0.5,-0.5) .. (ac2);
	\fill[color=red!80!black] (ac2) circle (3pt) node[color=red!80!black, above, font=\footnotesize] {$\triangleright_M$};
	\fill[color=red!80!black] (ac) circle (3pt) node[color=red!80!black, above, font=\footnotesize] {$\triangleright_M$};
	\fill[color=green!30!black] (md) circle (3pt) node[color=green!30!black, below, font=\footnotesize] {$\mu^\vee$};
	\fill[color=green!80!black] (u) circle (3pt) node[color=green!80!black, below, font=\footnotesize] {};
	\end{tikzpicture}}%%popende 
%%%%%%%%%%%%%%%%%%%%%%  
\arrow{r}{(\alpha_M^{\textrm{r}})^{-1}}
\arrow{r}[swap]{\cong}
& 
%%%%%%%%%%%%%%%%%%%%%%
\tikzzbox{\begin{tikzpicture}[ultra thick,scale=0.6,baseline=-0.1cm]
	\coordinate (ml) at (-3,0);
	\coordinate (mr) at (0,0);
	\coordinate (m) at (-2,-0.5);
	\coordinate (ac2) at (-2.5,0);
	\coordinate (md) at (-1,-0.5);
	\coordinate (u) at (-0.3,-0.5);
	%
	% lines: 
	\draw[color=orange!80,opacity=0.545] (ml) -- (mr);
	\draw[color=green!70!black,opacity=0.545] (m) .. controls +(0.25,-0.25) and +(-0.25,-0.25) .. (md);
	\draw[color=green!70!black,opacity=0.545] (m) .. controls +(0.25,0.25) and +(-0.25,0.25) .. (md);
	\draw[color=green!70!black,opacity=0.545] (md) -- (u);
	\draw[color=green!70!black,opacity=0.545] (m) -- (ac2);
	\fill[color=red!80!black] (ac2) circle (3pt) node[color=red!80!black, above, font=\footnotesize] {$\triangleright_M$};
	\fill[color=green!30!black] (m) circle (3pt) node[color=green!30!black, below, font=\footnotesize] {$\mu$};
	\fill[color=green!30!black] (md) circle (3pt) node[color=green!30!black, below, font=\footnotesize] {$\mu^\vee$};
	\fill[color=green!80!black] (u) circle (3pt) node[color=green!80!black, below, font=\footnotesize] {};
	\end{tikzpicture}}%%popende 
%%%%%%%%%%%%%%%%%%%%%%  
\arrow[r, "\tev_{\mu}", ->,shift left=1.5]
\arrow[r, "\coev_{\mu}", swap, <-,shift right=1.5]
& 
%%%%%%%%%%%%%%%%%%%%%%
\tikzzbox{\begin{tikzpicture}[ultra thick,scale=0.6,baseline=-0.1cm]
	\coordinate (ml) at (-3,0);
	\coordinate (mr) at (0,0);
	\coordinate (ac2) at (-2.5,0);
	\coordinate (u) at (-0.3,-0.5);
	%
	% lines: 
	\draw[color=orange!80,opacity=0.545] (ml) -- (mr);
	\draw[color=green!70!black,opacity=0.545] (u) .. controls +(-1,0) and +(0.5,-0.5) .. (ac2);
	\fill[color=red!80!black] (ac2) circle (3pt) node[color=red!80!black, above, font=\footnotesize] {$\triangleright_M$};
	\fill[color=green!80!black] (u) circle (3pt) node[color=green!80!black, below, font=\footnotesize] {};
	\end{tikzpicture}}%%popende 
%%%%%%%%%%%%%%%%%%%%%%  
\arrow{r}{u_M^{\textrm{r}}}
\arrow{r}[swap]{\cong}
& 
%%%%%%%%%%%%%%%%%%%%%%
\tikzzbox{\begin{tikzpicture}[ultra thick,scale=0.6,baseline=-0.1cm]
	\coordinate (ml) at (-3,0);
	\coordinate (mr) at (0,0);
	%
	% lines: 
	\draw[color=orange!80,opacity=0.545] (ml) -- (mr);
	\end{tikzpicture}}%%popende 
%%%%%%%%%%%%%%%%%%%%%%  
\end{tikzcd} 
\end{equation}
and 
\begin{equation}
\begin{tikzcd}[column sep=4em, row sep=2.5em]
%%%%%%%%%%%%%%%%%%%%%%
\tikzzbox{\begin{tikzpicture}[ultra thick,scale=0.6,baseline=-0.4cm]
	\coordinate (ml) at (-3,0);
	\coordinate (mr) at (0,0);
	\coordinate (aml) at (-3,-1);
	\coordinate (amr) at (0,-1);
	\coordinate (ac2) at (-0.8,0);
	\coordinate (u) at (-0.3,-0.5);
	%
	% lines: 
	\draw[color=orange!80,opacity=0.545] (ml) -- (mr);
	\draw[color=green!70!black,opacity=0.545] (aml) -- (amr);
	\end{tikzpicture}}%%popende 
%%%%%%%%%%%%%%%%%%%%%%  
\arrow{r}{(u_M^{\textrm{r}})^{-1}}
\arrow{r}[swap]{\cong}
& 
%%%%%%%%%%%%%%%%%%%%%%
\tikzzbox{\begin{tikzpicture}[ultra thick,scale=0.6,baseline=-0.4cm]
	\coordinate (ml) at (-3,0);
	\coordinate (mr) at (0,0);
	\coordinate (aml) at (-3,-1);
	\coordinate (amr) at (0,-1);
	\coordinate (ac2) at (-0.8,0);
	\coordinate (u) at (-0.3,-0.5);
	%
	% lines: 
	\draw[color=orange!80,opacity=0.545] (ml) -- (mr);
	\draw[color=green!70!black,opacity=0.545] (aml) -- (amr);
	\draw[color=green!70!black,opacity=0.545] (u) .. controls +(-0.25,0) and +(0.25,-0.25) .. (ac2);
	\fill[color=red!80!black] (ac2) circle (3pt) node[color=red!80!black, above, font=\footnotesize] {$\triangleright_M$};
	\fill[color=green!80!black] (u) circle (3pt) node[color=green!80!black, below, font=\footnotesize] {};
	\end{tikzpicture}}%%popende 
%%%%%%%%%%%%%%%%%%%%%%  
\arrow[r, "\tcoev_{\mu}", ->,shift left=1.5]
\arrow[r, "\ev_{\mu}", swap, <-,shift right=1.5]
& 
%%%%%%%%%%%%%%%%%%%%%%
\tikzzbox{\begin{tikzpicture}[ultra thick,scale=0.6,baseline=-0.4cm]
	\coordinate (ml) at (-3,0);
	\coordinate (mr) at (0,0);
	\coordinate (aml) at (-3,-1);
	\coordinate (amr) at (0,-1);
	\coordinate (ac2) at (-2.5,0);
	\coordinate (u) at (-0.3,-0.5);
	\coordinate (m) at (-0.8,-0.75);
	\coordinate (md) at (-1.8,-0.75);
	%
	% lines: 
	\draw[color=orange!80,opacity=0.545] (ml) -- (mr);
	\draw[color=green!70!black,opacity=0.545] (amr) -- (m) -- (md) -- (aml);
	\draw[color=green!70!black,opacity=0.545] (u) -- (m) -- (md) -- (ac2);
	\fill[color=red!80!black] (ac2) circle (3pt) node[color=red!80!black, above, font=\footnotesize] {$\triangleright_M$};
	\fill[color=green!80!black] (u) circle (3pt) node[color=green!80!black, below, font=\footnotesize] {};
	\fill[color=green!30!black] (md) circle (3pt) node[color=green!30!black, below, font=\footnotesize] {$\mu^\vee$};
	\fill[color=green!30!black] (m) circle (3pt) node[color=green!30!black, below, font=\footnotesize] {$\mu$};
	\end{tikzpicture}}%%popende 
%%%%%%%%%%%%%%%%%%%%%% 
\arrow{r}{\textrm{Frobenius}}
\arrow{r}[swap]{\textrm{\& unit}}
& 
%%%%%%%%%%%%%%%%%%%%%%
\tikzzbox{\begin{tikzpicture}[ultra thick,scale=0.6,baseline=-0.4cm]
	\coordinate (ml) at (-3,0);
	\coordinate (mr) at (0,0);
	\coordinate (aml) at (-3,-1);
	\coordinate (amr) at (0,-0.375);
	\coordinate (ac2) at (-2.5,0);
	\coordinate (m) at (-2,-0.375);
	\coordinate (u) at (-0.8,-0.75);
	\coordinate (md) at (-1.5,-0.75);
	%
	% lines: 
	\draw[color=orange!80,opacity=0.545] (ml) -- (mr);
	\draw[color=green!70!black,opacity=0.545] (amr) -- (m) -- (md) -- (aml);
	\draw[color=green!70!black,opacity=0.545] (u) -- (md) -- (ac2);
	\fill[color=red!80!black] (ac2) circle (3pt) node[color=red!80!black, above, font=\footnotesize] {$\triangleright_M$};
	\fill[color=green!80!black] (u) circle (3pt) node[color=green!80!black, below, font=\footnotesize] {};
	\fill[color=green!30!black] (md) circle (3pt) node[color=green!30!black, below, font=\footnotesize] {$\mu^\vee$};
	\fill[color=green!30!black] (m) circle (3pt) node[color=green!30!black, left, font=\footnotesize] {$\mu$};
	\end{tikzpicture}}%%popende 
%%%%%%%%%%%%%%%%%%%%%% 
\arrow{r}{(\alpha_M^{\textrm{r}})^{-1}}
\arrow{r}[swap]{\cong}
& 
%%%%%%%%%%%%%%%%%%%%%%
\tikzzbox{\begin{tikzpicture}[ultra thick,scale=0.6,baseline=-0.4cm]
	\coordinate (ml) at (-3,0);
	\coordinate (mr) at (0,0);
	\coordinate (aml) at (-3,-1);
	\coordinate (amr) at (0,-0.375);
	\coordinate (ac2) at (-2.5,0);
	\coordinate (ac1) at (-1.5,0);
	\coordinate (m) at (-2,-0.375);
	\coordinate (u) at (-0.8,-0.75);
	\coordinate (md) at (-1.5,-0.75);
	%
	% lines: 
	\draw[color=orange!80,opacity=0.545] (ml) -- (mr);
	\draw[color=green!70!black,opacity=0.545] (m) -- (md) -- (aml);
	\draw[color=green!70!black,opacity=0.545] (u) -- (md) -- (ac2);
	\draw[color=green!70!black,opacity=0.545] (amr) .. controls +(-1,0) and +(0.5,-0.5) .. (ac);
	\fill[color=red!80!black] (ac) circle (3pt) node[color=red!80!black, above, font=\footnotesize] {$\triangleright_M$};
	\fill[color=red!80!black] (ac2) circle (3pt) node[color=red!80!black, above, font=\footnotesize] {$\triangleright_M$};
	\fill[color=green!80!black] (u) circle (3pt) node[color=green!80!black, below, font=\footnotesize] {};
	\fill[color=green!30!black] (md) circle (3pt) node[color=green!30!black, below, font=\footnotesize] {$\mu´^\vee$};
	\end{tikzpicture}}%%popende 
%%%%%%%%%%%%%%%%%%%%%% 
\end{tikzcd}
.
\end{equation}
The data for $\triangleleft_M$ is completely analogous. 
Note that this does not imply that they agree with ${\triangleright_M}^\vee$ and ${\triangleleft_M}^\vee$ (which are part of the data of the Gray category with duals~$\tric$). 
Uniqueness of left and right adjoints gives two potentially different isomorphisms between ${}^\dagger$ and ${}^\vee$. 
Using this we can give a simple condition for a 1-morphism between orbifold data to be in $\orb{\tric}$:  

\begin{proposition}
	\label{prop:1MorphismsOrb}
Let $\mathcal{M}\colon \mathcal{A}\to \mathcal{B}$ be a 1-morphism in $\ealg(\tric)$ between orbifold data such that the above 3-isomorphisms ${\triangleright_M}^\vee \to {\triangleright_M}^\dagger$ and ${\triangleleft_M}^\vee \to {\triangleleft_M}^\dagger$ agree pairwise. 
Then $\mathcal{M}$ is also a 1-morphism in $\orb{\tric}$.   
\end{proposition}  

\begin{proof}
The assumption of the proposition allows us to replace ${\triangleright_M}^\vee$ and ${\triangleleft_M}^\vee$ with ${\triangleright_M}^\dagger$ and ${\triangleleft_M}^\dagger$ when verifying the analogues of \eqref{eq:O4}--\eqref{eq:O8} for $\cat{M}$. 
This allows us to use the defining properties of orbifold data \eqref{eq:O1}--\eqref{eq:O8} to prove the relations for~$\mathcal{M}$. 
The verification is now a straightforward computation. 
Especially, if one uses that the 3-morphisms associated to two different ``admissible'' stratifications of two stratified 3-balls with identical boundaries agree by the defining properties of orbifold data (see \cite{CRS1, CMRSS1} for a detailed discussion).    
\end{proof} 

In the remainder of Section~\ref{subsec:3catOrbData}, we will compute the relative product~$\btimes_{\mathcal A}$ in~$\orb{\tric}$. 
We stress that for this it is not necessary to impose the conditions of Figure~\ref{fig:2MorphismsInCorb} on 1-morphisms. 
Hence we will work in~$\ealg(\tric)$, but we exclusively consider objects $\mathcal A \in \ealg(\tric)$ that satisfy the conditions in Figure~\ref{fig:OrbifoldDatumAxioms}, i.\,e.\ $\mathcal A \in \orb{\tric}$. 
%arXiv_v2: 
	A similar result for bimodules over separable algebras is proven in~\cite[Thm.\,3.1.6]{Decoppet_Morita}. 
	By Remark~\ref{Rem: Connection to seperable algebras} it also applies here. 
	We give the full proof in the following to introduce the relevant conventions, and the description in terms of the cobalanced tensor product is new to the best of our knowledge.       
For a 1-morphism~$\mathcal N$ in $\ealg(\tric)$ with codomain~$\mathcal A$, we can set
\be 
\label{eq:Pdefined}
P \equiv P_{\mathcal M, \mathcal A, \mathcal N} := 
%%%%%%%%%%%%%%%%%%%%%%
\tikzzbox{\begin{tikzpicture}[ultra thick,scale=0.75, baseline=-0.1cm]
	\coordinate (r0) at (0,-0.75);
	\coordinate (r1) at (0,0.75);
	\coordinate (l0) at (-3,-0.75);
	\coordinate (l1) at (-3,0.75);
	\coordinate (acm) at (-2,-0.75);
	\coordinate (acn) at (-2,0.75);
	\coordinate (md) at (-1,0);
	\coordinate (u) at (-0.25,0);
	%
	% lines: 
	\draw[color=orange!50,opacity=0.545] (l0) -- (r0);
	\draw[color=orange!50,opacity=0.545] (l1) -- (r1);
	\draw[color=green!70!black,opacity=0.545] (acm) -- (md);
	\draw[color=green!70!black,opacity=0.545] (acn) -- (md);
	\draw[color=green!70!black,opacity=0.545] (u) -- (md);
	%
	% vertices: 
	\fill[color=red!30!black] (r0) circle (0pt) node[color=red!80!black, right, font=\footnotesize] {$N$};
	\fill[color=red!30!black] (r1) circle (0pt) node[color=red!80!black, right, font=\footnotesize] {$M$};
	\fill[color=red!30!black] (l0) circle (0pt) node[color=red!80!black, left, font=\footnotesize] {$N$};
	\fill[color=red!30!black] (l1) circle (0pt) node[color=red!80!black, left, font=\footnotesize] {$M$};
	\fill[color=red!80!black] (acm) circle (3pt) node[color=red!80!black, below, font=\footnotesize] {};
	\fill[color=red!80!black] (acn) circle (3pt) node[color=red!80!black, above, font=\footnotesize] {};
	\fill[color=green!30!black] (md) circle (3pt) node[color=green!30!black, below, font=\footnotesize] {$\mu^\vee_A$};
	\fill[color=green!80!black] (u) circle (3pt) node[color=green!80!black, above, font=\footnotesize] {$u_A$};
	\end{tikzpicture}}%%popende 
%%%%%%%%%%%%%%%%%%%%%%  
\cong 
%%%%%%%%%%%%%%%%%%%%%%
\tikzzbox{\begin{tikzpicture}[ultra thick,scale=0.75, baseline=-0.1cm]
	\coordinate (r0) at (0,-0.75);
	\coordinate (r1) at (0,0.75);
	\coordinate (l0) at (-3,-0.75);
	\coordinate (l1) at (-3,0.75);
	\coordinate (acm) at (-2,-0.75);
	\coordinate (acn) at (-1,0.75);
	%
	% lines: 
	\draw[color=orange!50,opacity=0.545] (l0) -- (r0);
	\draw[color=orange!50,opacity=0.545] (l1) -- (r1);
	\draw[color=green!70!black,opacity=0.545] (acm) -- (acn);
	%
	% vertices: 
	\fill[color=red!30!black] (r0) circle (0pt) node[color=red!80!black, right, font=\footnotesize] {$N$};
	\fill[color=red!30!black] (r1) circle (0pt) node[color=red!80!black, right, font=\footnotesize] {$M$};
	\fill[color=red!30!black] (l0) circle (0pt) node[color=red!80!black, left, font=\footnotesize] {$N$};
	\fill[color=red!30!black] (l1) circle (0pt) node[color=red!80!black, left, font=\footnotesize] {$M$};
	\fill[color=red!80!black] (acm) circle (3pt) node[color=red!80!black, below, font=\footnotesize] {};
	\fill[color=red!80!black] (acn) circle (3pt) node[color=red!80!black, above, font=\footnotesize] {};
	\end{tikzpicture}}%%popende 
%%%%%%%%%%%%%%%%%%%%%% 
\cong 
%%%%%%%%%%%%%%%%%%%%%%
\tikzzbox{\begin{tikzpicture}[ultra thick,scale=0.75, baseline=-0.1cm]
	\coordinate (r0) at (0,-0.75);
	\coordinate (r1) at (0,0.75);
	\coordinate (l0) at (-3,-0.75);
	\coordinate (l1) at (-3,0.75);
	\coordinate (acm) at (-1,-0.75);
	\coordinate (acn) at (-2,0.75);
	%
	% lines: 
	\draw[color=orange!50,opacity=0.545] (l0) -- (r0);
	\draw[color=orange!50,opacity=0.545] (l1) -- (r1);
	\draw[color=green!70!black,opacity=0.545] (acm) -- (acn);
	%
	% vertices: 
	\fill[color=red!30!black] (r0) circle (0pt) node[color=red!80!black, right, font=\footnotesize] {$N$};
	\fill[color=red!30!black] (r1) circle (0pt) node[color=red!80!black, right, font=\footnotesize] {$M$};
	\fill[color=red!30!black] (l0) circle (0pt) node[color=red!80!black, left, font=\footnotesize] {$N$};
	\fill[color=red!30!black] (l1) circle (0pt) node[color=red!80!black, left, font=\footnotesize] {$M$};
	\fill[color=red!80!black] (acm) circle (3pt) node[color=red!80!black, below, font=\footnotesize] {};
	\fill[color=red!80!black] (acn) circle (3pt) node[color=red!80!black, above, font=\footnotesize] {};
	\end{tikzpicture}}%%popende 
%%%%%%%%%%%%%%%%%%%%%%  
\ee 
where here and in similar situations, we usually drop indices like ``$\mathcal M, \mathcal A, \mathcal N$'' as indicated. 
Importantly, $P$ is a 2-cocone over 
$
M\btimes A \btimes A \btimes N 
\mathrel{\substack{\textstyle\rightarrow\\[-0.6ex]
		\textstyle\rightarrow \\[-0.6ex]
		\textstyle\rightarrow}} 
M \btimes A \btimes N 
\rightrightarrows
M \btimes N
$, 
i.\,e.\ it is balanced with respect to the $\mathcal A$-actions (recall the discussion around~\eqref{eq:BalancingMap}): 

\begin{lemma}
	\label{lem:PBalanced}
	$P$ together with  
	\be 
	\begin{tikzcd}[column sep=5em, row sep=2.5em]
	%%%%%%%%%%%%%%%%%%%%%%
	\tikzzbox{
}%%popende 
	%%%%%%%%%%%%%%%%%%%%%%   
	\end{tikzcd} 
	\ee  
	is a balanced 2-morphism $M\btimes N \lra M\btimes N$. 
\end{lemma}

\begin{proof}
	We have to verify the commutativity of the following pentagon diagram: 
	\be 
	\begin{tikzcd}[column sep=5em, row sep=2.5em]
	& 
	%%%%%%%%%%%%%%%%%%%%%%
	\tikzzbox{%
}%%popende 
	%%%%%%%%%%%%%%%%%%%%%%   
	\arrow{l}{}
	\end{tikzcd}
	\ee 
	The left and right subdiagrams commute thanks to the 2-3 move~$\eqref{eq:O1}$. 
	The two paths around the middle diagram can be represented, in the graphical calculus, as two different ``admissible'' stratifications of two stratified 3-balls with identical boundaries. 
	Hence by the defining properties of orbifold data, both paths are equal (see again \cite{CRS1, CMRSS1} for a detailed discussion). 
\end{proof}

By our assumption on~$\tric$, the 2-colimit $M \relpro_{\mathcal A} N$ 
%arXiv_v3:
	%exisits. 
	 exists. 
Part of its structure is a (balanced) 2-morphism
\be 
\label{eq:Pi}
\Pi \equiv \Pi_{\mathcal M, \mathcal A, \mathcal N} \colon M \btimes N 
	\lra 
	M \relpro_{\mathcal A} N
\ee 
such that the pre-composition functor 
\begin{align}
\Pi^* \colon \Hom_\tric \big( M\relpro_{\mathcal A} N, T \big) 
	& 
	\lra \textrm{Bal}\big( M\btimes N, T \big) 
	\nonumber
	\\ 
	\varphi 
	& 
	\lmt \varphi \circ \Pi 
	\label{eq:BalancedEquiv}
\end{align}
is an equivalence for every 1-morphism~$T$ parallel to $M\btimes N$. 
In particular, for $T = M\btimes N$ we obtain from this universal property an essentially unique 2-morphism 
\be 
\label{eq:Sigma}
\Sigma 
	:= 
	(\Pi^*)^{-1}(P) \colon M \relpro_{\mathcal A} N \lra M\btimes N \, , 
\ee 
from which it immediately follows that 
\be 
\label{eq:SigmaPiP}
\Sigma \circ \Pi \cong P \, . 
\ee 
Using the graphical notation introduced in~\eqref{eq:RelativeGraphicalNotation} we summarise~\eqref{eq:Pi}, \eqref{eq:Sigma} and~\eqref{eq:SigmaPiP} as 
\be 
\Pi = 
%%%%%%%%%%%%%%%%%%%%%%
\tikzzbox{\begin{tikzpicture}[ultra thick,scale=0.75, baseline=-0.1cm]
	\coordinate (r0) at (0,-0.75);
	\coordinate (r1) at (0,0.75);
	\coordinate (l0) at (-3,-0.75);
	\coordinate (l1) at (-3,0.75);
	\coordinate (acm) at (-1.5,-0.75);
	\coordinate (acn) at (-1.5,0.75);
	%
	% area: 
	\fill[color=orange!50,opacity=0.35] ($(l0)+(0,-0.8pt)$) -- ($(acm)+(0,-0.8pt)$) -- ($(acn)+(0,0.8pt)$) -- ($(l1)+(0,0.8pt)$);
	%
	% lines: 
	\draw[color=orange!50,opacity=0.545] (l0) -- (r0);
	\draw[color=orange!50,opacity=0.545] (l1) -- (r1);
	\draw[blue!50!green] ($(acm)+(0,-0.8pt)$) -- ($(acn)+(0,0.8pt)$);
	%
	% vertices: 
	\fill[color=red!30!black] (r0) circle (0pt) node[color=red!80!black, right, font=\footnotesize] {$N$};
	\fill[color=red!30!black] (r1) circle (0pt) node[color=red!80!black, right, font=\footnotesize] {$M$};
	\fill[color=red!30!black] (l0) circle (0pt) node[color=red!80!black, left, font=\footnotesize] {$N$};
	\fill[color=red!30!black] (l1) circle (0pt) node[color=red!80!black, left, font=\footnotesize] {$M$};
	\fill[color=red!30!black] ($(l1)-0.5*(l1)+0.5*(l0)+(-0.05,-0.05)$) circle (0pt) node[color=red!80!black, left, font=\footnotesize] {$\displaystyle{\relpro_{\mathcal A}}$};
	\fill[color=blue!50!green] (acn) circle (0pt) node[color=blue!50!green, above, font=\footnotesize] {$\Pi$};
	\end{tikzpicture}}%%popende 
%%%%%%%%%%%%%%%%%%%%%%   
\, , \quad 
\Sigma = 
%%%%%%%%%%%%%%%%%%%%%%
\tikzzbox{\begin{tikzpicture}[ultra thick,scale=0.75, baseline=-0.1cm]
	\coordinate (r0) at (0,-0.75);
	\coordinate (r1) at (0,0.75);
	\coordinate (l0) at (-3,-0.75);
	\coordinate (l1) at (-3,0.75);
	\coordinate (acm) at (-1.5,-0.75);
	\coordinate (acn) at (-1.5,0.75);
	%
	% area: 
	\fill[color=orange!50,opacity=0.35] ($(r0)+(0,-0.8pt)$) -- ($(acm)+(0,-0.8pt)$) -- ($(acn)+(0,0.8pt)$) -- ($(r1)+(0,0.8pt)$);
	%
	% lines: 
	\draw[color=orange!50,opacity=0.545] (l0) -- (r0);
	\draw[color=orange!50,opacity=0.545] (l1) -- (r1);
	\draw[blue!50!purple] ($(acm)+(0,-0.8pt)$) -- ($(acn)+(0,0.8pt)$);
	%
	% vertices: 
	\fill[color=red!30!black] (r0) circle (0pt) node[color=red!80!black, right, font=\footnotesize] {$N$};
	\fill[color=red!30!black] (r1) circle (0pt) node[color=red!80!black, right, font=\footnotesize] {$M$};
	\fill[color=red!30!black] (l0) circle (0pt) node[color=red!80!black, left, font=\footnotesize] {$N$};
	\fill[color=red!30!black] (l1) circle (0pt) node[color=red!80!black, left, font=\footnotesize] {$M$};
	\fill[color=red!30!black] ($(r1)-0.5*(r1)+0.5*(r0)+(0.05,-0.05)$) circle (0pt) node[color=red!80!black, right, font=\footnotesize] {$\displaystyle{\relpro_{\mathcal A}}$};
	\fill[color=blue!50!purple] (acn) circle (0pt) node[color=blue!50!purple, above, font=\footnotesize] {$\Sigma$};
	\end{tikzpicture}}%%popende 
%%%%%%%%%%%%%%%%%%%%%%    
\, , \quad 
P \cong 
%%%%%%%%%%%%%%%%%%%%%%
\tikzzbox{\begin{tikzpicture}[ultra thick,scale=0.75, baseline=-0.1cm]
	\coordinate (r0) at (0,-0.75);
	\coordinate (r1) at (0,0.75);
	\coordinate (l0) at (-3,-0.75);
	\coordinate (l1) at (-3,0.75);
	\coordinate (acm) at (-1.5,-0.75);
	\coordinate (acn) at (-1.5,0.75);
	%
	% area: 
	\fill[color=orange!50,opacity=0.35] ($(acm)+(0.75,-0.8pt)$) -- ($(acm)+(-0.75,-0.8pt)$) -- ($(acn)+(-0.75,0.8pt)$) -- ($(acn)+(0.75,0.8pt)$);
	%
	% lines: 
	\draw[color=orange!50,opacity=0.545] (l0) -- (r0);
	\draw[color=orange!50,opacity=0.545] (l1) -- (r1);
	\draw[blue!50!purple] ($(acm)+(-0.75,-0.8pt)$) -- ($(acn)+(-0.75,0.8pt)$);
	\draw[blue!50!green] ($(acm)+(+0.75,-0.8pt)$) -- ($(acn)+(+0.75,0.8pt)$);
	%
	% vertices: 
	\fill[color=red!30!black] (r0) circle (0pt) node[color=red!80!black, right, font=\footnotesize] {$N$};
	\fill[color=red!30!black] (r1) circle (0pt) node[color=red!80!black, right, font=\footnotesize] {$M$};
	\fill[color=red!30!black] (l0) circle (0pt) node[color=red!80!black, left, font=\footnotesize] {$N$};
	\fill[color=red!30!black] (l1) circle (0pt) node[color=red!80!black, left, font=\footnotesize] {$M$};
	\fill[color=blue!50!purple] ($(acn)+(-0.75,0.8pt)$) circle (0pt) node[color=blue!50!purple, above, font=\footnotesize] {$\Sigma$};
	\fill[color=blue!50!green] ($(acn)+(0.75,0.8pt)$) circle (0pt) node[color=blue!50!green, above, font=\footnotesize] {$\Pi$};
	\end{tikzpicture}}%%popende 
%%%%%%%%%%%%%%%%%%%%%%  
\, . 
\ee 
Below in Proposition~\ref{prop:RelProd} we will see that~$\Pi$ and~$\Sigma$ are part of the splitting data of the 
%arXiv_v2: 
	%higher idempotent~$P$. 
	condensation monad~$P$. 
Following \cite{DouglasReutter2018, GaiottoJohnsonFreyd}, we introduce the latter as follows: 

\begin{definition}
	\label{def:2IdempotentSplit}
	Let $p\colon V \lra V$ be a 2-morphism in a 3-category. 
	\begin{enumerate}
		\item 
		The structure of a 
		%arXiv_v2: 
			%\textsl{2-idempotent} 
			 \textsl{condensation 2-monad} 
		on~$p$ is given by 3-morphisms $f\colon p \circ p \lra p$ and $g \colon p \lra p \circ p$ such that $fg = 1_p$, and such that~$f$ and~$g$ give~$p$ the structure of a non-unital $\Delta$-separable Frobenius algebra. 
		\item 
		A \textsl{splitting} of a 
		%arixv_v2: 
			%2-idempotent
			 condensation 2-monad 
		$(p,f,g)$ is given by two 2-morphisms $\pi\colon V \lra U$, $\sigma\colon U \lra V$ together with two 3-morphisms $\varepsilon\colon \pi\circ\sigma \lra 1_U$, $\eta\colon 1_U \lra \pi\circ\sigma$ such that $\sigma\circ\pi \cong p$ and $\varepsilon\eta = 1_{1_U}$, and $f,g$ are induced by $\varepsilon, \eta$. 
		Then~$U$ is the \textsl{image} of the split 
		%arXiv_v2:
			%idempotent
			 condensation 2-monad 
		$(p,\pi,\sigma,\varepsilon,\eta)$. 
	\end{enumerate}
\end{definition}

\begin{proposition}
	\label{prop:RelProd}
	Let~$\mathcal A$ be an orbifold datum in~$\tric$, and let $\mathcal M, \mathcal N$ be 1-morphisms in $\ealg(\tric)$ such that $M\relpro_{\mathcal A} N$ exists. 
	Then $M\relpro_{\mathcal A} N$ is the image of the split 
	%arXiv_v2:
		%2-idempotent 
		 condensation 2-monad
	$(P,\Pi,\Sigma,\varepsilon,\eta)$, where $\varepsilon,\eta$ are induced by $\widetilde{\ev}_\mu, \coev_\mu$ (as explained around~\eqref{eq:Splitting} below), respectively. 
\end{proposition}

\begin{proof}
	We first observe that 
	\vspace{0.4cm}
	\be 
	\label{eq:Splitting}
	%%%%%%%%%%%%%%%%%%%%%%%%%%%%%%
	\begin{tikzpicture}[
	baseline=(current bounding box.base),
	%>=stealth,
	descr/.style={fill=white,inner sep=3.5pt},
	normal line/.style={->}
	]
	\matrix (m) [matrix of math nodes, row sep=2em, column sep=5em, text height=1.5ex, text depth=0.1ex] {%
		\Pi \circ \Sigma \circ \Pi 
		\cong 
		\Pi \circ P 
		= 
		%%%%%%%%%%%%%%%%%%%%%%
		\tikzzbox{\begin{tikzpicture}[ultra thick,scale=0.75, baseline=-0.1cm]
			\coordinate (r0) at (0,-0.75);
			\coordinate (r1) at (0,0.75);
			\coordinate (l0) at (-3,-0.75);
			\coordinate (l1) at (-3,0.75);
			\coordinate (acm) at (-2,-0.75);
			\coordinate (acn) at (-2,0.75);
			\coordinate (acmm) at (-1.5,-0.75);
			\coordinate (acnn) at (-1.5,0.75);
			\coordinate (md) at (-1,0);
			\coordinate (u) at (-0.25,0);
			%
			% area: 
			\fill[color=orange!50,opacity=0.35] ($(l0)+(0,-0.8pt)$) -- ($(acm)+(0,-0.8pt)$) -- ($(acn)+(0,0.8pt)$) -- ($(l1)+(0,0.8pt)$);
			%
			% lines: 
			\draw[color=orange!50,opacity=0.545] (l0) -- (r0);
			\draw[color=orange!50,opacity=0.545] (l1) -- (r1);
			\draw[color=green!70!black,opacity=0.545] (acmm) -- (md);
			\draw[color=green!70!black,opacity=0.545] (acnn) -- (md);
			\draw[color=green!70!black,opacity=0.545] (u) -- (md);
			\draw[blue!50!green] ($(acm)+(0,-0.8pt)$) -- ($(acn)+(0,0.8pt)$);
			%
			% vertices: 
			%		\fill[color=red!30!black] (r0) circle (0pt) node[color=red!80!black, right, font=\footnotesize] {$N$};
			%		\fill[color=red!30!black] (r1) circle (0pt) node[color=red!80!black, right, font=\footnotesize] {$M$};
			%		\fill[color=red!30!black] (l0) circle (0pt) node[color=red!80!black, left, font=\footnotesize] {$N$};
			%		\fill[color=red!30!black] (l1) circle (0pt) node[color=red!80!black, left, font=\footnotesize] {$M$};
			%		\fill[color=red!30!black] ($(l1)-0.5*(l1)+0.5*(l0)+(-0.05,-0.05)$) circle (0pt) node[color=red!80!black, left, font=\footnotesize] {$\displaystyle{\relpro_{\mathcal A}}$};
			%
			\fill[color=blue!50!green] (acn) circle (0pt) node[color=blue!50!green, above, font=\footnotesize] {$\Pi$};
			\fill[color=red!80!black] (acmm) circle (3pt) node[color=red!80!black, below, font=\footnotesize] {};
			\fill[color=red!80!black] (acnn) circle (3pt) node[color=red!80!black, above, font=\footnotesize] {};
			\fill[color=green!30!black] (md) circle (3pt) node[color=green!30!black, below, font=\footnotesize] {};
			\fill[color=green!80!black] (u) circle (3pt) node[color=green!80!black, above, font=\footnotesize] {};
			\end{tikzpicture}}%%popende 
		%%%%%%%%%%%%%%%%%%%%%%  
		\cong 
		%%%%%%%%%%%%%%%%%%%%%%
		\tikzzbox{\begin{tikzpicture}[ultra thick,scale=0.75, baseline=-0.1cm]
			\coordinate (r0) at (0,-0.75);
			\coordinate (r1) at (0,0.75);
			\coordinate (l0) at (-3,-0.75);
			\coordinate (l1) at (-3,0.75);
			\coordinate (acm) at (-2,-0.75);
			\coordinate (acn) at (-2,0.75);
			\coordinate (acnn) at (-1.5,0.75);
			\coordinate (md) at (-0.5,0.25);
			\coordinate (m) at (-1,0.5);
			\coordinate (u) at (-0.1,0.05);
			%
			% area: 
			\fill[color=orange!50,opacity=0.35] ($(l0)+(0,-0.8pt)$) -- ($(acm)+(0,-0.8pt)$) -- ($(acn)+(0,0.8pt)$) -- ($(l1)+(0,0.8pt)$);
			%
			% lines: 
			\draw[color=orange!50,opacity=0.545] (l0) -- (r0);
			\draw[color=orange!50,opacity=0.545] (l1) -- (r1);
			\draw[color=green!70!black,opacity=0.545] (acnn) -- (m);
			\draw[color=green!70!black,opacity=0.545] (m) .. controls +(0.4,0.2) and +(0.2,0.4) .. (md);
			\draw[color=green!70!black,opacity=0.545] (m) .. controls +(-0.2,-0.4) and +(-0.2,-0.4) .. (md);
			\draw[color=green!70!black,opacity=0.545] (u) -- (md);
			\draw[blue!50!green] ($(acm)+(0,-0.8pt)$) -- ($(acn)+(0,0.8pt)$);
			%
			% vertices: 
			%		\fill[color=red!30!black] (r0) circle (0pt) node[color=red!80!black, right, font=\footnotesize] {$N$};
			%		\fill[color=red!30!black] (r1) circle (0pt) node[color=red!80!black, right, font=\footnotesize] {$M$};
			%		\fill[color=red!30!black] (l0) circle (0pt) node[color=red!80!black, left, font=\footnotesize] {$N$};
			%		\fill[color=red!30!black] (l1) circle (0pt) node[color=red!80!black, left, font=\footnotesize] {$M$};
			%		\fill[color=red!30!black] ($(l1)-0.5*(l1)+0.5*(l0)+(-0.05,-0.05)$) circle (0pt) node[color=red!80!black, left, font=\footnotesize] {$\displaystyle{\relpro_{\mathcal A}}$};
			%
			\fill[color=blue!50!green] (acn) circle (0pt) node[color=blue!50!green, above, font=\footnotesize] {$\Pi$};
			\fill[color=red!80!black] (acnn) circle (3pt) node[color=red!80!black, right, font=\footnotesize] {};
			\fill[color=green!30!black] (m) circle (3pt) node[color=green!30!black, left, font=\footnotesize] {};
			\fill[color=green!30!black] ($(m)+(0.05,-0.1)$) circle (0pt) node[color=green!30!black, left, font=\footnotesize] {$\mu$};
			\fill[color=green!30!black] (md) circle (3pt) node[color=green!30!black, below, font=\footnotesize] {};
			\fill[color=green!80!black] (u) circle (3pt) node[color=green!80!black, above, font=\footnotesize] {};
			\end{tikzpicture}}%%popende 
		%%%%%%%%%%%%%%%%%%%%%%   
		& 
		%%%%%%%%%%%%%%%%%%%%%%
		\tikzzbox{\begin{tikzpicture}[ultra thick,scale=0.75, baseline=-0.1cm]
			\coordinate (r0) at (0,-0.75);
			\coordinate (r1) at (0,0.75);
			\coordinate (l0) at (-3,-0.75);
			\coordinate (l1) at (-3,0.75);
			\coordinate (acm) at (-1.5,-0.75);
			\coordinate (acn) at (-1.5,0.75);
			%
			% area: 
			\fill[color=orange!50,opacity=0.35] ($(l0)+(0,-0.8pt)$) -- ($(acm)+(0,-0.8pt)$) -- ($(acn)+(0,0.8pt)$) -- ($(l1)+(0,0.8pt)$);
			%
			% lines: 
			\draw[color=orange!50,opacity=0.545] (l0) -- (r0);
			\draw[color=orange!50,opacity=0.545] (l1) -- (r1);
			\draw[blue!50!green] ($(acm)+(0,-0.8pt)$) -- ($(acn)+(0,0.8pt)$);
			%
			% vertices: 
			%		\fill[color=red!30!black] (r0) circle (0pt) node[color=red!80!black, right, font=\footnotesize] {$N$};
			%		\fill[color=red!30!black] (r1) circle (0pt) node[color=red!80!black, right, font=\footnotesize] {$M$};
			%		\fill[color=red!30!black] (l0) circle (0pt) node[color=red!80!black, left, font=\footnotesize] {$N$};
			%		\fill[color=red!30!black] (l1) circle (0pt) node[color=red!80!black, left, font=\footnotesize] {$M$};
			%		\fill[color=red!30!black] ($(l1)-0.5*(l1)+0.5*(l0)+(-0.05,-0.05)$) circle (0pt) node[color=red!80!black, left, font=\footnotesize] {$\displaystyle{\relpro_{\mathcal A}}$};
			%
			\fill[color=blue!50!green] (acn) circle (0pt) node[color=blue!50!green, above, font=\footnotesize] {$\Pi$};
			\end{tikzpicture}}%%popende 
		%%%%%%%%%%%%%%%%%%%%%%    
		, 
		\\
	};
	\path[font=\footnotesize, transform canvas={yshift=+0.8mm}] (m-1-1) edge[->] node[above] {$\tev_\mu$} (m-1-2);
	\path[font=\footnotesize, transform canvas={yshift=-0.8mm}] (m-1-1) edge[<-] node[below] {$\coev_\mu$} (m-1-2);
	\end{tikzpicture}
	%%%%%%%%%%%%%%%%%%%%%%%%%%%%%%
	\ee 
	\vspace{0.3cm}
	
	\noindent 
	and that $\tev_\mu \coev_\mu = 1_{1_A}$ according to~\eqref{eq:O8}. 
 	Hence by the universal property of $M\relpro_{\mathcal A} N$ (setting $T= M\relpro_{\mathcal A} N$ in \eqref{eq:BalancedEquiv}) there are essentially unique 3-morphisms $\varepsilon\colon \Pi\circ\Sigma \lra 1_{M\relpro_{\mathcal A} N}$ and $\eta\colon 1_{M\relpro_{\mathcal A} N} \lra \Pi \circ \Sigma$ with $\varepsilon\eta = 1_{1_{M \relpro_{\mathcal A} N}}$. 
	Moreover, we have
	\be
	%%%%%%%%%%%%%%%%%%%%%%%%%%%%%%
	\begin{tikzpicture}[
	baseline=(current bounding box.base),
	%>=stealth,
	descr/.style={fill=white,inner sep=3.5pt},
	normal line/.style={->}
	]
	\matrix (m) [matrix of math nodes, row sep=2em, column sep=6em, text height=1.5ex, text depth=0.1ex] {%
		P^2 \cong \Sigma \circ \Pi \circ \Sigma \circ \Pi  
		& 
		\Sigma\circ\Pi \cong P  
		\, .  
		\\
	};
	\path[font=\footnotesize, transform canvas={yshift=+0.8mm}] (m-1-1) edge[->] node[above] {$1_\Sigma \circ \varepsilon \circ 1_\Pi$} (m-1-2);
	\path[font=\footnotesize, transform canvas={yshift=-0.8mm}] (m-1-1) edge[<-] node[below] {$1_\Sigma \circ \eta \circ 1_\Pi$} (m-1-2);
	\end{tikzpicture}
	%%%%%%%%%%%%%%%%%%%%%%%%%%%%%%
	\ee 
\end{proof}

In summary, the relative product $\relpro_{\mathcal A}$ can be computed by splitting the 
%arXiv_v2
	%2-idempotent~$P$
	 condensation 2-monad~$P$ 
from~\eqref{eq:Pdefined}. 

\medskip 

It will be useful to establish a diagrammatic representation of the functor $(\Pi^*)^{-1}$, and also to consider the dual of the relative product $\relpro_{\mathcal A}$. 
Concerning the former, recall the equivalence~$\Pi^*$ in~\eqref{eq:BalancedEquiv}. 
We write 
\begin{align}
(\Pi^*)^{-1} \colon \textrm{Bal}\big(M\btimes N, T \big) 
& 
\stackrel{\cong}{\lra} \Hom_\tric\big(M\relpro_{\mathcal A} N , T \big) 
\nonumber
\\
%%%%%%%%%%%%%%%%%%%%%%
\tikzzbox{\begin{tikzpicture}[ultra thick,scale=0.75, baseline=-0.1cm]
	\coordinate (r0) at (0,-0.75);
	\coordinate (r1) at (0,0.75);
	\coordinate (l) at (-3,0);
	\coordinate (psi) at (-1.5,0);
	%
	% lines: 
	\draw[color=orange!50,opacity=0.545] (r0) -- (psi);
	\draw[color=orange!50,opacity=0.545] (r1) -- (psi);
	\draw[color=orange!50,opacity=0.545] (l) -- (psi);
	%
	%
	% vertices: 
	\fill[color=red!30!black] (r0) circle (0pt) node[color=red!80!black, right, font=\footnotesize] {$N$};
	\fill[color=red!30!black] (r1) circle (0pt) node[color=red!80!black, right, font=\footnotesize] {$M$};
	\fill[color=red!30!black] (l) circle (0pt) node[color=red!80!black, left, font=\footnotesize] {$T$};
	\fill[color=green!50!purple] (psi) circle (3pt) node[color=green!50!purple, above, font=\footnotesize] {$\psi$};
	\end{tikzpicture}}%%popende 
%%%%%%%%%%%%%%%%%%%%%% 
& 
\longmapsto 
%%%%%%%%%%%%%%%%%%%%%%
\tikzzbox{\begin{tikzpicture}[ultra thick,scale=0.75, baseline=-0.1cm]
	\coordinate (r0) at (0,-0.75);
	\coordinate (r1) at (0,0.75);
	\coordinate (l) at (-3,0);
	\coordinate (psi) at (-1.5,0);
	%
	% area: 
	\fill[color=orange!50,opacity=0.35] ($(r0)+(0,-0.8pt)$) -- ($(psi)+(0,-0.8pt)$) -- ($(r1)+(0,0.8pt)$);
	%
	% lines: 
	\draw[color=orange!50,opacity=0.545] (r0) -- (psi);
	\draw[color=orange!50,opacity=0.545] (r1) -- (psi);
	\draw[color=orange!50,opacity=0.545] (l) -- (psi);
	%
	%
	% vertices: 
	\fill[color=red!30!black] (r0) circle (0pt) node[color=red!80!black, right, font=\footnotesize] {$N$};
	\fill[color=red!30!black] (r1) circle (0pt) node[color=red!80!black, right, font=\footnotesize] {$M$};
	\fill[color=red!30!black] (l) circle (0pt) node[color=red!80!black, left, font=\footnotesize] {$T$};
	\fill[color=red!30!black] ($(r1)-0.5*(r1)+0.5*(r0)+(0.05,-0.05)$) circle (0pt) node[color=red!80!black, right, font=\footnotesize] {$\displaystyle{\relpro_{\mathcal A}}$};
	\fill[color=green!50!purple] (psi) circle (3pt) node[color=green!50!purple, above, font=\footnotesize] {$\psi$};
	\end{tikzpicture}}%%popende 
%%%%%%%%%%%%%%%%%%%%%%  
\label{eq:PiStarInv}
\end{align}
for its weak inverse. 
Hence applying $\Pi^* = (-)\circ \Pi$ to the image of $(\Pi^*)^{-1}$ amounts to deleting the shading associated to the relative product. 

Given $\mathcal A, \mathcal M, \mathcal N$ as above, the categories $\textrm{Cobal}(S,M\btimes N)$ of \textsl{cobalanced} 2-morphisms $S\lra M\btimes N$ are defined dually to $\textrm{Bal}(M\btimes N, T)$. 
In particular, objects~$\xi$ of $\textrm{Cobal}(S,M\btimes N)$ come together with 3-isomorphisms 
\be 
%%%%%%%%%%%%%%%%%%%%%%
\tikzzbox{\begin{tikzpicture}[ultra thick,scale=0.75, baseline=-0.1cm]
	\coordinate (r0) at (0,-0.75);
	\coordinate (r1) at (0,0.75);
	\coordinate (l) at (3,0);
	\coordinate (A) at (0,0);
	\coordinate (psi) at (1.5,0);
	\coordinate (ac) at ($(r0)-0.5*(r0)+0.5*(psi)$);
	%
	% lines: 
	\draw[color=orange!50,opacity=0.545] (r0) -- (psi);
	\draw[color=orange!50,opacity=0.545] (r1) -- (psi);
	\draw[color=orange!50,opacity=0.545] (l) -- (psi);
	\draw[color=green!70!black,opacity=0.545] (A) -- (ac);
	%
	%
	% vertices: 
	\fill[color=green!30!black] (A) circle (0pt) node[color=green!30!black, left, font=\footnotesize] {$A$};
	\fill[color=red!30!black] (r0) circle (0pt) node[color=red!80!black, left, font=\footnotesize] {$N$};
	\fill[color=red!30!black] (r1) circle (0pt) node[color=red!80!black, left, font=\footnotesize] {$M$};
	\fill[color=red!30!black] (l) circle (0pt) node[color=red!80!black, right, font=\footnotesize] {$S$};
	\fill[color=green!50!purple] (psi) circle (3pt) node[color=green!50!purple, above, font=\footnotesize] {$\xi$};
	\fill[color=red!80!black] (ac) circle (3pt) node[color=red!80!black, below, font=\footnotesize] {};
	\end{tikzpicture}}%%popende 
%%%%%%%%%%%%%%%%%%%%%% 
\cong 
%%%%%%%%%%%%%%%%%%%%%%
\tikzzbox{\begin{tikzpicture}[ultra thick,scale=0.75, baseline=-0.1cm]
	\coordinate (r0) at (0,-0.75);
	\coordinate (r1) at (0,0.75);
	\coordinate (l) at (3,0);
	\coordinate (A) at (0,0);
	\coordinate (psi) at (1.5,0);
	\coordinate (ac) at ($(r1)-0.5*(r1)+0.5*(psi)$);
	%
	% lines: 
	\draw[color=orange!50,opacity=0.545] (r0) -- (psi);
	\draw[color=orange!50,opacity=0.545] (r1) -- (psi);
	\draw[color=orange!50,opacity=0.545] (l) -- (psi);
	\draw[color=green!70!black,opacity=0.545] (A) -- (ac);
	%
	%
	% vertices: 
	\fill[color=green!30!black] (A) circle (0pt) node[color=green!30!black, left, font=\footnotesize] {$A$};
	\fill[color=red!30!black] (r0) circle (0pt) node[color=red!80!black, left, font=\footnotesize] {$N$};
	\fill[color=red!30!black] (r1) circle (0pt) node[color=red!80!black, left, font=\footnotesize] {$M$};
	\fill[color=red!30!black] (l) circle (0pt) node[color=red!80!black, right, font=\footnotesize] {$S$};
	\fill[color=green!50!purple] (psi) circle (3pt) node[color=green!50!purple, above, font=\footnotesize] {$\xi$};
	\fill[color=red!80!black] (ac) circle (3pt) node[color=red!80!black, below, font=\footnotesize] {};
	\end{tikzpicture}}%%popende 
%%%%%%%%%%%%%%%%%%%%%% 
\ee
that are subject to analogous coherence constraints. 
Furthermore, the \textsl{corelative product} $M\relpro^{\mathcal A} N$ may be dually defined as a 2-limit involving a structure map $\widetilde\Sigma \colon M \relpro^{\mathcal A} N \lra M\btimes N$, or equivalently by the universal property that there are equivalences 
\begin{align}
\widetilde\Sigma_* \colon \Hom_\tric \big( S, M \relpro^{\mathcal A} N \big) 
	& 
	\stackrel{\cong}{\lra} \textrm{Cobal}\big( S, M\btimes N \big) 
	\nonumber
	\\ 
	\varphi 
	& 
	\lmt \widetilde\Sigma \circ \varphi \, . 
\end{align} 
Analogously to the case of $\relpro_{\mathcal A}$ one finds that there is a 2-morphism $\widetilde\Pi := (\widetilde\Sigma_*)^{-1}(P) \colon M\btimes N \lra M\relpro^{\mathcal A} N$ such that $\widetilde\Pi \circ \widetilde\Sigma \cong 1_{M\relpro^{\mathcal A} N}$, which is also part of a splitting of~$P$. 
Since such splittings are essentially unique, we must have $\Pi\cong\widetilde\Pi$ and $\Sigma \cong\widetilde\Sigma$, and hence 
\be 
M\relpro_{\mathcal A} N \cong M \relpro^{\mathcal A} N \, . 
\ee 
In fact we can, and do, choose $\Pi = \widetilde\Pi$, $\Sigma = \widetilde\Sigma$ and write 
\begin{align}
(\Sigma_*)^{-1} \colon \textrm{Cobal}\big(S, M\btimes N \big) 
& 
\stackrel{\cong}{\lra} \Hom_\tric\big(S, M\relpro^{\mathcal A} N\big) 
\nonumber
\\
%%%%%%%%%%%%%%%%%%%%%%
\tikzzbox{\begin{tikzpicture}[ultra thick,scale=0.75, baseline=-0.1cm]
	\coordinate (r0) at (0,-0.75);
	\coordinate (r1) at (0,0.75);
	\coordinate (l) at (3,0);
	\coordinate (A) at (0,0);
	\coordinate (psi) at (1.5,0);
	%
	% lines: 
	\draw[color=orange!50,opacity=0.545] (r0) -- (psi);
	\draw[color=orange!50,opacity=0.545] (r1) -- (psi);
	\draw[color=orange!50,opacity=0.545] (l) -- (psi);
	%
	%
	% vertices: 
	\fill[color=red!30!black] (r0) circle (0pt) node[color=red!80!black, left, font=\footnotesize] {$N$};
	\fill[color=red!30!black] (r1) circle (0pt) node[color=red!80!black, left, font=\footnotesize] {$M$};
	\fill[color=red!30!black] (l) circle (0pt) node[color=red!80!black, right, font=\footnotesize] {$S$};
	\fill[color=green!50!purple] (psi) circle (3pt) node[color=green!50!purple, above, font=\footnotesize] {$\xi$};
	\end{tikzpicture}}%%popende 
%%%%%%%%%%%%%%%%%%%%%%  
& 
\longmapsto 
%%%%%%%%%%%%%%%%%%%%%%
\tikzzbox{\begin{tikzpicture}[ultra thick,scale=0.75, baseline=-0.1cm]
	\coordinate (r0) at (0,-0.75);
	\coordinate (r1) at (0,0.75);
	\coordinate (l) at (3,0);
	\coordinate (psi) at (1.5,0);
	%
	% area: 
	\fill[color=orange!50,opacity=0.35] ($(r0)+(0,-0.8pt)$) -- ($(psi)+(0,-0.8pt)$) -- ($(r1)+(0,0.8pt)$);
	%
	% lines: 
	\draw[color=orange!50,opacity=0.545] (r0) -- (psi);
	\draw[color=orange!50,opacity=0.545] (r1) -- (psi);
	\draw[color=orange!50,opacity=0.545] (l) -- (psi);
	%
	%
	% vertices: 
	\fill[color=red!30!black] (r0) circle (0pt) node[color=red!80!black, left, font=\footnotesize] {$N$};
	\fill[color=red!30!black] (r1) circle (0pt) node[color=red!80!black, left, font=\footnotesize] {$M$};
	\fill[color=red!30!black] (l) circle (0pt) node[color=red!80!black, right, font=\footnotesize] {$S$};
	\fill[color=red!30!black] ($(r1)-0.5*(r1)+0.5*(r0)+(-0.05,0.05)$) circle (0pt) node[color=red!80!black, left, font=\footnotesize] {$\displaystyle{\relpro^{\mathcal A}}$};
	\fill[color=green!50!purple] (psi) circle (3pt) node[color=green!50!purple, above, font=\footnotesize] {$\xi$};
	\end{tikzpicture}}%%popende 
%%%%%%%%%%%%%%%%%%%%%%   
\, , 
\label{eq:SigmaStarInv}
\end{align} 
so that 
\be 
\label{eq:CobalEquiv}
\Pi \cong (\Sigma_*)^{-1}(P) \colon M\btimes N \lra M\relpro_{\mathcal A} N \, , 
\ee 
which corresponds to~\eqref{eq:PiStarInv} and~\eqref{eq:BalancedEquiv}, respectively. 

As an application that will be useful below in Section~\ref{subsec:properties}, we note that while $\mu\circ\mu^\vee \ncong 1_A$ in general, we have: 

\begin{lemma}
	\label{lem:Filled}
	For any orbifold datum~$\mathcal A$ in~$\tric$, it holds that 
	\be 
	\label{eq:Filled}
	%%%%%%%%%%%%%%%%%%%%%%
	\tikzzbox{\begin{tikzpicture}[ultra thick,scale=0.75, baseline=-0.1cm]
		\coordinate (l) at (-3,0);
		\coordinate (psi) at (-1.5,0);
		\coordinate (r0) at (0,-0.75);
		\coordinate (r1) at (0,0.75);
		\coordinate (rr0) at (0.5,-0.75);
		\coordinate (rr1) at (0.5,0.75);
		\coordinate (md) at (2,0);
		\coordinate (r) at (3.5,0);
		%
		% area: 
		\fill[color=green!50,opacity=0.35] (psi) .. controls +(1,0.5) and +(-0.5,0) .. (r1) -- (rr1) .. controls +(0.5,0) and +(-1,0.5) .. (md) 
		.. controls +(-1,-0.5) and +(0.5,0) .. (rr0) -- (r0)
		.. controls +(-0.5,0) and +(1,-0.5) .. (psi);	
		%
		% lines: 
		\draw[color=green!70!black,opacity=0.545] (psi) .. controls +(1,0.5) and +(-0.5,0) .. (r1) -- (rr1) .. controls +(0.5,0) and +(-1,0.5) .. (md) 
		.. controls +(-1,-0.5) and +(0.5,0) .. (rr0) -- (r0)
		.. controls +(-0.5,0) and +(1,-0.5) .. (psi);
		\draw[color=green!70!black,opacity=0.545] (l) -- (psi);
		\draw[color=green!70!black,opacity=0.545] (md) -- (r);
		%
		% vertices: 
		\fill[color=green!30!black] (psi) circle (3pt) node[color=green!30!black, above, font=\footnotesize] {$\mu$};
		\fill[color=green!30!black] (md) circle (3pt) node[color=green!30!black, above, font=\footnotesize] {$\mu^\vee$};
		\end{tikzpicture}}%%popende 
	%%%%%%%%%%%%%%%%%%%%%%   
	\cong 
	1_A \, . 
	\ee  
\end{lemma}

\begin{proof}
	Since by definition 
	\be 
	%%%%%%%%%%%%%%%%%%%%%%
	\tikzzbox{\begin{tikzpicture}[ultra thick,scale=0.75, baseline=-0.1cm]
		\coordinate (r0) at (0,-0.75);
		\coordinate (r1) at (0,0.75);
		\coordinate (l) at (3,0);
		\coordinate (psi) at (1.5,0);
		%
		% area: 
		\fill[color=green!50,opacity=0.35] ($(r0)+(0,-0.8pt)$) -- ($(psi)+(0,-0.8pt)$) -- ($(r1)+(0,0.8pt)$);
		%
		% lines: 
		\draw[color=green!70!black,opacity=0.545] (r0) -- (psi);
		\draw[color=green!70!black,opacity=0.545] (r1) -- (psi);
		\draw[color=green!70!black,opacity=0.545] (l) -- (psi);
		%
		% vertices: 
		\fill[color=green!30!black] (psi) circle (3pt) node[color=green!30!black, above, font=\footnotesize] {$\mu^\vee$};
		\end{tikzpicture}}%%popende 
	%%%%%%%%%%%%%%%%%%%%%%  
	= 
	(\Sigma_*)^{-1}(\mu^\vee) 
	\, , \quad 
	%%%%%%%%%%%%%%%%%%%%%%
	\tikzzbox{\begin{tikzpicture}[ultra thick,scale=0.75, baseline=-0.1cm]
		\coordinate (r0) at (0,-0.75);
		\coordinate (r1) at (0,0.75);
		\coordinate (l) at (-3,0);
		\coordinate (psi) at (-1.5,0);
		%
		% area: 
		\fill[color=green!50,opacity=0.35] ($(r0)+(0,-0.8pt)$) -- ($(psi)+(0,-0.8pt)$) -- ($(r1)+(0,0.8pt)$);
		%
		% lines: 
		\draw[color=green!70!black,opacity=0.545] (r0) -- (psi);
		\draw[color=green!70!black,opacity=0.545] (r1) -- (psi);
		\draw[color=green!70!black,opacity=0.545] (l) -- (psi);
		%
		% vertices: 
		\fill[color=green!30!black] (psi) circle (3pt) node[color=green!30!black, above, font=\footnotesize] {$\mu$};
		\end{tikzpicture}}%%popende 
	%%%%%%%%%%%%%%%%%%%%%%   
	= 
	(\Pi^*)^{-1}(\mu) \, , 
	\ee 
	the isomorphism~\eqref{eq:Filled} exists if 
	\be 
	(\Sigma_*)^{-1}(\mu^\vee) \cong 
	%%%%%%%%%%%%%%%%%%%%%%
	\tikzzbox{\begin{tikzpicture}[ultra thick,scale=0.75, baseline=-0.1cm]
		\coordinate (r0) at (0,-0.75);
		\coordinate (r1) at (-0.25,0.75);
		\coordinate (l0) at (-3,-0.75);
		\coordinate (l1) at (-3,0.75);
		\coordinate (acm) at (-1.5,-0.75);
		\coordinate (acn) at (-1.5,0.75);
		%
		% area: 
		\fill[color=green!50,opacity=0.35] (acm) -- (l0) -- (l1) -- (acn);
		%
		% lines: 
		\draw[color=green!70!black,opacity=0.545] (l0) -- (r0);
		\draw[color=green!70!black,opacity=0.545] (l1) -- (r1);
		\draw[blue!50!green] ($(acm)+(0,-0.8pt)$) -- ($(acn)+(0,0.8pt)$);
		%
		% vertices: 
		\fill[color=blue!50!green] ($(acn)+(0,0.8pt)$) circle (0pt) node[color=blue!50!green, above, font=\footnotesize] {$\Pi$};
		\fill[color=green!80!black] (r1) circle (3pt) node[color=green!80!black, above, font=\footnotesize] {};
		\end{tikzpicture}}%%popende 
	%%%%%%%%%%%%%%%%%%%%%% 
	\, . 
	\ee 
	But since~$\Sigma_*$ is an equivalence, this is equivalent to 
	\be 
	\mu^\vee \cong 
	%%%%%%%%%%%%%%%%%%%%%%
	\tikzzbox{\begin{tikzpicture}[ultra thick,scale=0.75, baseline=-0.1cm]
		\coordinate (r0) at (0,-0.75);
		\coordinate (r1) at (-0.25,0.75);
		\coordinate (l0) at (-3,-0.75);
		\coordinate (l1) at (-3,0.75);
		\coordinate (acm) at (-1.5,-0.75);
		\coordinate (acn) at (-1.5,0.75);
		%
		% area: 
		\fill[color=green!50,opacity=0.35] ($(acm)+(0.75,-0.8pt)$) -- ($(acm)+(-0.75,-0.8pt)$) -- ($(acn)+(-0.75,0.8pt)$) -- ($(acn)+(0.75,0.8pt)$);
		%
		% lines: 
		\draw[color=green!70!black,opacity=0.545] (l0) -- (r0);
		\draw[color=green!70!black,opacity=0.545] (l1) -- (r1);
		\draw[blue!50!purple] ($(acm)+(-0.75,-0.8pt)$) -- ($(acn)+(-0.75,0.8pt)$);
		\draw[blue!50!green] ($(acm)+(+0.75,-0.8pt)$) -- ($(acn)+(+0.75,0.8pt)$);
		%
		% vertices: 
		\fill[color=blue!50!purple] ($(acn)+(-0.75,0.8pt)$) circle (0pt) node[color=blue!50!purple, above, font=\footnotesize] {$\Sigma$};
		\fill[color=blue!50!green] ($(acn)+(0.75,0.8pt)$) circle (0pt) node[color=blue!50!green, above, font=\footnotesize] {$\Pi$};
		\fill[color=green!80!black] (r1) circle (3pt) node[color=green!80!black, above, font=\footnotesize] {};
		\end{tikzpicture}}%%popende 
	%%%%%%%%%%%%%%%%%%%%%%  
	\cong
	%%%%%%%%%%%%%%%%%%%%%%
	\tikzzbox{\begin{tikzpicture}[ultra thick,scale=0.75, baseline=-0.1cm]
		\coordinate (r0) at (0,-0.75);
		\coordinate (r1) at (-0.25,0.75);
		\coordinate (l0) at (-3,-0.75);
		\coordinate (l1) at (-3,0.75);
		\coordinate (acmm) at (-2,-0.75);
		\coordinate (acnn) at (-2,0.75);
		\coordinate (md) at (-1.5,0);
		\coordinate (u) at (-.75,0);
		%
		% lines: 
		\draw[color=green!70!black,opacity=0.545] (l0) -- (r0);
		\draw[color=green!70!black,opacity=0.545] (l1) -- (r1);
		\draw[color=green!70!black,opacity=0.545] (acmm) -- (md);
		\draw[color=green!70!black,opacity=0.545] (acnn) -- (md);
		\draw[color=green!70!black,opacity=0.545] (u) -- (md);
		%
		% vertices: 
		\fill[color=green!30!black] (acmm) circle (3pt) node[color=red!80!black, below, font=\footnotesize] {};
		\fill[color=green!30!black] (acnn) circle (3pt) node[color=red!80!black, above, font=\footnotesize] {};
		\fill[color=green!30!black] (md) circle (3pt) node[color=green!30!black, below, font=\footnotesize] {};
		\fill[color=green!80!black] (u) circle (3pt) node[color=green!80!black, above, font=\footnotesize] {};
		\fill[color=green!80!black] (r1) circle (3pt) node[color=green!80!black, above, font=\footnotesize] {};
		\end{tikzpicture}}%%popende 
	%%%%%%%%%%%%%%%%%%%%%%  
	\cong 
	%%%%%%%%%%%%%%%%%%%%%%
	\tikzzbox{\begin{tikzpicture}[ultra thick,scale=0.75, baseline=-0.1cm]
		\coordinate (r0) at (0,-0.75);
		\coordinate (r1) at (0,0.75);
		\coordinate (l) at (3,0);
		\coordinate (psi) at (1.5,0);
		%
		% lines: 
		\draw[color=green!70!black,opacity=0.545] (r0) -- (psi);
		\draw[color=green!70!black,opacity=0.545] (r1) -- (psi);
		\draw[color=green!70!black,opacity=0.545] (l) -- (psi);
		%
		% vertices: 
		\fill[color=green!30!black] (psi) circle (3pt) node[color=green!30!black, above, font=\footnotesize] {};
		\end{tikzpicture}}%%popende 
	%%%%%%%%%%%%%%%%%%%%%%  
	= 
	\mu^\vee
	\, . 
	\ee 
\end{proof}

\subsection{Properties}
\label{subsec:properties}

\subsubsection{Adjoints}
\label{subsubsec:Adjoints}

In this section we show that 1- and 2-morphisms in $\ealg(\tric)$ have adjoints. 

\medskip  

Let~$\tric$ be a Gray category with duals as before, and let $\mathcal F = (F,\triangleright_F,\triangleleft_F) \colon \mathcal M \lra \mathcal M'$ be a 2-morphism in $\orb{\tric}$ as in Definitions~\ref{def:2MorphismsEAlg} and~\ref{def:Corb}. 
We set 
\be 
\widetilde{\triangleright}_F := 
%%%%%%%%%%%%%%%%%%%%%% 
\tikzzbox{\begin{tikzpicture}[ultra thick,scale=2.5,color=blue!50!black, baseline=0.3cm, >=stealth, 
	style={x={(-0.6cm,-0.4cm)},y={(1cm,-0.2cm)},z={(0cm,0.9cm)}}]
	%: where to put leftmost T-line: 
	\pgfmathsetmacro{\yy}{0.2}
	\coordinate (T) at (0.5, 0.4, 0);
	\coordinate (L) at (0.5, 0, 0);
	\coordinate (R1) at (0.3, 1, 0);
	\coordinate (R2) at (0.7, 1, 0);
	% top vertices: 
	\coordinate (1T) at (0.5, 0.4, 1);
	\coordinate (1L) at (0.5, 0, 1);
	\coordinate (1R1) at (0.3, 1, );
	\coordinate (1R2) at (0.7, 1, );
	%
	% Poincare dual triangle: 
	\coordinate (p3) at (0.1, 0.1, 0.5);
	\coordinate (p2) at (0.5, 0.95, 0.5);
	\coordinate (p1) at (0.9, 0.1, 0.5);
	%
	% A-plane (front)
	\fill [orange!80,opacity=0.545] (L) -- (T) -- (1T) -- (1L);
	\fill [orange!80,opacity=0.545] (R1) -- (T) -- (1T) -- (1R1);
	%
	% F-line START %%%%%%%%%%%%%%%%%%%%%%
	\draw[postaction={decorate}, decoration={markings,mark=at position 0.5 with {\arrow[draw=blue!50!black]{<}}}] ($(T)+0.8*(L)-0.8*(T)$) -- ($(T)+0.8*(L)-0.8*(T)+(0,0,0.7)$);
	\draw[color=blue!50!black] 
	($(T)+0.8*(L)-0.8*(T)+(0,0,0.7)$) .. controls +(0,0,0.2) and +(0,0,0.2) .. ($(T)+0.2*(L)-0.2*(T)+(0,0,0.7)$);
	\draw[color=blue!50!black] 
	($(T)+0.2*(R1)-0.2*(T)+(0,0,0.3)$) .. controls +(0,0,0.2) and +(0,0,-0.2) .. ($(T)+0.2*(L)-0.2*(T)+(0,0,0.7)$);
	\draw[color=blue!50!black] 
	($(T)+0.2*(R1)-0.2*(T)+(0,0,0.3)$) .. controls +(0,0,-0.2) and +(0,0,-0.2) .. ($(T)+0.8*(R1)-0.8*(T)+(0,0,0.3)$);
	\draw[postaction={decorate}, decoration={markings,mark=at position 0.5 with {\arrow[draw=blue!50!black]{<}}}] ($(T)+0.8*(R1)-0.8*(T)+(0,0,0.3)$) -- ($(T)+0.8*(R1)-0.8*(T)+(0,0,1)$);
	% F-line END %%%%%%%%%%%%%%%%%%%%%%
	%
	\fill [green!50,opacity=0.6] (R2) -- (T) -- (1T) -- (1R2);
	%
	%T-line: 
	\draw[color=red!80!black] (T) -- (1T);
	%
	%labels: 
	\fill[color=black] (0.5,1.14,0.04) circle (0pt) node[left] (0up) { {\scriptsize$a$} };
	\fill[color=black] (0.7,0.5,0.05) circle (0pt) node[left] (0up) { {\scriptsize$a$} };
	\fill[color=black] (0.3,0.5,1.02) circle (0pt) node[left] (0up) { {\scriptsize$b$} };
	\fill[red!80!black] (0.5, 0.38, 0.08) circle (0pt) node[right] (0up) {{\scriptsize$\triangleright_{M'}$}};
	\fill[red!80!black] (0.5, 0.38, 0.8) circle (0pt) node[right] (0up) {{\scriptsize$\triangleright_{M}$}};
	\fill[color=red!80!black] (0.3,0.25,0.07) circle (0pt) node[left] (0up) { {\scriptsize$M'$} };
	\fill[color=red!80!black] (0.15,0.95,0.05) circle (0pt) node[left] (0up) { {\scriptsize$M'$} };
	\fill[color=green!50!black] (0.55,0.95,0.05) circle (0pt) node[left] (0up) { {\scriptsize$A$} };
	\fill[color=red!80!black] (0.55,0.95,1.05) circle (0pt) node[left] (0up) { {\scriptsize$M$} };
	\fill[color=red!80!black] (0.6, 0.15, 0.97) circle (0pt) node (0up) { {\scriptsize$M$} };
	\fill[color=blue!60!black] ($(T)+0.8*(L)-0.8*(T)$) circle (0pt) node[below] (0up) { {\scriptsize$F^\vee$} };
	\fill[color=blue!60!black] ($(T)+0.8*(R1)-0.8*(T)+(0,0,1)$) circle (0pt) node[above] (0up) { {\scriptsize$F^\vee$} };
	\fill[color=blue!60!black] ($(T)+0.53*(1T)-0.53*(T)$) circle (1pt) node[left] (0up) { {\scriptsize$\triangleright_F^{-1}$} };
	%
	% black boundaries: 
	\draw [black,opacity=1, very thin] (1T) -- (1L) -- (L) -- (T);
	\draw [black,opacity=1, very thin] (1T) -- (1R1) -- (R1) -- (T);
	\draw [black,opacity=1, very thin] (1T) -- (1R2) -- (R2) -- (T);
	\end{tikzpicture}}%%popende
%%%%%%%%%%%%%%%%%%%%%% 
\, , \qquad 
\widetilde{\triangleleft}_F := 
%%%%%%%%%%%%%%%%%%%%%% 
\tikzzbox{\begin{tikzpicture}[ultra thick,scale=2.5,color=blue!50!black, baseline=0.3cm, >=stealth, 
	style={x={(-0.6cm,-0.4cm)},y={(1cm,-0.2cm)},z={(0cm,0.9cm)}}]
	%: where to put leftmost T-line: 
	\pgfmathsetmacro{\yy}{0.2}
	\coordinate (T) at (0.5, 0.4, 0);
	\coordinate (L) at (0.5, 0, 0);
	\coordinate (R1) at (0.3, 1, 0);
	\coordinate (R2) at (0.7, 1, 0);
	% top vertices: 
	\coordinate (1T) at (0.5, 0.4, 1);
	\coordinate (1L) at (0.5, 0, 1);
	\coordinate (1R1) at (0.3, 1, );
	\coordinate (1R2) at (0.7, 1, );
	%
	% Poincare dual triangle: 
	\coordinate (p3) at (0.1, 0.1, 0.5);
	\coordinate (p2) at (0.5, 0.95, 0.5);
	\coordinate (p1) at (0.9, 0.1, 0.5);
	%
	% A-plane (front)
	\fill [orange!80,opacity=0.545] (L) -- (T) -- (1T) -- (1L);
	\fill [green!50,opacity=0.545] (R1) -- (T) -- (1T) -- (1R1);
	\fill [orange!80,opacity=0.545] (R2) -- (T) -- (1T) -- (1R2);
	%
	% F-line START %%%%%%%%%%%%%%%%%%%%%%
	\draw[postaction={decorate}, decoration={markings,mark=at position 0.5 with {\arrow[draw=blue!50!black]{<}}}] ($(T)+0.8*(L)-0.8*(T)$) -- ($(T)+0.8*(L)-0.8*(T)+(0,0,0.7)$);
	\draw[color=blue!50!black] 
	($(T)+0.8*(L)-0.8*(T)+(0,0,0.7)$) .. controls +(0,0,0.2) and +(0,0,0.2) .. ($(T)+0.2*(L)-0.2*(T)+(0,0,0.7)$);
	\draw[color=blue!50!black] 
	($(T)+0.2*(R2)-0.2*(T)+(0,0,0.3)$) .. controls +(0,0,0.2) and +(0,0,-0.2) .. ($(T)+0.2*(L)-0.2*(T)+(0,0,0.7)$);
	\draw[color=blue!50!black] 
	($(T)+0.2*(R2)-0.2*(T)+(0,0,0.3)$) .. controls +(0,0,-0.2) and +(0,0,-0.2) .. ($(T)+0.8*(R2)-0.8*(T)+(0,0,0.3)$);
	\draw[postaction={decorate}, decoration={markings,mark=at position 0.5 with {\arrow[draw=blue!50!black]{<}}}] ($(T)+0.8*(R2)-0.8*(T)+(0,0,0.3)$) -- ($(T)+0.8*(R2)-0.8*(T)+(0,0,1)$);
	% F-line END %%%%%%%%%%%%%%%%%%%%%%
	%
	%
	%T-line: 
	\draw[color=red!80!black] (T) -- (1T);
	%
	%labels: 
	\fill[color=black] (0.5,1.14,0.04) circle (0pt) node[left] (0up) { {\scriptsize$b$} };
	\fill[color=black] (0.7,0.5,0.05) circle (0pt) node[left] (0up) { {\scriptsize$a$} };
	\fill[color=black] (0.3,0.5,1.02) circle (0pt) node[left] (0up) { {\scriptsize$b$} };
	\fill[red!80!black] (0.5, 0.38, 0.1) circle (0pt) node[right] (0up) {\hspace{-0.6cm}{\scriptsize$\triangleleft_{M'}$}};
	\fill[red!80!black] (0.5, 0.38, 0.8) circle (0pt) node[right] (0up) {{\scriptsize$\triangleleft_{M}$}};
	\fill[color=red!80!black] (0.3,0.25,0.2) circle (0pt) node[left] (0up) { {\scriptsize$M'$} };
	\fill[color=green!50!black] (0.15,0.95,0.03) circle (0pt) node[left] (0up) { {\scriptsize$B$} };
	\fill[color=red!80!black] (0.55,0.95,0.05) circle (0pt) node[left] (0up) { {\scriptsize$M'$} };
	\fill[color=red!80!black] (0.55,0.8,0.7) circle (0pt) node[left] (0up) { {\scriptsize$M$} };
	\fill[color=red!80!black] (0.6, 0.15, 0.97) circle (0pt) node (0up) { {\scriptsize$M$} };
	\fill[color=blue!60!black] ($(T)+0.8*(L)-0.8*(T)$) circle (0pt) node[below] (0up) { {\scriptsize$F^\vee$} };
	\fill[color=blue!60!black] ($(T)+0.8*(R2)-0.8*(T)+(0,0,1)$) circle (0pt) node[above] (0up) { {\scriptsize$F^\vee$} };
	\fill[color=blue!60!black] ($(T)+0.5*(1T)-0.5*(T)$) circle (1pt) node[left] (0up) { {\scriptsize$\triangleleft_F^{-1}$} };
	%
	% black boundaries: 
	\draw [black,opacity=1, very thin] (1T) -- (1L) -- (L) -- (T);
	\draw [black,opacity=1, very thin] (1T) -- (1R1) -- (R1) -- (T);
	\draw [black,opacity=1, very thin] (1T) -- (1R2) -- (R2) -- (T);
	\end{tikzpicture}}%%popende
%%%%%%%%%%%%%%%%%%%%%% 
\ee 
and observe that 
\be 
\mathcal F^\vee 
= \big(F,\triangleright_F,\triangleleft_F \big)^\vee 
\defi \big( F^\vee, \widetilde{\triangleright}_F, \widetilde{\triangleleft}_F \big) 
\ee 
is a 2-morphism in $\orb{\tric}$. 
Indeed, the defining properties hold thanks to the properties of $\triangleright_F,\triangleleft,_F$ and the Zorro moves for~$F$ in~$\tric$.\footnote{%arXiv_v2: 
	A similar result for modules over separable algebras in a monoidal bicategory is proven in~\cite[Prop.\,2.2.7]{Decoppet_Rig}.}

\begin{lemma}
	\label{lem:Adj2Mor}
	Let $(F,\triangleright_F,\triangleleft_F)$ be a 2-morphism in $\orb{\tric}$ such that $\widetilde{\triangleright}_F, \widetilde{\triangleleft}_F$ are invertible and 
	\be
	\label{eq:CompDual}
	\big(\widetilde{\triangleright}_F\big)^{-1} = 
	%%%%%%%%%%%%%%%%%%%%%% 
	\tikzzbox{\begin{tikzpicture}[ultra thick,scale=2.5,color=blue!50!black, baseline=0.3cm, >=stealth, 
		style={x={(-0.6cm,-0.4cm)},y={(1cm,-0.2cm)},z={(0cm,0.9cm)}}]
		%: where to put leftmost T-line: 
		\pgfmathsetmacro{\yy}{0.2}
		\coordinate (T) at (0.5, 0.4, 0);
		\coordinate (L) at (0.5, 0, 0);
		\coordinate (R1) at (0.3, 1, 0);
		\coordinate (R2) at (0.7, 1, 0);
		% top vertices: 
		\coordinate (1T) at (0.5, 0.4, 1);
		\coordinate (1L) at (0.5, 0, 1);
		\coordinate (1R1) at (0.3, 1, );
		\coordinate (1R2) at (0.7, 1, );
		%
		% Poincare dual triangle: 
		\coordinate (p3) at (0.1, 0.1, 0.5);
		\coordinate (p2) at (0.5, 0.95, 0.5);
		\coordinate (p1) at (0.9, 0.1, 0.5);
		%
		% A-plane (front)
		\fill [orange!80,opacity=0.545] (L) -- (T) -- (1T) -- (1L);
		\fill [orange!80,opacity=0.545] (R1) -- (T) -- (1T) -- (1R1);
		%
		% F-line START %%%%%%%%%%%%%%%%%%%%%%
		\draw[postaction={decorate}, decoration={markings,mark=at position 0.5 with {\arrow[draw=blue!50!black]{>}}}] ($(T)+0.8*(L)-0.8*(T)+(0,0,1)$) -- ($(T)+0.8*(L)-0.8*(T)+(0,0,0.3)$);
		\draw[color=blue!50!black] 
		($(T)+0.8*(L)-0.8*(T)+(0,0,0.3)$) .. controls +(0,0,-0.2) and +(0,0,-0.2) .. ($(T)+0.2*(L)-0.2*(T)+(0,0,0.3)$);
		\draw[color=blue!50!black] 
		($(T)+0.2*(R1)-0.2*(T)+(0,0,0.7)$) .. controls +(0,0,-0.2) and +(0,0,0.2) .. ($(T)+0.2*(L)-0.2*(T)+(0,0,0.3)$);
		\draw[color=blue!50!black] 
		($(T)+0.2*(R1)-0.2*(T)+(0,0,0.7)$) .. controls +(0,0,0.2) and +(0,0,0.2) .. ($(T)+0.8*(R1)-0.8*(T)+(0,0,0.7)$);
		\draw[postaction={decorate}, decoration={markings,mark=at position 0.5 with {\arrow[draw=blue!50!black]{>}}}] ($(T)+0.8*(R1)-0.8*(T)+(0,0,0.7)$) -- ($(T)+0.8*(R1)-0.8*(T)+(0,0,0)$);
		% F-line END %%%%%%%%%%%%%%%%%%%%%%
		%
		\fill [green!50,opacity=0.6] (R2) -- (T) -- (1T) -- (1R2);
		%
		%T-line: 
		\draw[color=red!80!black] (T) -- (1T);
		%
		%labels: 
		%		\fill[color=black] (0.5,1.14,0.04) circle (0pt) node[left] (0up) { {\scriptsize$a$} };
		\fill[color=black] (0.7,0.5,0.05) circle (0pt) node[left] (0up) { {\scriptsize$a$} };
		\fill[color=black] (0.3,0.5,1.02) circle (0pt) node[left] (0up) { {\scriptsize$b$} };
		\fill[red!80!black] (0.5, 0.38, 0.08) circle (0pt) node[right] (0up) {{\scriptsize$\triangleright_{M}$}};
		\fill[red!80!black] (0.5, 0.38, 0.8) circle (0pt) node[right] (0up) {\hspace{-0.6cm}{\scriptsize$\triangleright_{M'}$}};
		\fill[color=red!80!black] (0.45,0.15,0.07) circle (0pt) node[left] (0up) { {\scriptsize$M$}};
		\fill[color=red!80!black] (0.15,0.95,0.85) circle (0pt) node[left] (0up) { {\scriptsize$M'$} };
		\fill[color=green!50!black] (0.55,0.95,0.05) circle (0pt) node[left] (0up) { {\scriptsize$A$} };
		\fill[color=blue!60!black] ($(T)+0.8*(L)-0.8*(T)+(0,0,1)$) circle (0pt) node[above] (0up) { {\scriptsize$F^\vee$} };
		\fill[color=blue!60!black] ($(T)+0.8*(R1)-0.8*(T)+(0,0,0)$) circle (0pt) node[below] (0up) { {\scriptsize$F^\vee$} };
		\fill[color=black] ($(T)+0.47*(1T)-0.47*(T)$) circle (1pt) node[left] (0up) { {\scriptsize$\triangleright_F$} };
		%
		% black boundaries: 
		\draw [black,opacity=1, very thin] (1T) -- (1L) -- (L) -- (T);
		\draw [black,opacity=1, very thin] (1T) -- (1R1) -- (R1) -- (T);
		\draw [black,opacity=1, very thin] (1T) -- (1R2) -- (R2) -- (T);
		\end{tikzpicture}}%%popende
	%%%%%%%%%%%%%%%%%%%%%% 
	\, , \qquad 
	\big(\widetilde{\triangleleft}_F\big)^{-1} = 
	%%%%%%%%%%%%%%%%%%%%%% 
	\tikzzbox{\begin{tikzpicture}[ultra thick,scale=2.5,color=blue!50!black, baseline=0.3cm, >=stealth, 
		style={x={(-0.6cm,-0.4cm)},y={(1cm,-0.2cm)},z={(0cm,0.9cm)}}]
		%: where to put leftmost T-line: 
		\pgfmathsetmacro{\yy}{0.2}
		\coordinate (T) at (0.5, 0.4, 0);
		\coordinate (L) at (0.5, 0, 0);
		\coordinate (R1) at (0.3, 1, 0);
		\coordinate (R2) at (0.7, 1, 0);
		% top vertices: 
		\coordinate (1T) at (0.5, 0.4, 1);
		\coordinate (1L) at (0.5, 0, 1);
		\coordinate (1R1) at (0.3, 1, );
		\coordinate (1R2) at (0.7, 1, );
		%
		% Poincare dual triangle: 
		\coordinate (p3) at (0.1, 0.1, 0.5);
		\coordinate (p2) at (0.5, 0.95, 0.5);
		\coordinate (p1) at (0.9, 0.1, 0.5);
		%
		% A-plane (front)
		\fill [orange!80,opacity=0.545] (L) -- (T) -- (1T) -- (1L);
		\fill [green!50,opacity=0.545] (R1) -- (T) -- (1T) -- (1R1);
		\fill [orange!80,opacity=0.545] (R2) -- (T) -- (1T) -- (1R2);
		%
		% F-line START %%%%%%%%%%%%%%%%%%%%%%
		\draw[postaction={decorate}, decoration={markings,mark=at position 0.5 with {\arrow[draw=blue!50!black]{>}}}] ($(T)+0.8*(L)-0.8*(T)+(0,0,1)$) -- ($(T)+0.8*(L)-0.8*(T)+(0,0,0.3)$);
		\draw[color=blue!50!black] 
		($(T)+0.8*(L)-0.8*(T)+(0,0,0.3)$) .. controls +(0,0,-0.2) and +(0,0,-0.2) .. ($(T)+0.2*(L)-0.2*(T)+(0,0,0.3)$);
		\draw[color=blue!50!black] 
		($(T)+0.2*(R2)-0.2*(T)+(0,0,0.7)$) .. controls +(0,0,-0.2) and +(0,0,0.2) .. ($(T)+0.2*(L)-0.2*(T)+(0,0,0.3)$);
		\draw[color=blue!50!black] 
		($(T)+0.2*(R2)-0.2*(T)+(0,0,0.7)$) .. controls +(0,0,0.2) and +(0,0,0.2) .. ($(T)+0.8*(R2)-0.8*(T)+(0,0,0.7)$);
		\draw[postaction={decorate}, decoration={markings,mark=at position 0.5 with {\arrow[draw=blue!50!black]{>}}}] ($(T)+0.8*(R2)-0.8*(T)+(0,0,0.7)$) -- ($(T)+0.8*(R2)-0.8*(T)+(0,0,0)$);
		% F-line END %%%%%%%%%%%%%%%%%%%%%%
		%
		%
		%T-line: 
		\draw[color=red!80!black] (T) -- (1T);
		%
		%labels: 
		\fill[color=black] (0.5,1.14,0.04) circle (0pt) node[left] (0up) { {\scriptsize$b$} };
		\fill[color=black] (0.7,0.5,0.05) circle (0pt) node[left] (0up) { {\scriptsize$a$} };
		\fill[color=black] (0.3,0.5,1.02) circle (0pt) node[left] (0up) { {\scriptsize$b$} };
		\fill[red!80!black] (0.5, 0.38, 0.1) circle (0pt) node[right] (0up) {\hspace{-0.55cm}{\scriptsize$\triangleleft_{M}$}};
		\fill[red!80!black] (0.5, 0.38, 0.8) circle (0pt) node[right] (0up) {\hspace{-0.6cm}{\scriptsize$\triangleleft_{M'}$}};
		\fill[color=green!50!black] (0.15,0.95,0.03) circle (0pt) node[left] (0up) { {\scriptsize$B$} };
		\fill[color=red!80!black] (0.55,0.95,0.9) circle (0pt) node[left] (0up) { {\scriptsize$M'$} };
		\fill[color=red!80!black] (0.6, 0.15, 0.13) circle (0pt) node (0up) { {\scriptsize$M$} };
		\fill[color=black] ($(T)+0.8*(L)-0.8*(T)+(0,0,1)$) circle (0pt) node[above] (0up) { {\scriptsize$F^\vee$} };
		\fill[color=black] ($(T)+0.8*(R2)-0.8*(T)+(0,0,0)$) circle (0pt) node[below] (0up) { {\scriptsize$F^\vee$} };
		\fill[color=black] ($(T)+0.47*(1T)-0.47*(T)$) circle (1pt) node[left] (0up) { {\scriptsize$\triangleleft_F$} };
		%
		% black boundaries: 
		\draw [black,opacity=1, very thin] (1T) -- (1L) -- (L) -- (T);
		\draw [black,opacity=1, very thin] (1T) -- (1R1) -- (R1) -- (T);
		\draw [black,opacity=1, very thin] (1T) -- (1R2) -- (R2) -- (T);
		\end{tikzpicture}}%%popende
	%%%%%%%%%%%%%%%%%%%%%% 
	\, . 
	\ee 
	Then $(F,\triangleright_F,\triangleleft_F)^\vee$ is left and right adjoint to $(F,\triangleright_F,\triangleleft_F)$ in $\orb{\tric}$, as witnessed by the adjunction 3-morphisms of~$F$ in~$\tric$. 
\end{lemma}

\begin{proof}
	To show that~$\mathcal F^\vee$ is a left adjoint we note that 
	\be 
	%%%%%%%%%%%%%%%%%%%%%% 
	\tikzzbox{\begin{tikzpicture}[ultra thick,scale=2.5,color=blue!50!black, baseline=0.3cm, >=stealth, 
		style={x={(-0.6cm,-0.4cm)},y={(1cm,-0.2cm)},z={(0cm,0.9cm)}}]
		%: where to put leftmost T-line: 
		\pgfmathsetmacro{\yy}{0.2}
		\coordinate (T) at (0.5, 0.4, 0);
		\coordinate (L) at (0.5, 0, 0);
		\coordinate (R1) at (0.3, 1, 0);
		\coordinate (R2) at (0.7, 1, 0);
		% top vertices: 
		\coordinate (1T) at (0.5, 0.4, 1);
		\coordinate (1L) at (0.5, 0, 1);
		\coordinate (1R1) at (0.3, 1, );
		\coordinate (1R2) at (0.7, 1, );
		%
		% Poincare dual triangle: 
		\coordinate (p3) at (0.1, 0.1, 0.5);
		\coordinate (p2) at (0.5, 0.95, 0.5);
		\coordinate (p1) at (0.9, 0.1, 0.5);
		%
		% A-plane (front)
		\fill [orange!80,opacity=0.545] (L) -- (T) -- (1T) -- (1L);
		\fill [green!50,opacity=0.545] (R1) -- (T) -- (1T) -- (1R1);
		\fill [orange!80,opacity=0.545] (R2) -- (T) -- (1T) -- (1R2);
		%
		% F-line START %%%%%%%%%%%%%%%%%%%%%%
		\draw ($(T)+0.2*(L)-0.2*(T)+(0,0,0)$) .. controls +(0,0,0.2) and +(0,0,-0.2) ..
		($(T)+0.8*(R2)-0.8*(T)+(0,0,0.7)$);
		\draw ($(T)+0.8*(R2)-0.8*(T)+(0,0,0.7)$) .. controls +(0,0,0.2) and +(0,0,0.2) .. ($(T)+0.2*(R2)-0.2*(T)+(0,0,0.7)$);
		\draw[postaction={decorate}, decoration={markings,mark=at position 0.7 with {\arrow[draw=blue!50!black]{>}}}] ($(T)+0.2*(R2)-0.2*(T)+(0,0,0.7)$) .. controls +(0,0,-0.2) and +(0,0,0.2) .. ($(T)+0.8*(L)-0.8*(T)+(0,0,0)$);
		% F-line END %%%%%%%%%%%%%%%%%%%%%%
		%
		%
		%T-line: 
		\draw[color=red!80!black] (T) -- (1T);
		%
		%labels: 
		\fill[color=black] ($(T)+0.8*(L)-0.8*(T)$) circle (0pt) node[below] (0up) { {\scriptsize$F^\vee$} };
		\fill[color=black] ($(T)+0.2*(L)-0.2*(T)$) circle (0pt) node[below] (0up) { {\scriptsize$F$} };
		\fill[color=black] ($(T)+0.47*(1T)-0.47*(T)$) circle (1pt) node[left] (0up) { {\scriptsize$\widetilde\triangleleft_F$} };
		\fill[color=black] ($(T)+0.16*(1T)-0.16*(T)$) circle (1pt) node[left] (0up) { {\scriptsize$\triangleleft_F$} };
		%
		% black boundaries: 
		\draw [black,opacity=1, very thin] (1T) -- (1L) -- (L) -- (T);
		\draw [black,opacity=1, very thin] (1T) -- (1R1) -- (R1) -- (T);
		\draw [black,opacity=1, very thin] (1T) -- (1R2) -- (R2) -- (T);
		\end{tikzpicture}}%%popende
	%%%%%%%%%%%%%%%%%%%%%% 
	=
	%%%%%%%%%%%%%%%%%%%%%% 
	\tikzzbox{\begin{tikzpicture}[ultra thick,scale=2.5,color=blue!50!black, baseline=0.3cm, >=stealth, 
		style={x={(-0.6cm,-0.4cm)},y={(1cm,-0.2cm)},z={(0cm,0.9cm)}}]
		%: where to put leftmost T-line: 
		\pgfmathsetmacro{\yy}{0.2}
		\coordinate (T) at (0.5, 0.4, 0);
		\coordinate (L) at (0.5, 0, 0);
		\coordinate (R1) at (0.3, 1, 0);
		\coordinate (R2) at (0.7, 1, 0);
		% top vertices: 
		\coordinate (1T) at (0.5, 0.4, 1);
		\coordinate (1L) at (0.5, 0, 1);
		\coordinate (1R1) at (0.3, 1, );
		\coordinate (1R2) at (0.7, 1, );
		%
		% Poincare dual triangle: 
		\coordinate (p3) at (0.1, 0.1, 0.5);
		\coordinate (p2) at (0.5, 0.95, 0.5);
		\coordinate (p1) at (0.9, 0.1, 0.5);
		%
		% A-plane (front)
		\fill [orange!80,opacity=0.545] (L) -- (T) -- (1T) -- (1L);
		\fill [green!50,opacity=0.545] (R1) -- (T) -- (1T) -- (1R1);
		\fill [orange!80,opacity=0.545] (R2) -- (T) -- (1T) -- (1R2);
		%
		% F-line START %%%%%%%%%%%%%%%%%%%%%%
		\draw ($(T)+0.2*(L)-0.2*(T)$) .. controls +(0,0,0.2) and +(0,0,-0.2) ..
		($(T)+0.8*(R2)-0.8*(T)+(0,0,0.3)$);
		\draw ($(T)+0.8*(R2)-0.8*(T)+(0,0,0.3)$) -- ($(T)+0.8*(R2)-0.8*(T)+(0,0,0.5)$);
		\draw ($(T)+0.8*(R2)-0.8*(T)+(0,0,0.5)$) .. controls +(0,0,0.2) and +(0,0,0.2) .. ($(T)+0.5*(R2)-0.5*(T)+(0,0,0.5)$);
		\draw ($(T)+0.5*(R2)-0.5*(T)+(0,0,0.5)$) .. controls +(0,0,-0.2) and +(0,0,-0.2) .. ($(T)+0.2*(R2)-0.2*(T)+(0,0,0.5)$);
		\draw ($(T)+0.2*(R2)-0.2*(T)+(0,0,0.5)$) .. controls +(0,0,0.2) and +(0,0,-0.2) .. ($(T)+0.2*(L)-0.2*(T)+(0,0,0.7)$);
		\draw ($(T)+0.2*(L)-0.2*(T)+(0,0,0.7)$) .. controls +(0,0,0.2) and +(0,0,0.2) .. ($(T)+0.8*(L)-0.8*(T)+(0,0,0.7)$);
		\draw[postaction={decorate}, decoration={markings,mark=at position 0.7 with {\arrow[draw=blue!50!black]{>}}}] ($(T)+0.8*(L)-0.8*(T)+(0,0,0.7)$) .. controls +(0,0,-0.2) and +(0,0,0.2) .. ($(T)+0.8*(L)-0.8*(T)+(0,0,0)$);
		% F-line END %%%%%%%%%%%%%%%%%%%%%%
		%
		%
		%T-line: 
		\draw[color=red!80!black] (T) -- (1T);
		%
		%labels: 
		\fill[color=black] ($(T)+0.8*(L)-0.8*(T)$) circle (0pt) node[below] (0up) { {\scriptsize$F^\vee$} };
		\fill[color=black] ($(T)+0.2*(L)-0.2*(T)$) circle (0pt) node[below] (0up) { {\scriptsize$F$} };
		\fill[color=black] ($(T)+0.59*(1T)-0.59*(T)$) circle (1pt) node[right] (0up) { {\scriptsize$\triangleleft_F^{-1}$} };
		\fill[color=black] ($(T)+0.09*(1T)-0.09*(T)$) circle (1pt) node[left] (0up) { {\scriptsize$\triangleleft_F$} };
		%
		% black boundaries: 
		\draw [black,opacity=1, very thin] (1T) -- (1L) -- (L) -- (T);
		\draw [black,opacity=1, very thin] (1T) -- (1R1) -- (R1) -- (T);
		\draw [black,opacity=1, very thin] (1T) -- (1R2) -- (R2) -- (T);
		\end{tikzpicture}}%%popende
	%%%%%%%%%%%%%%%%%%%%%% 
	=
	%%%%%%%%%%%%%%%%%%%%%% 
	\tikzzbox{\begin{tikzpicture}[ultra thick,scale=2.5,color=blue!50!black, baseline=0.3cm, >=stealth, 
		style={x={(-0.6cm,-0.4cm)},y={(1cm,-0.2cm)},z={(0cm,0.9cm)}}]
		%: where to put leftmost T-line: 
		\pgfmathsetmacro{\yy}{0.2}
		\coordinate (T) at (0.5, 0.4, 0);
		\coordinate (L) at (0.5, 0, 0);
		\coordinate (R1) at (0.3, 1, 0);
		\coordinate (R2) at (0.7, 1, 0);
		% top vertices: 
		\coordinate (1T) at (0.5, 0.4, 1);
		\coordinate (1L) at (0.5, 0, 1);
		\coordinate (1R1) at (0.3, 1, );
		\coordinate (1R2) at (0.7, 1, );
		%
		% Poincare dual triangle: 
		\coordinate (p3) at (0.1, 0.1, 0.5);
		\coordinate (p2) at (0.5, 0.95, 0.5);
		\coordinate (p1) at (0.9, 0.1, 0.5);
		%
		% A-plane (front)
		\fill [orange!80,opacity=0.545] (L) -- (T) -- (1T) -- (1L);
		\fill [green!50,opacity=0.545] (R1) -- (T) -- (1T) -- (1R1);
		\fill [orange!80,opacity=0.545] (R2) -- (T) -- (1T) -- (1R2);
		%
		% F-line START %%%%%%%%%%%%%%%%%%%%%%
		\draw ($(T)+0.2*(L)-0.2*(T)$) .. controls +(0,0,0.2) and +(0,0,-0.2) ..
		($(T)+0.2*(R2)-0.2*(T)+(0,0,0.3)$);
		\draw ($(T)+0.2*(R2)-0.2*(T)+(0,0,0.3)$) -- ($(T)+0.2*(R2)-0.2*(T)+(0,0,0.5)$);
		\draw ($(T)+0.2*(R2)-0.2*(T)+(0,0,0.5)$) .. controls +(0,0,0.2) and +(0,0,-0.2) .. ($(T)+0.2*(L)-0.2*(T)+(0,0,0.7)$);
		\draw ($(T)+0.2*(L)-0.2*(T)+(0,0,0.7)$) .. controls +(0,0,0.2) and +(0,0,0.2) .. ($(T)+0.8*(L)-0.8*(T)+(0,0,0.7)$);
		\draw[postaction={decorate}, decoration={markings,mark=at position 0.7 with {\arrow[draw=blue!50!black]{>}}}] ($(T)+0.8*(L)-0.8*(T)+(0,0,0.7)$) .. controls +(0,0,-0.2) and +(0,0,0.2) .. ($(T)+0.8*(L)-0.8*(T)+(0,0,0)$);
		% F-line END %%%%%%%%%%%%%%%%%%%%%%
		%
		%
		%T-line: 
		\draw[color=red!80!black] (T) -- (1T);
		%
		%labels: 
		\fill[color=black] ($(T)+0.8*(L)-0.8*(T)$) circle (0pt) node[below] (0up) { {\scriptsize$F^\vee$} };
		\fill[color=black] ($(T)+0.2*(L)-0.2*(T)$) circle (0pt) node[below] (0up) { {\scriptsize$F$} };
		\fill[color=black] ($(T)+0.59*(1T)-0.59*(T)$) circle (1pt) node[right] (0up) { {\scriptsize$\triangleleft_F^{-1}$} };
		\fill[color=black] ($(T)+0.14*(1T)-0.14*(T)$) circle (1pt) node[left] (0up) { {\scriptsize$\triangleleft_F$} };
		%
		% black boundaries: 
		\draw [black,opacity=1, very thin] (1T) -- (1L) -- (L) -- (T);
		\draw [black,opacity=1, very thin] (1T) -- (1R1) -- (R1) -- (T);
		\draw [black,opacity=1, very thin] (1T) -- (1R2) -- (R2) -- (T);
		\end{tikzpicture}}%%popende
	%%%%%%%%%%%%%%%%%%%%%% 
	=
	%%%%%%%%%%%%%%%%%%%%%% 
	\tikzzbox{\begin{tikzpicture}[ultra thick,scale=2.5,color=blue!50!black, baseline=0.3cm, >=stealth, 
		style={x={(-0.6cm,-0.4cm)},y={(1cm,-0.2cm)},z={(0cm,0.9cm)}}]
		%: where to put leftmost T-line: 
		\pgfmathsetmacro{\yy}{0.2}
		\coordinate (T) at (0.5, 0.4, 0);
		\coordinate (L) at (0.5, 0, 0);
		\coordinate (R1) at (0.3, 1, 0);
		\coordinate (R2) at (0.7, 1, 0);
		% top vertices: 
		\coordinate (1T) at (0.5, 0.4, 1);
		\coordinate (1L) at (0.5, 0, 1);
		\coordinate (1R1) at (0.3, 1, );
		\coordinate (1R2) at (0.7, 1, );
		%
		% Poincare dual triangle: 
		\coordinate (p3) at (0.1, 0.1, 0.5);
		\coordinate (p2) at (0.5, 0.95, 0.5);
		\coordinate (p1) at (0.9, 0.1, 0.5);
		%
		% A-plane (front)
		\fill [orange!80,opacity=0.545] (L) -- (T) -- (1T) -- (1L);
		\fill [green!50,opacity=0.545] (R1) -- (T) -- (1T) -- (1R1);
		\fill [orange!80,opacity=0.545] (R2) -- (T) -- (1T) -- (1R2);
		%
		% F-line START %%%%%%%%%%%%%%%%%%%%%%
		\draw ($(T)+0.2*(L)-0.2*(T)$) -- ($(T)+0.2*(L)-0.2*(T)+(0,0,0.4)$);
		\draw ($(T)+0.2*(L)-0.2*(T)+(0,0,0.4)$) .. controls +(0,0,0.2) and +(0,0,0.2) .. ($(T)+0.8*(L)-0.8*(T)+(0,0,0.4)$);
		\draw[postaction={decorate}, decoration={markings,mark=at position 0.7 with {\arrow[draw=blue!50!black]{>}}}] ($(T)+0.8*(L)-0.8*(T)+(0,0,0.4)$) -- ($(T)+0.8*(L)-0.8*(T)+(0,0,0)$);
		% F-line END %%%%%%%%%%%%%%%%%%%%%%
		%
		%
		%T-line: 
		\draw[color=red!80!black] (T) -- (1T);
		%
		%labels: 
		\fill[color=black] ($(T)+0.8*(L)-0.8*(T)$) circle (0pt) node[below] (0up) { {\scriptsize$F^\vee$} };
		\fill[color=black] ($(T)+0.2*(L)-0.2*(T)$) circle (0pt) node[below] (0up) { {\scriptsize$F$} };
		%
		% black boundaries: 
		\draw [black,opacity=1, very thin] (1T) -- (1L) -- (L) -- (T);
		\draw [black,opacity=1, very thin] (1T) -- (1R1) -- (R1) -- (T);
		\draw [black,opacity=1, very thin] (1T) -- (1R2) -- (R2) -- (T);
		\end{tikzpicture}}%%popende
	%%%%%%%%%%%%%%%%%%%%%% 
	\, , 
	\ee 
	and a similar identity holds for~$\widetilde{\triangleright}_F$. 
	Hence $\ev_F$ is indeed a 3-morphism in $\ealg(\tric)$. 
	Compatibility of $\coev_F$ and $\widetilde{\triangleleft}_F, \widetilde{\triangleright}_F$ is checked analogously, and the Zorro moves in $\orb{\tric}$ are identical to those in~$\tric$. 
	Note that we do not need the assumption~\eqref{eq:CompDual} for~$\mathcal F^\vee$ to be left adjoint to~$\mathcal F$. 
	
	To show that~$\mathcal F^\vee$ is also right adjoint to~$\mathcal F$, we compute 
	\be 
	%%%%%%%%%%%%%%%%%%%%%% 
	\tikzzbox{\begin{tikzpicture}[ultra thick,scale=2.5,color=blue!50!black, baseline=0.3cm, >=stealth, 
		style={x={(-0.6cm,-0.4cm)},y={(1cm,-0.2cm)},z={(0cm,0.9cm)}}]
		%: where to put leftmost T-line: 
		\pgfmathsetmacro{\yy}{0.2}
		\coordinate (T) at (0.5, 0.4, 0);
		\coordinate (L) at (0.5, 0, 0);
		\coordinate (R1) at (0.3, 1, 0);
		\coordinate (R2) at (0.7, 1, 0);
		% top vertices: 
		\coordinate (1T) at (0.5, 0.4, 1);
		\coordinate (1L) at (0.5, 0, 1);
		\coordinate (1R1) at (0.3, 1, );
		\coordinate (1R2) at (0.7, 1, );
		%
		% Poincare dual triangle: 
		\coordinate (p3) at (0.1, 0.1, 0.5);
		\coordinate (p2) at (0.5, 0.95, 0.5);
		\coordinate (p1) at (0.9, 0.1, 0.5);
		%
		% A-plane (front)
		\fill [orange!80,opacity=0.545] (L) -- (T) -- (1T) -- (1L);
		\fill [green!50,opacity=0.545] (R1) -- (T) -- (1T) -- (1R1);
		\fill [orange!80,opacity=0.545] (R2) -- (T) -- (1T) -- (1R2);
		%
		% F-line START %%%%%%%%%%%%%%%%%%%%%%
		\draw ($(T)+0.3*(R2)-0.3*(T)+(0,0,1)$) .. controls +(0,0,-0.2) and +(0,0,0.2) .. ($(T)+0.8*(L)-0.8*(T)+(0,0,0.3)$);
		\draw[color=blue!50!black] 
		($(T)+0.8*(L)-0.8*(T)+(0,0,0.3)$) .. controls +(0,0,-0.2) and +(0,0,-0.2) .. ($(T)+0.2*(L)-0.2*(T)+(0,0,0.3)$);
		\draw[postaction={decorate}, decoration={markings,mark=at position 0.7 with {\arrow[draw=blue!50!black]{>}}}]
		($(T)+0.2*(L)-0.2*(T)+(0,0,0.3)$) .. controls +(0,0,0.2) and +(0,0,-0.2) .. ($(T)+0.8*(R2)-0.8*(T)+(0,0,1)$);
		% F-line END %%%%%%%%%%%%%%%%%%%%%%
		%
		%
		%T-line: 
		\draw[color=red!80!black] (T) -- (1T);
		%
		%labels: 
		\fill[color=black] ($(T)+0.3*(R2)-0.3*(T)+(0,0,1)$) circle (0pt) node[above] (0up) { {\scriptsize$F^\vee$} };
		\fill[color=black] ($(T)+0.8*(R2)-0.8*(T)+(0,0,1)$) circle (0pt) node[above] (0up) { {\scriptsize$F$} };
		\fill[color=black] ($(T)+0.46*(1T)-0.46*(T)$) circle (1pt) node[left] (0up) { \hspace{0.1cm}{\scriptsize$\triangleleft_F$} };
		\fill[color=black] ($(T)+0.72*(1T)-0.72*(T)$) circle (1pt) node[left] (0up) { \hspace{0.1cm}{\scriptsize$\widetilde\triangleleft_F$} };
		%
		% black boundaries: 
		\draw [black,opacity=1, very thin] (1T) -- (1L) -- (L) -- (T);
		\draw [black,opacity=1, very thin] (1T) -- (1R1) -- (R1) -- (T);
		\draw [black,opacity=1, very thin] (1T) -- (1R2) -- (R2) -- (T);
		\end{tikzpicture}}%%popende
	%%%%%%%%%%%%%%%%%%%%%%  
	=
	%%%%%%%%%%%%%%%%%%%%%% 
	\tikzzbox{\begin{tikzpicture}[ultra thick,scale=2.5,color=blue!50!black, baseline=0.3cm, >=stealth, 
		style={x={(-0.6cm,-0.4cm)},y={(1cm,-0.2cm)},z={(0cm,0.9cm)}}]
		%: where to put leftmost T-line: 
		\pgfmathsetmacro{\yy}{0.2}
		\coordinate (T) at (0.5, 0.4, 0);
		\coordinate (L) at (0.5, 0, 0);
		\coordinate (R1) at (0.3, 1, 0);
		\coordinate (R2) at (0.7, 1, 0);
		% top vertices: 
		\coordinate (1T) at (0.5, 0.4, 1);
		\coordinate (1L) at (0.5, 0, 1);
		\coordinate (1R1) at (0.3, 1, );
		\coordinate (1R2) at (0.7, 1, );
		%
		% Poincare dual triangle: 
		\coordinate (p3) at (0.1, 0.1, 0.5);
		\coordinate (p2) at (0.5, 0.95, 0.5);
		\coordinate (p1) at (0.9, 0.1, 0.5);
		%
		% A-plane (front)
		\fill [orange!80,opacity=0.545] (L) -- (T) -- (1T) -- (1L);
		\fill [green!50,opacity=0.545] (R1) -- (T) -- (1T) -- (1R1);
		\fill [orange!80,opacity=0.545] (R2) -- (T) -- (1T) -- (1R2);
		%
		% F-line START %%%%%%%%%%%%%%%%%%%%%%
		\draw ($(T)+0.4*(R2)-0.4*(T)+(0,0,1)$) -- ($(T)+0.4*(R2)-0.4*(T)+(0,0,0.7)$);
		\draw ($(T)+0.4*(R2)-0.4*(T)+(0,0,0.7)$) .. controls +(0,0,-0.1) and +(0,0,-0.1) .. ($(T)+0.1*(R2)-0.1*(T)+(0,0,0.7)$);
		\draw ($(T)+0.1*(R2)-0.1*(T)+(0,0,0.7)$) .. controls +(0,0,0.1) and +(0,0,-0.1) .. ($(T)+0.1*(L)-0.1*(T)+(0,0,0.8)$);
		\draw ($(T)+0.1*(L)-0.1*(T)+(0,0,0.8)$) .. controls +(0,0,0.1) and +(0,0,0.1) .. ($(T)+0.5*(L)-0.5*(T)+(0,0,0.8)$) -- ($(T)+0.5*(L)-0.5*(T)+(0,0,0.2)$);
		\draw($(T)+0.5*(L)-0.5*(T)+(0,0,0.2)$) .. controls +(0,0,-0.1) and +(0,0,-0.1) .. ($(T)+0.1*(L)-0.1*(T)+(0,0,0.2)$);
		\draw ($(T)+0.1*(L)-0.1*(T)+(0,0,0.2)$) .. controls +(0,0,0.1) and +(0,0,-0.1) .. ($(T)+0.1*(R2)-0.1*(T)+(0,0,0.3)$);
		\draw ($(T)+0.1*(R2)-0.1*(T)+(0,0,0.3)$) .. controls +(0,0,0.1) and +(0,0,0.1) .. ($(T)+0.4*(R2)-0.4*(T)+(0,0,0.3)$);
		\draw ($(T)+0.4*(R2)-0.4*(T)+(0,0,0.3)$) .. controls +(0,0,-0.1) and +(0,0,-0.1) .. ($(T)+0.8*(R2)-0.8*(T)+(0,0,0.3)$);
		\draw[postaction={decorate}, decoration={markings,mark=at position 0.7 with {\arrow[draw=blue!50!black]{>}}}]
		($(T)+0.8*(R2)-0.8*(T)+(0,0,0.3)$) -- ($(T)+0.8*(R2)-0.8*(T)+(0,0,1)$);
		% F-line END %%%%%%%%%%%%%%%%%%%%%%
		%
		%
		%T-line: 
		\draw[color=red!80!black] (T) -- (1T);
		%
		%labels: 
		\fill[color=black] ($(T)+0.4*(R2)-0.4*(T)+(0,0,1)$) circle (0pt) node[above] (0up) { {\scriptsize$F^\vee$} };
		\fill[color=black] ($(T)+0.8*(R2)-0.8*(T)+(0,0,1)$) circle (0pt) node[above] (0up) { {\scriptsize$F$} };
		\fill[color=black] ($(T)+0.27*(1T)-0.27*(T)$) circle (1pt) node[right] (0up) { \hspace{-0.68cm} {\scriptsize$\triangleleft_F$} };
		\fill[color=black] ($(T)+0.74*(1T)-0.74*(T)$) circle (1pt) node[right] (0up) { };
		\fill[color=black] ($(T)+0.74*(1T)-0.74*(T)+(0,0,0.05)$) circle (0pt) node[right] (0up) { \hspace{-0.11cm}{\scriptsize$\triangleleft_F^{-1}$} };
		%
		% black boundaries: 
		\draw [black,opacity=1, very thin] (1T) -- (1L) -- (L) -- (T);
		\draw [black,opacity=1, very thin] (1T) -- (1R1) -- (R1) -- (T);
		\draw [black,opacity=1, very thin] (1T) -- (1R2) -- (R2) -- (T);
		\end{tikzpicture}}%%popende
	%%%%%%%%%%%%%%%%%%%%%%  
	=
	%%%%%%%%%%%%%%%%%%%%%% 
	\tikzzbox{\begin{tikzpicture}[ultra thick,scale=2.5,color=blue!50!black, baseline=0.3cm, >=stealth, 
		style={x={(-0.6cm,-0.4cm)},y={(1cm,-0.2cm)},z={(0cm,0.9cm)}}]
		%: where to put leftmost T-line: 
		\pgfmathsetmacro{\yy}{0.2}
		\coordinate (T) at (0.5, 0.4, 0);
		\coordinate (L) at (0.5, 0, 0);
		\coordinate (R1) at (0.3, 1, 0);
		\coordinate (R2) at (0.7, 1, 0);
		% top vertices: 
		\coordinate (1T) at (0.5, 0.4, 1);
		\coordinate (1L) at (0.5, 0, 1);
		\coordinate (1R1) at (0.3, 1, );
		\coordinate (1R2) at (0.7, 1, );
		%
		% Poincare dual triangle: 
		\coordinate (p3) at (0.1, 0.1, 0.5);
		\coordinate (p2) at (0.5, 0.95, 0.5);
		\coordinate (p1) at (0.9, 0.1, 0.5);
		%
		% A-plane (front)
		\fill [orange!80,opacity=0.545] (L) -- (T) -- (1T) -- (1L);
		\fill [green!50,opacity=0.545] (R1) -- (T) -- (1T) -- (1R1);
		\fill [orange!80,opacity=0.545] (R2) -- (T) -- (1T) -- (1R2);
		%
		% F-line START %%%%%%%%%%%%%%%%%%%%%%
		\draw ($(T)+0.4*(R2)-0.4*(T)+(0,0,1)$) -- ($(T)+0.4*(R2)-0.4*(T)+(0,0,0.7)$);
		\draw ($(T)+0.4*(R2)-0.4*(T)+(0,0,0.7)$) .. controls +(0,0,-0.2) and +(0,0,-0.2) .. ($(T)+0.8*(R2)-0.8*(T)+(0,0,0.7)$);
		\draw[postaction={decorate}, decoration={markings,mark=at position 0.7 with {\arrow[draw=blue!50!black]{>}}}]
		($(T)+0.8*(R2)-0.8*(T)+(0,0,0.7)$) -- ($(T)+0.8*(R2)-0.8*(T)+(0,0,1)$);
		% F-line END %%%%%%%%%%%%%%%%%%%%%%
		%
		%
		%T-line: 
		\draw[color=red!80!black] (T) -- (1T);
		%
		%labels: 
		\fill[color=black] ($(T)+0.4*(R2)-0.4*(T)+(0,0,1)$) circle (0pt) node[above] (0up) { {\scriptsize$F^\vee$} };
		\fill[color=black] ($(T)+0.8*(R2)-0.8*(T)+(0,0,1)$) circle (0pt) node[above] (0up) { {\scriptsize$F$} };
		%
		% black boundaries: 
		\draw [black,opacity=1, very thin] (1T) -- (1L) -- (L) -- (T);
		\draw [black,opacity=1, very thin] (1T) -- (1R1) -- (R1) -- (T);
		\draw [black,opacity=1, very thin] (1T) -- (1R2) -- (R2) -- (T);
		\end{tikzpicture}}%%popende
	%%%%%%%%%%%%%%%%%%%%%%  
	\ee 
	where we used a Zorro move and~\eqref{eq:CompDual} in the first step, and the fact that $\widetilde\triangleleft_F^{-1}$ is left inverse to~$\widetilde\triangleleft_F$ in the second step. 
	This proves compatibility of $\tcoev_F$ and~$\widetilde{\triangleleft}_F$, and the other three compatibility conditions are checked analogously. 
\end{proof}

We note that the conditions in~\eqref{eq:CompDual} are precisely the conditions~\eqref{eq:T6} and~\eqref{eq:T7}. 
These relations generalise \cite[(T6)\,\&\,(T7)]{CMRSS1}, which in turn have those in \cite[Fig.\,3.1]{MuleRunk} as a motivating special case. 

\medskip 

We now turn to adjoints of 1-morphisms 
$
\mathcal M 
	= (M, \triangleright_M, \triangleleft_M, u^{\textrm{r}}_M, u^{\textrm{r}}_M, \alpha^{\textrm{l}}_M,\alpha^{\textrm{m}}_M, \alpha^{\textrm{l}}_M)
	\colon \mathcal A \lra \mathcal B  
$ 
in $\orb{\tric}$, cf.\ Definitions~\ref{def:1MorphismsEAlg} and~\ref{def:Corb}. 
To lift the adjoint~$M^\#$ of~$M$ in~$\tric$ to an adjoint~$\mathcal M^\#$ of~$\mathcal M$ in $\orb{\tric}$, we first set 
\be 
\begin{tikzcd}[column sep=5em, row sep=0.5em]
%%%%%%%%%%%%%%%%%%%%%%
\tikzzbox{
}%%popende 
%%%%%%%%%%%%%%%%%%%%%%  
\, .
\end{tikzcd}
\ee 
Here the 3-isomorphisms $\mathcal P_{\mathcal A}, \mathcal P_{\mathcal B}$ are obtained from the 3-isomorphisms $\Theta_{\triangleleft_M}, \Theta_{\triangleright_M}$ (recall~\eqref{eq:ThetaThetaTheta}) and the cusp isomorphisms for $A,B,M$. 
It then follows from the relations~\eqref{eq:ThetaThetaTheta} that the diagram
\be 
\label{eq:PPP}
\begin{tikzcd}[column sep=5em, row sep=2em]
%%%%%%%%%%%%%%%%%%%%%%
\tikzzbox{%
}%%popende 
%%%%%%%%%%%%%%%%%%%%%%
\end{tikzcd}
\ee 
commutes, and similary for the $\mathcal A$-action. 

Next we define 
\be 
\begin{tikzcd}[column sep=5em, row sep=0.5em]
u_{M^\#}^{\textrm{l}} 
:= 
\Bigg(
%%%%%%%%%%%%%%%%%%%%%%
\tikzzbox{%
}%%popende 
%%%%%%%%%%%%%%%%%%%%%%   
\Bigg) , 
\end{tikzcd}
\ee
\be 
\begin{tikzcd}[column sep=5em, row sep=0.5em]
u_{M^\#}^{\textrm{r}} 
:= 
\Bigg(
%%%%%%%%%%%%%%%%%%%%%%
\tikzzbox{%
}%%popende 
%%%%%%%%%%%%%%%%%%%%%%   
\Bigg)  
\end{tikzcd}
\ee 
and 
\be 
\alpha_{M^\#}^{\textrm{l}} 
:= 
\left(
\begin{tikzcd}[column sep=4.5em, row sep=0.5em]
%%%%%%%%%%%%%%%%%%%%%%
\tikzzbox{%
}%%popende 
%%%%%%%%%%%%%%%%%%%%%%     
\end{tikzcd}
\right) .
\ee 
The map~$\alpha_{M^\#}^{\textrm{r}}$ is defined analogously, and we set 
\be 
\alpha_{M^\#}^{\textrm{m}} 
:= 
\left(
\begin{tikzcd}[column sep=3em, row sep=0.5em]
%%%%%%%%%%%%%%%%%%%%%%
\tikzzbox{%
}%%popende 
%%%%%%%%%%%%%%%%%%%%%%     
\end{tikzcd}
\right) 
\ee 
where here and below we do not indicate the use of cusp isomorphisms and tensorators unless we think it is helpful (as in the fourth step above, to produce something that is manifestly the domain of $\mathcal P_{\mathcal B}^{-1}$). 

\begin{lemma}
	\label{lem:AdjointsFor1Morphisms}
	Let 
	$
	\mathcal M 
	= (M, \triangleright_M, \triangleleft_M, u^{\textrm{l}}_M, u^{\textrm{r}}_M, \alpha^{\textrm{l}}_M,\alpha^{\textrm{m}}_M, \alpha^{\textrm{r}}_M)
		\colon \mathcal A \lra \mathcal B
	$ 
	be a 1-morphism in $\orb{\tric}$. 
	Then 
	\be 
	\mathcal M^\# 
		:= \big(M^\#, \triangleright_{M^\#}, \triangleleft_{M^\#}, u^{\textrm{l}}_{M^\#}, u^{\textrm{r}}_{M^\#}, \alpha^{\textrm{l}}_{M^\#},\alpha^{\textrm{m}}_{M^\#}, \alpha^{\textrm{r}}_{M^\#}\big)
			\colon \mathcal B \lra \mathcal A 
	\ee 
	is also a 1-morphism in $\orb{\tric}$. 
\end{lemma}

\begin{proof}
	The coherence axioms for $u^{\textrm{l}}_{M^\#}, u^{\textrm{r}}_{M^\#}$ follow from those for $u^{\textrm{l}}_{M}, u^{\textrm{r}}_{M}$ and the definitions, for example: 
	\be 
	\begin{tikzcd}[column sep=3em, row sep=4em]
	%%%%%%%%%%%%%%%%%%%%%%
	\tikzzbox{\begin{tikzpicture}[ultra thick,scale=0.75,baseline=-0.1cm]
		\coordinate (ml) at (-3,0);
		\coordinate (mr) at (0,0);
		\coordinate (ac) at (-1,0);
		\coordinate (B) at (-0.25,-0.4);
		\coordinate (ac2) at (-2,0);
		\coordinate (B2) at (0,-1);
		%
		% lines: 
		\draw[color=orange!80,opacity=0.545] (ml) -- (mr);
		\draw[color=green!70!black,opacity=0.545] (ac) -- (B);
		\draw[color=green!70!black,opacity=0.545] (ac2) -- (B2);
		%
		% vertices: 
		\fill[color=red!80!black] (ac) circle (3pt) node[color=red!80!black, above, font=\footnotesize] {$\triangleright_{M^\#}$};
		\fill[color=red!80!black] (ac2) circle (3pt) node[color=red!80!black, above, font=\footnotesize] {$\triangleright_{M^\#}$};
		\fill[color=green!80!black] (B) circle (3pt) node[color=green!80!black, below, font=\footnotesize] {};
		\end{tikzpicture}}%%popende 
	%%%%%%%%%%%%%%%%%%%%%%       
	\arrow[equal]{dr}
	\arrow{dddd}[swap]{\alpha_{M^\#}^{\textrm{l}}}
	%\arrow[bend left=55]{rrrdd}{u_{M^\#}^{\textrm{l}}}
	\arrow[in=90, out=20]{rrrdd}{u_{M^\#}^{\textrm{l}}}
	%\arrow{rrrdd}%[controls={+(1.5,0.5) and +(-1,0.8)}]
	&
	&
	&
	\\
	& 
	%%%%%%%%%%%%%%%%%%%%%%
	\tikzzbox{\begin{tikzpicture}[ultra thick,scale=0.75,baseline=0.6cm]
		\coordinate (ml) at (-1.5,-1);
		\coordinate (mr) at (2,1);
		\coordinate (ac) at (0,0);
		\coordinate (ac2) at (1,1.5);
		\coordinate (B) at (2.5,0.5);
		\coordinate (B2) at (2,1.75);
		\coordinate (ev) at (-1,0.5);
		\coordinate (coev) at (1,-0.5);
		\coordinate (ev2) at (0,2);
		\coordinate (mr2) at (2.5,2.5);
		%
		% lines: 
		\draw[color=orange!80,opacity=0.545] (ml) -- (coev) -- (ev) -- (mr) -- (ev2) -- (mr2);
		\draw[color=green!70!black,opacity=0.545] (ac) -- (B);
		\draw[color=green!70!black,opacity=0.545] (ac2) -- (B2);
		%
		% vertices: 
		\fill[color=red!80!black] (ac) circle (3pt) node[color=red!80!black, above, font=\footnotesize] {$\triangleleft_{M}$};
		\fill[color=red!80!black] (ac2) circle (3pt) node[color=red!80!black, above, font=\footnotesize] {$\triangleleft_{M}$};
		\fill[color=green!80!black] (B2) circle (3pt) node[color=green!80!black, below, font=\footnotesize] {};
		\end{tikzpicture}}%%popende 
	%%%%%%%%%%%%%%%%%%%%%%  
	\arrow{r}{u_M^{\textrm{r}}}
	\arrow{d}[sloped]{\textrm{tensorator}}
	\arrow{d}[swap, sloped]{\textrm{\& cusp}}
	& 
	%%%%%%%%%%%%%%%%%%%%%%
	\tikzzbox{\begin{tikzpicture}[ultra thick,scale=0.75,baseline=0.6cm]
		\coordinate (ml) at (-1.5,-1);
		\coordinate (mr) at (2,1);
		\coordinate (ac) at (0,0);
		\coordinate (ac2) at (1,1.5);
		\coordinate (B) at (2.5,0.5);
		\coordinate (ev) at (-1,0.5);
		\coordinate (coev) at (1,-0.5);
		\coordinate (ev2) at (0,2);
		\coordinate (mr2) at (2.5,2.5);
		%
		% lines: 
		\draw[color=orange!80,opacity=0.545] (ml) -- (coev) -- (ev) -- (mr) -- (ev2) -- (mr2);
		\draw[color=green!70!black,opacity=0.545] (ac) -- (B);
		%
		% vertices: 
		\fill[color=red!80!black] (ac) circle (3pt) node[color=red!80!black, above, font=\footnotesize] {$\triangleleft_{M}$};
		\end{tikzpicture}}%%popende 
	%%%%%%%%%%%%%%%%%%%%%% 
	\arrow{d}[sloped]{\textrm{tensorator}}
	\arrow{d}[swap, sloped]{\textrm{\& cusp}}
	& 
	\\
	& 
	%%%%%%%%%%%%%%%%%%%%%%
	\tikzzbox{\begin{tikzpicture}[ultra thick,scale=0.75,baseline=-0.1cm]
		\coordinate (ml) at (-2,-1);
		\coordinate (mr) at (1.5,1);
		\coordinate (ac) at (-0.5,0.25);
		\coordinate (ac2) at (0.5,-0.25);
		\coordinate (B) at (1,0.5);
		\coordinate (B2) at (1.5,0);
		\coordinate (ev) at (-1,0.5);
		\coordinate (coev) at (1,-0.5);
		%
		% lines: 
		\draw[color=orange!80,opacity=0.545] (ml) -- (coev) -- (ev) -- (mr);
		\draw[color=green!70!black,opacity=0.545] (ac) -- (B);
		\draw[color=green!70!black,opacity=0.545] (ac2) -- (B2);
		%
		% vertices: 
		\fill[color=red!80!black] (ac) circle (3pt) node[color=red!80!black, below, font=\footnotesize] {$\triangleleft_{M}$};
		\fill[color=red!80!black] (ac2) circle (3pt) node[color=red!80!black, above, font=\footnotesize] {$\triangleleft_{M}$};
		\fill[color=green!80!black] (B) circle (3pt) node[color=green!80!black, below, font=\footnotesize] {};
		\end{tikzpicture}}%%popende 
	%%%%%%%%%%%%%%%%%%%%%%   
	\arrow{r}[swap]{u_M^{\textrm{r}}}
	\arrow{d}[swap]{\alpha_M^{\textrm{r}}}
	& 
	%%%%%%%%%%%%%%%%%%%%%%
	\tikzzbox{\begin{tikzpicture}[ultra thick,scale=0.75,baseline=-0.1cm]
		\coordinate (ml) at (-2,-1);
		\coordinate (mr) at (1.5,1);
		\coordinate (ac2) at (0,0);
		\coordinate (B2) at (1.5,0.25);
		\coordinate (ev) at (-1,0.5);
		\coordinate (coev) at (1,-0.5);
		%
		% lines: 
		\draw[color=orange!80,opacity=0.545] (ml) -- (coev) -- (ev) -- (mr);
		\draw[color=green!70!black,opacity=0.545] (ac2) -- (B2);
		%
		% vertices: 
		\fill[color=red!80!black] (ac2) circle (3pt) node[color=red!80!black, above, font=\footnotesize] {$\triangleleft_{M}$};
		\end{tikzpicture}}%%popende 
	%%%%%%%%%%%%%%%%%%%%%%
	\arrow[equal]{r}
	& 
	%%%%%%%%%%%%%%%%%%%%%%
	\tikzzbox{\begin{tikzpicture}[ultra thick,scale=0.75,baseline=-0.1cm]
		\coordinate (ml) at (-3,0);
		\coordinate (mr) at (0,0);
		\coordinate (ac) at (-1.5,0);
		\coordinate (B) at (0,-1);
		%
		% lines: 
		\draw[color=orange!80,opacity=0.545] (ml) -- (mr);
		\draw[color=green!70!black,opacity=0.545] (ac) -- (B);
		%
		% vertices: 
		\fill[color=green!30!black] (B) circle (0pt) node[color=green!30!black, right, font=\footnotesize] {};
		green!30!black
		\fill[color=red!80!black] (ac) circle (3pt) node[color=red!80!black, above, font=\footnotesize] {$\triangleright_{M^\#}$};
		\end{tikzpicture}}%%popende 
	%%%%%%%%%%%%%%%%%%%%%%    
	\\
	& 
	%%%%%%%%%%%%%%%%%%%%%%
	\tikzzbox{\begin{tikzpicture}[ultra thick,scale=0.75,baseline=-0.1cm]
		\coordinate (ml) at (-1.5,-1);
		\coordinate (mr) at (2,1);
		\coordinate (ac) at (0,0);
		\coordinate (mu) at (1,0.125);
		\coordinate (B) at (1.75,0.5);
		\coordinate (B2) at (2,-0.25);
		\coordinate (ev) at (-1,0.5);
		\coordinate (coev) at (1,-0.5);
		%
		% lines: 
		\draw[color=orange!80,opacity=0.545] (ml) -- (coev) -- (ev) -- (mr);
		\draw[color=green!70!black,opacity=0.545] (ac) -- (mu);
		\draw[color=green!70!black,opacity=0.545] (B) -- (mu);
		\draw[color=green!70!black,opacity=0.545] (B2) -- (mu);
		%
		% vertices: 
		\fill[color=red!80!black] (ac) circle (3pt) node[color=red!80!black, above, font=\footnotesize] {$\triangleleft_{M}$};
		\fill[color=green!30!black] (mu) circle (3pt) node[color=green!30!black, above, font=\footnotesize] {};
		\fill[color=green!80!black] (B) circle (3pt) node[color=green!80!black, below, font=\footnotesize] {};
		\end{tikzpicture}}%%popende 
	%%%%%%%%%%%%%%%%%%%%%%  
	\arrow{ru}[swap]{u_B^{\textrm{r}}}
	&
	&
	\\ 
	%%%%%%%%%%%%%%%%%%%%%%
	\tikzzbox{\begin{tikzpicture}[ultra thick,scale=0.75,baseline=-0.1cm]
		\coordinate (ml) at (-3,0);
		\coordinate (mr) at (0,0);
		\coordinate (mu) at (-1,-0.5);
		\coordinate (B) at (-0.25,-0.4);
		\coordinate (ac2) at (-2,0);
		\coordinate (B2) at (0,-1);
		%
		% lines: 
		\draw[color=orange!80,opacity=0.545] (ml) -- (mr);
		\draw[color=green!70!black,opacity=0.545] (mu) -- (B);
		\draw[color=green!70!black,opacity=0.545] (ac2) -- (B2);
		%
		% vertices: 
		\fill[color=green!30!black] (mu) circle (3pt) node[color=green!30!black, below, font=\footnotesize] {};
		\fill[color=red!80!black] (ac2) circle (3pt) node[color=red!80!black, above, font=\footnotesize] {$\triangleright_{M^\#}$};
		\fill[color=green!80!black] (B) circle (3pt) node[color=green!80!black, below, font=\footnotesize] {};
		\end{tikzpicture}}%%popende 
	%%%%%%%%%%%%%%%%%%%%%%  
	%\arrow[bend right=55, swap]{rrruu}{u_{B}^{\textrm{r}}}
	\arrow[in=-90, out=-20, swap]{rrruu}{u_{B}^{\textrm{r}}}
	\arrow[equal]{ur}
	&
	&
	&
	\end{tikzcd}
	\ee 
	Here the outer subdiagrams commute by definition, the upper inner subdiagram commutes by naturality, and the commutativity of the lower inner subdiagram is a consequence of the coherence axioms for $u^{\textrm{l}}_{M}, u^{\textrm{r}}_{M}$. 
	
	The 2-3 move for~$\alpha^{\textrm{l}}_{M}$ implies the 2-3 move for~$\alpha^{\textrm{r}}_{M^\#}$, 
	\be 
	\begin{tikzcd}[column sep=2em, row sep=4em]
	& 
	& 
	& 
	%%%%%%%%%%%%%%%%%%%%%%
	\tikzzbox{
}%%popende 
	%%%%%%%%%%%%%%%%%%%%%%  
	\arrow[in=-35, out=-145]{llllll}
	\end{tikzcd}
	\ee 
	where all subdiagrams commute either by definition or due to naturality and coherence for cusp and tensorator isomorphisms. 
	Similarly, the 2-3 move for~$\alpha^{\textrm{l}}_{M^\#}$ follows from the one for~$\alpha^{\textrm{r}}_{M}$. 
	
	Finally, the two 2-3 moves~\eqref{eq:23MovesFor1Morphisms2} and~\eqref{eq:23MovesFor1Morphisms3} involving~$\alpha^{\textrm{m}}_{M^\#}$ hold because of those for~$\alpha^{\textrm{m}}_{M}$ as well as naturality of $\mathcal P_{\mathcal A}, \mathcal P_{\mathcal B}$ and the identities~\eqref{eq:PPP}. 
	We provide details for the 2-3 move involving~$\alpha^{\textrm{m}}_{M^\#}$ and~$\alpha^{\textrm{l}}_{M^\#}$ in Figure~\ref{fig:23moveMiddleRight}, where all subdiagrams commute by naturality and/or the relations indicated. 
	The 2-3 move involving~$\alpha^{\textrm{m}}_{M^\#}$ and~$\alpha^{\textrm{r}}_{M^\#}$ is checked analogously. 
\end{proof}

\begin{figure}[p]
	\begin{center}
		$$
		\!\!\!
		\begin{tikzcd}[column sep=2em, row sep=3.5em]
		%%%%%%%%%%%%%%%%%%%%%%
		\tikzzbox{
}%%popende 
		%%%%%%%%%%%%%%%%%%%%%% 
		\arrow{l}
		\end{tikzcd}
		$$
		\caption{The 2-3 move involving~$\alpha^{\textrm{m}}_{M^\#}$ and~$\alpha^{\textrm{l}}_{M^\#}$ states that the outer paths between the five shaded 2-morphisms are equal; non-labelled arrows are given by cusp and tensorator structure 3-morphisms.}
		\label{fig:23moveMiddleRight}
	\end{center}
\end{figure}

We want to identify~$\mathcal M^\#$ as a left adjoint of~$\mathcal M$ in $\orb{\tric}$. 
To do so, we first define the 3-isomorphism
\be 
\zeta_{\mathcal M}
:= 
\left(\hspace{-0.1cm}
\begin{tikzcd}[column sep=3em, row sep=0.5em]
%%%%%%%%%%%%%%%%%%%%%%
\tikzzbox{
}%%popende 
%%%%%%%%%%%%%%%%%%%%%%  
\end{tikzcd}
\right) 
\ee 
where we use 
\be 
\xi := 
%%%%%%%%%%%%%%%%%%%%%%
\tikzzbox{%
}%%popende 
%%%%%%%%%%%%%%%%%%%%%%   
\, . 
\ee 
With this, and using the graphical notation for relative products introduced in~\eqref{eq:PiStarInv} and~\eqref{eq:SigmaStarInv}, we can set 
\begin{align}
\label{eq:evMcal}
\ev_{\mathcal M} 
:= 
\left( 
%%%%%%%%%%%%%%%%%%%%%%
\tikzzbox{%
}%%popende 
%%%%%%%%%%%%%%%%%%%%%%  
\colon 
\begin{tikzcd}[column sep=3em, row sep=0.5em]
%%%%%%%%%%%%%%%%%%%%%%
\tikzzbox{%
}%%popende 
%%%%%%%%%%%%%%%%%%%%%%  
\colon 
\begin{tikzcd}[column sep=3em, row sep=0.5em]
%%%%%%%%%%%%%%%%%%%%%%
\tikzzbox{%
}%%popende 
%%%%%%%%%%%%%%%%%%%%%% 
\end{tikzcd}
\ee 
\be 
\triangleright_{\coev_{\mathcal M}} := 
\left( 
\begin{tikzcd}[column sep=3em, row sep=0.5em]
%%%%%%%%%%%%%%%%%%%%%%
\tikzzbox{%
}%%popende 
%%%%%%%%%%%%%%%%%%%%%%  
\end{tikzcd}
\right) 
\ee 
\be 
\triangleleft_{\coev_{\mathcal M}} := 
\left( 
\begin{tikzcd}[column sep=3em, row sep=0.5em]
%%%%%%%%%%%%%%%%%%%%%%
\tikzzbox{%
}%%popende 
%%%%%%%%%%%%%%%%%%%%%%  
\end{tikzcd}
\right) 
\ee 

\begin{proposition}
	\label{prop:Adj1Mor}
	Every 1-morphism~$\mathcal M$ in $\orb{\tric}$ has a left and right adjoint~$\mathcal M^\#$. 
\end{proposition}

\begin{proof}
	We have to verify that~\eqref{eq:evMcal} and its variants for right adjunction are indeed 2-morphisms in $\orb{\tric}$, and that they satisfy the Zorro moves up to 3-isomorphism. 
	To address the former, we have to check that $\triangleright_{\ev_{\mathcal M}}, \triangleleft_{\ev_{\mathcal M}}, \triangleright_{\coev_{\mathcal M}}, \triangleleft_{\coev_{\mathcal M}}$ satisfy the constraints of Definition~\ref{def:2MorphismsEAlg}. 
	This is straightforward, and we give details only for the most involved case of $\triangleleft_{\ev_{\mathcal M}}$. 
	Here the constraint to be checked is
	\be 
	%%%%%%%%%%%%%%%%%%%%%%%%%%%%
	\tikzzbox{
}%%popende 
	%%%%%%%%%%%%%%%%%%%%%%  
	\, . 
	\ee 
	Inserting the definitions and suppressing the contribution of $\ev_M$ to $\ev_{\mathcal M}$, this becomes
	\be
	\label{eq:interm}
	%%%%%%%%%%%%%%%%%%%%%%%%%%%%
	\tikzzbox{%
}%%popende 
	%%%%%%%%%%%%%%%%%%%%%%  
	\, . 
	\ee 
	On the left-hand side the middle~$\zeta_{\mathcal M}$ and~$\zeta_{\mathcal M}^{-1}$ cancel, after which a Zorro move for $\triangleleft_{M^\#}^\vee$ can be applied. 
	Pre-composition with $\coev_{\triangleright_M}$ and post-composition with $\ev_{\triangleright_M}$ then makes the identity~\eqref{eq:interm} equivalent to
	\be 
	%%%%%%%%%%%%%%%%%%%%%%%%%%%%
	\tikzzbox{%
}%%popende 
	%%%%%%%%%%%%%%%%%%%%%%  
	\ee 
	which is precisely the 2-3 move for~$\alpha_{M^\#}^{\textrm{r}}$. 
	
	It remains to verify the Zorro moves for $\ev_{\mathcal M}, \coev_{\mathcal M}$. 
	For this, we first compute 
	\begin{align}
	&
	\overbrace{%
		%%%%%%%%%%%%%%%%%%%%%%
		\tikzzbox{%
}%%popende 
	%%%%%%%%%%%%%%%%%%%%%% 
	\, , 
	\end{align}
	where we used~\eqref{eq:BalancedEquiv}, \eqref{eq:CobalEquiv} as well as $\Pi^* = (-)\circ \Pi$ and $\Sigma_* = \Sigma\circ (-)$. 
	The other Zorro move follows more directly with the help of~\eqref{lem:Filled}: 
	\be 
	\overbrace{%
		%%%%%%%%%%%%%%%%%%%%%%
		\tikzzbox{%
}%%popende 
	%%%%%%%%%%%%%%%%%%%%%%  
	\, . 
	\ee 
\end{proof}

\begin{remark}
	In the case that $\tric$ is symmetric monoidal we expect $\orb{\tric}$ to be symmetric monoidal as well. 
	If $\tric$ admits duals for objects we expect in addition that the objects in $\orb{\tric}$ also admit duals. 
	In the case that $\tric$ has only one object this follows from the fact that every algebra is 1-dualisable in the usual Morita category, see e.\,g.\ \cite{ClaudiaOwen}. 
	Together with Proposition~\ref{prop:Adj1Mor} and Lemma~\ref{lem:Adj2Mor} this would imply that every object in $\orb{\tric}$ is fully dualisable. 
	Based on the close connection to oriented theories we expect all objects in $\orb{\tric}$ to have a canonical $\operatorname{SO}(3)$-homotopy fixed-point structure, in the sense that they trivialise the $\operatorname{SO}(3)$-action on $\orb{\tric}^\times$. 
	Furthermore, the 1-morphisms should come with a relative $\operatorname{SO}(2)$-homotopy fixed-point structure as discussed in~\cite[Sect.\,4.3]{Lurie}.      
\end{remark}

\subsubsection{Universal property}
\label{subsubsec:UniversalProperty}

In this section we discuss a universal property for orbifold completion. 
Our discussion is rigorous in dimensions~1 and~2, but it is short on technical details in the 3- and higher-dimensional case. This can be understood as an extension of higher idempotent completions
%arXiv_v2:
	or condensation completion 
as developed in~\cite{GaiottoJohnsonFreyd} to pivotal higher categories. 

\medskip 

Recall that the \textsl{idempotent completion} of a 1-category~$\mathcal C$ is a fully faithful functor $\mathcal C \longhookrightarrow \overline{\mathcal C}$ to a category~$\overline{\mathcal C}$ in which every idempotent splits, such that the following universal property holds: 
for every functor $\mathcal C \lra \overline{\mathcal D}$ with idempotent complete codomain, there exists an essentially unique functor $\overline{\mathcal C} \lra 
\overline{\mathcal D}$ such that 
\be 
\begin{tikzcd}[column sep=3em, row sep=3em]
\mathcal C
\ar[r] 
\ar[d, hook] 
& 
\overline{\mathcal D}
\\ 
\overline{\mathcal C}
\ar[ur, dashed] 
& 
\end{tikzcd}
\ee 
commutes up to isomorphism. 
It is straightforward to check that the category whose objects are idempotents in~$\mathcal C$ and whose morphisms $e\lra e'$ are morphisms in~$\mathcal C$ that are invariant under pre- and post-composition with~$e$ and~$e'$, respectively, satisfies this universal property. 

\medskip 

The 2-dimensional case is completely analogous. 
The \textsl{idempotent completion} of a 2-category~$\mathcal B$ is a fully faithful 2-functor $\mathcal B \longhookrightarrow \overline{\mathcal B}$ to a 2-category~$\overline{\mathcal B}$ in which every 
%arXiv_v3: 
	%\textsl{2-idempotent} 
	 \textsl{condensation 2-monad} 
splits (in the sense of Definition~\ref{def:2IdempotentSplit} with 2-morphism replaced by 1-morphisms, and 3-morphisms replaced by 2-morphisms) and every Hom category is idempotent complete, such that the following universal property holds: 
for every 2-functor $\mathcal B \lra \overline{\mathcal D}$ with idempotent complete codomain, there exists an essentially unique 2-functor $\overline{\mathcal B} \lra 
\overline{\mathcal D}$ such that 
\be 
\begin{tikzcd}[column sep=3em, row sep=3em]
\mathcal B
\ar[r] 
\ar[d, hook] 
& 
\overline{\mathcal D}
\\ 
\overline{\mathcal B}
\ar[ur, dashed] 
& 
\end{tikzcd}
\ee 
commutes up to equivalence. 
%arXiv_v2: 
	We refer to \cite[App.\,A]{2EW} for a more detailed discussion. 

As explained in \cite{GaiottoJohnsonFreyd} (see \cite{FrauenbergerMasterThesis} for a detailed discussion), $\overline{\mathcal B}$ can be taken to be the 2-category of not necessarily (co)unital $\Delta$-separable 
%arXiv_v2: 
	Frobenius 
algebras, their bimodules and bimodule maps in~$\mathcal B$ (after idempotent completing the Hom categories of~$\mathcal B$). 
Apart from the structure of units and counits, this is precisely the ``equivariant completion''~$\mathcal B_{\textrm{eq}}$ introduced in~\cite{cr1210.6363}. 
As argued in \cite{GaiottoJohnsonFreyd}, this further generalises to a notion of idempotent completion of $n$-categories for arbitrary dimension~$n$. 
This notion of ``Karoubi envelope'' or ``condensation monads'', i.\,e.\ idempotents for $n=1$ and $\Delta$-separable Frobenius algebras for $n=2$, is motivated by the cobordism hypothesis for fully extended \textsl{framed} TQFTs. 

Our notion of orbifold completion arises in the context of \textsl{oriented} TQFTs. 
For $n=1$ this coincides with idempotent completion, in line with the fact that framings and orientations are equivalent for 1-dimensional manifolds. 
For $n=2$, we recalled the orbifold completion $\orb{\mathcal B}$ of a pivotal 2-category~$\mathcal B$ in Definition~\ref{def:2dOrbifoldCompletion}. 
To define a natural variant of split idempotents in such a context recall that for every 1-morphism $X\colon \alpha\to \beta$ in a pivotal 2-category $\mathcal{B}$ the morphism $X^\vee\circ X $ has a canonical Frobenius algebra structure~\cite[Prop.\,4.3]{cr1210.6363} where for example the multiplication is induced by the evaluation: 
\begin{align}
X^\vee \circ X \circ X^\vee \circ X \xrightarrow{\id_{X^\vee}\circ {\tev}_X\circ \id_X} X^\vee \circ X \, .
\end{align}  
Based on this observation we make the following:
\begin{definition}
	\label{def:OrbifoldCondensation}
	Let~$\mathcal B$ be a pivotal 2-category, and let $\alpha,\beta \in \mathcal B$. 
	\begin{enumerate}
		\item 
		An \textsl{orbifold condensation of~$\alpha$ onto~$\beta$} is a 1-morphism $X\colon \alpha \lra \beta$ such that $1_\beta \xrightarrow{\coev_X} X \circ X^\vee \xrightarrow{\widetilde{\ev}_X} 1_\beta$ is the identity 2-morphism on $\id_\beta$. 
		\item 
		An orbifold datum $A\colon \alpha \lra \alpha$ in~$\mathcal B$ \textsl{splits} if there exists an orbifold condensation~$X$ of~$\alpha$ such that $X^\vee \circ X \cong A$ as Frobenius algebras, where $X^\vee \circ X$ is equipped with the canonical Frobenius structure recalled above. 
	\end{enumerate}
\end{definition}

It is immediate from the definition that if~$X$ is an orbifold condensation, then $X^\vee\circ X$ carries the structure of an orbifold datum. 
The pair of maps $(\widetilde{\ev}_X,\coev_X)$ is a 1-condensation of $X^\vee\circ X$ onto~$1_\beta$ in as defined in~\cite{GaiottoJohnsonFreyd}.  
Orbifold completion in dimension two satisfies the property that every orbifold datum splits: 

\begin{proposition}
	\label{prop:EverythingSplitsInBorb}
	Let~$\mathcal B$ be a pivotal 2-category. 
	Then every orbifold datum in~$\orb{\mathcal B}$ splits. 
\end{proposition}

\begin{proof}
	Let $A\colon \alpha \lra \alpha$ be an object in $\orb{\mathcal B}$, and let $\mathcal A\colon A \lra A$ be an orbifold datum in $\orb{\mathcal B}$. 
	It is straightforward to check that the Frobenius algebra structure of~$\mathcal A$ in $\orb{\mathcal B}$ induces the structure of an orbifold datum on the underlying 1-morphism of~$\mathcal A$ in~$\mathcal B$, see the proof of \cite[Prop.\,4.2]{cr1210.6363} for details. 
	We denote this object of $\orb{\mathcal B}$ by $(\alpha,\mathcal A)$, and write $(\alpha,A) \in \orb{\mathcal B}$ for our chosen orbifold datum $A\colon \alpha\lra \alpha$. 
	
	The underlying 1-morphism~$\mathcal A$ in~$\mathcal B$ has a canonical $\mathcal A$-$A$-bimodule structure, which makes it a 1-morphism $X := {}_{\mathcal A}\mathcal A_A\colon (\alpha,A) \lra (\alpha,\mathcal A)$ in $\orb{\mathcal B}$. 
	We claim that~$X$ is an orbifold condensation of $(\alpha,A)$ onto $(\alpha,\mathcal A)$ in $\orb{\mathcal B}$. 
	To see this, we note that by composing the two middle arrows in 
	\be 
	\begin{tikzcd}[column sep=5em, row sep=3em]
	(\alpha,A)
	\ar[r, out=30, in=150, "X={}_{\mathcal A}\mathcal A_A"] 
	\ar[loop, out=225, in=135, looseness=8, "\mathcal A"] 
	& 
	(\alpha,\mathcal A)
	\ar[l, out=210, in=-30, "X^\vee={}_A\mathcal A_{\mathcal A}"] 
	\ar[loop, out=-45, in=45, looseness=8, "1_{(\alpha,\mathcal A)}", swap] 
	\end{tikzcd}
	\ee 
	to form an endomorphism of $(\alpha,A)$, we directly obtain the orbifold datum $\mathcal A\colon (\alpha,A) \lra (\alpha,A)$ in $\orb{\mathcal B}$, while the opposite composition $X\circ X^\vee$ condenses onto~$1_{(\alpha,\mathcal A)}$.
	Indeed, this splitting is witnessed by the adjunction 2-morphisms $\tev_X, \coev_X$ which are induced from the (co)multiplication of~$\mathcal A$. 
\end{proof}

Better still, $\orb{\mathcal B}$ is universal among all pivotal 2-categories in which all orbifold data split: 

\begin{proposition}
	\label{prop:BorbSatisfiesUniversalProperty}
	Let~$\mathcal B$ be a pivotal 2-category whose Hom categories are idempotent complete. Then the inclusion $\mathcal B \longhookrightarrow \orb{\mathcal B}$, $\alpha\lmt (\alpha,1_\alpha)$ satisfies the following universal property: 
	for every pivotal 2-functor $\mathcal B \lra \overline{\mathcal D}$ with pivotal codomain in which every orbifold condensation splits, there exists an essentially unique pivotal 2-functor $\orb{\mathcal B} \lra 
	\overline{\mathcal D}$ such that 
	\be 
	\begin{tikzcd}[column sep=3em, row sep=3em]
	\mathcal B
	\ar[r] 
	\ar[d, hook] 
	& 
	\overline{\mathcal D}
	\\ 
	\orb{\mathcal B}
	\ar[ur, dashed] 
	& 
	\end{tikzcd}
	\ee 
	commutes up to equivalence. 
\end{proposition}

\begin{proof} 
	Let $F\colon \mathcal B \lra \overline{\mathcal D}$ be a pivotal 2-functor. 
	Then its lift $\orb{\mathcal B} \lra \overline{\mathcal D}$ sends $(\alpha,A) \in \orb{\mathcal B}$ to the image of the splitting of the orbifold datum $F(A) \colon F(\alpha) \lra F(\alpha)$, which exists by assumption, and it is essentially unique by construction. 
\end{proof}

As a direct consequence, we recover the fact that $\orb{(\orb{\mathcal B})} \cong \orb{\mathcal B}$. 
Moreover, the universal property of $\orb{\mathcal B}$ should imply that various potential additional structures on~$\mathcal B$ are inherited by $\orb{\mathcal B}$, in particular braided monoidal structures. 

\medskip 

The above suggests the following universal property for an algebraic notion of orbifold completion in arbitrary dimension~$n$. 
Let~$\tric$ be an $n$-category with compatible adjoints for all morphisms and orbifold complete Hom $(n-1)$-categories. 
For a 1-morphism $X\colon \alpha\to \beta$ consider the endomorphism $A\coloneqq X^\vee\circ X \colon \alpha\to \alpha$. 
The adjunction data induce on $A$ the potential structure of an orbifold datum. 
In case this satisfies the orbifold conditions we call $X$ an \emph{orbifold condensation} of $\alpha$ onto $\beta$. An orbifold datum in $\tric$ \emph{splits} if there exists an orbifold condensation giving rise to it.    
Then the orbifold completion $\orb{\tric}$ satisfies the universal property that for every pivotal $n$-functor $F\colon \tric \lra \overline{\mathcal U}$ whose codomain has compatible adjoints and in which every orbifold datum splits, there exists an essentially unique lift of~$F$ to an $n$-functor $\orb{\mathcal T} \lra \overline{\mathcal U}$. 
As in the 2-categorical case the splitting of an orbifold datum $\mathcal{A}$ in the orbifold completion should be given by $\mathcal{A}$ seen as a 1-morphism in $\orb \tric$. 
%arXiv_v2: 
	Making this rigorous for $n=3$ amounts to working in the 4-category of 3-categories with compatible adjoints for 1- and 2-morphisms, as well as their structure-preserving 3-functors, 3-natural transformations, modifications, and perturbations. 
	To our knowledge, the formidable task of rigorously constructing this 4-category has not yet been fully completed.

\section{Examples of completions}
\label{sec:ExamppleCompletion}

In this section we study some natural examples of 3-dimensional orbifold completion and show how they relate to the existing literature. 
We start, in Section~\ref{subsec:StateSumModels}, by studying the delooping of 2-dimensional state sum models described by the pivotal 2-category of separable symmetric Frobenius algebras $\ssFrob$. 
Then in Section~\ref{subsec:DomainWallsForRT} we turn to a example involving modular fusion categories. 

The 3-categories~$\tric$ appearing in this section are not strict enough to be Gray categories with duals.
We expect that our construction of the 3-category $\EC$ can be adapted to work for less strict 3-categories. 
However, an appropriate weakening of the notion of a generic 3-category with duals has not appeared in the literature, to the best of our knowledge. 
Giving a general description of the subcategory $\orb\tric$ is not straightforward. 
In the examples we will consider there are however natural guesses on how to adapt the conditions. 
Another approach to the problem is to construct strictifications of the 3-categories which appear. 
One way of doing this is to first construct a 3-dimensional defect TQFT and then extract the corresponding Gray category with duals constructed in~\cite{CMS}.

\subsection{State sum models} 
\label{subsec:StateSumModels}

Here we discuss how state sum models fit within the theory of orbifold completion. 
In Section~\ref{subsec:2Dstatesums} we briefly recall the 2-dimensional case and formalise the notion of the ``Euler completion'' of a pivotal 2-category. 
In Section~\ref{subsec:EWor} we build on \cite{FEW} to construct a pivotal equivalence between the 2-categories of separable symmetric Frobenius $\K$-algebras and of semisimple $\K$-linear Calabi--Yau categories, 
%arXiv_v2: 
	where throughout~$\K$ some fixed algebraically closed field of characteristic~0.  
One of our main applications is then proved in Section~\ref{subsubsec:3dStateSumModels}, where we show that the 3-category of spherical fusion categories and bimodule categories with trace, constructed by Schaumann in \cite{Bimodtrace}, is a subcategory of the orbifold completion of the Euler completion of the delooping of the 2-categories considered in Section~\ref{subsec:EWor}.

\subsubsection{Recollection on 2-dimensional state sum models}
\label{subsec:2Dstatesums}

We start by briefly recalling the construction of 2-dimensional state sum models from separable symmetric Frobenius algebras, from the perspective of orbifold completion. 
The pivotal 2-category of defects in the trivial 2-dimensional field theory is the delooping $\Bar\vs$, where $\vs$ is the category of finite-dimensional vector spaces. 
Its orbifold completion is the pivotal 2-category $\Delta\ssFrob$ of $\Delta$-separable symmetric Frobenius algebras, bimodules and bimodule maps, see Section~\ref{subsec:2categories}. 
The category $\Delta\ssFrob$ can be used to construct 2-dimensional defect TQFTs of state sum type -- however not all of them. 
For example, among the invertible state sum models with underlying vector space $A=\C$, only those whose counit~$\varepsilon$ is multiplication with $\pm 1$ are $\Delta$-separable, while for every invertible $\lambda\in \C^\times$ there is a TQFT $\zz_{\lambda}^{\textrm{eu}}$ with counit $\varepsilon(c)=\lambda c$ that assigns the invariant $\lambda^{2g-2}$ to a closed surface of genus $g$. 
Note that here $2g-2$ is the Euler characteristic of the surface. 
The TQFT $\zz_{\lambda}^{\textrm{eu}}$ can be obtained from $\zz_{\pm 1}^{\textrm{eu}}$ by the procedure of ``Euler completion'' described in~\cite[Sect.\,2.5]{CRS1}. 
This example suggests that we should consider the Euler completion of the defect theory corresponding to $\Delta\ssFrob$ in general. 

We now describe the \textsl{Euler completion} $E(\mathcal B)$ of a pivotal 2-category $\mathcal{B}$. 
This follows closely \cite[Sect.\,4.2]{CRS1} and \cite[App.\,B]{CMRSS1}. An object of $E(\mathcal{B})$ is a pair of an object $b\in \mathcal{B}$ and an element $\psi_b\in \Aut (1_{b})$. 
1- and 2-morphisms are the same as those of $\mathcal{B}$ and their composition is unchanged. The pivotal structure is changed by setting the adjunction 2-morphisms of 1-morphisms $X\colon (b,\psi_b) \lra (b',\psi'_{b'})$ to be 
\be 
\label{eq:EulerAdjunctionMaps}
\tikzzbox{%
	%%%%%%%%%%%%%%%%%%%%%% 
	\begin{tikzpicture}[very thick,scale=1.0,color=blue!50!black, baseline=0.5cm]
	\coordinate (X) at (0.5,0);
	\coordinate (Xd) at (-0.5,0);
	\coordinate (d1) at (-1,0);
	\coordinate (d2) at (+1,0);
	\coordinate (u1) at (-1,1.25);
	\coordinate (u2) at (+1,1.25);
	%
	% colouring: 
	\fill [orange!40!white, opacity=0.7] (d1) -- (d2) -- (u2) -- (u1); 
	\draw[thin] (d1) -- (d2) -- (u2) -- (u1) -- (d1); 
	%
	% strings: 
	\draw[directed] (X) .. controls +(0,1) and +(0,1) .. (Xd);
	%
	% labels: 
	\fill (X) circle (0pt) node[below] {{\small $X\vphantom{X^\vee}$}};
	\fill (Xd) circle (0pt) node[below] {{\small ${}^\vee\!X$}};
	\fill[red!80!black] (0,0.15) circle (0pt) node {{\scriptsize $b'$}};
	\fill[red!80!black] (0.8,1) circle (0pt) node {{\scriptsize $b$}};
	\fill[black] (0.15,0.4) circle (1.5pt) node[left] {{\scriptsize $\psi'_{b'}\hspace{-0.3em}$}};
	\fill[black] (-0.8,1) circle (1.5pt) node[right] {{\scriptsize $\hspace{-0.1em}\psi_{b}^{-1}$}};
	\end{tikzpicture}
	%%%%%%%%%%%%%%%%%%%%%% 
}
\, , \quad 
\tikzzbox{%
	%%%%%%%%%%%%%%%%%%%%%% 
	\begin{tikzpicture}[very thick,scale=1.0,color=blue!50!black, baseline=0.5cm]
	\coordinate (X) at (-0.5,1.25);
	\coordinate (Xd) at (0.5,1.25);
	\coordinate (d1) at (-1,0);
	\coordinate (d2) at (+1,0);
	\coordinate (u1) at (-1,1.25);
	\coordinate (u2) at (+1,1.25);
	%
	% colouring: 
	\fill [orange!40!white, opacity=0.7] (d1) -- (d2) -- (u2) -- (u1); 
	\draw[thin] (d1) -- (d2) -- (u2) -- (u1) -- (d1);
	%
	% strings: 
	\draw[redirected] (X) .. controls +(0,-1) and +(0,-1) .. (Xd);
	%
	% labels: 
	\fill (X) circle (0pt) node[above] {{\small $X\vphantom{X^\vee}$}};
	\fill (Xd) circle (0pt) node[above] {{\small ${}^\vee\!X$}};
	\fill[red!80!black] (0,1.1) circle (0pt) node {{\scriptsize $b$}};
	\fill[red!80!black] (0.8,0.2) circle (0pt) node {{\scriptsize $b'$}};
	\fill[black] (0.15,0.8) circle (1.5pt) node[left] {{\scriptsize $\psi_{b}\hspace{-0.3em}$}};
	\fill[black] (-0.8,0.2) circle (1.5pt) node[right] {{\scriptsize $\hspace{-0.1em}{\psi'_{b'}}^{-1}$}};
	\end{tikzpicture}
	%%%%%%%%%%%%%%%%%%%%%% 
}
\, ,
\quad
\tikzzbox{%
	%%%%%%%%%%%%%%%%%%%%%% 
	\begin{tikzpicture}[very thick,scale=1.0,color=blue!50!black, baseline=0.5cm]
	\coordinate (X) at (0.5,0);
	\coordinate (Xd) at (-0.5,0);
	\coordinate (d1) at (-1,0);
	\coordinate (d2) at (+1,0);
	\coordinate (u1) at (-1,1.25);
	\coordinate (u2) at (+1,1.25);
	%
	% colouring: 
	\fill [orange!40!white, opacity=0.7] (d1) -- (d2) -- (u2) -- (u1); 
	\draw[thin] (d1) -- (d2) -- (u2) -- (u1) -- (d1);
	%
	% strings: 
	\draw[redirected] (X) .. controls +(0,1) and +(0,1) .. (Xd);
	%
	% labels: 
	\fill (Xd) circle (0pt) node[below] {{\small $X\vphantom{X^\vee}$}};
	\fill (X) circle (0pt) node[below] {{\small $X^\vee$}};
	\fill[red!80!black] (0,0.15) circle (0pt) node {{\scriptsize $b$}};
	\fill[red!80!black] (0.8,1) circle (0pt) node {{\scriptsize $b'$}};
	\fill[black] (0.15,0.4) circle (1.5pt) node[left] {{\scriptsize $\psi_{b}\hspace{-0.3em}$}};
	\fill[black] (-0.8,1) circle (1.5pt) node[right] {{\scriptsize $\hspace{-0.1em}{\psi'_{b'}}^{-1}$}};
	\end{tikzpicture}
	%%%%%%%%%%%%%%%%%%%%%% 
}
\, , \quad 
\tikzzbox{%
	%%%%%%%%%%%%%%%%%%%%%% 
	\begin{tikzpicture}[very thick,scale=1.0,color=blue!50!black, baseline=0.5cm]
	\coordinate (X) at (-0.5,1.25);
	\coordinate (Xd) at (0.5,1.25);
	\coordinate (d1) at (-1,0);
	\coordinate (d2) at (+1,0);
	\coordinate (u1) at (-1,1.25);
	\coordinate (u2) at (+1,1.25);
	%
	% colouring: 
	\fill [orange!40!white, opacity=0.7] (d1) -- (d2) -- (u2) -- (u1); 
	\draw[thin] (d1) -- (d2) -- (u2) -- (u1) -- (d1);
	%
	% strings: 
	\draw[directed] (X) .. controls +(0,-1) and +(0,-1) .. (Xd);
	%
	% labels: 
	\fill (Xd) circle (0pt) node[above] {{\small $X\vphantom{X^\vee}$}};
	\fill (X) circle (0pt) node[above] {{\small $X^\vee$}};
	\fill[red!80!black] (0,1.1) circle (0pt) node {{\scriptsize $b'$}};
	\fill[red!80!black] (0.8,0.2) circle (0pt) node {{\scriptsize $b$}};
	\fill[black] (0.15,0.8) circle (1.5pt) node[left] {{\scriptsize $\psi'_{b'}\hspace{-0.3em}$}};
	\fill[black] (-0.8,0.2) circle (1.5pt) node[right] {{\scriptsize $\hspace{-0.1em}{\psi_{b}}^{-1}$}};
	\end{tikzpicture}
	%%%%%%%%%%%%%%%%%%%%%% 
} 
\, .
\ee 

\begin{remark}
	Let~$\zz$ be a 2-dimensional defect TQFT, let~$\zz^\odot$ be its Euler completion as in \cite[Def.\,2.24]{CRS1}, and let $\mathcal B_\zz$ and $\mathcal B_{\zz^\odot}$ be the associated pivotal 2-categories as constructed in \cite{dkr1107.0495}. 
	Then it follows directly from the definitions that $E(\mathcal B_{\zz})$ is not exactly equal to $\mathcal B_{\zz^\odot}$, but they are pivotally equivalent by a rescaling of 2-morphisms as in~\cite[App.\,B]{CMRSS1}. 
\end{remark}  

Note that in the case $\mathcal{B}=\Delta\ssFrob$ the elements of $\Aut (1_{A})$ are the invertible elements in the centre of~$A$. 
There is a canonical equivalence of 2-categories $\Gamma\colon \ssFrob\to E(\Delta\ssFrob)$, where $\ssFrob$ is the 2-category of all separable symmetric Frobenius algebras. 
We use this equivalence to define a pivotal structure on $\ssFrob$. 
To write down the equivalence explicitly recall~\cite{Abrams} that for every separable symmetric Frobenius algebra $(A,1,\varepsilon,\mu, \Delta)$ in $\vs$, the characteristic element $\omega \coloneqq \mu \circ \Delta (1)$ is an invertible element of the centre $Z(A)$. 
The equivalence~$\Gamma$ sends~$A$ to $(A,1,\mu, \Delta'= \Delta \circ \mu(-\otimes\omega^{-1}), \varepsilon'= \varepsilon \circ \mu(-\otimes\omega) )$ together with $\psi_{\Gamma(A)}=\omega$. 

This shows that the 2-category of 2-dimensional defect state sum models can be identified with $\ssFrob$, in agreement with the answer expected from the stratified cobordism hypothesis.  
In the next section we show that $\ssFrob$ is equivalent to the 2-category of semisimple Calabi--Yau categories, which will be helpful when analysing the orbifold completion of $\Bar\ssFrob$.

\subsubsection{Oriented Eilenberg--Watts theorem}
\label{subsec:EWor}

Let $\K$ be an algebraically closed field.
The (semisimple) classical Eilenberg--Watts theorem is the statement that the 2-category  $\sAlg$ of semisimple $\K$-algebras, bimodules and bimodule maps is equivalent to the 
2-category $\Tvs$ of Kapranov--Voevodsky 2-vector spaces (i.\,e.\ semisimple $\K$-linear abelian categories with finitely many isomorphism classes of simple objects), linear functors and natural transformations. 
The symmetric monoidal equivalence $\operatorname{EW}\colon \sAlg \longrightarrow \Tvs$ sends a semisimple algebra $A$ to its category of finite-dimensional left modules $A\text{-Mod}$, a bimodule ${}_B Y_A$ to the linear functor 
\begin{align}
	\EW(Y) \colon A\text{-Mod} & \longrightarrow B\text{-Mod} \\ 
	M & \longmapsto Y \otimes_A M 
\end{align}   
and a bimodule map $f\colon Y \longrightarrow X$ to the corresponding natural transformation constructed from
the morphisms $\EW(f)_M \colon  Y \otimes_A M \xrightarrow{f \otimes_A \id_M } X \otimes_A M$.    

The maximal subgroupoids $(\sAlg)^\times$ and $\Tvs^{\times}$ on both sides of the equivalence classify the space of fully extended framed TQFTs with the corresponding target. 
Passing to the 2-groupoid of homotopy fixed-points for the $\SO(2)$-action coming from the cobordism hypothesis induces an equivalence ${\ssFrob}^\times \longrightarrow (\CYCat)^\times $ between the core of the 2-category of separable symmetric Frobenius 
algebras with compatible bimodules and the 2-category of Calabi--Yau categories and Calabi--Yau functors~\cite{FEW}. 
Both sides classify fully extended oriented TQFTs with values in $\Alg$ and $\Tvs$, respectively. 
We extend this to an equivalence of pivotal 2-categories, 
\begin{align}
	\EW^{\operatorname{or}} \colon {\ssFrob} \stackrel{\cong}{\lra} \CYCat \,,
\end{align}  
which we call the \emph{oriented Eilenberg--Watts equivalence} due to its relation to oriented field theories. 
While we do not make this precise here, 
this should be thought of as the induced equivalence on categories of fully extended TQFTs with defects induced by the ordinary Eilenberg--Watts equivalence. 
In the next section we use this equivalence to systematically study the relation between spherical fusion categories and orbifold data for the trivial 3-dimensional defect TQFT. 

We begin by defining the objects of interest in the semisimple case relevant for us. 
(Many interesting examples of Calabi--Yau categories are \textsl{not} semisimple, such as derived categories of coherent sheaves on Calabi--Yau manifolds, or homotopy categories of matrix factorisations.) 

\begin{definition}
	A \emph{semisimple Calabi--Yau category} is a 2-vector space $\Ca$ together with a linear 
	trace map $\operatorname{tr}_c \colon \End (c) \longrightarrow \K$ for each $c\in \Ca$ such that  
	\begin{itemize}
		\item 
		for all morphisms $f\colon c \longrightarrow c'$ and $g\colon c' \longrightarrow c$ in $\Ca$ the relation $\operatorname{tr}_{c'}(f\circ g)= \operatorname{tr}_c({g\circ f})$ holds, and
		\item 
		for all $c,c'\in \Ca$ the pairing $\Ca(c,c') \otimes \Ca(c',c) \xrightarrow{\  \circ  \ } \Ca(c,c) \xrightarrow{\operatorname{tr}_c} \K$ is non-degenerate. 
	\end{itemize}
	A \emph{Calabi--Yau functor} $\mathcal{F} \colon (\Ca, \operatorname{tr}) \longrightarrow (\Ca' , \operatorname{tr}')$ is a linear functor $\mathcal{F}$ such that $\operatorname{tr}_c(f)= \operatorname{tr}'_{\mathcal{F}(c)}(\mathcal{F}(f)) $ for all $c\in \Ca$ and $f\in \End(c)$. We denote by 
	$\CYCat$ the symmetric monoidal 2-category of semisimple Calabi--Yau categories, linear functors and natural transformations.
\end{definition}

\begin{remark}
Note that the 1-morphisms in $\CYCat$ are all linear functors and not just the Calabi--Yau functors. 
This implies that $\CYCat$ is equivalent to $\Tvs$ as a 2-category. 
The identity functor gives an isomorphism between different Calabi--Yau structures on the same category, and every semisimple category admits a Calabi--Yau structure.   
However, we will see that $\CYCat$ is equipped with an additional structure, namely that of a pivotal 2-category. 
%, see Section~\ref{subsec:2categories}.  
\end{remark}

\begin{example}\label{Ex: CY structure on AMod}
	Let $(A, \lambda)$ be a separable symmetric Frobenius algebra, where~$\lambda$ denotes the Frobenius trace. 
	The category $A\text{-Mod}$ of finite-dimensional $A$-modules has a natural Calabi--Yau structure: 
	We give an interpretation of the trace in terms of the pivotal 2-category $\ssFrob$. 
	A left $A$-module~$M$ is a 1-morphism from the monoidal unit~$\K$ to~$A$. 
	The trace of an endomorphism $f\colon M \longrightarrow M$ is given by computing the trace in the pivotal 2-category $\ssFrob$ which gives a 2-endomorphism of~$1_\K$, which can be canonically identified with a number. 
	In the language of defect TQFTs, $M$ describes a line defect between the trivial theory and the TQFT associated to~$A$, and $f$ describes a point defect. 
	The trace agrees with the partition function of ${S}^2$ with one defect line labelled by $M$ and one point insertion labelled by $f$. 

In~\cite{FEW} a slightly different Calabi--Yau structure on $ A\text{-Mod}$ was constructed:
For $M\in A\text{-Mod}$ consider the dual right module 
$M^* \coloneqq \Hom_A(M,A)$ which comes with a natural isomorphism $\End_A(M)\cong M^*\otimes_A M$ and an 
evaluation map $\ev_M \colon  M^* \otimes_A M \longrightarrow A $. 
%arXiv_v2:
	%The linear maps $\tr_M \colon \End_A(M)
	%\cong M^*\otimes_A M \xrightarrow{\ev_M} A \xrightarrow{ \ \lambda \ } \K$ define the Calabi--Yau structure on $A\text{-Mod}$ as shown in \cite[Sect.\,3.1]{FEW}. 
	%To highlight the difference between this Calabi--Yau structure and ours, we look at the Frobenius algebra~$\K$ with counit given by an invertible scalar~$\lambda$. 
	%A $\K$-module is just a vector space.
	%The trace constructed in \cite{FEW} is the usual trace times $\lambda$. 
	%arXiv_v3:
		%Ours instead is the usual trace times $\lambda^2$. 
	 In~\cite[Sect.\,3.1]{FEW} the linear map $\textrm{tr}_M\colon \End_A(M)
	 \cong M^*\otimes_A M \xrightarrow{\ev_M} A \xrightarrow{ \ \lambda \ } \K$ is used as a trace to define the Calabi--Yau structure on $A\text{-Mod}$. 
	 We instead use the trace defined by the adjunction data from Equation~\eqref{eq:EulerAdjunctionMaps}. 
	 %arXiv_v3:
	 	%To highlight the difference between $\textrm{tr}_M$ and our Calabi--Yau trace, we look at the Frobenius algebra~$\K$ with counit given by an invertible scalar~$\lambda$, so a $\K$-module is just a vector space.
	 	%The trace $\textrm{tr}_M$ of \cite{FEW} by construction is the usual trace times~$\lambda$. 
	 	%Ours instead is the usual trace times~$\lambda^2$, which follows from the expressions in~\eqref{eq:EulerAdjunctionMaps}. 
	 	It follows that in the example where~$A$ is the Frobenius algebra~$\Bbbk$ with counit $\lambda\in\Bbbk$, the trace of~\cite{FEW} is the usual trace times~$\lambda$, while our trace is the usual one times~$\lambda^2$. 
The justification for our choice is that it gives rise to a pivotal equivalence between $\ssFrob$ and $\CYCat$ as we show in Proposition~\ref{prop:EWor} below. 
\end{example} 

\begin{example}
	\label{ex:Spherical}
	Every spherical fusion category $\Ca$ is Calabi--Yau with trace induced by the pivotal structure. 
	The tensor product is a Calabi--Yau functor $\otimes \colon \Ca \boxtimes\Ca \longrightarrow \Ca$. 
	A module category $M$ with trace 
	over $\Ca$ as defined in~\cite{Bimodtrace} is a module category $M$ over $\Ca$ together with a compatible Calabi--Yau structure (see Definition~\ref{Def:Mod cat with trace})
which implies that the action functor is Calabi--Yau. 
\end{example}

One reason we do not restrict to Calabi--Yau functors is that the adjoint of a Calabi--Yau functor is not Calabi--Yau in general. 
%arXiv_v3: 
	For example, for a spherical fusion category $\Ca$ the monoidal product $\otimes \colon \Ca \boxtimes \Ca \to \Ca$ is Calabi--Yau, but its adjoint sends the monoidal unit to $\bigoplus_{[x]\in \pi_0(\Ca)} x\boxtimes x^*$, where the direct sum runs over isomorphism classes of simple objects. 
	For this functor to be Calabi--Yau the equation $1=\tr_{\Ca}(\mathbbm{1})= \tr_{\Ca \boxtimes \Ca}\left( \otimes^\vee (\mathbbm{1})\right) = \sum_{[x]\in \pi_0(\Ca)} \tr_{\Ca}(x)^2$ would have to hold. 
	However, this is not true for most spherical fusion categories; for example for $\Z_2$-graded vector spaces the right-hand side is~2.  

There is another interplay between Calabi--Yau structures and adjoints. 
Namely, it allows us to coherently identify left and right adjoints, leading to a pivotal 2-category. 

\begin{proposition}\label{prob:CYCatPivotal}
Let $F\colon (\Ca, \tr) \to (\Ca', \tr') $ be a functor between Calabi--Yau categories, and let $F^\vee \colon \Ca'\to \Ca$ its right adjoint. 
Then $F^\vee$ is also left adjoint to $F$. 
Furthermore, this gives a pivotal structure on the 2-category $\CYCat$.
\end{proposition}

\begin{proof}
We first show that $F^\vee$ is also left adjoint. 
This follows from the sequence of natural isomorphisms 
$\Ca(c,F(c')) \cong \Ca(F(c'),c)^* \cong \Ca'(c',F^\vee(c))^* \cong \Ca'(F^\vee(c),c') $, where the identification with the dual spaces is through the trace and the middle map is the dual inverse from the isomorphism which is part of the original adjunction. 

It should be straightforward to check the conditions required for this to define a pivotal structure directly. 
However, there is also a different way: in Remark~\ref{Rem: Inverse to EW} below we construct an equivalence $\CYCat\to \ssFrob$ which is compatible with the tentative pivotal structure considered here. 
But then the relations we want to check directly follow from the fact that $\ssFrob$ is a pivotal 2-category, and that any equivalence induces a bijection on 2-morphisms.  
\end{proof}

For a pair of objects $c,c'$ in a Calabi--Yau category $(\cat{C},\tr)$ the non-degenerate trace pairing $\cat{C}(c,c') \otimes \cat{C}(c',c) \xrightarrow{\circ } \cat{C}(c,c) \xrightarrow{\tr_c} \K$ induces an isomorphism 
\be 
\label{eq:tau}
\tau_c \colon \cat{C}(c',c) \stackrel{\cong}{\lra} \cat{C}(c,c')^* 
\ee 
which we will use in Lemma~\ref{lem:FormulaLeftAdjoint} below.
To make $F^\vee\colon \mathcal{C}'\to \mathcal{C}$ a right adjoint to $F\colon \mathcal{C}\to \mathcal{C}'$ we can equivalently equip it with a unit $\eta^R\colon \id_{\mathcal{C}}\Longrightarrow F^R \circ F$ and a counit $\epsilon^R\colon F \circ F^R\Longrightarrow \id_{\mathcal{C}'}$ satisfying the Zorro identities. 
Similarly the structure of a left adjoint can be described by $\eta^L\colon \id_{\mathcal{C}'}\Longrightarrow F \circ F^L$ and a counit $\epsilon^L\colon F^L \circ F\Longrightarrow \id_{\mathcal{C}}$. 

\begin{lemma}\label{lem:FormulaLeftAdjoint}
Let $F\colon \cat{C}\to \cat{C}'$ be a functor between Calabi--Yau categories, and let $(F^\vee,\eta^R,\epsilon^R)$ be a right adjoint to $F$. 
The structure of a left adjoint on $F^\vee$ introduced in Proposition~\ref{prob:CYCatPivotal} has unit and counit given by
\begin{align}
\epsilon^L_c &= \tau_c^{-1}[\tr_{F(c)}(\epsilon^R_{F(c)}\circ F(-))] \in \cat{C}(F^RF(c),c) \,, 
\\
\eta^L_{c'} &= \tau_{c'}^{-1}[\tr_{F^R(c')}(F^R(-)\circ \eta^R_{F^R(c')})] \in \cat{C}'(c',FF^R(c')) \,.
\end{align}
\end{lemma}
 
\begin{proof}
We will only prove the formula for $\epsilon^L_c$, the proof for $\eta^L$ is analogous. 
The map $\epsilon^L_c \colon F^R \circ F(c)\to c$ is the image of the identity on $F(c)$ under the isomorphism $\cat{C}'(F(c),F(c)) \cong \cat{C}(F^RF(c),c)$. 
We just have to chase it through the isomorphism constructed in the proof of Proposition~\ref{prob:CYCatPivotal}: 
\begin{align}
\cat{C}'(F(c),F(c)) \ni \id_{F(c)} 
	& \longmapsto \tr_{F(c)}(-)\in \cat{C}'(F(c),F(c))^* 
	\nonumber
	\\
	& \longmapsto \tr_{F(c)}(\epsilon^R \circ -) \in \cat{C}'(F(c),FF^RF(c))^* 
	\nonumber
	\\
	&\longmapsto \tr_{F(c)}(\epsilon^R \circ F(-)) \in \cat{C}(c,F^RF(c))^* 
	\nonumber
	\\
	&\longmapsto \tau_c^{-1}[\tr_{F(c)}(\epsilon^R_{F(c)}\circ F(-))]  \,.
\end{align}   
Here we used that the isomorphism $\mathcal{C}'(c,F^R(c'))\cong \mathcal{C}(F(c),c')$ can be expressed in terms of the counit as 
\begin{align}
\mathcal{C}(c,F^R(c')) \xrightarrow{ \; F\; } \mathcal{C}'(F(c),FF^R(c')) \xrightarrow{\epsilon^R \circ -} \mathcal{C}'(F(c),c')	\, . 
\end{align}
\end{proof}

Now we can relate the notion of pivotal equivalence introduced in Definition~\ref{def:PivotalEquivalence} to Calabi--Yau functors in $\CYCat$. 

\begin{lemma}
An equivalence $F\colon \cat{C}\to \cat{C}'$ between semisimple Calabi--Yau categories is 
	%arXiv_v3: 
	%a pivotal equivalence 
	 pivotal 
if and only if $F$ and $F^{-1}$ are Calabi--Yau functors.  
\end{lemma}

\begin{proof}
To prove the statement it is enough to check it on simple objects. 
For this we fix a set $I$ of representatives of simple objects $i\in \cat{C}$, and similarly a set~$I'$ of simple objects $i'\in \cat{C}'$. 
The Calabi--Yau structures can be described by a collection of non-zero numbers $\{\lambda_i\}_{i\in I}$ and $\{\lambda_{i'}\}_{i'\in I'}$ defining the trace on the endomorphisms of the simple object.    
Since $F$ is an equivalence it sends a simple object $i\in \cat{C}$ to a simple $i'\in \cat{C}'$ (up to isomorphism). 
We can fix a right adjoint $G$ by setting $G(i')=i$. 
The counit $\epsilon^R_i \colon FG(i')=i'\to i'$ and unit $\eta^R_i\colon i\to GF(i)=i$ can both be chosen to be the identity. 
From Lemma~\ref{lem:FormulaLeftAdjoint} we can directly read off that $\epsilon^L_i=\tfrac{\lambda_{i'}}{\lambda_i}\id_i$ and $\eta^L_i=\tfrac{\lambda_{i}}{\lambda_{i'}}\id_{i'}$. 
We conclude that the condition that $F$ is a pivotal equivalence is equivalent to $\lambda_{i}=\lambda_{i'}$. 
\end{proof}

Replacing the 2-category $\sAlg$ by $\ssFrob$, Example~\ref{Ex: CY structure on AMod} shows that the Eilenberg--Watts equivalence lifts to an equivalence $\EW^{\operatorname{or}}\colon \ssFrob \to \CYCat$. 
This is also compatible with the pivotal structures introduced above: 

\begin{proposition}
	\label{prop:EWor}
The 2-functor $\EW^{\operatorname{or}}\colon \ssFrob \to \CYCat$ is pivotal.
\end{proposition} 

\begin{proof}
Let ${}_B M_A$ be a 1-morphism from $A$ to $B$ in $\ssFrob$, with adjoint ${}_A M^*_B$. We can assume that
the choice of right adjoint to $\EW^{\operatorname{or}}(M)$ is $\EW^{\operatorname{or}}(M^*)$ with adjunction data 
\begin{align}
\Hom_B(M\otimes_A Y, Z) & \lra \Hom_A(Y, M^* \otimes_B Z) \nonumber
\\
\tikzzbox{%
	%%%%%%%%%%%%%%%%%%%%%% 
	\begin{tikzpicture}[very thick,scale=1.0,color=blue!50!black, baseline=0cm]
	\coordinate (phi) at (0,0);
	\coordinate (M) at (-0.2,-1);
	\coordinate (Y) at (+0.2,-1);
	\coordinate (Z) at (0,+1);
	\coordinate (d1) at (-1,-1);
	\coordinate (d2) at (+1,-1);
	\coordinate (u1) at (-1,+1);
	\coordinate (u2) at (+1,+1);
	%
	% colouring: 
	\fill [orange!40!white, opacity=0.7] (d1) -- (d2) -- (u2) -- (u1); 
	\draw[thin] (d1) -- (d2) -- (u2) -- (u1) -- (d1);
	%
	% relative tensor product: 
	\fill [red!40!white, opacity=0.8] (M) -- ($(M)+(0,1)$) -- ($(Y)+(0,1)$) -- (Y); 
	%
	% strings: 
	\draw (M) -- ($(M)+(0,1)$);
	\fill (M) circle (0pt) node[above left] {{\small $M\hspace{-0.3em}$}};
	\draw (Y) -- ($(Y)+(0,1)$);
	\fill (Y) circle (0pt) node[above right] {{\small $\hspace{-0.3em}Y$}};
	\draw (Z) -- ($(Z)+(0,-1)$);
	\fill (Z) circle (0pt) node[below right] {{\small $\hspace{-0.3em}Z$}};
	%
	% labels: 
	\fill[red!80!black] (0-0.8,0.7) circle (0pt) node {{\scriptsize $B$}};
	\fill[red!80!black] (0,-0.6) circle (0pt) node {{\scriptsize $A$}};
	%
	%% PHI: 
	\fill[color=white] (phi) node[inner sep=3.9pt,draw, rounded corners=1pt, fill, color=white] (R2) {{\scriptsize$\;\psi\;$}};
	\draw[line width=1pt, color=black] (phi) node[inner sep=4pt,draw, rounded corners=1pt] (R) {{\scriptsize$\;\psi\;$}};	 
	\end{tikzpicture}
	%%%%%%%%%%%%%%%%%%%%%% 
} 
& \longmapsto 
\tikzzbox{%
	%%%%%%%%%%%%%%%%%%%%%% 
	\begin{tikzpicture}[very thick,scale=1.0,color=blue!50!black, baseline=0cm]
	\coordinate (phi) at (0,0);
	\coordinate (M) at (-0.8,1);
	\coordinate (M2) at (-0.8,-0.2);
	\coordinate (Y) at (+0.2,-1);
	\coordinate (Z) at (0,+1);
	\coordinate (d1) at (-1,-1);
	\coordinate (d2) at (+1,-1);
	\coordinate (u1) at (-1,+1);
	\coordinate (u2) at (+1,+1);
	%
	% colouring: 
	\fill [orange!40!white, opacity=0.7] (d1) -- (d2) -- (u2) -- (u1); 
	\draw[thin] (d1) -- (d2) -- (u2) -- (u1) -- (d1);
	%
	% relative tensor product: 
	\fill [purple!40!white, opacity=0.8] ($(phi)+(-0.2,-0.2)$) .. controls +(0,-0.5) and +(0,-0.5) .. (M2) -- (M) -- (Z) -- (phi); 
	%
	% strings: 
	\draw[redirected] ($(phi)+(-0.2,-0.2)$) .. controls +(0,-0.5) and +(0,-0.5) .. (M2);
	\draw (M) -- (M2); 
	\fill (M) circle (0pt) node[below right] {{\small $\hspace{-0.3em}M$}};
	\draw (Y) -- ($(Y)+(0,1)$);
	\fill (Y) circle (0pt) node[above right] {{\small $\hspace{-0.3em}Y$}};
	\draw (Z) -- ($(Z)+(0,-1)$);
	\fill (Z) circle (0pt) node[below right] {{\small $\hspace{-0.3em}Z$}};
	%
	% labels: 
	\fill[red!80!black] (-0.55,0) circle (0pt) node {{\scriptsize $B$}};
	\fill[red!80!black] (0,-0.6) circle (0pt) node {{\scriptsize $A$}};
	%
	%% PHI: 
	\fill[color=white] (phi) node[inner sep=3.9pt,draw, rounded corners=1pt, fill, color=white] (R2) {{\scriptsize$\;\psi\;$}};
	\draw[line width=1pt, color=black] (phi) node[inner sep=4pt,draw, rounded corners=1pt] (R) {{\scriptsize$\;\psi\;$}};	 
	\end{tikzpicture}
	%%%%%%%%%%%%%%%%%%%%%% 
} 
\end{align} 
where here and below we graphically denote the relative tensor product by a stronger colouring. 
This means that the right vertical arrow in~\eqref{eq:ConditionPivotal} is the 
identity, so we need to show that the left vertical map is also the identity. 
For this first note that the adjunction data coming from $\ssFrob$ is 
\begin{align}
\Hom_B(Z,M\otimes_A Y) & \lra \Hom_A(M^* \otimes_B Z,Y) \nonumber
\\
\tikzzbox{%
	%%%%%%%%%%%%%%%%%%%%%% 
	\begin{tikzpicture}[very thick,scale=1.0,color=blue!50!black, baseline=0cm]
	\coordinate (phi) at (0,0);
	\coordinate (M) at (-0.2,1);
	\coordinate (Y) at (+0.2,1);
	\coordinate (Z) at (0,-1);
	\coordinate (d1) at (-1,-1);
	\coordinate (d2) at (+1,-1);
	\coordinate (u1) at (-1,+1);
	\coordinate (u2) at (+1,+1);
	%
	% colouring: 
	\fill [orange!40!white, opacity=0.7] (d1) -- (d2) -- (u2) -- (u1); 
	\draw[thin] (d1) -- (d2) -- (u2) -- (u1) -- (d1);
	%
	% relative tensor product: 
	\fill [red!40!white, opacity=0.8] (M) -- ($(M)+(0,-1)$) -- ($(Y)+(0,-1)$) -- (Y); 
	%
	% strings: 
	\draw (M) -- ($(M)+(0,-1)$);
	\fill (M) circle (0pt) node[below left] {{\small $M\hspace{-0.3em}$}};
	\draw (Y) -- ($(Y)+(0,-1)$);
	\fill (Y) circle (0pt) node[below right] {{\small $\hspace{-0.3em}Y$}};
	\draw (Z) -- ($(Z)+(0,1)$);
	\fill (Z) circle (0pt) node[above right] {{\small $\hspace{-0.3em}Z$}};
	%
	% labels: 
	\fill[red!80!black] (-0.8,0.7) circle (0pt) node {{\scriptsize $B$}};
	\fill[red!80!black] (0,0.6) circle (0pt) node {{\scriptsize $A$}};
	%
	%% PHI: 
	\fill[color=white] (phi) node[inner sep=3.9pt,draw, rounded corners=1pt, fill, color=white] (R2) {{\scriptsize$\;\varphi\;$}};
	\draw[line width=1pt, color=black] (phi) node[inner sep=4pt,draw, rounded corners=1pt] (R) {{\scriptsize$\;\varphi\;$}};	 
	\end{tikzpicture}
	%%%%%%%%%%%%%%%%%%%%%% 
} 
& \longmapsto 
\tikzzbox{%
	%%%%%%%%%%%%%%%%%%%%%% 
	\begin{tikzpicture}[very thick,scale=1.0,color=blue!50!black, baseline=0cm]
	\coordinate (phi) at (0,0);
	\coordinate (M) at (-0.8,0);
	\coordinate (M2) at (-0.8,-1);
	\coordinate (Y) at (+0.2,1);
	\coordinate (Z) at (0,-1);
	\coordinate (d1) at (-1,-1);
	\coordinate (d2) at (+1,-1);
	\coordinate (u1) at (-1,+1);
	\coordinate (u2) at (+1,+1);
	%
	% colouring: 
	\fill [orange!40!white, opacity=0.7] (d1) -- (d2) -- (u2) -- (u1); 
	\draw[thin] (d1) -- (d2) -- (u2) -- (u1) -- (d1);
	%
	% relative tensor product: 
	\fill [purple!40!white, opacity=0.8] (phi) -- (-0.2,0.2) .. controls +(0,0.5) and +(0,0.5) .. ($(M)+(0,0.2)$) -- (M2) -- (Z); 
	%
	% strings: 
	\draw[directed] (-0.2,0.2) .. controls +(0,0.5) and +(0,0.5) .. ($(M)+(0,0.2)$);
	\draw ($(M)+(0,0.2)$) -- (M2);
	\fill (M2) circle (0pt) node[above right] {{\small $\hspace{-0.3em}M$}};
	\draw (Y) -- ($(Y)+(0,-1)$);
	\fill (Y) circle (0pt) node[below right] {{\small $\hspace{-0.3em}Y$}};
	\draw (Z) -- ($(Z)+(0,1)$);
	\fill (Z) circle (0pt) node[above right] {{\small $\hspace{-0.3em}Z$}};
	%
	% labels: 
	\fill[red!80!black] (-0.5,0) circle (0pt) node {{\scriptsize $B$}};
	\fill[red!80!black] (0,0.6) circle (0pt) node {{\scriptsize $A$}};
	%
	%% PHI: 
	\fill[color=white] (phi) node[inner sep=3.9pt,draw, rounded corners=1pt, fill, color=white] (R2) {{\scriptsize$\;\varphi\;$}};
	\draw[line width=1pt, color=black] (phi) node[inner sep=4pt,draw, rounded corners=1pt] (R) {{\scriptsize$\;\varphi\;$}};	 
	\end{tikzpicture}
	%%%%%%%%%%%%%%%%%%%%%% 
} 
\end{align}
The statement now follows from the commutativity of the diagram 
\be 
\begin{tikzcd}[column sep=3em, row sep=2em]
\Hom_B(Z, M\otimes_A Y) 
\ar[rr]
\ar[ddd, out=-120, in=120, shift right=10] 
&&
\Hom_A(M^*\otimes_B Z,Y)
\ar[ddd, out=-60, in=60, shift left=10] 
\\[-1.5em]
\tikzzbox{%
	%%%%%%%%%%%%%%%%%%%%%% 
	\begin{tikzpicture}[very thick,scale=1.0,color=blue!50!black]
	\coordinate (phi) at (0,0);
	\coordinate (M) at (-0.2,1);
	\coordinate (Y) at (+0.2,1);
	\coordinate (Z) at (0,-1);
	\coordinate (d1) at (-1,-1);
	\coordinate (d2) at (+1,-1);
	\coordinate (u1) at (-1,+1);
	\coordinate (u2) at (+1,+1);
	%
	% colouring: 
	\fill [orange!40!white, opacity=0.7] (d1) -- (d2) -- (u2) -- (u1); 
	\draw[thin] (d1) -- (d2) -- (u2) -- (u1) -- (d1);
	%
	% relative tensor product: 
	\fill [red!40!white, opacity=0.8] (M) -- ($(M)+(0,-1)$) -- ($(Y)+(0,-1)$) -- (Y); 
	%
	% strings: 
	\draw (M) -- ($(M)+(0,-1)$);
	\fill (M) circle (0pt) node[below left] {{\small $M\hspace{-0.3em}$}};
	\draw (Y) -- ($(Y)+(0,-1)$);
	\fill (Y) circle (0pt) node[below right] {{\small $\hspace{-0.3em}Y$}};
	\draw (Z) -- ($(Z)+(0,1)$);
	\fill (Z) circle (0pt) node[above right] {{\small $\hspace{-0.3em}Z$}};
	%
	% labels: 
	\fill[red!80!black] (-0.8,0.7) circle (0pt) node {{\scriptsize $B$}};
	\fill[red!80!black] (0,0.6) circle (0pt) node {{\scriptsize $A$}};
	\fill[black] (0,1.5) circle (0pt) node {\rotatebox{90}{$\in$}};
	%
	%% PHI: 
	\fill[color=white] (phi) node[inner sep=3.9pt,draw, rounded corners=1pt, fill, color=white] (R2) {{\scriptsize$\;\varphi\;$}};
	\draw[line width=1pt, color=black] (phi) node[inner sep=4pt,draw, rounded corners=1pt] (R) {{\scriptsize$\;\varphi\;$}};	 
	\end{tikzpicture}
	%%%%%%%%%%%%%%%%%%%%%% 
} 
\ar[rr,mapsto]
\ar[d,mapsto]
&&
\tikzzbox{%
	%%%%%%%%%%%%%%%%%%%%%% 
	\begin{tikzpicture}[very thick,scale=1.0,color=blue!50!black]
	\coordinate (phi) at (0,0);
	\coordinate (M) at (-0.8,0);
	\coordinate (M2) at (-0.8,-1);
	\coordinate (Y) at (+0.2,1);
	\coordinate (Z) at (0,-1);
	\coordinate (d1) at (-1,-1);
	\coordinate (d2) at (+1,-1);
	\coordinate (u1) at (-1,+1);
	\coordinate (u2) at (+1,+1);
	%
	% colouring: 
	\fill [orange!40!white, opacity=0.7] (d1) -- (d2) -- (u2) -- (u1); 
	\draw[thin] (d1) -- (d2) -- (u2) -- (u1) -- (d1);
	%
	% relative tensor product: 
	\fill [purple!40!white, opacity=0.8] (phi) -- (-0.2,0.2) .. controls +(0,0.5) and +(0,0.5) .. ($(M)+(0,0.2)$) -- (M2) -- (Z); 
	%
	% strings: 
	\draw[directed] (-0.2,0.2) .. controls +(0,0.5) and +(0,0.5) .. ($(M)+(0,0.2)$);
	\draw ($(M)+(0,0.2)$) -- (M2);
	\fill (M2) circle (0pt) node[above right] {{\small $\hspace{-0.3em}M$}};
	\draw (Y) -- ($(Y)+(0,-1)$);
	\fill (Y) circle (0pt) node[below right] {{\small $\hspace{-0.3em}Y$}};
	\draw (Z) -- ($(Z)+(0,1)$);
	\fill (Z) circle (0pt) node[above right] {{\small $\hspace{-0.3em}Z$}};
	%
	% labels: 
	\fill[red!80!black] (-0.5,0) circle (0pt) node {{\scriptsize $B$}};
	\fill[red!80!black] (0,0.6) circle (0pt) node {{\scriptsize $A$}};
	\fill[black] (0,1.5) circle (0pt) node {\rotatebox{90}{$\in$}};
	%
	%% PHI: 
	\fill[color=white] (phi) node[inner sep=3.9pt,draw, rounded corners=1pt, fill, color=white] (R2) {{\scriptsize$\;\varphi\;$}};
	\draw[line width=1pt, color=black] (phi) node[inner sep=4pt,draw, rounded corners=1pt] (R) {{\scriptsize$\;\varphi\;$}};	 
	\end{tikzpicture}
	%%%%%%%%%%%%%%%%%%%%%% 
} 
\ar[d,mapsto]
\\
\tikzzbox{%
	%%%%%%%%%%%%%%%%%%%%%% 
	\begin{tikzpicture}[very thick,scale=1.0,color=blue!50!black]
	\coordinate (phi) at (0,-0.75);
	\coordinate (blank) at (0,0.75);
	\coordinate (M) at (-0.8,0);
	\coordinate (M2) at (-0.2,0);
	\coordinate (Y) at (+0.2,0);
	\coordinate (Z) at (0,-1);
	\coordinate (d1) at (-2,-2);
	\coordinate (d2) at (+1,-2);
	\coordinate (u1) at (-2,+2);
	\coordinate (u2) at (+1,+2);
	%
	% colouring: 
	\fill [orange!40!white, opacity=0.7] (d1) -- (d2) -- (u2) -- (u1); 
	\draw[thin] (d1) -- (d2) -- (u2) -- (u1) -- (d1);
	%
	% relative tensor product: 
	\fill [red!40!white, opacity=0.8] ($(phi)+(-0.2,0)$) -- ($(phi)+(0.2,0)$) -- ($(blank)+(0.2,0)$) -- ($(blank)+(-0.2,0)$); 
	%
	% strings: 
	\draw ($(phi)+(0.2,0)$) -- ($(blank)+(0.2,0)$);
	\draw ($(phi)+(-0.2,0)$) -- ($(blank)+(-0.2,0)$);
	\draw (blank) -- ($(blank)+(0,0.5)$);
	\draw[directed] ($(blank)+(0,0.5)$) .. controls +(0,0.75) and +(0,0.75) .. ($(blank)+(-1.5,0.5)$);
	\draw ($(phi)+(-1.5,-0.5)$) -- ($(blank)+(-1.5,0.5)$);
	\draw[redirected] ($(phi)+(0,-0.5)$) .. controls +(0,-0.75) and +(0,-0.75) .. ($(phi)+(-1.5,-0.5)$);
	\draw (phi) -- ($(phi)+(0,-0.5)$);
	%
	% labels: 
	\fill[red!80!black] (-1,0) circle (0pt) node {{\scriptsize $B$}};
	\fill[red!80!black] (0,0) circle (0pt) node {{\scriptsize $A$}};
	\fill (M2) circle (0pt) node[left] {{\small $M\hspace{-0.3em}$}};
	\fill (Y) circle (0pt) node[right] {{\small $\hspace{-0.3em}Y$}};
	\fill ($(phi)+(0,-0.5)$) circle (0pt) node[right] {{\small $\hspace{-0.3em}Z$}};
	\fill[black] (-0.5,-2.75) circle (0pt) node {\rotatebox{-90}{$\in$}};
	%
	%% BLANK: 
	\fill[color=white] (blank) node[inner sep=3.9pt,draw, rounded corners=1pt, fill, color=white] (R2) {{\scriptsize$\;-\;$}};
	\draw[line width=1pt, color=black] (blank) node[inner sep=4pt,draw, rounded corners=1pt] (R) {{\scriptsize$\;-\;$}};	 
	%
	%% PHI: 
	\fill[color=white] (phi) node[inner sep=3.9pt,draw, rounded corners=1pt, fill, color=white] (R2) {{\scriptsize$\;\varphi\;$}};
	\draw[line width=1pt, color=black] (phi) node[inner sep=4pt,draw, rounded corners=1pt] (R) {{\scriptsize$\;\varphi\;$}};	 
	\end{tikzpicture}
	%%%%%%%%%%%%%%%%%%%%%% 
} 
\ar[r,mapsto]
&
\tikzzbox{%
	%%%%%%%%%%%%%%%%%%%%%% 
	\begin{tikzpicture}[very thick,scale=1.0,color=blue!50!black]
	\coordinate (phi) at (0,-0.75);
	\coordinate (blank) at (0,0.75);
	\coordinate (M) at (-0.8,0);
	\coordinate (M2) at (-0.2,0);
	\coordinate (Y) at (+0.2,0);
	\coordinate (Z) at (0,-1);
	\coordinate (d1) at (-2,-2);
	\coordinate (d2) at (+1,-2);
	\coordinate (u1) at (-2,+2);
	\coordinate (u2) at (+1,+2);
	%
	% colouring: 
	\fill [orange!40!white, opacity=0.7] (d1) -- (d2) -- (u2) -- (u1); 
	\draw[thin] (d1) -- (d2) -- (u2) -- (u1) -- (d1);
	%
	% relative tensor product: 
	\fill [red!40!white, opacity=0.8] ($(phi)+(-0.2,0)$) .. controls +(0,0.75) and +(0,-0.75) .. ($(blank)+(-1,0.2)$) -- ($(blank)+(-1,0.2)$) .. controls +(0,0.5) and +(0,0.5) .. ($(blank)+(-0.2,0.2)$)
	-- ($(blank)+(0.2,0.2)$) -- ($(phi)+(0.2,0)$); 
	%
	% strings: 
	\draw ($(phi)+(0.2,0)$) -- ($(blank)+(0.2,0)$);
	\draw ($(phi)+(-0.2,0)$) .. controls +(0,0.75) and +(0,-0.75) .. ($(blank)+(-1,0.2)$);
	\draw[directed] ($(blank)+(-1,0.2)$) .. controls +(0,0.5) and +(0,0.5) .. ($(blank)+(-0.2,0.2)$);
	\draw (blank) -- ($(blank)+(0,0.5)$);
	\draw[directed] ($(blank)+(0,0.5)$) .. controls +(0,0.75) and +(0,0.75) .. ($(blank)+(-1.5,0.5)$);
	\draw ($(phi)+(-1.5,-0.5)$) -- ($(blank)+(-1.5,0.5)$);
	\draw[redirected] ($(phi)+(0,-0.5)$) .. controls +(0,-0.75) and +(0,-0.75) .. ($(phi)+(-1.5,-0.5)$);
	\draw (phi) -- ($(phi)+(0,-0.5)$);
	%
	% labels: 
	\fill[red!80!black] (-1,0) circle (0pt) node {{\scriptsize $B$}};
	\fill[red!80!black] (0,0) circle (0pt) node {{\scriptsize $A$}};
	\fill (M2) circle (0pt) node {{\small $\hspace{-0.5em}M$}};
	\fill (Y) circle (0pt) node[right] {{\small $\hspace{-0.3em}Y$}};
	\fill ($(phi)+(0,-0.5)$) circle (0pt) node[right] {{\small $\hspace{-0.3em}Z$}};
	%
	%% BLANK: 
	\fill[color=white] (blank) node[inner sep=3.9pt,draw, rounded corners=1pt, fill, color=white] (R2) {{\scriptsize$\;-\;$}};
	\draw[line width=1pt, color=black] (blank) node[inner sep=4pt,draw, rounded corners=1pt] (R) {{\scriptsize$\;-\;$}};	 
	%
	%% PHI: 
	\fill[color=white] (phi) node[inner sep=3.9pt,draw, rounded corners=1pt, fill, color=white] (R2) {{\scriptsize$\;\varphi\;$}};
	\draw[line width=1pt, color=black] (phi) node[inner sep=4pt,draw, rounded corners=1pt] (R) {{\scriptsize$\;\varphi\;$}};	 
	\end{tikzpicture}
	%%%%%%%%%%%%%%%%%%%%%% 
} 
\ar[r,equal] 
& 
\tikzzbox{%
	%%%%%%%%%%%%%%%%%%%%%% 
	\begin{tikzpicture}[very thick,scale=1.0,color=blue!50!black]
	\coordinate (phi) at (0,-0.75);
	\coordinate (blank) at (0,0.75);
	\coordinate (M) at (-0.8,0);
	\coordinate (M2) at (-0.2,0);
	\coordinate (Y) at (+0.2,0);
	\coordinate (Z) at (0,-1);
	\coordinate (d1) at (-2,-2);
	\coordinate (d2) at (+1,-2);
	\coordinate (u1) at (-2,+2);
	\coordinate (u2) at (+1,+2);
	%
	% colouring: 
	\fill [orange!40!white, opacity=0.7] (d1) -- (d2) -- (u2) -- (u1); 
	\draw[thin] (d1) -- (d2) -- (u2) -- (u1) -- (d1);
	%
	% relative tensor product: 
	\fill [purple!40!white, opacity=0.8] ($(phi)+(-0.2,0)$) .. controls +(0,0.75) and +(0,-0.75) .. ($(blank)+(-1,0.2)$) -- ($(blank)+(-1,0.2)$) .. controls +(0,0.5) and +(0,0.5) .. ($(blank)+(-0.2,0.2)$)
	-- ($(blank)+(0.2,0.2)$) -- ($(phi)+(0.2,0)$); 
	%
	% strings: 
	\draw ($(phi)+(0.2,0)$) -- ($(blank)+(0.2,0)$);
	\draw ($(phi)+(-0.2,0)$) .. controls +(0,0.75) and +(0,-0.75) .. ($(blank)+(-1,0.2)$);
	\draw[redirected] ($(blank)+(-1,0.2)$) .. controls +(0,0.5) and +(0,0.5) .. ($(blank)+(-0.2,0.2)$);
	\draw (blank) -- ($(blank)+(0,0.5)$);
	\draw[directed] ($(blank)+(0,0.5)$) .. controls +(0,0.75) and +(0,0.75) .. ($(blank)+(-1.5,0.5)$);
	\draw ($(phi)+(-1.5,-0.5)$) -- ($(blank)+(-1.5,0.5)$);
	\draw[redirected] ($(phi)+(0,-0.5)$) .. controls +(0,-0.75) and +(0,-0.75) .. ($(phi)+(-1.5,-0.5)$);
	\draw (phi) -- ($(phi)+(0,-0.5)$);
	%
	% labels: 
	\fill[red!80!black] (-1,0) circle (0pt) node {{\scriptsize $A$}};
	\fill[red!80!black] (0,0) circle (0pt) node {{\scriptsize $B$}};
	\fill (M2) circle (0pt) node {{\small $\hspace{-0.5em}M$}};
	\fill (Y) circle (0pt) node[right] {{\small $\hspace{-0.3em}Z$}};
	\fill ($(phi)+(0,-0.5)$) circle (0pt) node[right] {{\small $\hspace{-0.3em}Y$}};
	\fill[black] (-0.5,-2.75) circle (0pt) node {\rotatebox{-90}{$\in$}};
	%
	%% BLANK: 
	\fill[color=white] (blank) node[inner sep=3.9pt,draw, rounded corners=1pt, fill, color=white] (R2) {{\scriptsize$\;\psi\;$}};
	\draw[line width=1pt, color=black] (blank) node[inner sep=4pt,draw, rounded corners=1pt] (R) {{\scriptsize$\;\psi\;$}};	 
	%
	%% PHI: 
	\fill[color=white] (phi) node[inner sep=3.9pt,draw, rounded corners=1pt, fill, color=white] (R2) {{\scriptsize$\;-\;$}};
	\draw[line width=1pt, color=black] (phi) node[inner sep=4pt,draw, rounded corners=1pt] (R) {{\scriptsize$\;-\;$}};	 
	\end{tikzpicture}
	%%%%%%%%%%%%%%%%%%%%%% 
} 
\\
\Hom_B(M\otimes_A Y, Z)^*
\ar[rr]
&&
\Hom_A(Y,M^*\otimes_B Z)^*
\end{tikzcd}
\ee 
where the anti-clockwise composition starting in the top left corner and ending in the top right corner is the left adjoint structure coming from the Calabi--Yau structure. 
\end{proof}

\begin{remark}\label{Rem: Inverse to EW}
Note that Proposition~\ref{prop:EWor} directly implies that there exist an inverse to $\EW^{\operatorname{or}}$ which is also a pivotal functor.
A pivotal inverse $(\EW^{\operatorname{or}})^{-1}\colon \CYCat \longrightarrow \ssFrob$ can be described explicitly by picking for every Calabi--Yau category $\mathcal{C}\in \CYCat$ a set of representatives $I$ for all isomorphism classes of simple objects in $\mathcal{C}$. 
The 2-functor $(\EW^{\operatorname{or}})^{-1}$ sends a Calabi--Yau category $(\mathcal{C},\tr)\in \CYCat$ to the algebra $\End(\bigoplus_{i\in I} i )$ which has a Frobenius structure coming from $\tr$. 
However, this is not the Frobenius structure on $(\EW^{\operatorname{or}})^{-1}$ which can be constructed by picking a square root for the window element~$\omega$ and changing the Frobenius structure by its inverse as we did in Section~\ref{subsec:2Dstatesums}. 
That such square roots exist was shown in~\cite{Mule1}.  
A functor $F\colon \mathcal{C}\longrightarrow \mathcal{C}'$ is sent to the bimodule $\Hom_{\mathcal{C}'}(\bigoplus_{j\in I'} j, F(\bigoplus_{i\in I} i ) ) $. A natural transformation gets mapped to the post-composition with its component at $\bigoplus_{i\in I} i$. 
That this is pivotal follows because $\EW^{\operatorname{or}}(\End(\bigoplus_{i\in I} i ))$ is pivotally equivalent to $\mathcal{C}$, and because of Proposition~\ref{Prop: Whitehead}.  
\end{remark}

\begin{remark}
We expect that the above extends to higher dimensions.
%arXiv_v2: 
	%In \cite{LukasNilsVincentas} we will consider the case $n=3$, i.\,e.\ an equivalence $2\EW \colon \mFus \longrightarrow \TTvs$ between multifusion categories and 3-vector spaces from~\cite{2EW}, and then provide an oriented version, i.\,e.\ an equivalence $2\EW^{\operatorname{or}}\colon \msFus \longrightarrow \CYTCat$ of 3-categories between the 3-category of spherical multifusion categories and Calabi--Yau 2-categories, to study the relation between spherical fusion 2-categories and orbifold data for the trivial 4-dimensional defect TQFT. 
	 In the case $n=3$, there is an equivalence $2\EW \colon \mFus \longrightarrow \TTvs$ between multifusion categories and 3-vector spaces proven in~\cite{2EW}, for which one expects also an oriented version, i.\,e.\ an equivalence $2\EW^{\operatorname{or}}\colon \msFus \longrightarrow \CYTCat$ of 3-categories between the 3-category of spherical multifusion categories and appropriate Calabi--Yau 2-categories. 
	 This is relevant for the relation between spherical fusion 2-categories and orbifold data for the trivial 4-dimensional defect TQFT, cf.\ \cite{LukasNilsVincentas}. 
\end{remark}

\subsubsection{3-dimensional state sum models} 
\label{subsubsec:3dStateSumModels}

As we have seen in Section~\ref{subsec:2Dstatesums} the 2-category describing 2-dimensional defect state sum models is $\ssFrob$. 
Following the ideas outlined in Section~\ref{sec:introduction} its delooping describes defects of the trivial 3-dimensional TQFT. 
%arXiv_v3: 
	%The 3-category describing 3-dimensional defect state sum models is (the Euler completion of) its orbifold completion. 
	 The 3-category describing 3-dimensional defect state sum models is the orbifold completion of its Euler completion. 
To analyse the orbifold completion we can employ the oriented Eilenberg--Watts equivalence to work within $\CYCat$ as explained in the previous section. 

To cover more interesting orbifold data, we have to work in an appropriate Euler completion. This means we consider the 3-category $E(\Bar\ssFrob)$ whose objects are given by pairs consisting of an object of $\Bar\ssFrob$ (there is only one) and an invertible 3-morphism $\phi\colon 1_{1_*}\longrightarrow 1_{1_*}$. 
%arXiv_v3: 
	%Explicitly this 3-category can be constructed by considering the pivotal 3-category~\cite{CMS} associated to the Euler completion of the corresponding defect field theory constructed in~\cite{CRS1}. 
	 We identify the 3-category $E(\Bar\ssFrob)$ with the one constructed as the Gray category with duals associated in~\cite{CMS} to the Euler completion~\cite{CRS1} of the ``trivial'' 3-dimensional defect TQFT constructed in \cite[Sect.\,4.2]{CRS3}; we do so in the spirit of the second paragraph of Section~\ref{sec:ExamppleCompletion}, noting that $\Bar\ssFrob$ is not a Gray category (but equivalent to one).
This construction also involves adding Euler defects for surface defects. 
Here we can ignore them because they are already taken care of by the Euler completion of $\Delta\ssFrob$ in Section~\ref{subsec:2Dstatesums}. 
In practice this leads to including $\phi$-insertions into the condition on orbifold data.     

Our first goal is to relate $\orb{E(\Bar\ssFrob)}$ to spherical fusion categories and their bimodule categories with trace. 
In~\cite{CRS3} it was already shown that spherical fusion categories give rise to orbifold data. 
From the perspective of Section~\ref{subsec:EWor} this construction can be understood as follows. 
Recall from Example~\ref{ex:Spherical} that spherical fusion categories~$S$ gives rise to $E_1$-algebras in $\CYCat$. 
Using the square root of the global dimension $\phi = (\dim S)^{-1/2}$ as the Euler defects upgrades this to an orbifold datum in $E(\Bar\CYCat)$. 
The paper~\cite{CRS3} works with the image of $(\EW^{\operatorname{or}})^{-1}$ in $\ssFrob$ when doing computations. 

Before restricting to the subcategory of orbifold data corresponding to spherical fusion categories we make a few comments on the general structure of orbifold data in $E(\Bar\CYCat)$. 
For this we fix an object $(\cat{C},\otimes, 1, \alpha)\in \mathcal{E}(\Bar\CYCat)$. 
This is a monoidal category~$\mathcal C$ which in addition is equipped with a Calabi--Yau structure. 
For this to be an orbifold datum, $(\cat{C},\otimes, 1, \alpha)$ has to satisfy the conditions in Figure~\ref{fig:OrbifoldDatumAxioms}. 
These conditions give rise to additional properties of the monoidal structure. One of the conditions is that~$\cat{C}$ is rigid. 
To show this we use an equivalent reformulation of rigidity: 

\begin{proposition}
	\label{prop:RigidStrong}
	Let $A \in \Tvs$ be endowed with the structure of a monoidal category. The following are equivalent: 
	\begin{enumerate}
		\item $A$ is rigid.
		\item the canonical lax bimodule structure on the right adjoint 
%		$\Delta_{\textrm{R}} \colon A \longrightarrow A \boxtimes A $ 
		of $\otimes \colon A \boxtimes A \longrightarrow A $ is strong,  
		\item the canonical op-lax bimodule structure on the left adjoint 
%		$\Delta_{\textrm{L}} \colon A \longrightarrow A \boxtimes A $ 
		of $\otimes \colon A \boxtimes A \longrightarrow A $ is strong. 
	\end{enumerate}
\end{proposition}  

\begin{proof}
	This is essentially~\cite[Def.-Prop.\,4.1]{BJS} applied to the semisimple case, where adjoints always exist. 
\end{proof}

\begin{remark}
	The equivalence breaks down as soon as one leaves the realm of 2-vector spaces. 
	For example for compactly generated presentable categories condition (ii) is equivalent to the condition that every compact projective 
	object has a left and right dual, see \cite[Def.-Prop.\,4.1]{BJS}.  
\end{remark} 

\begin{remark}\label{Rem: rigid}
The natural lax and op-lax bimodule structures are exactly given by $\alpha', \alpha''$, $\bar{\alpha}'$ and $\bar{\alpha}''$ as in~\eqref{eq:alpha}. 
The orbifold conditions enforce more than just the fact that these are invertible; they give specific inverses.  
As in every pivotal 2-category, the left and right adjoints of 1-morphisms between Calabi--Yau categories agree. 
Thanks to Proposition~\ref{prop:RigidStrong}, the relations in Figure~\ref{fig:OrbifoldDatumAxioms} except the bubble cancellation can now be reformulated as the statement that the two bimodule structures on the left and right adjoint agree under this identification (up to a scalar since we work in the Euler completion).    
	%Let us explain how to construct the canonical lax and op-lax bimodule structures. As a right adjoint to $\otimes $, $\Delta_R$ comes with natural transformations $\epsilon \colon \otimes \circ \Delta_R \longrightarrow \id_A$ and $\eta \colon \id_{A \boxtimes A} \longrightarrow \Delta_R \circ \otimes $ which 
%	we represent graphically as:
%	\begin{center}
%		\includegraphics{adjunction.pdf}
%	\end{center}
%	satisfying the usual snake type identities: 
%	\begin{center}
%		\includegraphics{relation_Adjunction.pdf}
%	\end{center}
%	Now the canonical lax left module functor structure for $\Delta_R$ is given by:
%	\begin{align}
%	(\id_A \boxtimes \Delta)
%	\end{align}
%	\begin{center}
%		\includegraphics{lax.pdf}
%	\end{center} 
%	The lax structure for the right action and the ob-lax versions for the left adjoint are obtained analogously 
%	and correspond to $\alpha', \alpha''$ and $\bar{\alpha}'$ from Equation~\eqref{eq:alpha}.   
\end{remark} 

%If $A$ is rigid we have a concrete description for the left and right adjoint. 
%\begin{lemma}
%	Let $(A,\otimes)$ be a rigid monoidal 2-vector space. There are canonical isomorphisms of bimodule functors
%	\begin{align}
%		\Delta_R(a) &\cong \int^{c\in C} a\otimes c^\vee \boxtimes c \cong \int^{c\in C}  c\boxtimes {}^\vee c \otimes a \\  
%		\Delta_L(a) &\cong \int^{c\in C} a\otimes {}^\vee c \boxtimes c \cong \int^{c\in C}  c\boxtimes c^\vee \otimes a
%	\end{align}
%\end{lemma}
%\begin{proof}
%	Check the universal property... Be careful with the bimodule functor structure.... \LM{I have not check this. atm its not really a statement and more like a conjecture...}
%\end{proof}
%
%\begin{remark}
%	In the CP-rigid case the coend over $A$ is replaced by a coend over 
%	the compact projective objects of $A$~\cite[Section 3.2]{BJSS}.
%\end{remark}

\begin{proposition} 
Let $A$ be an object of $\orb{E(\Bar\CYCat)}$. 
Then $A$ is a rigid linear monoidal category. 
\end{proposition}

\begin{proof}
This is a direct consequence of Proposition~\ref{prop:RigidStrong} and Remark~\ref{Rem: rigid}. 
\end{proof}
%arXiv_v2:
	%We
	Based on this proposition one might expect that every orbifold datum arises from a spherical 
	%arXiv_v3: 
		%(multi-)fusion 
		 \mbox{(multi-)}fusion 
	category. 
	However, it seems hard to deduce the existence and concrete choice of a pivotal structure directly from the orbifold conditions. 
	For this reason, we 
now restrict to orbifold data coming from spherical fusion categories and turn to 1-morphisms in $\orb{E(\Bar\CYCat)}$. 
According to Definition~\ref{def:Corb}, a 1-morphism $(\cat{C},\otimes, 1 , \alpha) \lra (\cat{C}',\otimes', 1' , \alpha')$ is given by a $\cat{C}$-$\cat{C'}$-bimodule category $\cat{M}$ which is equipped with a Calabi--Yau structure, subject to the conditions from Definition~\ref{def:Corb}.
For orbifold data coming from spherical fusion categories this is related to the notion of a trace on a bimodule category, introduced in \cite[Def.\,3.7]{Bimodtrace}.

\begin{definition}
	\label{Def:Mod cat with trace}
Let $\mathcal{C},\mathcal{C}'$ be spherical fusion categories and $\mathcal{M}$ a $\mathcal{C}$-$\mathcal{C}'$-bimodule category. 
A \textsl{trace} on $\mathcal{M}$ is a choice of Calabi--Yau structure $\tr_\mathcal{M}$ on $\mathcal{M}$ such that 
\be 
\operatorname{tr}_{c\triangleright m} \left( 
\,
\tikzzbox{%
	%%%%%%%%%%%%%%%%%%%%%% 
	\begin{tikzpicture}[very thick,scale=0.85,color=blue!50!black, baseline=-0.1cm]
	\coordinate (phi) at (0,0);
	\coordinate (c) at (-0.2,-1);
	\coordinate (m) at (+0.2,-1);
	\coordinate (c2) at (-0.2,1);
	\coordinate (m2) at (+0.2,1);
	\coordinate (d1) at (-1,-1);
	\coordinate (d2) at (+1,-1);
	\coordinate (u1) at (-1,+1);
	\coordinate (u2) at (+1,+1);
	%
	% colouring: 
	\fill [orange!40!white, opacity=0.7] (d1) -- (d2) -- (u2) -- (u1); 
	\draw[thin] (d1) -- (d2) -- (u2) -- (u1) -- (d1);
	%
	% strings: 
	\draw[color=green!50!black] (c) -- ($(c)+(0,1)$);
	\draw (m) -- ($(m)+(0,1)$);
	\draw[color=green!50!black] (c2) -- ($(c)+(0,1)$);
	\draw (m2) -- ($(m)+(0,1)$);
	%
	% labels: 
	\fill[color=green!50!black] (c) circle (0pt) node[above left] {{\small $c\hspace{-0.2em}$}};
	\fill (m) circle (0pt) node[above right] {{\small $\hspace{-0.2em}m$}};	
	\fill[color=green!50!black] (c2) circle (0pt) node[below left] {{\small $c\hspace{-0.2em}$}};
	\fill (m2) circle (0pt) node[below right] {{\small $\hspace{-0.2em}m$}};	
	%
	%% PHI: 
	\fill[color=white] (phi) node[inner sep=3.9pt,draw, rounded corners=1pt, fill, color=white] (R2) {{\scriptsize$\;\xi\;$}};
	\draw[line width=1pt, color=black] (phi) node[inner sep=4pt,draw, rounded corners=1pt] (R) {{\scriptsize$\;\xi\;$}};	 
	\end{tikzpicture}
	%%%%%%%%%%%%%%%%%%%%%% 
} 
\,
\right) 
= 
\operatorname{tr}_{m} \left( 
\,
\tikzzbox{%
	%%%%%%%%%%%%%%%%%%%%%% 
	\begin{tikzpicture}[very thick,scale=0.85,color=blue!50!black, baseline=-0.1cm]
	\coordinate (phi) at (0,0);
	\coordinate (c) at (-0.2,-0.2);
	\coordinate (c') at (-0.8,-0.2);
	\coordinate (m) at (+0.2,-1);
	\coordinate (c2) at (-0.2,0.2);
	\coordinate (c2') at (-0.8,0.2);
	\coordinate (m2) at (+0.2,1);
	\coordinate (d1) at (-1,-1);
	\coordinate (d2) at (+1,-1);
	\coordinate (u1) at (-1,+1);
	\coordinate (u2) at (+1,+1);
	%
	% colouring: 
	\fill [orange!40!white, opacity=0.7] (d1) -- (d2) -- (u2) -- (u1); 
	\draw[thin] (d1) -- (d2) -- (u2) -- (u1) -- (d1);
	%
	% strings: 
	\draw[redirectedgreen, color=green!50!black] (c) .. controls +(0,-0.5) and +(0,-0.5) .. (c');
	\draw[directedgreen, color=green!50!black] (c2) .. controls +(0,0.5) and +(0,0.5) .. (c2');
	\draw[color=green!50!black] (c') -- (c2');
	\draw (m) -- ($(m)+(0,1)$);
	\draw (m2) -- ($(m)+(0,1)$);
	%
	% labels: 
	\fill (m) circle (0pt) node[above right] {{\small $\hspace{-0.2em}m$}};	
	\fill (m2) circle (0pt) node[below right] {{\small $\hspace{-0.2em}m$}};	
	\fill[color=green!50!black] (c2) circle (0pt) node[above right] {{\small $\hspace{-0.35em}c$}};
	%
	%% PHI: 
	\fill[color=white] (phi) node[inner sep=3.9pt,draw, rounded corners=1pt, fill, color=white] (R2) {{\scriptsize$\;\xi\;$}};
	\draw[line width=1pt, color=black] (phi) node[inner sep=4pt,draw, rounded corners=1pt] (R) {{\scriptsize$\;\xi\;$}};	 
	\end{tikzpicture}
	%%%%%%%%%%%%%%%%%%%%%% 
} 
\,
\right) 
\, , \quad 
\operatorname{tr}_{m\triangleleft c'} \left( 
\,
\tikzzbox{%
	%%%%%%%%%%%%%%%%%%%%%% 
	\begin{tikzpicture}[very thick,scale=0.85,color=blue!50!black, baseline=-0.1cm]
	\coordinate (phi) at (0,0);
	\coordinate (c) at (-0.2,-1);
	\coordinate (m) at (+0.2,-1);
	\coordinate (c2) at (-0.2,1);
	\coordinate (m2) at (+0.2,1);
	\coordinate (d1) at (-1,-1);
	\coordinate (d2) at (+1,-1);
	\coordinate (u1) at (-1,+1);
	\coordinate (u2) at (+1,+1);
	%
	% colouring: 
	\fill [orange!40!white, opacity=0.7] (d1) -- (d2) -- (u2) -- (u1); 
	\draw[thin] (d1) -- (d2) -- (u2) -- (u1) -- (d1);
	%
	% strings: 
	\draw[color=blue!50!black] (c) -- ($(c)+(0,1)$);
	\draw[color=green!50!black] (m) -- ($(m)+(0,1)$);
	\draw[color=blue!50!black] (c2) -- ($(c)+(0,1)$);
	\draw[color=green!50!black] (m2) -- ($(m)+(0,1)$);
	%
	% labels: 
	\fill[color=blue!50!black] (c) circle (0pt) node[above left] {{\small $m\hspace{-0.2em}$}};
	\fill[color=green!50!black] (m) circle (0pt) node[above right] {{\small $\hspace{-0.2em}c'$}};	
	\fill[color=blue!50!black] ($(c2)+(0,0.1)$) circle (0pt) node[below left] {{\small $m\hspace{-0.2em}\vphantom{c'}$}};
	\fill[color=green!50!black] ($(m2)+(0,0.1)$) circle (0pt) node[below right] {{\small $\hspace{-0.2em}c'$}};	
	%
	%% PHI: 
	\fill[color=white] (phi) node[inner sep=3.9pt,draw, rounded corners=1pt, fill, color=white] (R2) {{\scriptsize$\;\zeta\;$}};
	\draw[line width=1pt, color=black] (phi) node[inner sep=4pt,draw, rounded corners=1pt] (R) {{\scriptsize$\;\zeta\;$}};	 
	\end{tikzpicture}
	%%%%%%%%%%%%%%%%%%%%%% 
} 
\,
\right) 
= 
\operatorname{tr}_{m} \left( 
\,
\tikzzbox{%
	%%%%%%%%%%%%%%%%%%%%%% 
	\begin{tikzpicture}[very thick,scale=0.85,color=blue!50!black, baseline=-0.1cm, xscale=-1]
	\coordinate (phi) at (0,0);
	\coordinate (c) at (-0.2,-0.2);
	\coordinate (c') at (-0.8,-0.2);
	\coordinate (m) at (+0.2,-1);
	\coordinate (c2) at (-0.2,0.2);
	\coordinate (c2') at (-0.8,0.2);
	\coordinate (m2) at (+0.2,1);
	\coordinate (d1) at (-1,-1);
	\coordinate (d2) at (+1,-1);
	\coordinate (u1) at (-1,+1);
	\coordinate (u2) at (+1,+1);
	%
	% colouring: 
	\fill [orange!40!white, opacity=0.7] (d1) -- (d2) -- (u2) -- (u1); 
	\draw[thin] (d1) -- (d2) -- (u2) -- (u1) -- (d1);
	%
	% strings: 
	\draw[redirectedgreen, color=green!50!black] (c) .. controls +(0,-0.5) and +(0,-0.5) .. (c');
	\draw[directedgreen, color=green!50!black] (c2) .. controls +(0,0.5) and +(0,0.5) .. (c2');
	\draw[color=green!50!black] (c') -- (c2');
	\draw (m) -- ($(m)+(0,1)$);
	\draw (m2) -- ($(m)+(0,1)$);
	%
	% labels: 
	\fill (m) circle (0pt) node[above left] {{\small $m\hspace{-0.2em}$}};	
	\fill (m2) circle (0pt) node[below left] {{\small $m\hspace{-0.2em}$}};	
	\fill[color=green!50!black] (c2) circle (0pt) node[above left] {{\small $c'\hspace{-0.55em}$}};
	%
	%% PHI: 
	\fill[color=white] (phi) node[inner sep=3.9pt,draw, rounded corners=1pt, fill, color=white] (R2) {{\scriptsize$\;\zeta\;$}};
	\draw[line width=1pt, color=black] (phi) node[inner sep=4pt,draw, rounded corners=1pt] (R) {{\scriptsize$\;\zeta\;$}};	 
	\end{tikzpicture}
	%%%%%%%%%%%%%%%%%%%%%% 
} 
\,
\right) 
\ee 
for all $c\in\cat{C}, c'\in\cat{C}'$ and $m\in\mathcal M$. 
Here we use the graphical notation of putting lines next to each other for the left and right action, see~\cite{Bimodtrace} for more details.
\end{definition} 

\begin{proposition}
If $\cat{M}$ is a bimodule category with trace as above, then it gives rise to a 1-morphism in $\orb{E(\Bar\CYCat)}$.  
\end{proposition} 

\begin{proof}
We use Proposition~\ref{prop:1MorphismsOrb} and show the statement only for the left action~$\triangleright$ of~$\cat{C}$ on~$\cat{M}$. 
The proof for $\triangleleft$ is completely analogous. 
Applied to the situation at hand, Proposition~\ref{prop:1MorphismsOrb} implies that we have to show that two isomorphisms between the left and right adjoint of $\triangleright$ agree. 
One of these isomorphisms involves the Calabi--Yau structure on $\cat{C}\boxtimes \cat{M}$, the other one only that on $\cat{C}$. 

We proceed by explicit computation. 
For this we denote by~$I$ a set of representatives for simple objects in~$\cat{C}$. 
The right adjoint is given by 
\begin{align}
\triangleright^R\colon \cat{M} & \to \cat{C}\boxtimes \cat{M} 
\\ 
m & \longmapsto \bigoplus_{i\in I} i\boxtimes i^* \triangleright m 
\end{align}
with unit and counit induced from those for the tensor product $\mu:=\otimes$ of~$\cat{C}$, as explained just above Proposition~\ref{prop:1MorphismsOrb}. 
Concretely we find for the unit 
\begin{align}
\eta^R \colon c \boxtimes m \lmt \bigoplus_{i\in I} i \boxtimes  i^* \triangleright (c \triangleright m) \coloneqq (\id_{\cat{C}}\boxtimes \  \triangleright  ) [\eta_{c,1}\boxtimes 1_m] \, ,
\end{align} 
suppressing coherence isomorphisms related to the unit of~$\cat{C}$.
We compute the unit and counit for the structure of a left adjoint on $\triangleright^R$ using the formulas from Lemma~\ref{lem:FormulaLeftAdjoint}. 
It is enough to check the desired relation for the unit, since this fixes the counit automatically. 
Spelling out the formula from Lemma~\ref{lem:FormulaLeftAdjoint} we find for the component of $\eta^L\colon \id_{\mathcal{M}}\Longrightarrow \triangleright \circ \triangleright^R$ at $m\in \mathcal{M}$ that
\begin{align}
\eta^L_m 
	= 
	\tau_m^{-1}\left[\tr_{\bigoplus_{i\in I} i\boxtimes i^* \triangleright m}\left(( i\boxtimes i^* \triangleright(-))\circ \eta^R_{ \triangleright^R (m)}\right)\right] 
	= 
	\tau_m^{-1}\left[\tr_{\bigoplus_{i\in I} i\otimes i^* \triangleright m}\left(( i\otimes i^* \triangleright(-))\circ \triangleright^R(\eta^R_{ \triangleright^R (m)})\right)\right]
\end{align}
where~$\tau_m$ is the isomorphism introduced in~\eqref{eq:tau}, and the equality uses the assumption on the compatibility of the action with the trace on $\mathcal{M}$. 
Using the graphical calculus the last expressions is given by 
\be 
\tau_m^{-1} \operatorname{tr}_{\bigoplus_i i\otimes i^* \triangleright m} \left( 
\,
\tikzzbox{%
	%%%%%%%%%%%%%%%%%%%%%% 
	\begin{tikzpicture}[very thick,scale=0.85,color=blue!50!black, baseline=-0.1cm]
	\coordinate (eta) at (0,-0.5);
	\coordinate (blank) at (1,1);
	\coordinate (d1) at (-2.2,-2);
	\coordinate (d2) at (+2.2,-2);
	\coordinate (u1) at (-2.2,+2);
	\coordinate (u2) at (+2.2,+2);
	%
	% colouring: 
	\fill [orange!40!white, opacity=0.7] (d1) -- (d2) -- (u2) -- (u1); 
	\draw[thin] (d1) -- (d2) -- (u2) -- (u1) -- (d1);
	%
	% strings: 
	\draw ($(blank)+(0.5,-3)$) -- ($(blank)+(0.5,+1)$);
	\draw[directedgreen] ($(eta)+(0.3,0)$) -- ($(eta)+(0.3,1.5)$);
	\draw[redirectedgreen] ($(eta)+(0.7,0)$) -- ($(eta)+(0.7,1.5)$);
	\draw[directedgreen] ($(eta)+(-0.7,0)$) -- ($(eta)+(-0.7,2.5)$);
	\draw[redirectedgreen] ($(eta)+(-0.3,0)$) -- ($(eta)+(-0.3,2.5)$);
	\draw[directedgreen] ($(eta)+(-0.7,-1.5)$) -- ($(eta)+(-0.7,0)$);
	\draw[redirectedgreen] ($(eta)+(-0.3,-1.5)$) -- ($(eta)+(-0.3,0)$);
	%
	%
	% labels: 
	\fill (1.25,-1.8) circle (0pt) node {{\small $m$}};
	%
	%% ETA: 
	\fill[color=white] (eta) node[inner sep=3.9pt,draw, rounded corners=1pt, fill, color=white] (R2) {{\scriptsize$\;\eta_\mu^R \big|_{i\otimes i^*,1}\;$}};
	\draw[line width=1pt, color=black] (eta) node[inner sep=4pt,draw, rounded corners=1pt] (R) {{\scriptsize$\;\eta_\mu^R \big|_{i\otimes i^*,1}\;$}};	 
	%
	%% BLANK: 
	\fill[color=white] (blank) node[inner sep=3.9pt,draw, rounded corners=1pt, fill, color=white] (R2) {{\scriptsize$\;\eta_\mu^R \big|_{i\otimes i^*,1}\;$}};
	\draw[line width=1pt, color=black] (blank) node[inner sep=4pt,draw, rounded corners=1pt] (R) {{\scriptsize$\;\hspace{1.8em}-\vphantom{\eta_\mu^R \big|_{i\otimes i^*,1_A}}\hspace{1.8em}\;$}};	 
	\end{tikzpicture}
	%%%%%%%%%%%%%%%%%%%%%% 
} 
\,
\right) 
=
\tau_m^{-1} \operatorname{tr}_{m} \left( 
\,
\tikzzbox{%
	%%%%%%%%%%%%%%%%%%%%%% 
	\begin{tikzpicture}[very thick,scale=0.85,color=blue!50!black, baseline=-0.1cm]
	\coordinate (eta) at (0,-0.5);
	\coordinate (blank) at (1,1);
	\coordinate (d1) at (-2.2,-2);
	\coordinate (d2) at (+2.2,-2);
	\coordinate (u1) at (-2.2,+2);
	\coordinate (u2) at (+2.2,+2);
	%
	% colouring: 
	\fill [orange!40!white, opacity=0.7] (d1) -- (d2) -- (u2) -- (u1); 
	\draw[thin] (d1) -- (d2) -- (u2) -- (u1) -- (d1);
	%
	% strings: 
	\draw ($(blank)+(0.5,-3)$) -- ($(blank)+(0.5,+1)$);
	\draw[directedgreen] ($(eta)+(0.3,0)$) -- ($(eta)+(0.3,1.5)$);
	\draw[redirectedgreen] ($(eta)+(0.7,0)$) -- ($(eta)+(0.7,1.5)$);
	\draw[color=green!50!black] ($(eta)+(-0.7,0)$) -- ($(eta)+(-0.7,0.5)$);
	\draw[color=green!50!black] ($(eta)+(-0.3,0)$) -- ($(eta)+(-0.3,0.5)$);
	\draw[directedgreen] ($(eta)+(-0.7,0.5)$).. controls +(0,0.5) and +(0,0.5) .. ($(eta)+(-1.4,0.5)$);
	\draw[redirectedgreen] ($(eta)+(-0.3,0.5)$).. controls +(0,1) and +(0,1) .. ($(eta)+(-1.8,0.5)$);
	\draw[color=green!50!black] ($(eta)+(-0.7,-0.5)$) -- ($(eta)+(-0.7,0)$);
	\draw[color=green!50!black] ($(eta)+(-0.3,-0.5)$) -- ($(eta)+(-0.3,0)$);
	\draw[redirectedgreen] ($(eta)+(-0.7,-0.5)$).. controls +(0,-0.5) and +(0,-0.5) .. ($(eta)+(-1.4,-0.5)$);
	\draw[directedgreen] ($(eta)+(-0.3,-0.5)$).. controls +(0,-1) and +(0,-1) .. ($(eta)+(-1.8,-0.5)$);
	\draw[color=green!50!black] ($(eta)+(-1.8,-0.5)$) -- ($(eta)+(-1.8,+0.5)$);
	\draw[color=green!50!black] ($(eta)+(-1.4,-0.5)$) -- ($(eta)+(-1.4,+0.5)$);
	%
	%
	% labels: 
	\fill (1.25,-1.8) circle (0pt) node {{\small $m$}};
	%
	%% ETA: 
	\fill[color=white] (eta) node[inner sep=3.9pt,draw, rounded corners=1pt, fill, color=white] (R2) {{\scriptsize$\;\hspace{1.62em}\eta_\mu^R\hspace{1.62em}\;$}};
	\draw[line width=1pt, color=black] (eta) node[inner sep=4pt,draw, rounded corners=1pt] (R) {{\scriptsize$\;\hspace{1.62em}\eta_\mu^R\hspace{1.62em}\;$}};	 
	%
	%% BLANK: 
	\fill[color=white] (blank) node[inner sep=3.9pt,draw, rounded corners=1pt, fill, color=white] (R2) {{\scriptsize$\;\eta_\mu^R \big|_{i\otimes i^*,1}\;$}};
	\draw[line width=1pt, color=black] (blank) node[inner sep=4pt,draw, rounded corners=1pt] (R) {{\scriptsize$\;\hspace{1.8em}-\vphantom{\eta_\mu^R \big|_{i\otimes i^*,1_A}}\hspace{1.8em}\;$}};	 
	\end{tikzpicture}
	%%%%%%%%%%%%%%%%%%%%%% 
} 
\,
\right) 
\ee 
from which we can read off 
\be 
\eta^L_m = 
\tikzzbox{%
	%%%%%%%%%%%%%%%%%%%%%% 
	\begin{tikzpicture}[very thick,scale=0.85,color=blue!50!black, baseline=-0.1cm]
	\coordinate (eta) at (0,-0.5);
	\coordinate (blank) at (1,1);
	\coordinate (d1) at (-2.2,-2);
	\coordinate (d2) at (+2.2,-2);
	\coordinate (u1) at (-2.2,+2);
	\coordinate (u2) at (+2.2,+2);
	%
	% colouring: 
	\fill [orange!40!white, opacity=0.7] (d1) -- (d2) -- (u2) -- (u1); 
	\draw[thin] (d1) -- (d2) -- (u2) -- (u1) -- (d1);
	%
	% strings: 
	\draw ($(blank)+(0.5,-3)$) -- ($(blank)+(0.5,+1)$);
	\draw[directedgreen] ($(eta)+(0.3,0)$) -- ($(eta)+(0.3,2.5)$);
	\draw[redirectedgreen] ($(eta)+(0.7,0)$) -- ($(eta)+(0.7,2.5)$);
	\draw[color=green!50!black] ($(eta)+(-0.7,0)$) -- ($(eta)+(-0.7,0.5)$);
	\draw[color=green!50!black] ($(eta)+(-0.3,0)$) -- ($(eta)+(-0.3,0.5)$);
	\draw[directedgreen] ($(eta)+(-0.7,0.5)$).. controls +(0,0.5) and +(0,0.5) .. ($(eta)+(-1.4,0.5)$);
	\draw[redirectedgreen] ($(eta)+(-0.3,0.5)$).. controls +(0,1) and +(0,1) .. ($(eta)+(-1.8,0.5)$);
	\draw[color=green!50!black] ($(eta)+(-0.7,-0.5)$) -- ($(eta)+(-0.7,0)$);
	\draw[color=green!50!black] ($(eta)+(-0.3,-0.5)$) -- ($(eta)+(-0.3,0)$);
	\draw[redirectedgreen] ($(eta)+(-0.7,-0.5)$).. controls +(0,-0.5) and +(0,-0.5) .. ($(eta)+(-1.4,-0.5)$);
	\draw[directedgreen] ($(eta)+(-0.3,-0.5)$).. controls +(0,-1) and +(0,-1) .. ($(eta)+(-1.8,-0.5)$);
	\draw[color=green!50!black] ($(eta)+(-1.8,-0.5)$) -- ($(eta)+(-1.8,+0.5)$);
	\draw[color=green!50!black] ($(eta)+(-1.4,-0.5)$) -- ($(eta)+(-1.4,+0.5)$);
	%
	%
	% labels: 
	\fill (1.25,-1.8) circle (0pt) node {{\small $m$}};
	%
	%% ETA: 
	\fill[color=white] (eta) node[inner sep=3.9pt,draw, rounded corners=1pt, fill, color=white] (R2) {{\scriptsize$\;\eta_\mu^R \big|_{i\otimes i^*,1}\;$}};
	\draw[line width=1pt, color=black] (eta) node[inner sep=4pt,draw, rounded corners=1pt] (R) {{\scriptsize$\;\eta_\mu^R \big|_{i\otimes i^*,1}\;$}};	 
	\end{tikzpicture}
	%%%%%%%%%%%%%%%%%%%%%% 
} 
\, . 
\ee 
We can do exactly the same computation with $\mathcal{M}$ replaced by $\mathcal{C}$ to find a formula for $\eta^L_\mu$. 
This shows that the structure of a left adjoint on $\triangleright^R$ coming from the Calabi--Yau structure on $\mathcal{C}\boxtimes \mathcal{M}$ agrees with the one induced by $\mu^R$. 
This finishes the proof.    
\end{proof}

From the definition it is clear that 2- and 3-morphism in $\orb{E(\Bar\CYCat)}$ are given by bimodule functors
and natural transformations. The discussion in this section can hence be summarised as follows. 

\begin{theorem}
	\label{thm:OrbOfTriv}
The 3-category $\catf{sFus}$ of spherical fusion categories, bimodule categories with trace, bimodule functors, 
and bimodule natural transformations defined in~\cite{Bimodtrace} is a subcategory of $\orb{E(\Bar\CYCat)}$, or equivalently of $\orb{E(\Bar\ssFrob)}$. 
\end{theorem}

\begin{remark}
A consequence of Theorem~\ref{thm:OrbOfTriv} is that every spherical fusion category
%arXiv_v2:
	(of non-zero global dimension, which is guaranteed by our assumption on the ground field)
is 3-dualisable as an object of the symmetric monoidal 3-category $\Alg_1(\Tvs)$. 
Indeed, as already mentioned it is always 1-dualisable. 
The higher dualisability follows from Lemma~\ref{lem:Adj2Mor} and Proposition~\ref{prop:Adj1Mor}. 
This recovers a special case of the main result of~\cite{DuaTen}. 

Spherical fusion categories are expected to also have a canonical $\operatorname{SO}(3)$-homotopy fixed-point structure. 
We hope that the techniques developed in the present paper are helpful in proving this.     
\end{remark}

\subsection{Domain walls between Reshetikhin--Turaev theories}
\label{subsec:DomainWallsForRT}

Any 
%arXiv_v3: 
	%semisimple 
modular fusion category $\cat{M}$ can be used to construct a 3-dimensional TQFT known as the Reshetikhin--Turaev theory based on $\cat{M}$, see~\cite{BookTuraev} for a textbook account. 
In \cite{ks1012.0911, CRS2} it was shown that surface defects within a Reshetikhin--Turaev theory can be constructed from $\Delta$-separable symmetric Frobenius algebras in $\cat{M}$. 
Furthermore, line and point defects can be constructed from bimodules and bimodule maps, respectively. 
This leads us to consider the 3-category with duals $\Bar\Delta\ssFrob(\mathcal{M})$, and in this section we explain how the study of $\orb{E(\Bar\Delta\ssFrob(\mathcal{M}))}$ can be used to reproduce and contextualise the results \cite{KMRS} of Koppen--Mulevi\v{c}ius--Runkel--Schweigert on defect TQFTs with different Reshetikhin--Turaev phases. 

\medskip 

Orbifold data in $\Bar\Delta\ssFrob(\mathcal{M})$ are extensively studied in~\cite{CRS3}. 
One source for orbifold data are commutative $\Delta$-separable Frobenius algebras $(A,\mu, \eta, \Delta, \varepsilon)$ in $\cat{M}$. 
The $A$-$(A\otimes A)$-bimodule which is part of the orbifold datum is~$A$ where the action of $A\otimes A$ uses the multiplication $\mu$. 
This can be understood as follows. 
Consider the 2-category $\Delta\ssFrob^{\operatorname{pt}}(\mathcal{M})$ with $\Delta$-separable symmetric Frobenius algebras in $\cat{M}$ as objects, algebra homomorphisms as 1-morphisms, and only identities as 2-morphisms. 
There is a 2-functor 
%arXiv_v2: 
	%$\Delta\ssFrob^{\operatorname{pt}}(\mathcal{M}) \to \Delta\ssFrob(\mathcal{M}) $ 
	 $\Phi \colon \Delta\ssFrob^{\operatorname{pt}}(\mathcal{M}) \to \Delta\ssFrob(\mathcal{M}) $ 
sending an algebra homomorphism to the bimodule it induces. 
Dunn additivity 
%arXiv_v2: 
	for ordinary $E_1$-algebras\footnote{Concretely, we use the equivalence between the category of algebras in algebras in a symmetric monoidal category $\Ca$ and the category of commutative algebras in $\Ca$.}
implies that
%arXiv_v2: 
	%objects~$\mathcal A$ in $\Delta\ssFrob^{\operatorname{pt}}(\mathcal{M})$ 
	 objects~$A$ in the orbifold completion of the delooping of the image of~$\Phi$ are \textsl{commutative} $\Delta$-separable symmetric Frobenius algebras.\footnote{Note that it does not make sense to look for orbifold data in $\Bar\Delta\ssFrob^{\operatorname{pt}}(\mathcal{M})$ because its 1-morphisms do not admit adjoints.}  
%arXiv_v2: 
	%By the results of \cite{CRS3}, such algebras~$\mathcal A$ always have the structure of orbifold data in $\Bar\Delta\ssFrob(\cat{M})$.\footnote{Note that it does not make sense to look for orbifold data in $\Delta\ssFrob^{\operatorname{pt}}(\mathcal{M})$ because its morphisms do not admit adjoints.} 
	%We can now use a similar approach to construct 1-morphisms in $\orb{(\Bar\Delta\ssFrob(\cat{M}))}$ which lie in the image of $\Delta\ssFrob^{\operatorname{pt}}(\mathcal{M})$. 
	%Concretely, a 1-morphism from~$B$ to~$A$ in $\orb{(\Bar\Delta\ssFrob(\cat{M}))}$ 
	 For such orbifold data~$A$ and~$B$, a 1-morphism from~$B$ to~$A$ in $\orb{(\Bar\Delta\ssFrob(\cat{M}))}$ 
is given by a $\Delta$-separable Frobenius algebra~$F$ in~$\cat{M}$ which at the same time is an $A$-$B$-bimodule where the algebra actions are morphisms of Frobenius algebras. 
The following result explains the relation to~\cite{KMRS}.

\begin{proposition}
Let $A$ and $B$ be objects of $\Delta\ssFrob (\mathcal{M})$ and $F$ an $A$-$B$-bimodule additionally equipped with the structure of a $\Delta$-separable symmetric Frobenius algebra. 
Then the following are equivalent: 
\begin{itemize}
	\item 
	The structure maps encoding the right and left action are morphisms of Frobenius algebras. 
	\item 
	The conditions 
	\begin{align}
	&
	


	\end{align}
of \cite[Def.\,2.10]{KMRS} hold.\footnote{We thank Vincentas Mulevi\v{c}ius for sharing these diagrams from \cite{KMRS}, and for producing the diagrams in \eqref{eq:LeftActionBimodule1}--\eqref{eq:LeftActionBimodule2}.} 
\end{itemize}    	
\end{proposition}  

\begin{proof}
It is enough to show the relation for the multiplication, because those imply the relations for the comultiplication. 
We show the equivalence for the left $A$-action. 
The proof for the right $B$-action is completely analogous. 

The condition that the left action is a bimodule map is 
\begin{equation}
\label{eq:LeftActionBimodule1}

	\begin{tikzpicture}[baseline={([yshift=-.5ex]current bounding box.center)}]
	\node at (0,0) {\def\svgscale{1.0} %% Creator: Inkscape inkscape 0.92.3, www.inkscape.org
%% Accompanies image file 'FA_bimodule_cond_lhs.pdf' (pdf, eps, ps)
\begingroup%
  \makeatletter%
  \providecommand\color[2][]{%
    \errmessage{(Inkscape) Color is used for the text in Inkscape, but the package 'color.sty' is not loaded}%
    \renewcommand\color[2][]{}%
  }%
  \providecommand\transparent[1]{%
    \errmessage{(Inkscape) Transparency is used (non-zero) for the text in Inkscape, but the package 'transparent.sty' is not loaded}%
    \renewcommand\transparent[1]{}%
  }%
  \providecommand\rotatebox[2]{#2}%
  \newcommand*\fsize{\dimexpr\f@size pt\relax}%
  \newcommand*\lineheight[1]{\fontsize{\fsize}{#1\fsize}\selectfont}%
  \ifx\svgwidth\undefined%
    \setlength{\unitlength}{48.52142059bp}%
    \ifx\svgscale\undefined%
      \relax%
    \else%
      \setlength{\unitlength}{\unitlength * \real{\svgscale}}%
    \fi%
  \else%
    \setlength{\unitlength}{\svgwidth}%
  \fi%
  \global\let\svgwidth\undefined%
  \global\let\svgscale\undefined%
  \makeatother%
  \begin{picture}(1,1.25613066)%
    \lineheight{1}%
    \setlength\tabcolsep{0pt}%
    \put(0,0){\includegraphics[width=\unitlength,page=1]{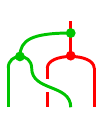}}%
    \put(-0.00313972,0){\color[rgb]{0,0.76862745,0}\makebox(0,0)[lt]{\smash{\begin{tabular}[t]{l}$A$\end{tabular}}}}%
    \put(0.39729423,0){\color[rgb]{1,0,0}\makebox(0,0)[lt]{\smash{\begin{tabular}[t]{l}$F$\end{tabular}}}}%
    \put(0.61514272,0){\color[rgb]{0,0.76862745,0}\makebox(0,0)[lt]{\smash{\begin{tabular}[t]{l}$A$\end{tabular}}}}%
    \put(0.86100695,0){\color[rgb]{1,0,0}\makebox(0,0)[lt]{\smash{\begin{tabular}[t]{l}$F$\end{tabular}}}}%
    \put(0.62915059,1.07632123){\color[rgb]{1,0,0}\makebox(0,0)[lt]{\smash{\begin{tabular}[t]{l}$F$\end{tabular}}}}%
  \end{picture}%
\endgroup%
 };
	\end{tikzpicture}
 
=

	\begin{tikzpicture}[baseline={([yshift=-.5ex]current bounding box.center)}]
	\node at (0,0) {\def\svgscale{1.0} %% Creator: Inkscape inkscape 0.92.3, www.inkscape.org
%% Accompanies image file 'FA_bimodule_cond_rhs.pdf' (pdf, eps, ps)
\begingroup%
  \makeatletter%
  \providecommand\color[2][]{%
    \errmessage{(Inkscape) Color is used for the text in Inkscape, but the package 'color.sty' is not loaded}%
    \renewcommand\color[2][]{}%
  }%
  \providecommand\transparent[1]{%
    \errmessage{(Inkscape) Transparency is used (non-zero) for the text in Inkscape, but the package 'transparent.sty' is not loaded}%
    \renewcommand\transparent[1]{}%
  }%
  \providecommand\rotatebox[2]{#2}%
  \newcommand*\fsize{\dimexpr\f@size pt\relax}%
  \newcommand*\lineheight[1]{\fontsize{\fsize}{#1\fsize}\selectfont}%
  \ifx\svgwidth\undefined%
    \setlength{\unitlength}{41.02147826bp}%
    \ifx\svgscale\undefined%
      \relax%
    \else%
      \setlength{\unitlength}{\unitlength * \real{\svgscale}}%
    \fi%
  \else%
    \setlength{\unitlength}{\svgwidth}%
  \fi%
  \global\let\svgwidth\undefined%
  \global\let\svgscale\undefined%
  \makeatother%
  \begin{picture}(1,1.48578859)%
    \lineheight{1}%
    \setlength\tabcolsep{0pt}%
    \put(0,0){\includegraphics[width=\unitlength,page=1]{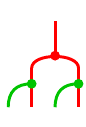}}%
    \put(-0.00371376,0){\color[rgb]{0,0.76862745,0}\makebox(0,0)[lt]{\smash{\begin{tabular}[t]{l}$A$\end{tabular}}}}%
    \put(0.28710174,0){\color[rgb]{1,0,0}\makebox(0,0)[lt]{\smash{\begin{tabular}[t]{l}$F$\end{tabular}}}}%
    \put(0.5447794,0){\color[rgb]{0,0.76862745,0}\makebox(0,0)[lt]{\smash{\begin{tabular}[t]{l}$A$\end{tabular}}}}%
    \put(0.8355949,0){\color[rgb]{1,0,0}\makebox(0,0)[lt]{\smash{\begin{tabular}[t]{l}$F$\end{tabular}}}}%
    \put(0.56134832,1.27310466){\color[rgb]{1,0,0}\makebox(0,0)[lt]{\smash{\begin{tabular}[t]{l}$F$\end{tabular}}}}%
  \end{picture}%
\endgroup%
 };
	\end{tikzpicture}
 \, .
\end{equation}
By inserting the unit of $A$ into one of the green strands we receive the desired relations. The other direction is also straightforward: 
\begin{equation}
\label{eq:LeftActionBimodule2}

	\begin{tikzpicture}[baseline={([yshift=-.5ex]current bounding box.center)}]
	\node at (0,0) {\def\svgscale{1.0}  };
	\end{tikzpicture}
 
=

	\begin{tikzpicture}[baseline={([yshift=-.5ex]current bounding box.center)}]
	\node at (0,0) {\def\svgscale{1.0} %% Creator: Inkscape inkscape 0.92.3, www.inkscape.org
%% Accompanies image file 'FA_bimodule_cond_calc_1.pdf' (pdf, eps, ps)
\begingroup%
  \makeatletter%
  \providecommand\color[2][]{%
    \errmessage{(Inkscape) Color is used for the text in Inkscape, but the package 'color.sty' is not loaded}%
    \renewcommand\color[2][]{}%
  }%
  \providecommand\transparent[1]{%
    \errmessage{(Inkscape) Transparency is used (non-zero) for the text in Inkscape, but the package 'transparent.sty' is not loaded}%
    \renewcommand\transparent[1]{}%
  }%
  \providecommand\rotatebox[2]{#2}%
  \newcommand*\fsize{\dimexpr\f@size pt\relax}%
  \newcommand*\lineheight[1]{\fontsize{\fsize}{#1\fsize}\selectfont}%
  \ifx\svgwidth\undefined%
    \setlength{\unitlength}{48.52142059bp}%
    \ifx\svgscale\undefined%
      \relax%
    \else%
      \setlength{\unitlength}{\unitlength * \real{\svgscale}}%
    \fi%
  \else%
    \setlength{\unitlength}{\svgwidth}%
  \fi%
  \global\let\svgwidth\undefined%
  \global\let\svgscale\undefined%
  \makeatother%
  \begin{picture}(1,1.25613066)%
    \lineheight{1}%
    \setlength\tabcolsep{0pt}%
    \put(0,0){\includegraphics[width=\unitlength,page=1]{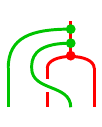}}%
    \put(-0.00313972,0){\color[rgb]{0,0.76862745,0}\makebox(0,0)[lt]{\smash{\begin{tabular}[t]{l}$A$\end{tabular}}}}%
    \put(0.39729423,0){\color[rgb]{1,0,0}\makebox(0,0)[lt]{\smash{\begin{tabular}[t]{l}$F$\end{tabular}}}}%
    \put(0.61514272,0){\color[rgb]{0,0.76862745,0}\makebox(0,0)[lt]{\smash{\begin{tabular}[t]{l}$A$\end{tabular}}}}%
    \put(0.86100695,0){\color[rgb]{1,0,0}\makebox(0,0)[lt]{\smash{\begin{tabular}[t]{l}$F$\end{tabular}}}}%
    \put(0.62915059,1.07632123){\color[rgb]{1,0,0}\makebox(0,0)[lt]{\smash{\begin{tabular}[t]{l}$F$\end{tabular}}}}%
  \end{picture}%
\endgroup%
 };
	\end{tikzpicture}
 
=

	\begin{tikzpicture}[baseline={([yshift=-.5ex]current bounding box.center)}]
	\node at (0,0) {\def\svgscale{1.0} %% Creator: Inkscape inkscape 0.92.3, www.inkscape.org
%% Accompanies image file 'FA_bimodule_cond_calc_2.pdf' (pdf, eps, ps)
\begingroup%
  \makeatletter%
  \providecommand\color[2][]{%
    \errmessage{(Inkscape) Color is used for the text in Inkscape, but the package 'color.sty' is not loaded}%
    \renewcommand\color[2][]{}%
  }%
  \providecommand\transparent[1]{%
    \errmessage{(Inkscape) Transparency is used (non-zero) for the text in Inkscape, but the package 'transparent.sty' is not loaded}%
    \renewcommand\transparent[1]{}%
  }%
  \providecommand\rotatebox[2]{#2}%
  \newcommand*\fsize{\dimexpr\f@size pt\relax}%
  \newcommand*\lineheight[1]{\fontsize{\fsize}{#1\fsize}\selectfont}%
  \ifx\svgwidth\undefined%
    \setlength{\unitlength}{48.52142059bp}%
    \ifx\svgscale\undefined%
      \relax%
    \else%
      \setlength{\unitlength}{\unitlength * \real{\svgscale}}%
    \fi%
  \else%
    \setlength{\unitlength}{\svgwidth}%
  \fi%
  \global\let\svgwidth\undefined%
  \global\let\svgscale\undefined%
  \makeatother%
  \begin{picture}(1,1.25613066)%
    \lineheight{1}%
    \setlength\tabcolsep{0pt}%
    \put(0,0){\includegraphics[width=\unitlength,page=1]{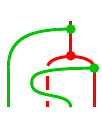}}%
    \put(-0.00313972,0){\color[rgb]{0,0.76862745,0}\makebox(0,0)[lt]{\smash{\begin{tabular}[t]{l}$A$\end{tabular}}}}%
    \put(0.39729423,0){\color[rgb]{1,0,0}\makebox(0,0)[lt]{\smash{\begin{tabular}[t]{l}$F$\end{tabular}}}}%
    \put(0.61514272,0){\color[rgb]{0,0.76862745,0}\makebox(0,0)[lt]{\smash{\begin{tabular}[t]{l}$A$\end{tabular}}}}%
    \put(0.86100695,0){\color[rgb]{1,0,0}\makebox(0,0)[lt]{\smash{\begin{tabular}[t]{l}$F$\end{tabular}}}}%
    \put(0.62915059,1.07632123){\color[rgb]{1,0,0}\makebox(0,0)[lt]{\smash{\begin{tabular}[t]{l}$F$\end{tabular}}}}%
  \end{picture}%
\endgroup%
 };
	\end{tikzpicture}
 
=

	\begin{tikzpicture}[baseline={([yshift=-.5ex]current bounding box.center)}]
	\node at (0,0) {\def\svgscale{1.0}  };
	\end{tikzpicture}

\end{equation}
where we only used the conditions from \cite[Def.\,2.10]{KMRS} as well as the bimodule condition. 
\end{proof}  

A Frobenius-algebra-cum-bimodule~$F$ satisfying either of the equivalent definitions above is called a \emph{$\Delta$-separable symmetric Frobenius algebra over $(A,B)$} in~\cite{KMRS}.      
To check that these define 1-morphisms in the orbifold completion, note that the 3-morphisms we want to check to be equal can be computed by inserting networks of  
$F$, $A$ and $B$. The relations follow exactly as in~\cite{CRS3} from using $\Delta$-separability and the bimodule property.    

To construct 2-morphism in $\orb{(\Bar\Delta\ssFrob)}$ we cannot work within $\Delta\ssFrob^{\operatorname{pt}}(\mathcal{M})$. Instead we just work in $\Bar\Delta\ssFrob$. 
For this let us fix two $\Delta$-separable symmetric Frobenius algebras $F,G$ over a pair of commutative $\Delta$-separable Frobenius algebras $(A,B)$ and an $G$-$F$-bimodule~$M$. 
We want to equip $M$ with the structure of a 2-morphism in $\orb{(\Bar\Delta\ssFrob)}$. 
According to Definition~\ref{def:2MorphismsEAlg}, for this we have to specify a  $G$-$(F\otimes B)$-bimodule isomorphism\footnote{Our choice of notation in this section are chosen so that they agree with the paper~\cite{KMRS}. This leads to some potential confusion when comparing to Definition~\ref{def:2MorphismsEAlg}: $A$ and $B$ are swapped, what is called $F$ here is called $M$ in Definition~\ref{def:2MorphismsEAlg}, and what is called $M$ here is called $F$ in Definition~\ref{def:2MorphismsEAlg}.}
\begin{align}\label{eq: different modules 1}
M \otimes_F F \longrightarrow  G \otimes_{G\otimes B } (M \otimes B )
\end{align} 
(corresponding to the third diagram in~\eqref{eq:2MorphismStructureMaps}), and a $G$-$(A \otimes F)$-bimodule isomorphism
\begin{align}\label{eq: different modules 2}
	M \otimes_F F \longrightarrow  G \otimes_{A\otimes G } (A \otimes M )
\end{align}
(corresponding to the second diagram in~\eqref{eq:2MorphismStructureMaps}), where the commutativity of $A$ is used to turn the right action on $G$ into a left action.
All bimodules appearing here are canonically isomorphic to $M$ as a vector space. 
The isomorphism from~$M$ to the corresponding bimodule inserts units.  
However there is no reason for the $G$-$(A\otimes F)$- or $G$-$(F\otimes B)$-bimodule structures to agree. 
In case they agree the identity will automatically give the bimodule  
isomorphisms we are looking for. 
The left actions in \eqref{eq: different modules 1} are the same under the identification with~$M$. 
For the right actions to agree we find the condition
\begin{align}
& \big( M\otimes B \xrightarrow{u_{F}} M\otimes F\otimes B \xrightarrow{\triangleleft_{F}} M\otimes F \xrightarrow{\triangleleft_{M}} M 
\big) 
\nonumber 
\\ 
& 
\qquad 
= 
\big( M\otimes B \xrightarrow{u_{G}} G \otimes M\otimes B  \xrightarrow{c_{M,B}} G\otimes B \otimes M \xrightarrow{\triangleleft_{G}} G \otimes M \xrightarrow{\triangleright_{M}} M \big) \, .
\end{align}
For the bimodule structures appearing in~\eqref{eq: different modules 2} to be equal we find the condition 
\begin{align}
& \big( M\otimes A \xrightarrow{u_{F}} M\otimes A \otimes F \xrightarrow{\triangleright_{F}} M\otimes F M \big) 
\nonumber 
\\ 
& \qquad 
= \big( M \xrightarrow{\triangleleft_{M}} \otimes A \xrightarrow{u_{G}} G \otimes M\otimes  A  \xrightarrow{c^{-1}_{G\otimes M,A}} A \otimes G \otimes M \xrightarrow{\triangleright_{G}} 
G \otimes M \xrightarrow{\triangleright_{M}} M \big) \,.
\end{align} 
A bimodule satisfying the conditions we found here is equivalent to what is called a \emph{$G$-$F$-bimodule over $(A,B)$} in~\cite[Def.\,2.12]{KMRS}. 

Finally 3-morphisms in $\orb{(\Bar\Delta\ssFrob)}$ are given by $G$-$D$-bimodule maps. 
Hence in summary we have shown: 

\begin{theorem}\label{thm:OrbOfssFrobM}
Let~$\mathcal M$ be a modular fusion category. 
The 3-category in which
\begin{itemize}
	\item 
	objects are commutative $\Delta$-separable Frobenius algebras in $\cat{M}$, 
	\item 
	1-morphisms from $B$ to $A$ are $\Delta$-separable symmetric Frobenius algebras $F$ over $(A,B)$, 
	\item 
	2-morphisms from $F$ to $G$ are $G$-$F$-bimodules $M$ over $(A,B)$, and
	\item 
	3-morphisms are bimodule maps 
\end{itemize} 
is a subcategory of $\orb{(\Bar\Delta\ssFrob(\mathcal{M}))}$.
\end{theorem} 

This gives a more systematic perspective on the defect field theory constructed in~\cite{KMRS} which should agree with the one constructed from the orbifold completion as we describe in Section~\ref{sec:TQFTs}.

\section{Orbifold construction of defect TQFTs}
\label{sec:TQFTs}

In this section the 3-dimensional orbifold completion is used to construct new defect TQFTs from old ones. 
Most of the techniques used have become fairly standard, hence our discussion is brief.  
After a short reminder on 3-dimensional defect TQFTs and their associated 3-categories in Section~\ref{subsec:DefectTQFTs}, in Section~\ref{subsec:DefectTQFTfromOrbifoldCompletion} we apply orbifold completion to the latter and recover defect state sum models and defects in Reshetikhin--Turaev theory as special cases.

\subsection{Defect TQFTs}
\label{subsec:DefectTQFTs}

Topological quantum field theories can be axiomatised as symmetric monoidal functors $\mathcal{Z} \colon \catf{Bord}_{n,n-1}\to \vs$ from the $n$-dimensional bordism category to vector spaces
%arXiv_v2: 
	over an algebraically closed field of characteristic~0, 
see e.\,g.\ \cite{AtiyahTQFT, Kockbook}. 
The objects of $\catf{Bord}_{n,n-1}$ are closed oriented $(n-1)$-dimensional manifolds. 
Morphism are 
%arXiv_v3: 
	%diffeomorphisms 
	 diffeomorphism 
classes of bordisms, i.\,e.\ a morphism $\Sigma\to \Sigma'$ is an equivalence class of compact oriented manifolds $Y$ with boundary together with an orientation-preserving diffeomorphism $\overline{\Sigma}\sqcup \Sigma' \to \partial Y$ (here $\overline{\Sigma}$ denotes the orientation reversal of $\Sigma$). 
Two such manifolds belong to the same equivalence class if there exists an orientation-preserving diffeomorphism between them which is compatible with the boundary parametrisation. 
Composition is defined by gluing manifolds along boundaries, illustrated for $n=2$ by
\be 
\label{eq:VanillaBord}
\vcenter{\hbox{
\begin{overpic}[scale=0.35]{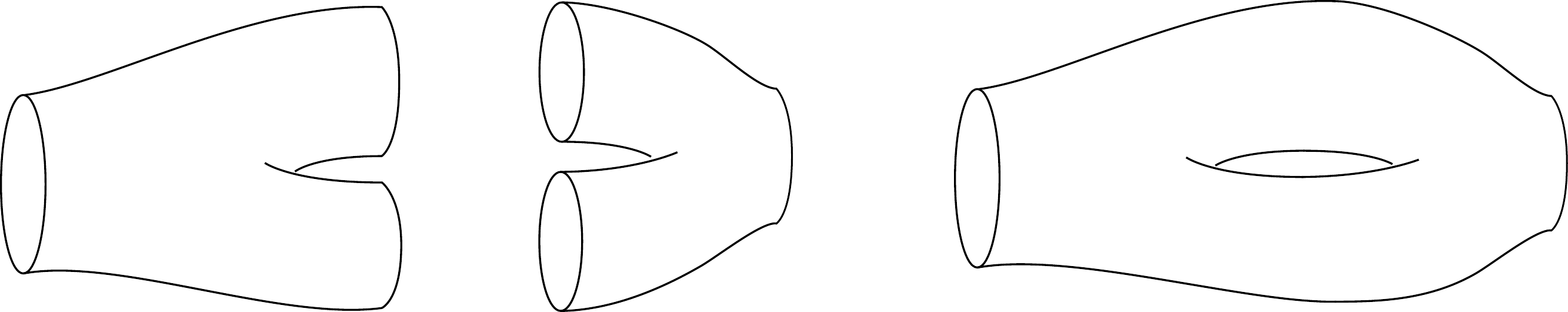}
\put(29.5,8.5){$\circ$} 
\put(55,8.5){$=$}
\end{overpic}}}
\ee 
The monoidal structure of $\catf{Bord}_{n,n-1}$ is given by taking disjoint unions.  

Since their introduction various generalisations of the above \textsl{closed} oriented TQFTs have been considered, among them variants with different tangential structures, different target categories, and extended TQFTs, see e.\,g.\ \cite{Lurie}.
To describe ``defects'' in such theories an extension of the framework involving stratified manifolds is used. 
In comparison with~\eqref{eq:VanillaBord}, the main idea of this refinement is captured in the example
\be 
\label{eq:DefectBord}
\vcenter{\hbox{
\begin{overpic}[scale=0.35]{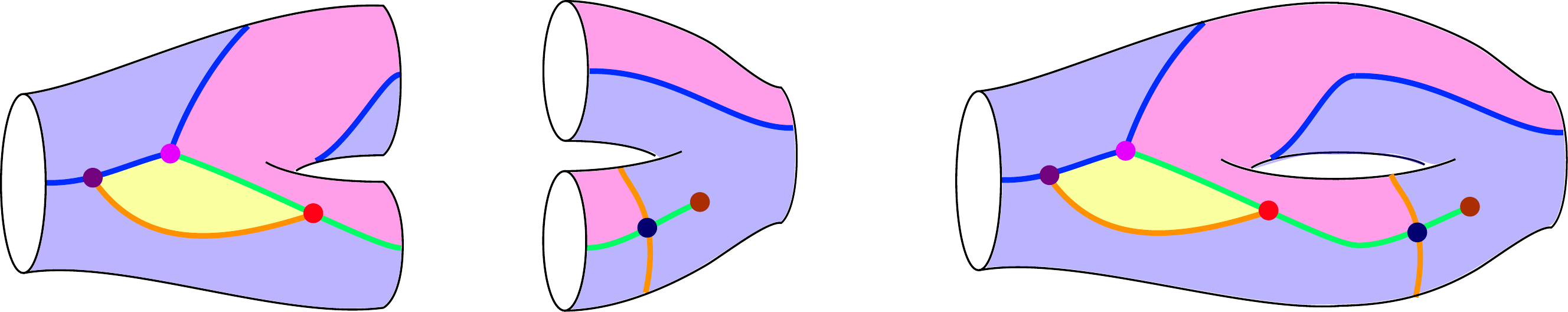}
	\put(29.5,8.5){$\circ$} 
	\put(55,8.5){$=$}
	\put(19,1){$u_1$} 
	\put(5,12){$u_1$}
	\put(23.5,10.5){$u_1$}
	\put(11,7){$u_2$} 
	\put(17,12){$u_3$}
	\put(42,16.5){$u_3$}
	\put(38,6.5){$u_3$}
	\put(38,1){$u_1$}
	\put(48,8){$u_1$} 
	\put(82,1){$u_1$} 
	\put(65,12){$u_1$}
	\put(97,8){$u_1$}
	\put(72,7){$u_2$} 
	\put(85,17.5){$u_3$}
\end{overpic}}}
\ee 
where all the coloured submanifolds come with prescribed orientations (which we suppress in the picture; orientations of top-dimensional strata by definition are induced by the orientation of the underlying manifold) and labels. 
We only display the labels of 2-strata, but also the 1- and 0-strata come with prescribed labels, which we suppress to avoid clutter.  

To make this picture precise we recall from \cite{CMS, CRS1} that by a \emph{3-dimensional oriented stratified manifold (without boundary)} we mean an oriented 3-manifold~$Y$ together with a filtration $Y=F_3\supset F_2 \supset F_1 \supset F_0 \supset F_{-1} = \varnothing$ such that 
\begin{itemize}
	\item 
	$F_{i}\setminus F_{i-1}$ is an $i$-dimensional submanifold of~$Y$ which is furthermore equipped with an orientation. 
	The connected components of $F_{i}\setminus F_{i-1}$ are called \emph{$i$-strata}, and the orientations of 3-strata are induced by the given orientation of~$Y$. 
	\item 
	There are only finitely many $i$-strata for every $i\in\{0,1,2,3\}$. 
	\item 
	If a stratum~$\sigma$ intersects non-trivially with the closure of another stratum~$\sigma'$ (i.\,e.\ $\sigma \cap \overline{\sigma'}\neq \varnothing$), then $\sigma\subset\overline{\sigma'}$.       
\end{itemize}
The definition can be extended to manifolds with boundary which in particular entails a stratification of the boundary. 
For many purposes (including our applications to defect TQFTs) the topology of stratified manifolds defined in this way is too wild. 
Conically stratified manifolds solve this by requiring that locally the stratified manifold is diffeomorphic to a (cylinder over a stratified $(n-1)$-ball, or a) cone over an $(n-1)$-dimensional stratified sphere. 
The type of allowed spheres is defined inductively; illustrative examples of such cylinders and cones for $n=2$ and $n=3$ are  
\be 
\label{eq:NeighbourhoodExamples}
%%%%%%%%%%%%%%%%%%%%%% 
\tikzzbox{\begin{tikzpicture}[thick,scale=2.5,color=black, baseline=1.2cm, >=stealth]
	% 2-strata:
	\fill [blue!60,opacity=0.4] (0,0) -- (0,1) -- (1,1) -- (1,0);
	% 1-strata_
	\draw[
	color=red!90!black, 
	very thick,
	>=stealth, 
	postaction={decorate}, decoration={markings,mark=at position .5 with {\arrow[draw]{>}}}
	] 
	(0.5,0) --  (0.5,1);
	\end{tikzpicture}}%%popende
%%%%%%%%%%%%%%%%%%%%%% 
\quad
%%%%%%%%%%%%%%%%%%%%%%
\tikzzbox{\begin{tikzpicture}[very thick,scale=0.7,color=blue!50!black, baseline=-0.1cm]
	\fill [blue!60,opacity=0.4] (0,0) circle (1.75);
	%%1-strata:  
	\draw[
	color=red!70!black, 
	>=stealth,
	decoration={markings, mark=at position 0.5 with {\arrow{>}},
	}, postaction={decorate}
	] 
	(0,0) -- (0:1.75);
	\draw[
	color=red!70!black, 
	>=stealth,
	decoration={markings, mark=at position 0.5 with {\arrow{>}},
	}, postaction={decorate} 
	] 
	(0,0) -- (120:1.75);
	\draw[
	color=red!70!black, 
	>=stealth,
	decoration={markings, mark=at position 0.5 with {\arrow{<}},
	}, postaction={decorate}
	] 
	(0,0) -- (240:1.75); 
	\fill[color=green!50!black] (0,0) circle (3.5pt) node[above] {};
	\end{tikzpicture}}%%popende%
%%%%%%%%%%%%%%%%%%%%%% 
\qquad \textrm{and} \qquad 
%%%%%%%%%%%%%%%%%%%%%% 
\tikzzbox{\begin{tikzpicture}[very thick,scale=1.0,color=blue!50!black, baseline=-1.9cm]
	%cylinder:
	\fill [blue!15,
	opacity=0.5, 
	left color=blue!15, 
	right color=white] 
	(-1.25,0) -- (-1.25,-3) arc (180:360:1.25 and 0.5) -- (1.25,0) arc (0:180:1.25 and -0.5);
	\fill [blue!35,opacity=0.1] (-1.25,-3) arc (180:360:1.25 and 0.5) -- (1.25,-3) arc (0:180:1.25 and 0.5);
	\fill [blue!25,opacity=0.5] (-1.25,0) arc (180:360:1.25 and 0.5) -- (1.25,0) arc (0:180:1.25 and 0.5);
	%
	%2-strata: 
	\fill [red,opacity=0.4] (0,0) -- (0,-3) -- (1.25,-3) -- (1.25,0);
	\fill [opacity=0.3] (0,0) -- ($(0,0)+(120:1.25 and 0.5)$) -- ($(0,-3)+(120:1.25 and 0.5)$) -- (0,-3) -- (0,0);
	\fill [red,opacity=0.4] (0,0) -- ($(0,0)+(120:1.25 and 0.5)$) -- ($(0,-3)+(120:1.25 and 0.5)$) -- (0,-3) -- (0,0);
	\fill [red,opacity=0.4] (0,0) -- ($(0,0)+(245:1.25 and 0.5)$) -- ($(0,-3)+(245:1.25 and 0.5)$) -- (0,-3) -- (0,0);
	\draw[
	color=green!50!black, 
	>=stealth,
	decoration={markings, mark=at position 0.5 with {\arrow{>}},
	}, postaction={decorate}
	] 
	(0,-3) -- (0,0);
	\end{tikzpicture}}%%popende
%%%%%%%%%%%%%%%%%%%%%% 
\quad
%%%%%%%%%%%%%%%%%%%%%%
\tikzzbox{\begin{tikzpicture}[very thick,scale=1.2,color=green!60!black=-0.1cm, >=stealth, baseline=0]
	\fill[ball color=white!95!blue] (0,0) circle (0.95 cm);
	\coordinate (v1) at (-0.4,-0.6);
	\coordinate (v2) at (0.4,-0.6);
	\coordinate (v3) at (0.4,0.6);
	\coordinate (v4) at (-0.4,0.6);
	%
	% 2-strata: 
	\fill [red!80!black,opacity=0.3] (0,0) -- (v2) .. controls +(0,-0.25) and +(0,-0.25) .. (v1);
	\fill [red!80!black,opacity=0.3] (0,0) -- (v4) .. controls +(0,0.15) and +(0,0.15) .. (v3);
	\fill [red!80!black,opacity=0.3] (0,0) -- (v4) .. controls +(0.25,-0.1) and +(-0.05,0.5) .. (v2);
	\fill [red!80!black,opacity=0.3] (0,0) -- (v3) .. controls +(-0.9,0.99) and +(-0.75,0.4) .. (v1);
	\fill [red!80!black,opacity=0.3] (0,0) -- (v1) .. controls +(-0.15,0.5) and +(-0.15,-0.5) .. (v4);
	\fill [red!80!black,opacity=0.3] (0,0) -- (v3) .. controls +(0.25,-0.5) and +(0.25,0.5) .. (v2);
	%
	% 1-strata: 
	\draw[thick, opacity=0.6, postaction={decorate}, decoration={markings,mark=at position .6 with {\arrow[draw=green!60!black]{<}}}] (0,0) -- (v1);
	\draw[thick, opacity=0.6, postaction={decorate}, decoration={markings,mark=at position .6 with {\arrow[draw=green!60!black]{<}}}] (0,0) -- (v2);
	\draw[thick, opacity=0.6, postaction={decorate}, decoration={markings,mark=at position .7 with {\arrow[draw=green!60!black]{>}}}] (0,0) -- (v3);
	\draw[thick, opacity=0.6, postaction={decorate}, decoration={markings,mark=at position .5 with {\arrow[draw=green!60!black]{>}}}] (0,0) -- (v4);
	\draw[color=red!80!black, very thin, rounded corners=0.5mm, postaction={decorate}, decoration={markings,mark=at position .5 with {\arrow[draw=red!80!black]{>}}}] 
	(v2) .. controls +(0,-0.25) and +(0,-0.25) .. (v1);
	\draw[color=red!80!black, very thin, rounded corners=0.5mm, postaction={decorate}, decoration={markings,mark=at position .62 with {\arrow[draw=red!80!black]{>}}}] 
	(v4) .. controls +(0,0.15) and +(0,0.15) .. (v3);
	\draw[color=red!80!black, very thin, rounded corners=0.5mm, postaction={decorate}, decoration={markings,mark=at position .5 with {\arrow[draw=red!80!black]{>}}}] 
	(v4) .. controls +(0.25,-0.1) and +(-0.05,0.5) .. (v2);
	\draw[color=red!80!black, very thin, rounded corners=0.5mm, postaction={decorate}, decoration={markings,mark=at position .58 with {\arrow[draw=red!80!black]{>}}}] 
	(v3) .. controls +(-0.9,0.99) and +(-0.75,0.4) .. (v1);
	\draw[color=red!80!black, very thin, rounded corners=0.5mm, postaction={decorate}, decoration={markings,mark=at position .5 with {\arrow[draw=red!80!black]{>}}}] 
	(v1) .. controls +(-0.15,0.5) and +(-0.15,-0.5) .. (v4);
	\draw[color=red!80!black, very thin, rounded corners=0.5mm, postaction={decorate}, decoration={markings,mark=at position .5 with {\arrow[draw=red!80!black]{>}}}] 
	(v3) .. controls +(0.25,-0.5) and +(0.25,0.5) .. (v2);
	\fill[magenta!10!black] (0,0) circle (1.6pt) node[black, opacity=0.6, right, font=\tiny, left] {};
	\end{tikzpicture}}%%popende 
%%%%%%%%%%%%%%%%%%%%%%
\, . 
\ee 
We refer to~\cite[Sect.\,2]{CRS1} for a detailed exposition. 

In the setting of defect TQFTs we additionally assign labels to all strata. 
The 3-strata are labelled by elements of a prescribed set~$D_3$ (that is part of the data of a defect TQFT, see below) which we interpret as ``bulk theories''. 
Similarly 2-, 1- and 0-strata are labelled by elements of prescribed sets $D_2, D_1$ and~$D_0$ of ``surface defects'', ``line defects'' and ``point defects'', respectively. 
A $j$-stratum labelled by an element in~$D_j$ is a $j$-dimensional \textsl{defect}. 
The label sets~$D_j$ are completed to a \textsl{set of defect data}~$\mathds{D}$ by supplying adjacency data, i.\,e.\ prescribed rules which $D_2$-labels are allowed for the two 2-strata adjacent to a 3-stratum labelled by a given element of~$D_3$, and similarly how defects are allowed to meet in lines and points. 
We refer to \cite[Sect.\,2.3]{CMS} and \cite[Sect.\,2.2.2]{CRS1} for a formalisation of the notion of~$\mathds{D}$, as well as the precise definition of the \textsl{symmetric monoidal category of defect bordisms} $\DBord(\mathds{D})$. 
Its morphisms are basically equivalence classes of $\mathds{D}$-labelled stratified 3-bordisms, objects are stratified closed 2-manifolds whose $j$-strata are compatibly labelled with elements of~$D_{j+1}$, and glueing is a 3-dimensional version of the case illustrated in~\eqref{eq:DefectBord}. 

\begin{definition}
Let $\mathds{D}$ be a chosen set of 3-dimensional defect data, and let~$\mathcal C$ be a symmetric monoidal category. 
A \emph{3-dimensional defect TQFT} for~$\mathds{D}$ with values in~$\mathcal C$ is a symmetric monoidal functor  
\begin{align}
\cat{Z}\colon \DBord(\mathds{D}) \to \mathcal C \, .
\end{align}
\end{definition}

\begin{remark}
	\label{rem:Defect3category}
	\begin{enumerate}[label={(\roman*)}]
		\item 
		Usually one takes the target~$\mathcal C$ to be the category of (super) vector spaces. 
		To any such defect TQFT~$\zz$ one naturally associates a Gray category with duals~$\tric_{\zz}$, which we think of as an important invariant of~$\zz$. 
		This is the ``categorification'' of the analogous construction of pivotal 2-categories in one dimension lower, first explained in \cite{dkr1107.0495}. 
		Objects of~$\tric_{\zz}$ are elements in~$D_3$, 1- and 2-morphisms are (lists of) elements in~$D_2$ and~$D_1$, respectively, and the vector spaces of 3-morphisms are obtained by applying~$\zz$ to defect spheres (which are objects in $\DBord(\mathds{D})$) as illustrated by the boundary of the defect ball in~\eqref{eq:NeighbourhoodExamples}. 
		Hence $i$-cells in~$\tric_{\zz}$ correspond to $(3-i)$-dimensional defects. 
		For more details we refer to \cite[Sect.\,3]{CMS}. 
		\item 
		\label{item:DefectDataFrom3category}
		Conversely, given a Gray category with duals~$\tric$, one immediately obtains a set of defect data~$\mathds{D}^\tric$. 
		Its sets~$D_j^\tric$ are made up of the $(3-j)$-cells of~$\tric$, and its adjacency data are supplied by the source and target maps of~$\tric$ as well as its composition rules. 
	\end{enumerate}
\end{remark}

\subsection{Defect TQFTs associated to orbifold completion}
\label{subsec:DefectTQFTfromOrbifoldCompletion}

Let $\mathcal{Z}\colon \DBord (\mathds{D}) \to \vs$ be a defect TQFT and $\cat{T}_{\mathcal{Z}}$ the associated 3-category. 
The goal of this section is to construct a new defect TQFT 
\be 
\label{eq:Zorb}
\mathcal{Z}_{\textrm{orb}} \colon \DBord (\mathds{D}_{\textrm{orb}}) \to \vs
\ee  
where $\mathds{D}_{\operatorname{orb}} := \mathds{D}^{(\tric_\zz)_{\textrm{orb}}}$ are the defect data obtained from the orbifold completion of~$\tric_\zz$. 
The construction uses triangulations of stratified manifolds adapted to the stratification.
They are 
%arXiv_v3:  
	%a 
natural extensions of the transversal triangulations considered in~\cite{Meusburger3dDefectStateSumModels} to the slightly more general type of stratified manifolds considered in the present paper (where we allow more than two 2-strata to meet along a 1-stratum). 
In the following we use ``triangulations with total order''~$t$ (on the 0-simplices) as discussed e.\,g.\ in \cite[Sect.\,3.1]{CRS1}, to merge the original stratification of a bordism with the Poincar\'e dual of~$t$: 

\begin{definition}
Let~$Y$ be a defect bordism in $\DBord (\mathds{D})$.
A \emph{stratified triangulation of~$Y$} is a triangulation with total order~$t$ of the underlying manifold of~$Y$ such that: 
\begin{enumerate}[label={(\roman*)}]
	\item 
	\label{item:TransversalSimplex}
	All intersections of $i$-simplices of~$t$ with $j$-strata of~$Y$ are transversal\footnote{Transversality means that the fibres of the tangent spaces of the $i$-simplex and the $j$-stratum together span the fibres of $TY$.} in~$Y$.  
	\item 
	\label{item:NewStrata}
	We construct a proper substratification~$S^t_Y$ of~$Y$ from a stratification which is Poincar\'e dual to~$t$ as follows, where we demand that every $j$-stratum~$\sigma$ of~$Y$ is properly substratified in~$S_Y^t$, and that in every $(j+1)$-stratum in~$Y$ adjacent to~$\sigma$ there are (new) $j$-strata in $S_Y^t \setminus Y$ which meet~$\sigma$ (where here by ``$\tau \in S_Y^t \setminus Y$'' we mean ``$\tau \in S_Y^t$ is a $j$-stratum that is not contained in an (original) $j$-stratum in~$Y$''): 
	\begin{itemize}
		\item 
		For every 3-simplex $\delta_3\in t$, there is a 0-stratum $\sigma_0 \in S^t_Y \setminus Y$ in the interior of~$\delta_3$, unless every 2-face of~$\delta_3$ intersects with an original 2-stratum of~$Y$. 
		If there are original 2-strata~$\sigma_2$ of~$Y$ which intersect~$\delta_3$, then~$\sigma_0$ is placed in 
		%arXiv_v3:  
			%one such~$\delta_2$ 
			 one such~$\sigma_2$ 
		(becoming part of a substratification of~$\sigma_2$); otherwise~$\sigma_0$ is likewise placed in an original 3-stratum of~$Y$. 
		\item 
		For every 2-simplex $\delta_2\in t$, there is a 1-stratum $\sigma_1 \in S^t_Y$. 
		If there are original 1-strata in~$Y$ which intersect~$\delta_2$, $\sigma_1$ is taken to be (part of a substratification of) one of them; otherwise $\sigma_1 \in S^t_Y\setminus Y$ is a new 1-stratum which is Poincar\'e dual to~$\delta_2$. 
		If in the former case~$\sigma_1$ intersects with a 1-stratum~$\ell$ of~$Y$, there is another new 0-stratum at $\ell\cap\sigma_1$, substratifying both 1-strata. 
		\item 
		For every 1-simplex $\delta_1\in t$, there is a 2-stratum $\sigma_2\in S^t_Y$. 
		If there are original 2-strata in~$Y$ which intersect~$\delta_1$, $\sigma_2$ is taken to be (a substratification of) one of them; otherwise $\sigma_2 \in S^t_Y\setminus Y$ is a new 2-stratum which is Poincar\'e dual to~$\delta_1$. 
	\end{itemize}
	\item 
	The new strata in~$S_Y^t$ which meet the boundary of~$Y$ do so transversally. 
	\item 
	The new $j$-strata $\widetilde{\sigma} \in S_Y^t \setminus Y$ carry orientations such that the orientation of~$\widetilde{\sigma}$ together with the order-induced orientation of the Poincar\'e dual $(3-j)$-simplex in~$t$ (in that order) give the orientation of~$Y$. 
	The new sub-$j$-strata of $j$-strata of~$Y$ carry the induced orientation. 
\end{enumerate}
\label{def:StratifiedTriangulation}
\end{definition}
In the special case that~$Y$ is trivially stratified we have that the stratification of~$S^t_Y$ is precisely Poincar\'e dual to~$t$. 
For non-trivial original stratifications, Definition~\ref{def:StratifiedTriangulation}\,\ref{item:NewStrata} involves choices about how new Poincar\'e dual strata meet with old strata. 
In the construction of the TQFT $\mathcal{Z}_{\operatorname{orb}}$ below these choices will be immaterial by design of the orbifold completion $\orb{(\tric_{\zz})}$. 

Any sufficiently fine triangulation of~$Y$ can be deformed into a stratified triangulation by slightly moving its cells. 
Furthermore, any two stratified triangulations of a closed stratified manifold can be related by a sequence of ordinary Pachner moves between stratified triangulations, which can be seen by choosing a stratified triangulation of $Y \times [0,1]$ that restricts to the stratified triangulations on the boundary. 

Using the notion of stratified triangulation we can now describe the new defect TQFT $\mathcal{Z}_{\textrm{orb}}$ in~\eqref{eq:Zorb} in terms of the set of defect data
\be 
\mathds{D}_{\operatorname{orb}} := \mathds{D}^{(\tric_\zz)_{\textrm{orb}}}
\ee 
associated to $(\tric_\zz)_{\textrm{orb}}$ by the procedure outlined in Remark~\ref{rem:Defect3category}\,\ref{item:DefectDataFrom3category}. 
On a closed defect manifold $Y\colon \varnothing \lra \varnothing$ in $\DBord (\mathds{D}_{\textrm{orb}})$, we extend the construction of \cite[Sect.\,3.2]{CRS1} as follows to define the evaluation of the TQFT $\mathcal{Z}_{\textrm{orb}}$: 

\begin{figure}[hbt]
	\begin{center}
		\begin{overpic}[scale=0.7]{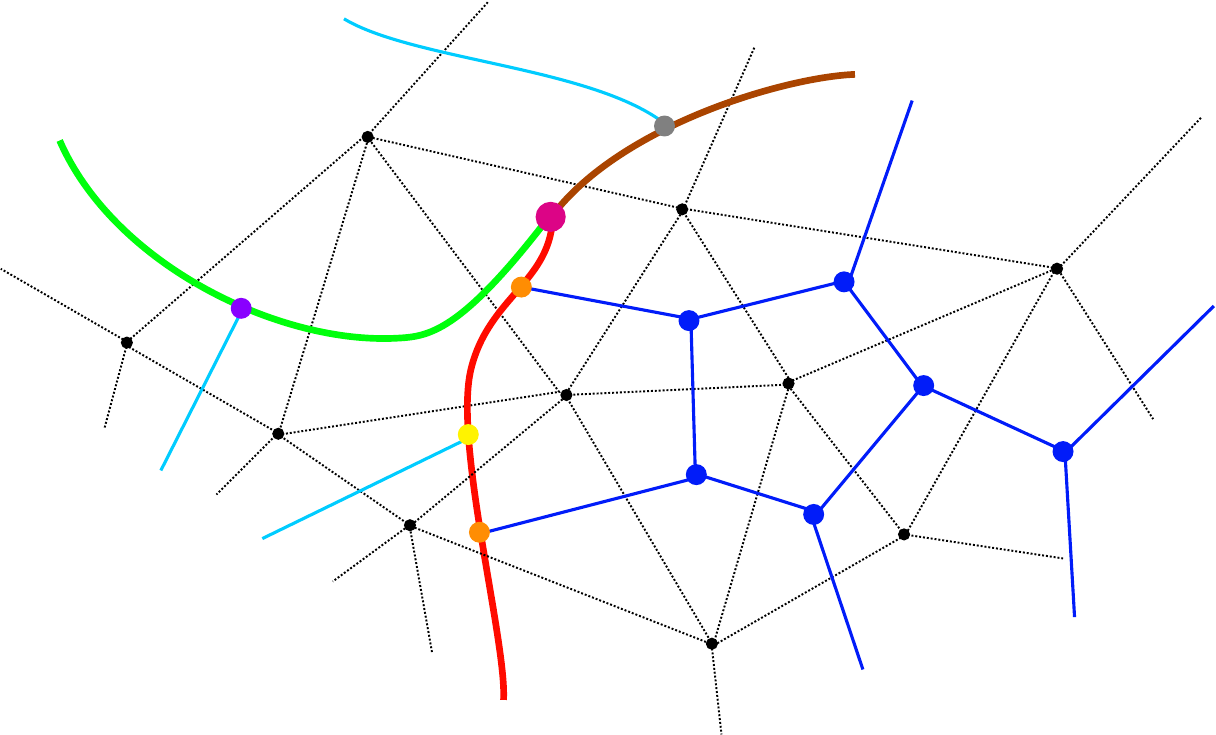} 
			\put(49,20){\color{blue} $A$}
			\put(50,36){\color{blue} $A$}
			\put(57,28){\color{blue} $A$}
			\put(62,17){\color{blue} $A$}
			\put(63,37){\color{blue} $A$}
			\put(69.5,12){\color{blue} $A$}
			\put(72,22){\color{blue} $A$}
			\put(70,45){\color{blue} $A$}
			\put(73,33){\color{blue} $A$}
			\put(80,24){\color{blue} $A$}
			\put(89,14){\color{blue} $A$}
			\put(94.5,32){\color{blue} $A$}
			\put(56.5,18.5){\color{blue} $\mu^\vee$ }
			\put(56.5,36){\color{blue} $\mu$ }
			\put(66,19.5){\color{blue} $\mu$ }
			\put(72,28){\color{blue} $\mu^\vee$ }
			\put(67.5,34){\color{blue} $\mu^\vee$ }
			\put(89,23){\color{blue} $\mu$ }
			\put(42,6){\color{red} $M$ }
			\put(40,20){\color{red} $M$ }
			\put(40,29.5){\color{red} $M$ }
			\put(45,38){\color{red} $M$ }
			\put(10,42){\color{green} $M'$ }
			\put(28,33.5){\color{green} $M'$ }
			\put(48,48){\color{brown} $M''$ }
			\put(62,54){\color{brown} $M''$ }
			\put(28,21){\color{cyan} $A'$}
			\put(12,25){\color{cyan} $A'$}
			\put(40,56){\color{cyan} $A''$}
			\put(40.5,14.5){\color{orange}$\triangleleft_M$}
			\put(43.5,34.5){\color{orange}$\triangleleft_M$} 
			\put(34.5,24.5){\color{yellow}$\triangleright_M$}
			\put(19.5,36.5){\color{violet}$\triangleleft_{M'}$}
			\put(54,52){\color{gray}$\triangleright_{M''}$}
			\put(44,44){\color{purple}$F$}
		\end{overpic}
	\end{center}
	\caption{A sketch of a 2-dimensional slice for the construction of a local patch of a defect bordism~$Y^t$ in $\DBord(\mathds{D})$ from a bordism~$Y$ in $\DBord (\mathds{D}_{\textrm{orb}})$. 
	The original stratification consists of the fat lines (corresponding to 2-strata in~$Y$) labelled by the ingredients $M,M',M''$ of $\mathcal M,\mathcal M',\mathcal M'' \in (\mathds{D}_{\operatorname{orb}})_2$, as well as the $F$-labelled point (corresponding to a 1-stratum in~$Y$). 
	The stratified triangulation is drawn in black. 
	The strata of the original stratification of~$Y$ are substratified in~$Y^t$; for example, the $M$-labelled line is substratified into four lines and three points. 
	The light and dark blue lines are labelled with the 1-morphisms $A,A',A''$ corresponding to the orbifold data $\mathcal A, \mathcal A', \mathcal A''$ assigned in~$Y$ the top-dimensional strata. 
	The blue vertices are labelled by the corresponding multiplication 2-morphisms $\mu,\mu',\dots$ or their adjoints, as dictated by the chosen total order of 0-simplices. 
	In addition, the new vertices (which are part of substratifications of lines) are labelled with the 2-morphisms~$\triangleright$ and~$\triangleleft$ which are part of the 1- morphisms $\mathcal M, \mathcal M', \mathcal M''$ in $(\tric_\zz)_{\textrm{orb}}$. 
	Note that we suppress orientations to avoid clutter.}
	\normalsize
%	\newline 
	\label{Fig: stratified triangulation}
\end{figure}

\begin{enumerate}
	\item 
	Pick a stratified triangulation~$t$ of~$Y$. 
	\item 
	Pass to a substratification~$S^t_Y$ of~$Y$ as in Definition~\ref{def:StratifiedTriangulation}. 
	In the next steps we describe how to label the new strata in $S^t_Y \subset Y$. 
	\item 
	Let~$\widetilde\sigma_3$ be a new 3-stratum in $S^t_Y \setminus Y$; because of the transversality condition in Definition~\ref{def:StratifiedTriangulation}\,\ref{item:TransversalSimplex}, $\widetilde\sigma_3$ is part of a substratification of a 3-stratum~$\sigma_3$ in~$Y$. 
	Since~$Y$ is a $\mathds{D}_{\operatorname{orb}}$-labelled defect bordism, $\sigma_3$ is labelled by an element $\mathcal A = (a,A,\mu_A, \alpha_A, u_A, u^{\textrm{l}}_A, u^{\textrm{r}}_A) \in ({D}_{\operatorname{orb}})_3 = \operatorname{Ob}((\tric_\zz)_{\textrm{orb}})$. 
	We assign the label~$a$ to~$\widetilde\sigma_3$. 
	\item 
	Let~$\widetilde\sigma_2$ be a new 2-stratum in $S^t_Y \setminus Y$ which is part of a substratification of an $\mathcal A$-labelled 3-stratum~$\sigma_3$ in~$Y$. 
	We assign the label~$A$ to~$\widetilde\sigma_2$. 
	\item 
	Let~$\widetilde\sigma_2$ be a new 2-stratum in $S^t_Y \setminus Y$ which is part of a substratification of a 2-stratum~$\sigma_3$ in~$Y$ which is labelled by $\mathcal M = (M, \triangleright_M, \triangleleft_M, u_M^{\textrm{l}}, u_M^{\textrm{r}}, \alpha_M^{\textrm{l}},\alpha_M^{\textrm{m}}, \alpha_M^{\textrm{r}}) \in ({D}_{\operatorname{orb}})_2$. 
	We assign the label~$M$ to~$\widetilde\sigma_2$. 
	\item 
	Let~$\widetilde\sigma_1$ be a new 1-stratum in $S^t_Y \setminus Y$ which is Poincar\'e dual to a 2-simplex in~$t$, and hence part of a substratification of an $\mathcal A$-labelled 3-stratum~$\sigma_3$ in~$Y$. 
	We assign the label~$\mu_A$ to~$\widetilde\sigma_1$. 
	We also observe that due to the transversality condition in Definition~\ref{def:StratifiedTriangulation}\,\ref{item:TransversalSimplex}, $\sigma_1$ cannot be adjacent to an old 2-stratum in~$Y$. 
	\item 
	Let~$\widetilde\sigma_1$ be a new 1-stratum in $S^t_Y \setminus Y$ which is part of a substratification of 2-stratum~$\sigma_2$ in~$Y$ which is labeled by $\mathcal M = (M, \triangleright_M, \triangleleft_M, u_M^{\textrm{l}}, u_M^{\textrm{r}}, \alpha_M^{\textrm{l}},\alpha_M^{\textrm{m}}, \alpha_M^{\textrm{r}}) \in ({D}_{\operatorname{orb}})_2$. 
	We assign the label~$\triangleright_M$ or~$\triangleleft_M$ to~$\widetilde\sigma_1$. 
	Which action of an algebra (whose underlying 1-morphism in~$\tric_\zz$ is assigned to the new 2-stratum which meets~$\widetilde\sigma_1$ and which is part of a substratification of an old 3-stratum) is to be used is uniquely determined by the orientations of the strata involved as well as the graphical calculus in Definition~\ref{def:1MorphismsEAlg}. 
	\item 
	Let~$\widetilde\sigma_0$ be a new 0-stratum in $S^t_Y \setminus Y$ which is part of a substratification of an $\mathcal A$-labelled 3-stratum~$\sigma_3$ in~$Y$. 
	We assign the label~$\alpha_A$ or its inverse to~$\widetilde\sigma_0$, uniquely determined by the orientations of the strata involved as well as the graphical calculus in Definition~\ref{def:ObjEAlg}. 
	\item 
	Let~$\widetilde\sigma_0$ be a new 0-stratum in $S^t_Y \setminus Y$ which is part of a substratification of an $\mathcal M$-labelled 2-stratum~$\sigma_2$ in~$Y$. 
	We assign the label $\alpha_M^{\textrm{l}},\alpha_M^{\textrm{m}}$ or~$\alpha_M^{\textrm{r}}$ (or their inverse) to~$\widetilde\sigma_0$, uniquely determined by the orientations of the strata involved as well as the graphical calculus in Definition~\ref{def:1MorphismsEAlg}. 
	\item 
	Let~$\widetilde\sigma_0$ be a new 0-stratum in $S^t_Y \setminus Y$ which is part of a substratification of a 1-stratum~$\sigma_1$ in~$Y$ which is labelled by $\mathcal F = (F,\triangleleft_F, \triangleright_F)$. 
	We assign the label $\triangleleft_F$ or~$\triangleright_F$ (or their inverse) to~$\widetilde\sigma_0$, uniquely determined by the orientations of the strata involved as well as the graphical calculus in Definition~\ref{def:2MorphismsEAlg}. 
	\item 
	We have thus lifted~$Y$ to a morphism $Y^t\colon \varnothing \lra \varnothing$ in $\DBord(\mathds{D})$, and we define 
	\be 
	\mathcal{Z}_{\operatorname{orb}}(Y) := \zz(Y^t) \, . 
	\ee 
\end{enumerate} 
The relations imposed in the construction of orbifold completion imply that $\mathcal{Z}_{\operatorname{orb}}(Y)$ does not depend on the specific choice of the stratified triangulation~$t$ and the induced substratification~$S^t_Y$, as follows directly from the discussion in~\cite[Sect.\,3.4]{CRS1} and Definition~\ref{def:Corb}. 

\medskip 

Extending the invariant $\mathcal{Z}_{\operatorname{orb}}(\varnothing\stackrel{Y}{\lra}\varnothing)$ to a full defect TQFT can be done in a standard way following \cite[Constr.\,3.8\,\&\,3.9]{CRS1}. 
We first describe the evaluation of~$\mathcal{Z}_{\operatorname{orb}}$ on objects. 
The state space $\zz_{\operatorname{orb}}(\Sigma)$ associated to a stratified $\mathds{D}_{\operatorname{orb}}$-labelled 2-manifold $\Sigma \in \DBord (\mathds{D}_{\textrm{orb}})$ is constructed in two steps. 
For any stratified triangulation~$\tau$ of~$\Sigma$, its Poincar\'e dual stratification~$S^\tau$ refines the original stratification of~$\Sigma$ to a stratification~$S^\tau_\Sigma$, whose new strata we label with elements of~$\mathds{D}_{\textrm{orb}}$ in analogy to the above discussion. 
This produces a new defect surface $\Sigma^\tau \in \DBord (\mathds{D}_{\textrm{orb}})$, and we can consider the vector space $\zz(\Sigma^\tau)$. 

For two distinct triangulations~$\tau$ and $\tau'$, the vector spaces $\zz(\Sigma^\tau)$ and $\zz(\Sigma^{\tau'})$ are not necessarily isomorphic or even equal. 
Any stratified triangulation~$t$ of $\Sigma\times [0,1]$ restricting to~$\tau$ and~$\tau'$ on the boundary gives rise to a linear map $C_\Sigma^t := \cat{Z}((\Sigma)^t \times [0,1])\colon \mathcal{Z}(\Sigma^{\tau}) \to \mathcal{Z}(\Sigma^{\tau'})$ which is independent of~$t$. 
The value of $\mathcal{Z}_{\operatorname{orb}}$ on~$\Sigma$ is defined to be the colimit over all stratified triangulations of~$\Sigma$ related by the maps~$C_\Sigma^t$. 
This vector space $\mathcal{Z}_{\operatorname{orb}}(\Sigma)$ is isomorphic to the image of the idempotent corresponding to the cylinder for any fixed stratified triangulation $\tau'=\tau$ of~$\Sigma$. 

Finally to construct the value on a bordism $Y\colon \Sigma\to \Sigma' \in \DBord(\mathds{D}_{\operatorname{orb}})$ we fix a stratified triangulation~$t$ of~$Y$ extending fixed triangulations~$\tau$ and~$\tau'$ of its source and target, respectively. 
As above we can construct a labelling of the dual triangulation by $\mathds{D}$ and hence get a morphism $Y^t\colon \Sigma^\tau \to (\Sigma')^{\tau'}$ in $\DBord(\mathds{D})$. 
By evaluating the original defect TQFT, this gives a linear map 
\begin{align}
\cat{Z}(Y^t)\colon \cat{Z}(\Sigma^\tau) \longrightarrow \cat{Z}((\Sigma')^{\tau'}) 
\end{align} 
which only depends on~$\tau$ and~$\tau'$, and which induces a well-defined map 
\begin{align}
	\cat{Z}_{\operatorname{orb}}(Y)\colon \cat{Z}_{\operatorname{orb}}(\Sigma) \longrightarrow \cat{Z}_{\operatorname{orb}}(\Sigma') 
\end{align} 
between the colimits. 
It follows directly from the construction that $\cat{Z}_{\operatorname{orb}}$ is a defect TQFT: 

\begin{theorem}\label{thm:OrbifoldDefectTQFT}
Let $\cat{Z}\colon \DBord (\mathds{D})\to \vs$ be a defect TQFT. 
The orbifold construction above gives a defect TQFT $\cat{Z}_{\operatorname{orb}}\colon \DBord (\mathds{D}_{\operatorname{orb}})\to \vs$. 
\end{theorem}

We call $\cat{Z}_{\operatorname{orb}}$ the \textsl{orbifold (TQFT) of~$\zz$}. 
Note that the definite article is appropriate; all the orbifold \textsl{closed} TQFTs~$\zz_{\mathcal A}$ constructed for some $\mathcal A \in (\tric_\zz)_{\textrm{orb}}$ in \cite{CRS1} are obtained from $\cat{Z}_{\operatorname{orb}}$ by restricting to trivially stratified $\mathcal A$-labelled bordisms. 
 
\begin{example}
The 3-category $E(\Bar\ssFrob)$ describes defects between 3-dimensional Euler theories of \cite{CRS3}, cf.\ Section~\ref{subsec:2Dstatesums}. 
The 3-category $\catf{sFus}$ of spherical fusion categories is a subcategory of the orbifold completion of $E(\Bar\ssFrob)$ by Theorem~\ref{thm:OrbOfTriv}. 
Hence from Theorem~\ref{thm:OrbifoldDefectTQFT} we get an associated defect TQFT. 
As shown in~\cite{CRS3}, the closed TQFT corresponding to an object in this orbifold construction is the one constructed by Turaev--Viro--Barrett--Westbury. 
Hence Theorem~\ref{thm:OrbifoldDefectTQFT} constructs defects between 3-dimensional state sum models from bimodule categories with trace and their compatible functors. 
We expect them to agree with the defects constructed in~\cite{Meusburger3dDefectStateSumModels} when considering only line defects with precisely two adjacent surface defects. 
We leave a detailed comparison for future work.  
\end{example} 

\begin{example}
Let $\cat{M}$ be a modular fusion category.
Defects in the corresponding Reshetikhin--Turaev theory are described by $\Bar\Delta\ssFrob(\cat{M})$. 
In Theorem~\ref{thm:OrbOfssFrobM} we identified a subcategory of its orbifold completion. 
The corresponding defect TQFT agrees with the one constructed in~\cite{KMRS}. 
Our construction extends theirs in the sense that we also describe the state spaces and linear maps associated to bordism with non-empty boundary. 

A special case of the constraints on 2-morphisms of $\orb{\tric}$ (recall Definition~\ref{def:Corb}) was first considered in the context of (a 1-categorical description of) Reshetikhin--Turaev theory in \cite[Fig.\,1.1\,\&\,3.1]{MuleRunk}, and used further in \cite{CMRSS1}. 
It is worth noting that from the 3-categorical Morita theoretic perspective adopted here, most of these constraints are identified as instances of universal coherence relations. 
\end{example}

\appendix

\section{2-(co)limits}
\label{app:2colimits}

In this appendix we recall some definitions related to 2-colimits relevant for the main part of the paper, see for example~\cite{JohnsonYauBook} for a text book treatment. 
As a simple consequence of the definitions we derive a coherence result for the computation of colimits over products. 

Before defining 2-colimits we recall that a 2-functor $\cat{F}\colon \cat{C}\to \Cat$ is \emph{representable} if there exists an object $c_\cat{F}\in \cat{C}$
and a natural isomorphism $\cat{C}(c_{\cat{F}},-) \longrightarrow \cat{F}(-) $ of
functors $\cat{C}\to \Cat$. The \emph{representing object} $c_\cat{F}$, if it exists, is unique
up to a contractible choice. 
In the 2-categorical setting this means that the 2-category $\operatorname{Rep}(\cat{F})$ with 
\begin{itemize}
	\item 
	objects pairs of an object
	$c_\cat{F}$ together with a natural isomorphism $\varphi_{c_\cat{F}}\colon \cat{C}(c_{\cat{F}},-) \longrightarrow \cat{F}(-)$, 
	\item 1-morphisms given by a 1-morphism $c_{\cat{F}}\xrightarrow{\  f  \ } c'_{\cat{F}}$ in $\cat{C}$ together with a 2-isomorphism 
	\begin{equation*}
	\begin{tikzcd} 
	& & \cat{C}(c_{\cat{F}}', -)\ar[dd, "f^*"] \\
	\cat{F}(-) \ar[rru,"\varphi_{c_\cat{F}}'"] \ar[rrd,"\varphi_{c_\cat{F}}",swap]  & & \ \\ 
 	& \ar[ruu,Rightarrow, "\vartheta_{f}" , shorten <= 20, shorten >= 1] & \cat{C}(c_{\cat{F}}, -)
	\end{tikzcd}
	\end{equation*} 
	\item and 2-morphisms given by 2-morphisms $\omega \colon f_1\to f_2$ in $\cat{C}$ which are compatible with $\vartheta_f$ in the sense that   
	\begin{equation*}
	\begin{tikzcd} 
	& & \cat{C}(c_{\cat{F}}', -)\ar[dd, "f_2^*"] \\
	\cat{F}(-) \ar[rru,"\varphi_{c_\cat{F}}'"] \ar[rrd,"\varphi_{c_\cat{F}}",swap]  &  &\ \\ 
	 & \ar[ruu,Rightarrow, "\vartheta_{f_2}" , shorten <= 20, shorten >= 1] & \cat{C}(c_{\cat{F}}, -)
	\end{tikzcd}
	= 
	\begin{tikzcd} 
	& & \cat{C}(c_{\cat{F}}', -)\ar[dd, "f_1^*",swap] \ar[dd, "f_2^*", bend left= 80] \\
	\cat{F}(-) \ar[rru,"\varphi_{c_\cat{F}}'"] \ar[rrd,"\varphi_{c_\cat{F}}",swap] & & \ \ar[r,Rightarrow, "\omega^*" , shorten <= 5, shorten >= 12] & \ \\ 
	 & \ar[ruu,Rightarrow, "\vartheta_{f_1}" , shorten <= 20, shorten >= 1] & \cat{C}(c_{\cat{F}}, -)
	\end{tikzcd}
	\end{equation*} 
	holds.
\end{itemize}
is either empty or equivalent to the 2-category with one object and only identity morphisms. This implies that
there is an isomorphism between any two objects $c_{\cat{F}}$ and $c'_{\cat{F}}$ and a unique 2-isomorphism between any pair of 1-morphisms. This in particular implies that if two compositions of 2-morphisms have the same source and target they must be equal.  
\begin{definition}
	Let $\cat{F} \colon {D} \to \cat{B}$ be a diagram in a 2-category $\cat{B}$. The \emph{2-colimit of $\cat{F}$}, if it exists, is an object $\colim_D \cat{F} \in \cat{B}$ which represents the functor 
	\begin{align*}
	\cat{B} & \to \Cat \\
	b & \longmapsto \operatorname{Nat}(\cat{F}, \Delta_b) \ \ ,   
	\end{align*}     
	where $\Delta_b$ denotes the constant diagram at $b$. 
\end{definition}
From the definition it is clear that if the 2-colimit of $\cat{F}$ exists then it is unique up to a contractible choice. 
Objects of the category $\operatorname{Nat}(\cat{F}, \Delta_b)$ are called \emph{cocones with apex $b$}.  
Evaluating the natural isomorphism $  \cat{B}(\colim_D \cat{F},-) \longrightarrow \operatorname{Nat}(\cat{F}, \Delta_-)$ at the identity $\colim_D \cat{F} \to \colim_D \cat{F}$ gives
rise to a cocone $\mathcal{F}(-)\longrightarrow \Delta_{\colim_D} \in \operatorname{Nat}(\cat{F}, \Delta_{\colim_D \cat{F}})$. The definition of a 2-colimit given above is equivalent to the condition 
that this cocone is universal, in the sense that 
the map \begin{align*}
(-)_*\colon \cat{B} (\colim_D \cat{F}, b)\to \operatorname{Nat}(\cat{F}, \Delta_{b})
\end{align*} induced by composition with the cocone is an equivalence. The contractibility of representing objects
translates to the contractibility of the 2-category of universal cocones. 

A 2-category $\cat{B}$ is called \emph{(finitely) cocomplete} if all (finite) 2-colimits exist in $\cat{B}$. In this case the weak uniqueness of colimits implies that any choice of 2-colimits will 
lead to a functor $\colim_D  \colon \cat{B}^D \to \cat{B}$. 
Let $\mathcal{G}\colon \cat{B}\to \cat{B}'$ be a functor and $\cat{F}\colon D \to \cat{B} $ a diagram. Note that applying $\mathcal{G}$ to the universal cocone $\cat{F}(-)\to \colim_D \cat{F}$ gives rise to a cocone over $\cat{G} \circ \cat{F}$ and hence up to contractible choice there is a 1-morphism $\colim_D \cat{G} \circ \cat{F} \to \cat{G}(\colim_D \cat{F})$. We say $\cat{G}$ \emph{preserves the colimit of $\cat{F}$} if this map is an equivalence. 
We say that $\cat{G}$ \emph{preserves (sifted) 2-colimits} if this map is an equivalence
for all (sifted) diagram categories. 

\begin{proposition}[Fubini theorem for 2-colimits]
	\label{prop:Fubini2}
	Let $\cat{F}\colon D_1\times D_2 \to \cat{B} $ be a diagram. There are 
	equivalences 
	\be 
		\colim_{d_1\in D_1}(\colim_{D_2} F(d_1,-)) \cong \colim_{D_1\times D_2} \cat{F} \cong \colim_{d_2\in D_2}(\colim_{D_1} F(-,d_2)) \, , 
	\ee 
	unique up to contractible choice.        
\end{proposition}
\begin{proof}
	We will only construct the equivalence for $\colim_{d_1\in D_1}(\colim_{D_2} \cat{F}(d_1,-))$. The second equivalence can be constructed analogously. For this
	it is enough to note that the map $F(d_1,d_2) \longrightarrow \colim_{D_2} \cat{F}(d_1,-) \longrightarrow \colim_{d_2\in D_2}(\colim_{D_1} F(-,d_2)) $ composing the maps which are part of the universal cocones, has the structure of a universal cocone
	for $\cat{F}(-,-)$. The statement now follows from the uniqueness of universal cocones.  
\end{proof}
Let us pick the equivalences described in the previous proposition. 
Then for a functor $\cat{F}\colon D_1\times D_2\times D_3 \to \cat{B}$ we can iterate the result of the proposition to get two equivalences 
\begin{equation*}
\begin{tikzcd} 
& \colim_{D_1\times D_2} (\colim_{ D_3}\cat{F}(d_1,d_2,-)) \ar[rd,"\cong" ]& 
\\
\colim_{D_1\times D_2\times D_3} \cat{F} \ar[ru, "\cong"] \ar[rd, "\cong",swap ] & & \colim_{D_1}(\colim_{D_2} (\colim_{D_3}\cat{F}(d_1,d_2,-)) \, .
\\ 
& \colim_{D_1} (\colim_{D_2\times D_3}\cat{F}(d_1,-,-)) \ar[ru, "\cong", swap]& 	
\end{tikzcd}
\end{equation*} 
In general there is no reason for these maps to agree. Furthermore, there is even a third isomorphism between $\colim_{D_1\times D_2\times D_3} \cat{F}$ and $ \colim_{D_1}(\colim_{D_2} (\colim_{D_3}\cat{F}(d_1,d_2,-))$ which we can construct by noting, similarly to the proof of the Fubini theorem, that $\colim_{D_1}(\colim_{D_2} (\colim_{D_3}\cat{F}(d_1,d_2,-))$ defines a universal cocone for $\cat{F}$. However, all these maps are related by unique 2-isomorphisms, since they all are part of maps of universal cocones.  

If we now consider diagrams indexed by the product of four categories there 
are even more ways to relate the different orders in which we can compute the 2-colimits. Again all these will be related by unique 2-isomorphisms. 
Furthermore, since the 2-category of universal cocones is contractible, two 2-isomorphisms which we can build by composition of those will agree as soon as they have the same source and target. 
This pattern continues to finite products of diagram categories and we record it in the following proposition.

\begin{proposition}[Coherence of the Fubini theorem for 2-colimits]
	\label{prop:Fubini2}
	Let $\cat{F}\colon D_1\times \dots \times D_n \to \cat{B}$ be a functor such that all the colimits appearing in the statement exist. The 2-category with 
	\begin{itemize}
		\item objects given by chosen representatives for different iterative ways and orders of computing the 2-colimit of $\cat{F}$ (note that all of them are canonically
		a universal cocone for $\cat{F}$),  
		\item 1-morphisms given by morphisms of universal cocones 
		\item 2-morphisms given by 2-morphisms compatible with the universal cocone structure   
	\end{itemize} 
	is contractible. 
\end{proposition}   
\begin{proof}
	The statement follows directly from the fact that the 2-category described in the proposition is a full subcategory of the 2-category of all universal cocones, which is contractible.  
\end{proof}


\begin{thebibliography}{CMRSS2}
	
\bibitem[Ab]{Abrams}
L.~Abrams,
\textsl{The quantum Euler class and the quantum cohomology of the Grassmannians}, 
\doi{10.1007/BF02773576}{Israel Journal of Mathematics \textbf{117} (2000), 335--352},
\href{https://arxiv.org/abs/q-alg/9712025}{arXiv:q-alg/9712025}.

\bibitem[At]{AtiyahTQFT}
M.~Atiyah, 
\textsl{Topological quantum field theories}, 
\href{http://www.numdam.org/item?id=PMIHES_1988__68__175_0}{Inst. Hautes \'{E}tudes Sci. Publ. Math. \textbf{68} (1988), 175--186}.

\bibitem[BCP]{BCP2}
I.~Brunner, N.~Carqueville, and D.~Plencner, 
\textsl{Discrete torsion defects}, 
\doi{10.1007/s00220-015-2297-9}{Comm. Math. Phys. \textbf{337} (2015), 429--453}, 
\href{https://arxiv.org/abs/1404.7497}{arXiv:1404.7497 [hep-th]}.

\bibitem[BJS]{BJS}
A.~Brochier, D.~Jordan, and N.~Snyder, 
\textsl{On dualizability of braided tensor categories}, 
\doi{10.1112/S0010437X20007630}{Compositio Mathematica \textbf{157}:3 (2021), 435--483}, 
\href{https://arxiv.org/abs/1804.07538}{arXiv:1804.07538 [math.CT]}.

\bibitem[BJSS]{BJSS}
A.~Brochier, D.~Jordan, P.~Safronov, and N.~Snyder, 
\textsl{Invertible braided tensor categories}, 
\doi{10.2140/agt.2021.21.2107}{Algebraic \& Geometric Topology \textbf{21} (2021), 2107--2140}, 
\href{https://arxiv.org/abs/2003.13812}{arXiv:2003.13812 [math.QA]}.

\bibitem[BMS]{BMS}
J.~Barrett, C.~Meusburger, and G.~Schaumann, 
\textsl{Gray categories with duals and their diagrams},
\href{https://arxiv.org/abs/1211.0529}{arXiv:1211.0529v3 [math.QA]}.

\bibitem[Ca]{2dDefectTQFTLectureNotes}
N.~Carqueville, 
\textsl{Lecture notes on 2-dimensional defect TQFT}, 
\doi{10.4064/bc114-2}{Banach Center Publications \textbf{114} (2018), 49--84}, 
\href{https://arxiv.org/abs/1607.05747}{\mbox{arXiv:}1607.05747 [math.QA]}. 

\bibitem[CDIS]{snow}
C.~Cordova, T.~T.~Dumitrescu, K.~Intriligator, and S.-H.~Shao, 
\textsl{Snowmass White Paper: Generalized Symmetries in Quantum Field Theory and Beyond}, 
\href{https://arxiv.org/abs/2205.09545}{arXiv:2205.09545 [hep-th]}.

\bibitem[CMRSS1]{CMRSS1}
N.~Carqueville, V.~Mulevi\v{c}ius, I.~Runkel, D.~Scherl, and G.~Schaumann, 
\textsl{Orbifold graph TQFTs}, 
\href{https://arxiv.org/abs/2101.02482}{arXiv:2101.02482 [math.QA]}.

\bibitem[CMRSS2]{CMRSS2}
N.~Carqueville, V.~Mulevi\v{c}ius, I.~Runkel, D.~Scherl, and G.~Schaumann, 
\textsl{Reshetikhin--Turaev TQFTs close under generalised orbifolds}, 
%arXiv_v3:
	\doi{10.1007/s00220-024-05068-6}{Communications in Mathematical Physics \textbf{405} (2024), 242}, 
\href{https://arxiv.org/abs/2109.04754}{arXiv:2109.04754 [math.QA]}.

\bibitem[CMM]{LukasNilsVincentas}
N.~Carqueville, V.~Mulevi\v{c}ius, and L.~Müller, 
\textsl{in preparation}.

\bibitem[CMS]{CMS}
N.~Carqueville, C.~Meusburger, and G.~Schaumann, 
\textsl{3-dimensional defect TQFTs and their tricategories}, 
\doi{10.1016/j.aim.2020.107024}{Adv. Math. \textbf{364} (2020) 107024},
\href{https://arxiv.org/abs/1603.01171}{arXiv:1603.01171 [math.QA]}.

\bibitem[Co]{Cooke2019}
J.~Cooke, 
\textsl{Excision of Skein Categories and Factorisation Homology}, 
\href{https://arxiv.org/abs/1910.02630}{arXiv:1910.02630 [math.QA]}.

\bibitem[CR]{cr1210.6363}
N.~Carqueville and I.~Runkel, 
\textsl{Orbifold completion of defect bicategories}, 
\doi{10.4171/QT/76}{Quantum Topology \textbf{7}:2 (2016) 203--279}, 
\href{https://arxiv.org/abs/1210.6363}{arXiv:1210.6363 [math.QA]}.

\bibitem[CRCR]{CRCR}
N.~Carqueville, A.~Ros Camacho, and I.~Runkel, 
\textsl{Orbifold equivalent potentials}, 
\doi{10.1016/j.jpaa.2015.07.015}{Journal of Pure and Applied Algebra \textbf{220} (2016), 759--781}, 
\href{https://arxiv.org/abs/1311.3354}{arXiv:1311.3354 [math.QA]}.

\bibitem[CRS1]{CRS1}
N.~Carqueville, I.~Runkel, and G.~Schaumann, 
\textsl{Orbifolds of $n$-dimensional defect TQFTs}, 
\doi{10.2140/gt.2019.23.781}{Geometry \& Topology \textbf{23} (2019), 781--864},  
\href{https://arxiv.org/abs/1705.06085}{arXiv:1705.06085 [math.QA]}.

\bibitem[CRS2]{CRS2}
N.~Carqueville, I.~Runkel, and G.~Schaumann, 
\textsl{Line and surface defects in Reshetikhin--Turaev TQFT},  
\doi{10.4171/QT/121}{Quantum Topology \textbf{10} (2019), 399--439}, 
\href{https://arxiv.org/abs/1710.10214}{arXiv:1710.10214 [math.QA]}.

\bibitem[CRS3]{CRS3}
N.~Carqueville, I.~Runkel, and G.~Schaumann, 
\textsl{Orbifolds of Reshetikhin--Turaev TQFTs}, 
\href{http://www.tac.mta.ca/tac/volumes/35/15/35-15abs.html}{Theory and Applications of Categories \textbf{35} (2020), 513--561}, 
\href{https://arxiv.org/abs/1809.01483}{arXiv:1809.01483 [math.QA]}.

\bibitem[D{\'e1}]{2EW}
T.~D.~D{\'e}coppet, 
\textsl{Multifusion Categories and Finite Semisimple 2-Categories},  
\doi{10.1016/j.jpaa.2022.107029}{J.\ Pure Appl.\ Algebra \textbf{226} (2022), 107029}, 
\href{https://arxiv.org/pdf/2012.15774}{arXiv:2012.15774 [math.QA]}.

%arXiv_v2:
	\bibitem[D{\'e}2]{Decoppet_Rig}
	T.~D.~D{\'e}coppet, 
	\textsl{Rigid and separable algebras in fusion 2-categories},  
	\doi{10.1016/j.aim.2023.108967}{Adv. Math. \textbf{419} (2023), 108967}, 
	\href{https://arxiv.org/abs/2205.06453}{arXiv:2205.06453 [math.QA]}.
	
%arXiv_v2:
	\bibitem[D{\'e}3]{Decoppet_Morita}
	T.~D.~D{\'e}coppet, 
	\textsl{The Morita Theory of Fusion 2-Categories},  
	\doi{10.21136/HS.2023.07}{Higher Structures \textbf{7}:1 (2023), 234--292}, 
	\href{https://arxiv.org/abs/2208.08722}{arXiv:2208.08722 [math.CT]}.

\bibitem[DSPS]{DuaTen}
C.~Douglas, C.~Schommer-Pries, and N.~Snyder, 
\textsl{Dualizable tensor categories}, 
\doi{10.1090/memo/1308}{Memoirs of the American Mathematical Society \textbf{268}, Number 1308, (2020)},
\href{https://www.arxiv.org/abs/1312.7188}{arXiv:1312.7188 [math.QA]}.

\bibitem[DKR]{dkr1107.0495}
A.~Davydov, L.~Kong, and I.~Runkel, 
\textsl{Field theories with defects and the centre functor}, \href{https://www.ams.org/bookstore?fn=20&arg1=pspumseries&ikey=PSPUM-83}{Mathematical Foundations of Quantum Field Theory and Perturbative String Theory, Proceedings of Symposia in Pure Mathematics, AMS, 2011}, 
\href{https://arxiv.org/abs/1107.0495}{arXiv:1107.0495 [math.QA]}.

\bibitem[DR]{DouglasReutter2018}
C.~Douglas and D.~Reutter, 
\textsl{Fusion 2-categories and a state-sum invariant for 4-manifolds}, 
\href{https://arxiv.org/abs/1812.11933}{arXiv:1812.11933 [math.QA]}. 

\bibitem[FFRS]{ffrs0909.5013}
J.~Fr\"ohlich, J.~Fuchs, I.~Runkel, and C.~Schweigert,
\textsl{Defect lines, dualities, and generalised orbifolds}, 
\doi{10.1142/9789814304634_0056}{Proceedings of the XVI International Congress on Mathematical Physics, Prague, August 3--8, 2009}, \href{https://arxiv.org/abs/0909.5013}{arXiv:0909.5013 [math-ph]}.

\bibitem[Fr]{FrauenbergerMasterThesis}
M.~Frauenberger, 
\textsl{Condensation monads in low dimension}, 
\doi{10.25365/thesis.72024}{Master thesis, University of Vienna (2022)}, 
\href{http://nbn-resolving.de/urn/resolver.pl?urn:nbn:at:at-ubw:1-17977.14847.768365-2}{urn:nbn:at:at-ubw:1-17977.14847.768365-2}.

\bibitem[FSV]{fsv1203.4568}
J.~Fuchs, C.~Schweigert, and A.~Valentino, 
\textsl{Bicategories for boundary conditions and for surface defects in 3-d TFT}, 
\doi{10.1007/s00220-013-1723-0}{Communications in Mathematical Physics \textbf{321}:2 (2013), 543--575}, 
\href{http://arxiv.org/abs/1203.4568}{arXiv:1203.4568 [hep-th]}.

\bibitem[GJF]{GaiottoJohnsonFreyd}
D.~Gaiotto and T.~Johnson-Freyd,
\textsl{Condensations in higher categories},
\arxiv{1905.09566}{arXiv:1905.09566 [math.CT]}.

\bibitem[GPS]{GPS}
R.~Gordon, A.~J.~Power, and R.~Street, 
\textsl{Coherence for Tricategories}, 
Memoirs of the American Mathematical Society \textbf{117}, 
American Mathematical Society, 1995.

\bibitem[GS]{ClaudiaOwen}
O.~Gwilliam and C.~Scheimbauer,
\textsl{Duals and adjoints in higher Morita categories},
\arxiv{1804.10924}{arXiv:1804.10924 [math.CT]}.

\bibitem[Gu]{Gurskibook}
N.~Gurski, 
\textsl{\doi{10.1017/CBO9781139542333}{Coherence in Three-Dimensional Category Theory}}, 
\textsl{Cambridge Tracts in Mathematics} \textbf{201}, Cambridge University Press, 2013.

\bibitem[He]{FEW}
J.~Hesse, 
\textsl{An equivalence between semisimple symmetric Frobenius algebras and Calabi-Yau categories}, 
\doi{10.1007/s40062-017-0181-3}{Journal of Homotopy and Related Structures \textbf{13} (2018), 251--272},  
\href{https://arxiv.org/abs/1609.06475}{arXiv:1609.06475 [math.QA]}.

\bibitem[JFS]{TheoClaudia}
T.~Johnson-Freyd and C.~Scheimbauer, 
\textsl{(Op)lax natural transformations, twisted quantum field theories, and ``even higher'' Morita categories}, 
\doi{10.1016/j.aim.2016.11.014}{Advances in Mathematics \textbf{307} (2017), 147--223},  
\href{https://arxiv.org/abs/1502.06526}{arXiv:1502.06526 [math.CT]}.

\bibitem[JY]{JohnsonYauBook}
N.~Johnson and D.~Yau, 
\textsl{2-Dimensional Categories}, 
\doi{10.1093/oso/9780198871378.001.0001}{Oxford University Press (2021)}, 
\href{https://arxiv.org/abs/2002.06055}{arXiv:2002.06055 [math.CT]}.

\bibitem[KMRS]{KMRS}
V.~Koppen, V.~Mulevi\v{c}ius, I.~Runkel, and C.~Schweigert, 
\textsl{Domain walls between 3d phases of Reshetikhin-Turaev TQFTs}, 
\doi{10.1007/s00220-022-04489-5}{Communications in Mathematical Physics \textbf{396} (2022), 1187--1220}, 
\href{https://arxiv.org/abs/2105.04613}{arXiv:2105.04613 [hep-th]}.

\bibitem[Ko]{Kockbook}
J.~Kock, 
\textsl{Frobenius algebras and 2D topological quantum field theories}, 
\textsl{London Mathematical Society Student Texts} \textbf{59}, Cambridge University Press, 2003.

\bibitem[KS]{ks1012.0911}
A.~Kapustin and N.~Saulina, 
\textsl{Surface operators in 3d Topological Field Theory and 2d Rational Conformal Field Theory}, 
Mathematical Foundations of Quantum Field Theory and Perturbative String Theory, 
Proceedings of Symposia in Pure Mathematics \textbf{83}, 175--198, 
American Mathematical Society, 2011, 
\href{https://arxiv.org/abs/1012.0911}{arXiv:1012.0911 [hep-th]}.

\bibitem[Lu]{Lurie}
J.~Lurie, 
\textsl{On the Classification of Topological Field Theories},
\href{https://projecteuclid.org/euclid.cdm/1254748657}{Current Developments in Mathematics \textbf{2008} (2009), 129--280}, 
\arxiv{0905.0465}{arXiv:0905.0465 [math.CT]}.

%arXiv_v3: 
	\bibitem[NRC]{NRC1}
	A.~Ros Camacho and R.~Newton, 
	\textsl{Strangely dual orbifold equivalence I}, 
	\doi{10.5427/jsing.2016.14c}{Journal of Singularities \textbf{14} (2016), 34--51}, 
	\href{https://arxiv.org/abs/1509.08069}{arXiv:1509.08069 [math.QA]}.

\bibitem[Me]{Meusburger3dDefectStateSumModels}
C.~Meusburger, 
\textsl{State sum models with defects based on spherical fusion categories}, 
\href{https://arxiv.org/abs/2205.06874}{arXiv:2205.06874 [math.QA]}.

\bibitem[MR1]{MuleRunk}
V.~Mulevi\v{c}ius and I.~Runkel, 
\textsl{Constructing modular categories from orbifold data}, 
\doi{10.4171/QT/170}{Quantum Topology \textbf{13}:3 (2022), 459--523}, 
\href{https://arxiv.org/abs/2002.00663}{arXiv:2002.00663 [math.QA]}.

\bibitem[MR2]{MuleRunk2}
V.~Mulevi\v{c}ius and I.~Runkel, 
\textsl{Fibonacci-type orbifold data in Ising modular categories}, 
\doi{10.1016/j.jpaa.2022.107301}{Journal of Pure and Applied Algebra \textbf{227} (2022)},
\href{https://arxiv.org/abs/2010.00932}{arXiv:2010.00932 [math.QA]}.

\bibitem[Mu]{Mule1}
V.~Mulevi\v{c}ius, 
\textsl{Condensation inversion and Witt equivalence via generalised orbifolds}, 
\href{https://arxiv.org/abs/2206.02611}{arXiv:2206.02611 [math.QA]}.

\bibitem[RW]{OEReck}
A.~Recknagel and P.~Weinreb, 
\textsl{Orbifold equivalence: structure and new examples}, 
\doi{10.5427/jsing.2018.17j}{Journal of Singularities \textbf{17} (2018), 216--244}, 
\href{https://arxiv.org/abs/1708.08359}{arXiv:1708.08359 [math.QA]}.

\bibitem[Sc1]{GregorDiss}
G.~Schaumann, 
\textsl{Duals in tricategories and in the tricategory of bimodule categories}, 
PhD thesis, 
University of Erlangen-N\"urnberg (2013), 
\href{http://nbn-resolving.de/urn/resolver.pl?urn:nbn:de:bvb:29-opus4-37321}{urn:nbn:de:bvb:29-opus4-37321}.

\bibitem[Sc2]{Bimodtrace}
G.~Schaumann, 
\textsl{Traces on module categories over fusion categories},
\doi{10.1016/j.jalgebra.2013.01.013}{Journal of Algebra \textbf{379} (2013), 382--425}, 
\href{https://arxiv.org/abs/1206.5716}{arXiv:1206.5716 [math.QA]}.

\bibitem[St]{street}
R.~Street, 
\textsl{Frobenius monads and pseudomonoids}, 
\doi{10.1063/1.1788852}{J.~Math.~Phys.\ \textbf{37945} (2004), 3930}. 

\bibitem[Tu]{BookTuraev}
V.~Turaev, 
\doi{10.1515/9783110435221}{\textsl{Quantum invariants of knots and 3-manifolds}, 
de Gruyter, 
2010,
2nd edition}.

\bibitem[Tr]{TrimbleSurfaceDiagrams}
T.~Trimble, 
\textsl{Surface diagrams}, 
\url{https://ncatlab.org/toddtrimble/published/Surface+diagrams}.

\end{thebibliography}
\end{document}